\documentclass[a4paper,14pt,oneside,openany,oldfontcommands]{memoir}

\usepackage[T1,T2A]{fontenc}                    
\usepackage[utf8]{inputenc}[2014/04/30]         
\usepackage[english, russian]{babel}[2014/03/24]
\usepackage{graphicx}

\usepackage{longtable,ltcaption}                    
\usepackage{multirow,makecell}                      
\usepackage{tabu, tabulary}  

\usepackage{setspace}
\usepackage{titlesec}
\usepackage{lipsum}
\usepackage{indentfirst}
\usepackage{xpatch}

\usepackage{geometry}  
\usepackage{hyperref}
\usepackage{ragged2e}

\usepackage{amsthm,amsmath,amscd}   
\usepackage{amsfonts,amssymb}            
\usepackage{mathtools}                           

\usepackage{csquotes}
\usepackage[
    backend=biber, 
    citestyle=numeric-comp,
    language=auto,
    autolang=other,
    sorting=nyt,
    style=gost-numeric,
    natbib=true
]{biblatex}
\vfuzz \hfuzz
\renewcommand{\Large}{\fontsize{16}{16pt}\selectfont} 
\renewenvironment{proof}{{\bfseries Доказательство.}}

\newtheorem{theorem}{\hspace{1.25cm}Теорема}[section]
\newtheorem{lemma}{\hspace{1.25cm}Лемма}[section]
\newtheorem{consequence}{\hspace{1.25cm}Следствие}[section]
\newtheorem{remark}{\hspace{1.25cm}Замечание}[section]
\newtheorem{propose}{\hspace{1.25cm}Предположение}[section]

\OnehalfSpacing*    
\setSpacing{1.42}   
\xpatchcmd{\proof}{\hskip\labelsep}{\hskip7\labelsep}{}{}  
\titlespacing*{\section}{0pt}{10pt}{10pt}
\titlespacing*{\chapter}{0pt}{0pt}{0pt}

\usepackage{xcolor}

\hypersetup{
    linktocpage=true,  
    plainpages=false,  
    colorlinks=true,
    pdflang={ru},
    citecolor=[rgb]{0,0,0},
    urlcolor=[rgb]{0,0,0},
    linkcolor=[rgb]{0,0,0}
}

\setrmarg{2.55em plus1fil}  

\setsecheadstyle{\bfseries\centering} 

\makepagestyle{plain}
\makeevenfoot{plain}{}{\thepage}{}
\makeoddfoot{plain}{}{\thepage}{}
\begin{document}
\selectlanguage{russian}

\thispagestyle{empty}
\begin{center}
\textbf{КАЗАНСКИЙ ФЕДЕРАЛЬНЫЙ УНИВЕРСИТЕТ}
\end{center}

\vspace{5cm}

{\Large
\begin{center}
\textbf{И.Я. ЗАБОТИН, Р.С. ЯРУЛЛИН }
\end{center}

\begin{center}
\textbf{М\hspace*{-0.4mm}Е\hspace*{-0.4mm}Т\hspace*{-0.4mm}О\hspace*{-0.4mm}Д\hspace*{-0.4mm}Ы О\hspace*{-0.4mm}Т\hspace*{-0.4mm}С\hspace*{-0.4mm}Е\hspace*{-0.4mm}Ч\hspace*{-0.4mm}Е\hspace*{-0.4mm}Н\hspace*{-0.4mm}И\hspace*{-0.4mm}Й Б\hspace*{-0.4mm}Е\hspace*{-0.4mm}З В\hspace*{-0.4mm}Л\hspace*{-0.4mm}О\hspace*{-0.4mm}Ж\hspace*{-0.4mm}Е\hspace*{-0.4mm}Н\hspace*{-0.4mm}И\hspace*{-0.4mm}Я
\mbox{А\hspace*{-0.4mm}П\hspace*{-0.4mm}П\hspace*{-0.4mm}Р\hspace*{-0.4mm}О\hspace*{-0.4mm}К\hspace*{-0.4mm}С\hspace*{-0.4mm}И\hspace*{-0.4mm}М\hspace*{-0.4mm}И\hspace*{-0.4mm}Р\hspace*{-0.4mm}У\hspace*{-0.4mm}Ю\hspace*{-0.4mm}Щ\hspace*{-0.4mm}И\hspace*{-0.4mm}Х М\hspace*{-0.4mm}Н\hspace*{-0.4mm}О\hspace*{-0.4mm}Ж\hspace*{-0.4mm}Е\hspace*{-0.4mm}С\hspace*{-0.4mm}Т\hspace*{-0.4mm}В Д\hspace*{-0.4mm}Л\hspace*{-0.4mm}Я З\hspace*{-0.4mm}А\hspace*{-0.4mm}Д\hspace*{-0.4mm}А\hspace*{-0.4mm}Ч}
М\hspace*{-0.4mm}А\hspace*{-0.4mm}Т\hspace*{-0.4mm}Е\hspace*{-0.4mm}М\hspace*{-0.4mm}А\hspace*{-0.4mm}Т\hspace*{-0.4mm}И\hspace*{-0.4mm}Ч\hspace*{-0.4mm}Е\hspace*{-0.4mm}С\hspace*{-0.4mm}К\hspace*{-0.4mm}О\hspace*{-0.4mm}Г\hspace*{-0.4mm}О П\hspace*{-0.4mm}Р\hspace*{-0.4mm}О\hspace*{-0.4mm}Г\hspace*{-0.4mm}Р\hspace*{-0.4mm}А\hspace*{-0.4mm}М\hspace*{-0.4mm}М\hspace*{-0.4mm}И\hspace*{-0.4mm}Р\hspace*{-0.4mm}О\hspace*{-0.4mm}В\hspace*{-0.4mm}А\hspace*{-0.4mm}Н\hspace*{-0.4mm}И\hspace*{-0.4mm}Я}
\end{center}
}

\vspace{0pt plus1fill}
\vspace{1.5cm}

\begin{center}
\includegraphics{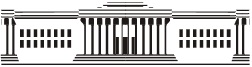}
\end{center}

\begin{center}
\textbf{КАЗАНЬ}\linebreak
\textbf{2019}
\end{center}

\setSpacing{1.2}
\thispagestyle{empty}
\begin{flushleft}
\textbf{УДК 519.85}\linebreak
\textbf{ББК 22.18 }\linebreak
\hspace*{1.32cm} \textbf{З-12 } 
\end{flushleft}
\begin{center}
\textit{Печатается по решению \linebreak
Редакционно-издательского совета\linebreak
Казанского федерального университета}
\end{center}
\begin{center}
\textbf{Научный редактор  }\linebreak
доктор физико-математических наук, профессор кафедры\linebreak
вычислительной математики КФУ \textbf{И.Б. Бадриев }
\end{center}
\begin{center}
\textbf{ Рецензенты:}\linebreak
доктор технических наук, профессор кафедры прикладной\linebreak
 математики и информатики КНИТУ-КАИ \textbf{Ш.И. Галиев;}\linebreak
кандидат физико-математических наук, доцент кафедры системного\linebreak
анализа и информационных технологий КФУ \textbf{А.А. Андрианова}
\end{center}
{\raggedleft
\hspace*{1.25cm}\textbf{Заботин И.Я.}}\hspace*{15.25cm}\linebreak
\textbf{З-12}\hspace*{0.165cm}\textbf{\mbox{М\hspace{-0.1mm}е\hspace{-0.1mm}т\hspace{-0.1mm}о\hspace{-0.1mm}д\hspace{-0.1mm}ы\hspace{0.16cm}о\hspace{-0.1mm}т\hspace{-0.1mm}с\hspace{-0.1mm}е\hspace{-0.1mm}ч\hspace{-0.1mm}е\hspace{-0.1mm}н\hspace{-0.1mm}и\hspace{-0.1mm}й\hspace{0.16cm}б\hspace{-0.1mm}е\hspace{-0.1mm}з\hspace{0.16cm}в\hspace{-0.1mm}л\hspace{-0.1mm}о\hspace{-0.1mm}ж\hspace{-0.1mm}е\hspace{-0.1mm}н\hspace{-0.1mm}и\hspace{-0.1mm}я\hspace{0.15cm}а\hspace{-0.1mm}п\hspace{-0.1mm}п\hspace{-0.1mm}р\hspace{-0.1mm}о\hspace{-0.1mm}к\hspace{-0.1mm}с\hspace{-0.1mm}и\hspace{-0.1mm}м\hspace{-0.1mm}и\hspace{-0.1mm}р\hspace{-0.1mm}у\hspace{-0.1mm}ю\hspace{-0.1mm}щ\hspace{-0.1mm}и\hspace{-0.1mm}х\hspace{0.15cm}м\hspace{-0.1mm}н\hspace{-0.1mm}о\hspace{-0.1mm}ж\hspace{-0.1mm}е\hspace{-0.1mm}с\hspace{-0.1mm}т\hspace{-0.1mm}в}}
\linebreak\hspace*{1.25cm}\mbox{\textbf{для задач  математического  программирования}\hspace{2.2mm}/\hspace{2.2mm}И.Я. Заботин,}
\linebreak\hspace*{1.25cm}Р.С. Яруллин. –  Казань:  Издательство Казанского университета, 2019. –
\linebreak\hspace*{1.25cm}182 с.
\begin{flushleft}
\textbf{ISBN  978-5-00130-223-0}
\end{flushleft}

\mbox{В связи с\hspace{0.2cm}потребностями\hspace{0.148cm}решения\hspace{0.198cm}оптимизационных задач~разработка} 
\mbox{методов\hspace{0.89mm}условной\hspace{0.89mm}минимизации\hspace{0.99mm}с\hspace{0.9mm}удобной численной\hspace{0.89mm}реализацией\hspace{0.89mm}продолжает} 
\mbox{пользоваться\hspace{0.7mm}вниманием\hspace{0.5mm}математиков.\hspace{0.7mm}В\hspace{0.6mm}настоящей\hspace{0.6mm}монографии\hspace{0.6mm}предлагаются}
\mbox{методы\hspace{2mm}решения\hspace{2.2mm}задач математического\hspace{2.3mm}программирования,\hspace{2mm}относящиеся\hspace{1.5mm}к}
\mbox{классу\hspace{2mm}методов\hspace{2mm}отсечений.\hspace{2.1mm}Предлагаемые\hspace{2mm}методы характерны,\hspace{2.1mm}в\hspace{2mm}частности,}
\mbox{тем,\hspace{1.6mm}что\hspace{1.6mm}в\hspace{1.6mm}них\hspace{1.6mm}заложена\hspace{1.6mm}возможность\hspace{1.6mm}периодического\hspace{1.6mm}обновления\hspace{1.6mm}аппрокси-}
\mbox{мирующих\hspace{3.7mm}множеств\hspace{3.7mm}за\hspace{3.5mm}счет\hspace{3.7mm}отбрасывания\hspace{3.7mm}накапливающихся\hspace{3.5mm}секущих}
\mbox{плоскостей, что делает эти методы удобными с практической точки зрения.}

\mbox{Монография\hspace{5mm}может\hspace{5mm}быть\hspace{4.5mm}использована\hspace{5mm}специалистами в\hspace{2mm}области}
\mbox{математического\hspace{2.8mm}программирования,\hspace{2.85mm}а\hspace{2.85mm}также\hspace{2.8mm}студентами,\hspace{2.8mm}обучающимися}
\mbox{по\hspace{2mm}специальностям\hspace{2.3mm}«Прикладная\hspace{2mm}математика»,\hspace{2.1mm}«Прикладная математика\hspace{2.1mm}и}
\mbox{информатика» и др.}

\begin{flushright}
\textbf{УДК 519.85}\linebreak
\textbf{ББК 22.18 }
\end{flushright}
\begin{flushleft}
\textbf{ISBN  978-5-00130-223-0}
\end{flushleft}
\begin{flushright}
\textbf{© Заботин И.Я., Яруллин Р.С., 2019\linebreak
© Издательство Казанского университета, 2019} 
\end{flushright}

\OnehalfSpacing*
\setSpacing{1.42}
{\Large
\tableofcontents*
\chapter*{Введение}
\addcontentsline{toc}{chapter}{Введение}

Довольно часто в различных областях прикладной математики возникают оптимизационные задачи, которые решаются методами нелинейного программирования. В связи с этим разработка методов условной минимизации с удобной численной реализацией продолжает пользоваться вниманием специалистов в области математического программирования. К настоящему времени накоплено значительное число методов решения задач нелинейного программирования. Различные подходы к построению таких  методов были предложены в 
\cite{n_andrianova_a_a_1,
n_vasillev_f_p_1,
n_golikov_a_i_1,
n_demyanov_v_f_2,
n_eremin_i_i_2,
n_karmanov_v_g_1,
n_konnov_i_v_2,
n_nesterov_yu_e_3,
n_polayk_b_t_1,
n_strekalovskii_a_s_1,
n_pshenichnii_b_n_2},
 и этот перечень работ можно существенно расширить. Среди упомянутых $-$ работы, посвященные построению методов минимизации выпуклых и невыпуклых, дифференцируемых и недифференцируемых функций. Сложилась и определенная классификация этих методов. Например, в отдельные классы выделились  методы штрафных и барьерных функций, методы возможных направлений, методы типа линеаризации, квазиньютоновские методы,  субградиентные и $\varepsilon$-субградиентные методы,  методы типа приведенных градиентов, методы центров и другие.

Известный и довольно широкий класс образуют методы отсечений, которые используются как для  задач целочисленного программирования, так и для задач непрерывной оптимизации. Построению методов отсечений, предназначенных для решения задач целочисленного программирования, посвящены, например, работы
\cite{n_kivistick_k_1,
n_kolokolov_a_a_3,
n_kolokolov_a_a_2,
n_saiko_i_v_1,
n_hamisov_o_v_3,
n_chervak_yu_yu_2,
n_shevchenko_v_n_1,
n_balas_e_1,
n_gomory_r_e_1}.  Пожалуй, большее число работ посвящено методам отсечений для задач нелинейного программирования. Разные подходы к построению таких методов можно найти в
\cite{n_anciferov_e_g_1,
n_anciferov_e_g_2,
n_anciferov_e_g_3,
n_apekina_e_v_1,
n_belih_t_i_2,
n_belih_t_i_1,
n_belih_t_i_3,
n_bulatov_v_p_7,
n_bulatov_v_p_5,
n_bulatov_v_p_4,
n_bulatov_v_p_6,
n_levin_yu_ya_1,
n_levitin_e_s_1,
n_veinot_a_f_1,
n_zoutendijk_g_1}.

Значительное число методов отсечений для задач нелинейного программирования основывается на  следующей идеи.  \mbox{При построении} итерационных точек в этих методах  используется последовательное погружение либо множества ограничений исходной задачи, либо надграфика целевой функции в некоторые множества более простой структуры. Каждое из аппроксимирующих множеств строится на основе предыдущего путем отсечения от него, как правило плоскостями, некоторого подмножества, содержащего текущую итерационную точку.

Одна из основных проблем методов отсечений такого типа при их практическом применении заключается в том, что \mbox{от шага} к шагу, как правило, неограниченно растет число отсекающих плоскостей, формирующих аппроксимирующие множества. Поэтому с увеличением числа шагов возрастает и трудоемкость решения задач нахождения итерационных точек. Указанный недостаток названных методов оказывается существенным  \mbox{при решении}   задач с высокой точностью даже небольшой размерности, поскольку  возникают  сложности с переполнением текущей информацией памяти  компьютера.  Таким образом, исследования, связанные с модификацией известных методов  отсечений \mbox{и разработкой} новых методов указанного класса, свободных от этого недостатка, представляются вполне обоснованными.  В работе предлагаются принципы обновления погружающих множеств в алгоритмах отсечений с целью периодического освобождения от накапливающихся дополнительных ограничений, \mbox{и на этих} принципах строятся новые методы указанного класса с аппроксимацией как допустимой области, так и надграфика целевой функции.

Ранее в некоторых  работах (см., напр., 
\cite{n_bulatov_v_p_4,
n_topkis_d_m_1}) было уделено определенное внимание  решению проблемы отбрасывания накапливающихся отсечений в методах  с аппроксимацией допустимой области. Так, в методе \cite{n_bulatov_v_p_4}, где для построения аппроксимирующих множеств используются конусы обобщенно-опорных векторов, заложена возможность отбрасывания  отсечений. Однако, сходимость этого метода  обоснована только \mbox{для сильно} выпуклой целевой функции. 

В настоящей работе предлагается  такой принцип  обновления аппроксимирующих область ограничений множеств, который можно применять в методах отсечений для задач с практически любыми, а именно непрерывными, целевыми функциями. Этот принцип  обновления   основывается на  введенном в работе критерии качества аппроксимирующих множеств в окрестности текущих итерационных точек.  Данный принцип позволяет включать в методы такие процедуры обновления, при   которых можно периодически отбрасывать произвольное число любых построенных ранее отсечений.  А именно, допустимо как полное, так и частничное обновление погружающих множеств. В случае частичного обновления можно оставлять, например, только <<активные>> на данном шаге отсекающие плоскости или $n+1$ последнее отсечение.

Наряду с уже  упомянутыми существует значительное число методов отсечений,  которые используют аппроксимацию \mbox{не области} ограничений, а надграфика целевой функции (напр., 
\cite{n_bulatov_v_p_5,
n_bulatov_v_p_4,
n_nurminskii_e_a_3,
n_kelley_j_e_1,
n_kiwiel_k_c_3,
n_kiwiel_k_c_2,
n_kiwiel_k_c_1,
n_lemarechal_c_1,
n_lemarechal_c_2}). Этим методам присуща та же проблема накопления отсекающих плоскостей. В данной работе  построены   методы отсечений с аппроксимацией надграфика, в которых главной особенностью является заложенная в них возможность периодического обновления погружающих множеств. Эти  обновления тоже заключаются в отбрасывании отсекающих плоскостей, и основаны  на двух введенных принципах, отличных от вышеупомянутого принципа. На тех итерациях методов, где качество аппроксимации становится достаточно хорошим,  фиксируется итерационная точка и при желании происходит отбрасывание любых построенных ранее отсекающих плоскостей. 

Далее, обычно известные методы отсечений применяются на практике только для двух следующих типов задач выпуклого программирования. Если допустимая область исходной задачи  имеет  вид многогранника, то применимы методы \mbox{с  аппроксимацией} надграфика целевой функции многогранными множествами. Если же область ограничений --  не  многогранник, а целевая функция является линейной или выпуклой квадратичной, то естественно применять методы отсечений \mbox{с аппроксимацией} допустимой области. Поскольку в обоих случаях вспомогательные задачи построения итерационных точек \mbox{в методах} являются задачами линейного или квадратичного  программирования, то для соответствующих типов задач эти методы являются практически реализуемыми.  Однако, если целевая функция нелинейна и допустимая область также задается с помощью нелинейных функций, то  известные методы отсечений на практике трудноприменимы, так как задачи построения приближений в этих методах близки по трудоемкости к исходной задаче.

В работе предлагаются, кроме уже названных, такие методы отсечений, в которых  независимо от вида задачи выпуклого программирования  при построении итерационных точек могут решаться на каждом шаге вспомогательные задачи линейного программирования. Дело в том, что в этих методах происходит одновременная аппроксимация многогранными множествами как надграфика целевой функции, так и допустимой области. Кроме того, подчеркнем, что для одного из таких методов разработана процедура обновления  множеств, аппроксимирующих и надграфик,  и  допустимую область. Указанные обновления  происходят за счет отбрасывания секущих плоскостей  \mbox{на основе}  одновременного оценивания качества аппроксимации надграфика и качества аппроксимации допустимой области.

При решении оптимизационных задач на практике довольно часто используются так называемые смешанные алгоритмы минимизации (напр., \cite[119]{n_zangwill_y_i_1}). Эти алгоритмы получаются \mbox{на основе} двух (а иногда и более) известных методов следующим образом. В смешанном алгоритме первые несколько итерационных точек строятся одним из выбранных методов. \mbox{Далее определенное} число приближений находится другим методом. Затем снова для построения точек привлекается первый метод, и т. д. \mbox{При этом} для каждого смешанного алгоритма возникает вопрос о сходимости, несмотря на то, что сходимость образующих его методов обоснована.

В работе показано, как можно строить смешанные алгоритмы условной минимизации на основе предлагаемых  методов отсечений и любых релаксационных методов выпуклого программирования. На тех итерациях предлагаемых методов, где происходят обновления аппроксимирующих множеств, фиксируются с условием релаксации  и основные итерационные точки. Для построения этих точек допустимо использование любых, в том числе новых, релаксационных методов минимизации. За счет этого и обеспечивается указанная выше возможность смешивания предлагаемых методов отсечений с другими методами.   При этом сходимость всех таких смешанных алгоритмов будет гарантирована сходимостью разработанных здесь методов отсечений, причем даже в том случае, когда упомянутые релаксационные методы, включенные в смешанные алгоритмы, являются эвристическими.  Заметим также, что заложенные в методах отсечений процедуры обновления погружающих множеств будут перенесены таким образом и в смешанные алгоритмы.

Далее, в последнее время появилась возможность в объединении группы компьютеров в единую мощную вычислительную систему  с целью  использования параллельных вычислений, которые  значительно сокращают время решения  задач. Возможность использования параллельных вычислений в некоторых методах  математического программирования обоснована, например, в \cite{n_garanga_v_a_1,
n_evtushenko_yu_g_4,
n_popov_l_d_1}. В данной работе показано, что параллельные вычисления можно применять и в методах отсечений. Большинство предлагаемых здесь методов  позволяют использовать параллельные вычисления  при построении как вспомогательных, так и основных итерационных точек. Например, распараллеливание вычислений при нахождении  основных итерационных точек   в этих методах может  происходить следующим образом. Различными алгоритмами параллельно строится некоторая совокупность вспомогательных точек, а затем из них в качестве итерационной точки выбирается рекордная.

Номера утверждений, замечаний, формул в книге задаются тремя числами, и для ссылок на них в тексте используется та же тройная нумерация. Первое из тройки чисел означает главу, в которой изложен материал, второе число -- параграф в этой главе, а третье -- номер, соответственно, утверждения, замечания или формулы в этом параграфе.

Опишем коротко содержание каждой из трех глав книги. \mbox{В первой} главе разрабатываются методы отсечений для решения задач нелинейного программирования и задачи отыскания проекции точки на выпуклое множество, основанные на аппроксимации допустимой области. Построение последовательности итерационных точек в предлагаемых методах происходит  с использованием операции частичного погружения допустимой области в аппроксимирующие ее многогранные множества. \mbox{В предложеных} методах из  $\S\S\ $1.1 $-$ 1.3  не требуется вложения каждого из аппроксимирующих множеств в предыдущее. \mbox{Такая особенность} методов позволяет  периодическое отбрасывание любых полученных в процессе решения дополнительных ограничений. В $\S$ 1.4 предлагается метод отсечений, допускающий параллельные вычисления. Обосновывается сходимость разработанных методов, обсуждаются  их реализации, предлагаются  процедуры обновлений аппроксимирующих множеств, приводятся оценки точности решения и скорости сходимости.

В $\S \ 1.1$ предлагается общий метод  решения задачи
\begin{equation}\label{e0.0.1}
\min \{f(x) : x\in D\},
\end{equation}
где $f(x)$ -- непрерывная в $n$-мерном евклидовом пространстве ${\rm R_n}$ функция, а множество $D\subset {\rm  R_n}$ выпукло и замкнуто. 

Предлагаемый метод отсечений вырабатывает вспомогательную последовательность приближений $\{y_i\}$, $i\in K=\{0,1,\ldots\}$, и фиксирует основную итерационную последовательность $\{x_k\}$, $k\in K$, в виде подпоследовательности $\{y_{i_k}\}$ последовательности $\{y_i\}$. На $i$-ом шаге $(i\ge 0)$ для построения точки $y_i$  приближенно решается  вспомогательная задача минимизации функции $f(x)$ на аппроксимирующем область ограничений  множестве $M_i$. Если точка $y_i$ в определенном смысле оказалась близка к множеству $D$, то фиксируются номер $i=i_k$ и точка $x_k=y_{i_k}$, и происходит обновление аппроксимирующего множества $M_i=M_{i_k}$. В противном случае считается, что качество аппроксимации области $D$ множеством $M_i$ является  неудовлетворительным, и на основе $M_i$ строится очередное аппроксимирующее множество $M_{i+1}$ путем отсечения от $M_i$ точки $y_i$  с использованием обобщенно-опорных элементов. Далее в множестве $M_{i+1}$ указанным выше способом  находится очередное приближение $y_{i+1}$ вспомогательной последовательности.

В том же параграфе приведено сравнение предложенного метода с идейно близким методом опорных множеств \cite{n_bulatov_v_p_4}, построенным для задачи (\ref{e0.0.1}) c выпуклой целевой функцией $f(x)$. Принципиальное отличие предложенного метода от метода \cite{n_bulatov_v_p_4} состоит в том, что в нем заложен критерий оценки качества аппроксимирующего множества, позволяющий вводить процедуры обновления погружающих множеств за счет периодического  отбрасывания любого числа любых отсекающих плоскостей.

Далее, в $\S$ 1.2 предлагается  метод отсечений для решения задачи (\ref{e0.0.1}), в которой область ограничений имеет более общий вид:
\begin{equation}\label{e0.0.2}
D=D^\prime \bigcap D^{\prime\prime},
\end{equation}
где $D^\prime$, $D^{\prime\prime}$ -- выпуклые замкнутые в ${\rm R_n}$ множества, причем второе  из них может заведомо иметь пустую внутренность. 

Метод отличается от разработанного в $\S \ 1.1$ метода, в частности, способами построения аппроксимирующих множеств. \mbox{А именно}, погружающее множество $M_{i+1}$ строится путем отсечений от $M_i$  точки $y_i$ плоскостями, построенными с помощью субградиентов функции, определяющей множество $D^\prime$.  Кроме того, в этом методе по сравнению с предыдущим предусмотрены более широкие возможности для выбора итерационных точек основной последовательности. За счет этих возможностей допустимо построение на основе этого метода отсечений некоторых смешанных алгоритмов. Сходимость таких смешанных алгоритмов будет гарантирована сходимостью предложенного метода отсечений. Метод из $\S \ 1.2$   характерен тем, что также предусматривает возможность  периодического обновления аппроксимирующих множеств на тех итерациях, где качество аппроксимирующего множества становится удовлетворительным.

В  $\S$ 1.3 главы 1 предлагается метод отсечений для нахождения проекции точки $y\in {\rm R_n}$ на выпуклое множество $D$, заданное в виде (\ref{e0.0.2}). На каждой итерации метода определенным образом строится многогранное множество, аппроксимирующее область $D^\prime$, и приближение получается путем проектирования точки $y$ на пересечение этого аппроксимирующего множества с множеством $D^{\prime\prime}$. В случае, если $D^{\prime\prime}$ задано системой линейных уравнений или неравенств, то построение приближения не вызывает большого труда, так как задача проектирования решается известными алгоритмами за конечное число шагов. Как и в методах из предыдущих параграфов, в этом проекционном методе тоже заложена возможность обновления аппроксимирующих область $D^\prime$ множеств за счет отбрасывания дополнительных ограничений.

В $\S \ 1.4$ для задачи (\ref{e0.0.1}) строится еще один метод отсечений  с аппроксимацией области ограничений, который отличается от методов из $\S\S\ 1.1, 1.2$ тем, что допускает возможность распараллеливания  процесса отыскания итерационных точек. В этом методе для построения основных итерационных точек сначала некоторым способом строится определенное число многогранных множеств, аппроксимирующих область ограничений исходной задачи. Далее в каждом из этих множеств находится вспомогательная точка приближенного минимума целевой функции. Наконец, в качестве основной итерационной точки выбирается рекордная из найденных вспомогательных точек. Заметим, что эти вспомогательные точки могут находится параллельно различными алгоритмами минимизации.


Во второй  главе предлагаются методы отсечений для решения задачи выпуклого программирования, которые используют операцию погружения не области ограничений, а  надграфика целевой функции исходной задачи. Как и методы главы 1, они характерны тем, что допускают  периодические обновления погружающих множеств за счет отбрасывания отсекающих плоскостей.  Все методы данной главы позволяют  на каждой итерации оценивать близость значения целевой функции в текущей итерационной точке к оптимальному значению. В главе предложены два критерия оценки качества аппроксимирующих множеств в окрестности  итерационных точек. Критерии позволяют включать в предлагаемые методы процедуры обновления погружающих множеств. Первый из этих критериев основан на упомянутых выше оценках близости текущих итерационых значений к оптимальному значению, а второй критерий опирается на оценку близости итерационных точек к надграфику целевой функции.  Доказывается сходимость предлагаемых в данной главе  методов, обсуждаются  их реализации, предлагаются  приемы периодических отбрасываний произвольного числа любых полученных в процессе решения отсекающих плоскостей, приводятся оценки точности решения и  скорости сходимости.

Первый параграф  второй главы посвящен решению  задачи  минимизации выпуклой в ${\rm R_n}$ дискретной функции максимума $f(x)$ на выпуклом замкнутом множестве $D\subset {\rm R_n}$. Предлагаемый для решения задачи метод заключается в следующем.  \mbox{На $i$-ом} шаге $(i\ge 0)$ находится решение $(y_i, \gamma_i)$  следующей задачи:
$$\gamma\to \min,$$
$$(x, \gamma)\in M_i,\quad  x\in G_i,\quad \gamma \ge \bar{\gamma}_i,$$
где $M_i$ $-$ множество, аппроксимирующее надграфик целевой функции, $G_i$ $-$ некоторое подмножество области $D$, 
$$\bar{\gamma}_i \le f^* = \min \{f(x) : x\in D\}.$$ 
Если разность  $f(y_i) - \gamma_i$ достаточно мала, то качество аппроксимации надграфика множеством $M_i$  считается удовлетворительным. В таком случае $y_i$ фиксируется в качестве точки основной  последовательности приближений, и происходит обновление аппроксимирующего множества $M_i$. В противном случае строится очередное аппроксимирующее множество $M_{i+1}$ путем отсечения от $M_i$ точки $(y_i, \gamma_i)$. 

Заметим, что метод обладает большими возможностями в выборе начального аппроксимирующего множества $M_0$, множеств $G_i$ и чисел $\bar{\gamma}_i$. За счет вида целевой функции задачи  у метода есть значительные возможности  и для построения отсекающих плоскостей,  которые формируют аппроксимирующие множества.

В $\S$ 2.2  предлагается  метод отсечений для решения задачи (\ref{e0.0.1}) с выпуклой функцией $f(x)$, который также использует аппроксимацию надграфика целевой функции. Главные   отличия  этого метода от метода из $\S$  $2.1$   заключаются  в следующем. \mbox{Во-первых},  в нем заложена возможность построения на его основе  смешанных алгоритмов с привлечением любых релаксационных сходящихся или эвристических методов условной минимизации без дополнительного обоснования сходимости таких смешанных алгоритмов. Во-вторых, метод позволяет привлекать параллельные вычисления при построении приближений. Кроме того, отметим, что в данном методе применяется  способ построения отсекающих плоскостей, отличный от используемого в предыдущем методе. 

В $\S$ $2.3$ предлагается модификация известного метода уровней (\cite{n_nesterov_yu_e_3,n_lemarechal_c_1}) для решения задачи (\ref{e0.0.1}), где $f(x)$ $-$ выпуклая  функция, а $D$ $-$ выпуклое ограниченное замкнутое множество.  Главное отличие предлагаемого метода от известного заключается в том, что в нем заложена возможность периодического отбрасывания накапливающихся в процессе счета отсекающих плоскостей, что делает метод привлекательным с практической точки зрения.   Другое немаловажное отличие этого метода от (\cite{n_nesterov_yu_e_3,n_lemarechal_c_1}) заключается в более общем способе задания  итерационных точек.  В частности, в предлагаемом методе при построении приближений можно привлекать параллельные вычисления.

Далее в $\S$ $2.4$, основываясь на способе построения аппроксимирующих  множеств из $\S \ 2.2$,  разрабатывается еще один метод отсечений для задачи (\ref{e0.0.1}). В методе при построении итерационных точек происходит аппроксимация множествами $M_i$  надграфика не целевой функции $f(x)$, а  некоторой вспомогательной функции  $g(x)$.
Существенным отличием предлагаемого метода от методов, построенных в предыдущих параграфах, является заложенный в нем новый критерий  качества аппроксимирующих  множеств $M_i$, а значит, и критерий отбрасывания отсекающих плоскостей. 
 А именно, как только  итерационная точка из $M_i$ становится достаточно близка к  надграфику   функции $g(x)$, качество аппроксимации  считается приемлимым, и происходит обновление множества $M_i$. 

Данный метод за конечное число шагов  позволяет найти $\varepsilon$-решение задачи (\ref{e0.0.1}). Если же после нахождения $\varepsilon$-решения процесс построения приближений продолжается, то любая предельная точка итерационной последовательности будет принадлежать множеству решений. 
Отметим, что данный метод, как и многие предыдущие,  позволяет строить смешанные алгоритмы с привлечением других известных методов условной минимизации. Кроме того, при построении приближений в методе можно  использовать и  параллельные вычисления.


В главах 1, 2 разработаны методы отсечений, в которых при построении итерационных точек,  используется аппроксимация либо допустимой области,  либо надграфика  функции. В третьей главе строятся методы отсечений,  в которых при нахождении приближений применяется аппроксимация как  допустимой области, так и надграфика целевой функции. В случае  аппроксимации указанных множеств многогранными множествами предлагаемые  методы позволяют независимо от вида функций ограничений и целевой функции исходной задачи решать на каждом шаге при построении приближения вспомогательную задачу линейного программирования. Обсуждаются реализации, обосновывается  сходимость разработанных методов.

В $\S$ 3.1 предлагается общий  метод отсечений для  задачи (\ref{e0.0.1}), где $f(x)$ $-$ выпуклая  функция. Метод строит последовательность приближений   следующим образом.
На $i$-ом шаге метода находится точка $(y_i, \gamma_i)$ как решение задачи
$$\gamma \to \min,$$
$$(x, \gamma)\in M_i, \quad x\in G_i, \quad \gamma \ge \alpha,$$
где $M_i\subset {\rm R_{n+1}}$ -- множество, аппроксимирующее надграфик целевой функции, $G_i\subset {\rm R_n}$ -- множество, аппроксимирующее область ограничений в окрестности решения исходной задачи, $\alpha$ -- нижняя оценка значения $f^*$. Компонента $y_i$, найденной точки, принимается за текущую итерационную точку. Далее строятся очередные аппроксимирующие  множества $M_{i+1}$ и $G_{i+1}$. Первое из них $-$ путем отсечения    от  множества $M_i$ точки $(y_i, \gamma_i)$, а второе $-$ на основе отсечения от $G_i$ приближения $y_i$ в случае $y_i\notin D$. Предлагаются различные варианты  выбора начальных погружающих множеств, а также различные способы построения отсекающих плоскостей.


Далее, в $\S$ 3.2 для минимизации выпуклой функции $f(x)$ на выпуклом множестве $D$  предлагается еще один метод отсечений, в котором происходит  аппроксимация многогранными множествами как области ограничений, так и надграфика целевой функции. Этот метод характерен тем, что построение каждой точки основной последовательности происходит в два этапа. На первом этапе фиксируется множество, аппроксимирующее область ограничений, и на основе некоторых вспомогательных точек последовательно строятся множества, аппроксимирующие надграфик. Когда качество аппроксимации надграфика становится в определенном смысле
приемлемым, первый этап завершается. На втором этапе находится основная итерационная точка, и путем ее отсечения строится очередное множество, аппроксимирующее область ограничений. Отметим, что после нахождения основной итерационной точки предусмотрена возможность обновления аппроксимирующего надграфик множества за счет отбрасывания любого числа ранее построенных секущих плоскостей.

На основе идеи одновременной аппроксимации надграфика целевой функции  и области ограничений в $\S \ 3.3$ предлагается метод отсечений для задачи выпуклого программирования (\ref{e0.0.1}) с областью $D$,  заданной в виде (\ref{e0.0.2}).  Принципиальная особенность предлагаемого метода, выделяющая его среди всех методов   работы, заключается в следующем. В метод заложены два критерия оценки качества аппроксимирующих множеств. Первый из них основан на оценке близости итерационной точки к множеству $D^\prime$ и позволяет отслеживать качество аппроксимации допустимой области. Второй из указанных критериев основан на оценке близости текущих значений целевой  функции $f(x)$ к ее оптимальным значениям  на  множествах,  аппроксимирующих область $D^\prime$, и характерезует качество аппроксимации надграфика. Обновление аппроксимирующих  множеств происходит на тех итерациях метода, когда одновременно качество аппроксимации    допустимой области  и   качество аппроксимации надграфика  являются удовлетворительными. 

В  параграфах $1.5$, $2.5$ и $3.4$ представлены результаты численных экспериментов, которые проводились с целью исследования методов, разработанных в $\S\S$  1.1, 1.2, 2.1, 2.2, 3.3.
Указанные методы программно реализованы на языке программирования C++ в среде Microsoft Visual Studio 2008. Методы  из $\S\S$  1.1, 1.2,  3.3 тестировались на задачах, в  которых допустимая область задавалась нелинейными функциями. Остальные из указанных методов исследовались на примерах, в  которых область ограничений являлась многогранником, а целевая функция выбиралась нелинейной.
Исследуемыми методами решено значительно число тестовых примеров с различным количеством переменных (до 50) и ограничений (до 50). Каждая  задача решалась  тестируемым методом как с  использованием процедур отбрасывания отсечений, так и без привлечения процедур обновления аппроксимирующих множеств. Рассчеты показали, что применение  в каждом из исследованных методов определенных процедур  обновления  способствовало значительному ускорению решению всех тестовых задач. Кроме того, результаты численных экспериментов позволили выработать рекомендации по заданию некоторых  параметров  методов.


\chapter[Методы отсечений с использованием процедур аппроксимации допустимой области]
{
\quad{} \quad{}\hspace{1mm}Методы отсечений
\newline с использованием процедур аппроксимации допустимой области }


Глава посвящена методам отсечений для решения задач нелинейного программирования и задачи отыскания проекции точки на выпуклое множество, основанным на аппроксимации допустимой области. Последовательности приближений в предлагаемых методах строятся с использованием операции частичного погружения допустимой области в аппроксимирующие ее многогранные множества. \mbox{В разработанных} методах, которые описаны в $\S\S$ 1.1 -- 1.3, не требуется вложения каждого из аппроксимирующих множеств в предыдущее. Такая особенность методов дает возможность периодического отбрасывания любых полученных в процессе решения дополнительных ограничений. \mbox{В $\S$ 1.4} предлагается метод отсечений, допускающий параллельные вычисления. Обосновываются сходимость методов, обсуждаются их реализации, приводятся оценки точности решения и скорости сходимости. Некоторые из результатов этой главы можно найти в \cite{
n_yarullin_r_s_3,
n_yarullin_r_s_6,
n_yarullin_r_s_12,
n_yarullin_r_s_17,
n_yarullin_r_s_10,
n_yarullin_r_s_13}.

\section
{Метод отсечений на основе аппроксимации допустимой области семейством опорных плоскостей без вложения погружающих множеств}

Предлагается метод для решения задачи нелинейного программирования с обновлением аппроксимирующих множеств. 

Пусть $f_j(x)$, $j\in J=\{1,\ldots,m\}$, -- выпуклые в $n$-мерном евклидовом пространстве ${\rm R_n}$ функции,
\begin{equation}\label{e1.1.0}
D=\{x\in {\rm R_n} : f_j(x)\le0, \ j\in J\},
\end{equation}
$f(x)$ -- непрерывная достигающая на $D$ минимального значения функция, и для всех $j\in J$ множества 
\begin{equation}\label{e1.1.00}
D_j=\{x\in {\rm R_n} : f_j(x)\le 0\}
\end{equation}
имеют непустую внутренность ${\rm int \,} D_j$. Решается задача 
\begin{equation}\label{e1.1.1}
\min \{f(x) : x\in D\}.
\end{equation}

Положим $F(x)=\max \limits_{j\in J} f_j(x)$, $D_\varepsilon=\{x\in {\rm R_n}  : F(x)\le \varepsilon\}$, где $\varepsilon\ge 0$, $f^*=\min \{f(x) : x\in D\}$, $X^*=\{x\in D : f(x) =f^*\}$, $E^*=\{x\in {\rm R_n} : f(x) \le f^*\}$, $W^1(x, D_j)=\{a\in {\rm R_n} : \langle a, z-x\rangle \le 0 \ \forall z\in D_j, \ \|a\|=1\}$ -- множество нормированных обобщенно-опорных векторов для множества $D_j$ в точке $x\in {\rm R_n}$, $K=\{0,1,\ldots\}$.

Сразу отметим, что в (\ref{e1.1.1}) внутренность множества $D$ может быть пустой.

Предлагаемый метод решения задачи (\ref{e1.1.1}) вырабатывает последовательности приближений $\{y_i\}$, $i\in K$, $\{x_k\}$, $k\in K$, и заключается в следующем.

\textbf{Метод 1.1.} Строится выпуклое замкнутое множество $M_0\subset {\rm  R_n}$, содержащее хотя бы одну точку множества $X^*$, например, точку $x^*$. Выбираются точки $v^j\in {\rm int\,} D_j$ для всех $j\in J$. Задается число $\varepsilon_0\ge 0$. Полагается $k=0$, $i=0$.

\textbf{1.} Находится точка
\begin{equation}\label{e1.1.2}
y_i \in M_i \bigcap E^*.
\end{equation}

\textbf{2.} Формируется множество 
$$J_i=\{j\in J : y_i \not \in D_j\}.$$
Если $J_i=\emptyset$, то $y_i\in X^*$, и процесс завершается.

\textbf{3}. Если $y_i\notin D_{\varepsilon_k}$, то выбирается выпуклое замкнутое множество $G_i\subset {\rm R_n}$, содержащее точку $x^*$, полагается 
\begin{equation}\label{e1.1.3}
Q_i=M_i\bigcap G_i,
\end{equation}
и следует переход к п. 4. В противном случае полагается $i_k=i$,
\begin{equation}\label{e1.1.4}
x_k=y_{i_k},
\end{equation}
и выбирается выпуклое замкнутое множество $Q_i=Q_{i_k}$ такое, что 
\begin{equation}\label{e1.1.5}
x^*\in Q_i.
\end{equation}
Задается число $\varepsilon_{k+1} \ge 0$, значение $k$ увеличивается на единицу, и следует переход к очередному пункту.

\textbf{4.} Для каждого $j\in J_i$ в интервале $(v^j, y_i)$ выбирается точка $z^j_i$ так, чтобы $z_i^j\notin {\rm int\,} D_j$ и при некотором 
$q_i^j\in [1, q]$, $q<+\infty$, для точки 
\begin{equation}\label{e1.1.6}
\bar{y}_i^j=y_i+q_i^j(z_i^j-y_i) 
\end{equation}
выполнялось включение $\bar{y}_i^j\in D_j$. Для всех  $j\in J\setminus J_i$ полагается $z_i^j=\bar{y}_i^j=y_i$.

\textbf{5.} Выбирается множество $H_i \subset J_i$ так, чтобы выполнялось включение $j_i \in H_i$, где номер $j_i$ удовлетворяет условию 
\begin{equation}\label{e1.1.7}
\|y_i-z_i^{j_i}\|=\max \limits_{j\in J_i}\|y_i-z_i^j\|.
\end{equation}

\textbf{6.} Для каждого $j\in H_i$ выбирается конечное множество $A_i^j \subset W^1(z_i^j, D_j)$,  полагается 
\begin{equation}\label{e1.1.8}
M_{i+1}=Q_i\bigcap \limits_{j\in H_i}\{x\in {\rm R_n} : \langle a, x- z_i^j\rangle \le 0 \ \forall a \in A^j_i\}, 
\end{equation}
и следует переход к п. 1 при $i$, увеличенном на единицу.

Сделаем некоторые  замечания, касающиеся данного метода.

\begin{remark}\label{r1.1.1}
Если на шаге 2 метода выполняется условие $J_i=\emptyset$, то согласно (\ref{e1.1.2}) точка $y_i$, действительно, является решением исходной задачи.
\end{remark}

\begin{remark}\label{r1.1.2}
Множество $M_i\bigcap E^*$ содержит, по крайней мере, точку $x^*$, а значит, выбор $y_i$ из условия (\ref{e1.1.2}) возможен. В частности, точку $y_i$ можно находить согласно равенству
\begin{equation}\label{e1.1.9}
f(y_i) =\min \{ f(x) : x\in M_i\}.
\end{equation}
\end{remark}

\begin{remark}\label{r1.1.3}
В случае, когда $f(x)$ выпукла, $J=\{1\}$, $v^1=v\in {\rm int\,} D$, причем ${\rm int\,} D \ne \emptyset$, $\varepsilon_0=0$ и точки $y_i$ выбраны согласно (\ref{e1.1.9}), то предлагаемый метод идейно близок к известному методу погружений В.\,П.~Булатова (напр., \cite{n_bulatov_v_p_4}) для задачи выпуклого программирования.
\end{remark}

\begin{remark}\label{r1.1.4}
В  упомянутом методе \cite{n_bulatov_v_p_4} отсечения строятся с помощью точек пересечения отрезка $[v, y_i]$ с границей множества $D$.  В предложенном здесь методе точки $z_i^j$, $j\in J_i$, используемые для построения отсечений, можно выбирать отличными от точек пересечения отрезков $[v^j, y_i]$ с границами множеств $D_j$, что существенно с точки зрения решения практических задач.
\end{remark}

\begin{remark}\label{r1.1.5}
Если в п. 5 положить $H_i=J_i$, то отпадает потребность в фиксировании согласно (\ref{e1.1.7}) номера $j_i$ "самого глубокого" \  на $i$-ом шаге отсечения.
\end{remark}

\begin{remark}\label{r1.4}
Для выбора начального аппроксимирующего множества $M_0$ имеется много возможностей. Если, например, положить
\begin{equation}\label{e1.1.10}
M_0=\bigcap \limits_{j\in J^\prime}D_j, 
\end{equation}
где $J^\prime\subset J$, то нет необходимости в задании точек $v^j\in{\rm int\,}D_j$, $j\in J^\prime$, поскольку тогда для всех $i\in K$ и $j\in J^\prime$ выполняются включения $y_i\in D_j$, и точки $v^j$, $j\in J^\prime$, в построении отсекающих плоскостей не участвуют. Отметим также, что выбор $M_0$ в виде (\ref{e1.1.10}) удобен в том случае, когда  $\bigcap \limits_{j\in J^\prime}D_j$ -- выпуклый многогранник.

Если $M_0={\rm R_n}$, а функция $f(x)$ достигает в ${\rm R_n}$ своего абсолютного минимального значения, то можно положить $y_0 = \underset{x\in {\rm R_n}}{\rm{arg \ min}} f(x)$. 

При выборе $M_0=D$ метод завершает процесс отыскания решения задачи (\ref{e1.1.1}) на нулевом шаге. 

Отметим, что если $\rm{int \,} D \ne \emptyset$ и известна точка $y\in \rm{int \,}D$, то в  методе удобно положить $v^j = v$ для всех $j\in J$.
\end{remark}

Приведем теперь принципиальное замечание для метода, касающееся возможности периодического обновления погружающих множеств $M_i$ за счет отбрасывания любого  числа построенных к $i$-ому шагу отсекающих плоскостей.

\begin{remark}\label{r1.1.7}
Пусть в (\ref{e1.1.3}) $G_i={\rm R_n}$ для всех $i\in K$, и для точки $y_i$ выполняется включение 
\begin{equation}\label{e1.1.11}
y_i\in D_{\varepsilon_k}.
\end{equation}
Положим 
\begin{equation}\label{e1.1.12}
Q_i=M_{r_i},
\end{equation}
где $0\le r_i\le i=i_k$. Ясно, что при всех $r_i=0,\ldots,i$ условие (\ref{e1.1.5}) выполняется. Согласно (\ref{e1.1.12}) в качестве $Q_i$ можно выбирать любое из построенных к $i$-ому шагу погружающих множеств $M_0,\ldots,M_i$. В частности, для всех $i\in K$, для которых выполняется (\ref{e1.1.11}), можно положить
$$Q_i=M_0.$$
Тем самым при каждом $i=i_k$, $k\in K$, произойдет отбрасывание всех накопившихся к шагу $i_k$ отсекающих плоскостей.
\end{remark}

Понятно, что множества $Q_i$ при выполнении условия (\ref{e1.1.11}) можно выбирать и не следуя этому замечанию. А именно, множества $Q_i$ удобно задавать следующим образом:
$$Q_i=M_i\bigcap \limits_{j\in H_{i-1}} \{x\in {\rm R_n} : \langle a, x-z^j_{i-1} \rangle \le 0 \ \forall a\in A^j_{i-1}\}, \quad i\ge 1.$$
Кроме того, отметим, что согласно (\ref{e1.1.5}) множества $Q_i$ при условии (\ref{e1.1.11}), как и множество $M_0$, можно выбирать совпадающими с ${\rm R_n}$.

Далее, отметим, что при всех $i\in K$, независимо от выполнения условия (\ref{e1.1.11}), множества $Q_i$ допустимо задавать в виде (\ref{e1.1.3}), поскольку включения (\ref{e1.1.5}) для них выполняются. В таком случае ввиду (\ref{e1.1.8})
$$M_{i+1}=M_i\bigcap G_i\bigcap \limits_{j\in H_i} \{x\in {\rm R_n} : \langle a, x- z_i^j \rangle \le 0 \quad \forall a\in A_i^j\} \ \forall i\in K,$$
то есть от шага к шагу происходит накопление отсекающих плоскостей и никаких обновлений погружающих множеств не происходит. Свойства последовательности $\{y_i\}$, построенной с условием (\ref{e1.1.3}) выбора множеств $Q_i$ при всех $i\in K$, будут обсуждаться ниже.

Перейдем к исследованию сходимости предложенного метода. Сначала изучим некоторые свойства построенной согласно методу последовательности $\{y_i\}$, $i\in K$. Везде далее в этом параграфе будем предполагать, что она ограничена.

\begin{lemma}\label{lemma1.1.1}
Пусть $U\subset {\rm R_n}$ -- выпуклое множество, $L$ -- его несущее подпространство, а множество $Q$ из аффинной оболочки $U$ ограничено и не содержится в ${\rm ri\,} U$ -- относительной внутренности множества $U$. Если точка $u\in {\rm R_n}$ такова, что выполняется включение $u\in {\rm ri \,}U$, то найдется такое число $\delta > 0$, что для всех $z\in Q\setminus {\rm ri\,}U$ и всех $a\in L\bigcap W^1(z,U)$ справедливо неравенство $\langle a, u-z\rangle\le -\delta$.
\end{lemma}

\begin{dokaz} утверждения приведено в \cite{n_zabotin_i_ya_50}.
\end{dokaz}

C учетом сделанного выше замечания о возможности выбора множеств $Q_i$ в виде (\ref{e1.1.3}), независимо от принадлежности точек $y_i$ множествам $D_{\varepsilon_k}$, сформулируем следующее вспомогательное утверждение.

\begin{lemma}\label{lemma1.1.2}
 Пусть последовательность $\{y_i\}$, $i\in K$, построена предложенным методом с условием, что для всех $i\in K$, начиная с некоторого номера $i^\prime \ge 0$, множества $Q_i$ выбраны согласно (\ref{e1.1.3}). Тогда любая предельная точка последовательности $\{y_i\}$, $ i\in K$, $i\ge i^\prime$, принадлежит множеству $D$.
\end{lemma}

\begin{proof}
Пусть $\{y_i\}, \ i\in K^\prime \subset K,$ -- любая сходящаяся подпоследовательность, выделенная из последовательности  точек $y_i$, $ i\in K$, $i\ge i^\prime$, и $\bar{y}$ -- ее предельная точка. Покажем справедливость включения
\begin{equation}\label{e1.1.13}
\bar{y}\in D.
\end{equation}

Пусть индекс $l\in J$ такой, что номер $j_i$, удовлетворяющий условию (\ref{e1.1.7}), совпадает с $l$ для бесконечного числа номеров $i\in K^\prime$. Положим 
$$K_l=\{i\in K^\prime : j_i=l\} $$
 и  докажем сначала равенство 
\begin{equation}\label{e1.1.14}
\lim \limits_{i\in K_l}\|z_i^l-y_i\|=0,
\end{equation} учитывая, что вместе с последовательностью $\{y_i\}, \ i\in K_l$, для каждого $j\in J$ построены и последовательности  $\{z_i^j\}, \ \{\bar{y}_i^j\}, \ i\in K_l.$

Заметим, что для всех $i\in K$ и $j\in J$ 
\begin{equation}\label{e1.1.15}
z_i^j=y_i+\gamma_i^j(v^j-y_i),
\end{equation} где $\gamma_i^j\in [0,1)$, причем $\gamma_i^l>0$ для всех $ i\in K_l$. Для произвольного $i\in K_l$ зафиксируем номер $p_i\in K_l$ такой, что $p_i>i$. \mbox{В силу} (\ref{e1.1.3}), (\ref{e1.1.8}) выполняется включение $M_{p_i}\subset M_i$. \mbox{Кроме того}, ввиду (\ref{e1.1.2}) $y_{p_i}\in M_{p_i}$, а любой элемент множества $A_i^l$  является обобщенно-опорным  и для множества $M_{p_i}$ в точке $z_i^l$. Следовательно,
 $$\langle a, y_{p_i}-z_i^l \rangle \le 0 \quad \forall a\in A_i^l.$$
 Отсюда с учетом (\ref{e1.1.15}) при $j=l$ для всех $a\in A_i^l$ имеем
$$\langle a, y_i - y_{p_i} \rangle \ge \gamma_i^l \langle a, y_i - v^l \rangle .$$
По лемме \ref{lemma1.1.1} найдется такое число $\delta_l>0$, что $\langle a, y_i-v^l\rangle \ge \delta_l$ для всех $i\in K_l, \ a\in A_i^l.$ Значит, $\langle a, y_i-y_{p_i}\rangle \ge \gamma_i^l\delta_l$ для всех $a\in A_i^l$, а поскольку $\|a\|=1$ для всех $a\in A_i^l$,  то 
\begin{equation}\label{e1.1.16}
\|y_i-y_{p_i}\|\ge \gamma_i^l\delta_l \quad \forall i,p_i\in K_l, \quad p_i>i.
\end{equation} Так как последовательность $\{y_i\}$, $ i\in K_l$, является сходящейся, то согласно (\ref{e1.1.16}) $\gamma_i^l \rightarrow 0, \  i \rightarrow \infty, \ i\in K_l$. Поэтому из (\ref{e1.1.15}) при $j=l$ с учетом ограниченности последовательности $\{\|v^l-y_i\|\}$, $i\in K_l$, следует равенство (\ref{e1.1.14}).

Далее, в силу условия (\ref{e1.1.7}) для всех $i\in K_l$ справедливы неравенства $\|z_i^l-y_i\|\ge \|z_i^j-y_i\|, \ j\in J_i$. Кроме того, согласно п. 4 метода $\|z_i^j-y_i\|=0$ для всех $j\in J\setminus J_i, \ i\in K_l.$ Следовательно, при любых $i \in K_l, \ j\in J$ $$\|z_i^l-y_i\|\ge \|z_i^j-y_i\|.$$Тогда ввиду (\ref{e1.1.14}) 
\begin{equation}\label{e1.1.17}
\lim \limits_{i\in K_l} \|z_i^j-y_i\|=0 \quad \forall j\in J. 
\end{equation}

Согласно п. 4 алгоритма для каждого $i\in K_l$ и $j\in J$ точка $\bar{y}_i^j$ либо совпадает с $y_i$, либо имеет вид (\ref{e1.1.6}). \mbox{Так как} последовательность $\{y_i\}, \ i\in K_l$, ограничена, то отсюда с учетом (\ref{e1.1.15}) следует ограниченность последовательностей $\{\bar{y}_i^j\}, \ i\in K_l$, для всех $j\in J$. Выделим теперь для каждого $j\in J$ из последовательности $\{\bar{y}_i^j\}$, $i\in K_l$, сходящуюся подпоследовательность  $\{\bar{y}_i^j\}$, $i\in K_l^j\subset K_l$, и пусть $u_j$ -- ее предельная точка. Заметим, что $u_j\in D_j$, $j\in J$,  в силу замкнутости множеств $D_j$. Положим для каждого $j\in J$ 
$$P_1^j=\{ i\in K_l^j : j\in J_i\}, \quad P_2^j=K_l^j\setminus P_1^j.$$ 
Хотя  бы одно из множеств $P_1^j$ или $P_2^j$ для каждого $j\in J$ состоит из бесконечного числа номеров. Перейдем теперь при каждом фиксированном $j\in J$ к пределу в равенствах (\ref{e1.1.6}) по $i\rightarrow \infty$, $i\in P_1^j$, c учетом (\ref{e1.1.17}), если множество $P_1^j$ бесконечно, или в равенствах $\bar{y}_i^j=y_i $ по $i\rightarrow \infty, \ i\in P_2^j$, если бесконечным является множество $P_2^j$. Тогда получим равенства $\bar{y}=u_j$ для всех $j\in J$, из которых следует включение (\ref{e1.1.13}). Лемма доказана.
\end{proof}

\begin{remark}\label{r1.1.8}
Если в методе положить
\begin{equation}\label{e1.1.18}
\varepsilon_0=0,
\end{equation}
то $D_{\varepsilon_0}=D$. Поскольку $y_i\notin D$ для всех $i\in K$ по построению, то согласно п. 3 метода ни одна из точек $x_k$ не будет зафиксирована ввиду условия (\ref{e1.1.4}). Кроме того, в этом случае не будут выбраны и числа $\varepsilon_1$, $\varepsilon_2$, $\ldots$.  В связи с этим отметим, что при выполнении (\ref{e1.1.18}) множества $Q_i$ для всех $i\in K$ будут иметь вид (\ref{e1.1.3}). 
\end{remark}

C учетом замечания \ref{r1.1.8} приведем утверждение, касающееся свойств последовательности $\{y_i\}$, $i\in K$,  построенной с условием (\ref{e1.1.18}).

\begin{theorem}\label{theorem1.1.1}
Пусть число $\varepsilon_0$ в методе  выбрано согласно  (\ref{e1.1.18}). Тогда любая предельная точка последовательности $\{y_i\}$, $i\in K$, принадлежит $X^*$, а если при этом для всех $i\in K$ выполняется  (\ref{e1.1.9}), то вся последовательность $\{y_i\}$, $i\in K$, сходится к множеству $X^*$.
\end{theorem}
\begin{proof}
Пусть $\{y_i\}, \ i\in K_1\subset K,$ -- любая сходящаяся подпоследовательность последовательности $\{y_i\}$, $ i\in K$, и $\bar{y}$ -- ее предельная точка. Как отмечено выше, для всех $i\in K$ выполняется (\ref{e1.1.3}).  Тогда по лемме \ref{lemma1.1.2} справедливо включение $\bar{y}\in D$, а значит, $f(\bar{y}) \ge f^*$.  С другой стороны, ввиду (\ref{e1.1.2}) $f(y_i)\le f^*$ для всех $i\in K$. Переходя в последнем неравенстве к пределу по $i\in K_1$, получим $f(\bar{y})\le f^*$. Таким образом, $f(\bar{y})=f^*$, и первое утверждение теоремы доказано.

Пусть теперь все точки $y_i$, $i\in K$, удовлетворяют условию (\ref{e1.1.9}). В силу (\ref{e1.1.3}), (\ref{e1.1.8}) $M_{i+1}\subset M_i, \ i\in K$, а значит, $f(y_{i+1})\ge f(y_i), \ i\in K$. Отсюда с учетом ограниченности $\{y_i\}, \ i\in K$, следует, что $\{f(y_i)\}, \ i\in K$, сходится. Тогда в силу уже доказанного первого утверждения последовательность $\{f(y_i)\}, \ i\in K$, является минимизирующей, и по теореме  (\cite[c. 62]{n_vasillev_f_p_1}) второе утверждение теоремы тоже доказано.
\end{proof}
Перейдем, наконец, к исследованию сходимости последовательности $\{x_k\}$, $k\in K$. Прежде всего покажем, что при условии положительности всех чисел $\varepsilon_k$, $k\in K$, наряду с последовательностью $\{y_i\}, \ i\in K$, будет построена соглано (\ref{e1.1.4}) и последовательность $\{x_k\}, \ k\in K$.
\begin{lemma}\label{lemma1.1.3}
Пусть последовательность $\{y_i\},\  i\in K$, построена предложенным методом, и при этом числа $\varepsilon_k$, $k\in K$, были выбраны так, что 
\begin{equation}\label{e1.1.19}
\varepsilon_k>0 \quad \forall k\in K. 
\end{equation}Тогда для каждого $k\in K$ существует такой номер $i=i_k$, что выполняется равенство  (\ref{e1.1.4}).
\end{lemma}

\begin{proof}
1) Пусть $k=0$. Если $y_0\in D_{\varepsilon_0}$, то согласно п. 3 метода $i_0=0, \ x_0=y_{i_0}=y_0$, и равенство (\ref{e1.1.4}) для $k=0$ выполняется. Поэтому будем считать, что $y_0\notin D_{\varepsilon_0}$. Покажем тогда существование номера $i=i_0>0$, для которого справедливо включение 
\begin{equation}\label{e1.1.20}
y_{i_0}\in D_{\varepsilon_0}. 
\end{equation}

Допустим противное, то есть
\begin{equation}\label{e1.1.21}
y_i\notin D_{\varepsilon_0} \quad \forall i\in K, \quad  i>0.
\end{equation}Выделим из ограниченной последовательности $\{y_i\}$, $i\in K$, $i>0$, сходящуюся подпоследовательность $\{y_i\}$, $i\in K^\prime$, и пусть $y^\prime$ -- ее предельная точка. Согласно п. 3 метода с учетом сделанных допущений для всех $i\in K$ множество $Q_i$ имеет вид (\ref{e1.1.3}). Тогда по лемме \ref{lemma1.1.2} выполняется включение $y^\prime \in D$, и, следовательно, 
\begin{equation}\label{e1.1.22}
F(y^\prime)\le 0.
\end{equation}
C другой стороны, в силу (\ref{e1.1.21}) $F(y_i)>\varepsilon_0$ для всех $i\in K^\prime$. Переходя в последнем неравенстве к пределу по $i\in K^\prime$, получим $F(y^\prime)\ge \varepsilon_0>0$, что противоречит (\ref{e1.1.22}). Таким образом, существование номера $i_0$, для которого справедливо (\ref{e1.1.20}), доказано, и равенство (\ref{e1.1.4}) имеет место при $k=0$.

2) Допустим теперь, что (\ref{e1.1.4}) выполняется при некотором фиксированном $k\ge 0$, то есть  $x_k=y_{i_k}$ при выбранном значении $k$. Покажем существование такого номера $i_{k+1}>i_k$, что 
\begin{equation}\label{e1.1.23}
y_{i_{k+1}}\in D_{\varepsilon_{k+1}}, 
\end{equation}тогда $x_{k+1}=y_{i_{k+1}}$, и лемма будет доказана.

Предположим противное, то есть 
\begin{equation}\label{e1.1.24}
y_i\notin D_{\varepsilon_{k+1}} \quad \forall i>i_k.
\end{equation}
Выберем среди точек $y_i, \ i>i_k$, сходящуюся подпоследовательность $\{y_i\}, \ i\in K^{\prime \prime}$, и пусть $y^{\prime \prime}$ -- ее предельная точка. Согласно (\ref{e1.1.24}) для всех $i\in K$, $i>i_k$, множества $Q_i$ по алгоритму задаются в виде (\ref{e1.1.3}). Значит, по лемме \ref{lemma1.1.2} выполняется неравенство $F(y^{\prime \prime})\le 0$. Как и в первой части доказательства, из условия (\ref{e1.1.24}) с учетом (\ref{e1.1.19}) легко получается неравенство $F(y^{\prime \prime})>0$, противоречащее предыдущему. Таким образом, показано существование номера $i_{k+1}$, удовлетворяющего (\ref{e1.1.23}). Лемма доказана.
\end{proof}

\begin{theorem}\label{theorem1.1.2}
Пусть последовательность $\{x_k\}, \ k\in K$, построена предложенным методом с условием, что числа $\varepsilon_k$, $k\in K$, выбраны согласно  (\ref{e1.1.19}), и
\begin{equation}\label{e1.1.25}
\varepsilon_k\rightarrow 0, \quad k\rightarrow \infty.
\end{equation}Тогда любая предельная точка этой последовательности принадлежит множеству $X^*$, а если для всех $k\in K$ выполняются неравенства 
\begin{equation}\label{e1.1.26}
f(x_{k+1})\ge f(x_k), 
\end{equation}
то вся последовательность $\{x_k\}$, $k\in K$, сходится к множеству $X^*$.
\end{theorem}

\begin{proof}
Так как в силу сделанного предположения последовательность $\{y_i\}$, $i\in K$, ограничена, то $\{x_k\}$, $k\in K$,  также ограничена. Пусть $\{x_k\}$, $k\in K^\prime \subset K$, -- любая сходящаяся подпоследовательность последовательности $\{x_k\}, \ k\in K$, и $\bar{x}$ -- ее предельная точка. Покажем, что 
\begin{equation}\label{e1.1.27}
\bar{x}\in X^*,
\end{equation}
тогда первое утверждение будет доказано. Действительно, для всех $k\in K$ справедливы неравенства
$$0<F(x_k)\le \varepsilon_k.$$
Отсюда с учетом (\ref{e1.1.25}) следует, что $\lim \limits_{k\in K}F(x_k)=0$. Тогда
 $$\lim \limits_{k\in K}F(x_k)=\lim \limits_{k\in K^\prime}F(x_k)=F(\bar{x})=0,$$
и $\bar{x}\in D$, следовательно,
$$f(\bar{x})\ge f^*.$$
C другой стороны, согласно (\ref{e1.1.2}) $f(x_k)\le f^*, \ k\in K^\prime$, а значит, $f(\bar{x})\le f^*$. Следовательно, $f(\bar{x})=f^*$, и включение (\ref{e1.1.27}) доказано.

Пусть теперь для последовательности $\{x_k\}$, $k\in K$,  выполняется условие (\ref{e1.1.26}). Тогда ввиду ограниченности $\{x_k\}$, $k\in K$, последовательность $\{f(x_k)\}$, $k\in K$, является сходящейся, и по доказанному выше $\lim \limits_{k\in K}f(x_k)=f^*$. Поэтому в силу вышеупомянутой  известной теоремы   (\cite{n_vasillev_f_p_1}, c. 62)  второе утверждение тоже доказано.
\end{proof}

Заметим, что условие (\ref{e1.1.26}) для построенной методом  последовательности $\{x_k\}$, $k\in K$, имеет место, если, например, для всех $k\in K$ выполняются включения $M_{i_{k+1}}\subset M_{i_k}$, а точки $y_{i_k}$ найдены согласно равенству (\ref{e1.1.9}), в котором $i=i_k$.

Отметим также, что для выделения из последовательности $\{x_k\}$, $k\in K$, минимизирующей подпоcледовательности $\{x_{k_l}\}, \  l\in K$, достаточно номера $k_l$ выбрать из условия $f(x_{k_{l+1}})\ge f(x_{k_l})$ для всех $ l\in K$.

При исследовании сходимости метода использовалось предположение об ограниченности последовательности $\{y_i\}$, $i\in K$. Ясно, что ограниченность $\{y_i\}$, $i\in K$,  можно обеспечить за счет соответствующего выбора множеств $M_0$ и $Q_i$.

Получим теперь для предложенного метода оценки точности решения задачи  и на их основе приведем оценки скорости сходимости последовательности $\{x_k\}$, $k\in K$.

Опишем сначала такой алгоритм метода, где на каждой итерации легко оценить значение $f^*$.

Пусть в (\ref{e1.1.1})  ${\rm int\,} D \ne \emptyset$, и известна точка $v\in{\rm int\,} D$. Положим в методе $v^j=v$ для всех $j\in J$. Пусть на $i$-ой итерации найдены согласно (\ref{e1.1.6}) точки $\bar{y}_i^j \in D_j$, $j\in J_i$, и зафиксирован такой номер $r_i\in J_i$, что $\bar{y}_i^{r_i}\in D$. Тогда, положив $\bar{y}_i=\bar{y}_i^{r_i}$, имеем оценки
$$f(y_i)\le f^*\le f(\bar{y}_i).$$

Покажем далее, что при некоторых дополнительных условиях на исходную задачу оценки близости значений $f(x_k)$ и приближений $x_k$ к значению $f^*$ и множеству $X^*$ можно получить и для последовательности $\{x_k\}$, $k\in K$, построенной общим методом.

Везде далее в данном параграфе будем считать, что для функций $f(x)$, $f_j(x)$, $j\in J$, и множества $D$ справедливо следующее

\begin{propose}\label{propose1.1.1}
Пусть функция $f(x)$ выпукла, а функции $f_j(x)$, $j\in J$, сильно выпуклы с константами сильной выпуклости $\mu_j$ соответственно. Кроме того, пусть множество $D$ удовлетворяет условию Слейтера, и любая точка абсолютного минимума функции $f(x)$, при наличии таковой, не принадлежит множеству ${\rm int\,} D$.
\end{propose}

\begin{theorem}\label{theorem1.1.3}
Пусть выполняется предположение \ref{propose1.1.1}, тогда решение задачи (\ref{e1.1.1}) единственно.
\end{theorem}

\begin{proof}
Допустим, что это утверждение неверно. Выберем  $x_1^*$, $x_2^*\in X^*$ так, что $x_1^* \ne x_2^*$, и положим $z=\alpha x_1^*+(1-\alpha)x_2^*$, где $0<\alpha<1$. Тогда с учетом строгой выпуклости фунции $F(x)$ и равенств $F(x_1^*)=F(x_2^*)=0$ имеем $F(z)<\alpha F(x_1^*)+(1-\alpha)F(x_2^*)=0$, а значит, $f(z)>f^*$. С другой стороны, $f(z)\le \alpha f(x_1^*)+(1-\alpha)f(x_2^*)=f^*$. Полученное противоречие доказывает утверждение.
\end{proof}

В связи с доказанным будем считать далее, что $X^*=\{x^*\}$.

Исследуем некоторые свойства функции $F(x)$.

\begin{lemma}\label{lemma1.1.4}
Функция $F(x)$ сильно выпуклая с константой сильной выпуклости $\mu = \min \limits_{j\in J} \mu_j$.
\end{lemma}

\begin{proof} Зафиксируем точки $u_0$, $u_1\in R_n$  и число $\alpha\in [0,1]$ такие, что $u_1 \ne u_2$.
С учетом условия сильной выпуклости функций $f_j(x)$, $j\in J$, справедливы следующие соотношения:
$F(\alpha u_0 + (1 - \alpha) u_1) =\max \limits_{j\in J} \{f_j(\alpha u_0 + (1 - \alpha)u_1)\}\le \max \limits_{j\in J} \{\alpha f_j(u_0) + (1 - \alpha) f_j(u_1) - \alpha (1 - \alpha) \mu_j\|u_0 - u_1\|^2\}\le \alpha \max \limits_{j\in J} f_j(u_0) + (1 - \alpha) \max \limits_{j\in J} f_j(u_1) - \alpha (1-\alpha) \min \limits_{j\in J} \mu_j \|u_0 - u_1\|^2.$ Отсюда следует справедливость утверждения.
\end{proof}

\begin{lemma}\label{lemma1.1.5}(\cite[c. 221]{n_vasillev_f_p_2}).
Пусть $U$ -- открытое выпуклое множество из ${\rm R_n}$, пусть функция $J(u)$ равномерно выпукла на $U$ с модулем выпуклости $\delta(t)$, $\partial F(u)$ -- субдифференциал функции $J(u)$ в точке $u\in {\rm R_n}$. Тогда необходимо выполняются неравенства
$$J(u)\ge J(v) + \langle c(v), u - v\rangle + \delta(|u-v|),$$
$$\langle c(u) - c(v), u - v\rangle \ge 2\delta (|u-v|)$$
при всех $c(v)\in \partial J(v)$, $c(u) \in \partial J(u)$ и всех $u$, $v\in U$.
\end{lemma}

Положим 
\begin{equation}\label{e1.1.28}
X_1^*(\delta) =\{x\in {\rm R_n} : \|x-x^*\|\le \delta\},
\end{equation}
\begin{equation}\label{e1.1.29}
X_2^*(\gamma)=\{x\in {\rm  R_n} : |f(x)-f^*|\le \gamma \},
\end{equation}
 где $\gamma, \ \delta \ge 0$. Пользуясь методикой, близкой к \cite{n_zabotin_ya_i_7}, докажем следующее вспомогательное утверждение.

\begin{lemma}\label{lemma1.6}
Пусть $\varepsilon >0$, точка $y\in D_\varepsilon$ такова, что $y\ne x^*$ и при всех $\alpha \in (0, 1]$ для точки $z=\alpha y + (1-\alpha)x^*$ выполняется  неравенство
\begin{equation}\label{e1.1.30}
F(z) >0.
\end{equation}
Тогда справедливо включение
\begin{equation}\label{e1.1.31}
y\in X_1^*(\delta),
\end{equation}
где $\delta =\displaystyle \sqrt{\frac{\varepsilon}{\mu}}$.
\end{lemma}

\begin{proof}
Согласно (\ref{e1.1.30}) и равенству $F(x^*)=0$ имеем $F(z)-F(x^*)=F(x^*+\alpha (y-x^*))-F(x^*)>0$, а значит,
$$
\frac{F(x^*+\alpha(y-x^*))-F(x^*)}{\alpha}>0
$$
для любого $\alpha \in (0,1]$. Переходя в последнем неравенстве к пределу по $\alpha \to +0$, получим для производной $F^\prime(x^*, y-x^*)$ функции $F(x)$ в точке $x^*$ по направлению $y-x^*$ следующее неравенство:
\begin{equation}\label{e1.1.32}
F^\prime(x^*, y-x^*) \ge 0.
\end{equation}
Известно (см., например, \cite{n_pshenichnii_b_n_0}, c. 74), что
$$F^\prime(x^*, y-x^*)=\max \{\langle c(x^*), y-x^*\rangle : c(x^*)\in \partial F(x^*)\},$$
где $\partial F(x^*)$ -- субдифференциал функции $F(x)$ в точке $x^*$. Тогда отсюда и из (\ref{e1.1.32}) вытекает существование такого $c^\prime(x^*) \in \partial F(x^*)$, что 
\begin{equation}\label{e1.1.33}
\langle c^\prime (x^*), y-x^*\rangle \ge 0.
\end{equation}

 Поскольку функция $F(x)$ сильно выпукла, то согласно лемме \ref{lemma1.1.5}  для вектора $c^\prime(x^*)$ имеет место неравенство 
$$F(y)\ge F(x^*) +\langle c^\prime (x^*), y-x^*\rangle + \mu \|y-x^*\|^2,$$
из которого с учетом (\ref{e1.1.33}) и равенства $F(x^*)=0$ следует, что $F(y)\ge \mu\|y-x^*\|^2$. Но по условию $F(y)\le \varepsilon$. Значит,
\begin{equation}\label{e1.1.34}
\|y-x^*\|\le \sqrt{\frac{\varepsilon}{\mu}},
\end{equation}
и утверждение леммы доказано.
\end{proof}

\begin{lemma}\label{lemma1.1.7}
 Пусть выполняются условия леммы \ref{lemma1.6}, и, кроме того, функция $f(x)$ удовлетворяет условию Липшица с константой $L$. Тогда 
\begin{equation}\label{e1.1.35}
y\in X^*_1(\delta)\bigcap X^*_2(\gamma),
\end{equation}
где $\delta$ определено в  (\ref{e1.1.31}), а $\gamma=L\delta$.
\end{lemma}

\begin{dokaz}  утверждения следует из (\ref{e1.1.31}), а также неравенства
\begin{equation}\label{e1.1.36}
|f(y)-f(x^*)|\le L\|y-x^*\|
\end{equation}
и оценки (\ref{e1.1.34}).
\end{dokaz}

\begin{lemma} \label{lemma1.1.8}
 Пусть $\varepsilon>0$. Тогда для любой точки $y\in D_\varepsilon \bigcap E^*$, $y\ne x^*$, выполняется  (\ref{e1.1.31}), а если функция $f(x)$ удовлетворяет условию Липшица, то справедливо и включение  (\ref{e1.1.35}).
\end{lemma}

\begin{proof}
Положим $z=\alpha y+(1-\alpha)x^*$, где $\alpha\in (0, 1]$. Докажем справедливость неравенства (\ref{e1.1.30}). Тогда для точки $y$ будут выполнены условия леммы \ref{e1.1.6}, и утверждение (\ref{e1.1.31}) будет доказано.

Предположим противное, допустим, что $F(z)\le 0$. Тогда $z\in D$ и $f(z)\ge f^*$. С другой стороны, поскольку $y\in E^*$, $x^*\in E^*$, а множество $E^*$ выпукло, то $z\in E^*$, и $f(z)\le f^*$. Следовательно,
\begin{equation}\label{e1.1.37}
f(z)=f^*.
\end{equation}
Однако, согласно теореме \ref{theorem1.1.3} решение задачи (\ref{e1.1.1}) единственно, а $z\ne x^*$,  так как  $y\ne x^*$ и $\alpha>0$. Полученное с равенством (\ref{e1.1.37}) противоречие доказывает неравенство (\ref{e1.1.30}), а значит, и включение (\ref{e1.1.31}).

Вторая часть утверждения леммы вытекает из (\ref{e1.1.31}), (\ref{e1.1.36}), (\ref{e1.1.34}).
\end{proof}

\begin{lemma}\label{lemma1.1.9}
 Пусть последовательность $\{y_i\}$, $i\in K$, построена предложенным методом. Если для некоторого номера $i\in K$ и $r>0$ выполняется $y_i\in D_r$, то $y_i \in X^*_1(\delta)$ при $\delta =\displaystyle \sqrt{\frac{r}{\mu}}$, а если к тому же $f(x)$ удовлетворяет условию Липшица, то $y_i \in X^*_1(\delta)\bigcap X^*_2(\gamma)$ при $\delta=\displaystyle \sqrt{\frac{r}{\mu}}$, $\gamma=L\delta$.
\end{lemma}

\begin{dokaz} утверждений непосредственно следует из леммы \ref{lemma1.1.8}, так как согласно (\ref{e1.1.2}) $y_i\in E^*$, и $y_i \ne x^*$ для всех $i\in K$.
\end{dokaz}

Лемма \ref{lemma1.1.9} дает возможность оценивать на каждой итерации близость точки $y_i$ к решению задачи. А именно, поскольку для всех $i\in K$ выполняется включение $y_i\in D_{r_i}$, где $r_i=F(y_i)$, то по  лемме \ref{lemma1.1.9}
$$\|y_i-x^*\|\le \delta_i,$$
где $\delta_i=\displaystyle \sqrt{\frac{F(y_i)}{\mu}}$. Кроме того, если функция $f(x)$ удовлетворяет условию Липшица, то для всех $i\in K$ справедливо также неравенство
$$|f(y_i)-f^*|\le L\delta_i.$$

\begin{theorem}\label{theorem1.1.4}
Пусть последовательность $\{x_k\}$, $k\in K$,  построена предложенным методом. Тогда для каждого $k\in K$ имеет место следующая оценка:
\begin{equation}\label{e1.1.38}
\|x_k-x^*\|\le \delta_k,
\end{equation}
где $\delta_k=\displaystyle \sqrt{\frac{\varepsilon_k}{\mu}}$. Если же $f(x)$ удовлетворяет условию Липшица, то, кроме  (\ref{e1.1.38}), при каждом $k\in K$ справедлива также оценка
\begin{equation}\label{e1.1.39}
|f(x_k)-f^*|\le L\delta_k.
\end{equation}
\end{theorem}

\begin{proof}
Согласно (\ref{e1.1.4}) $x_k=y_{i_k}$ для всех $k\in K$, причем $y_{i_k}\in D_{\varepsilon_k}$ и $y_{i_k}\ne x^*$. Тогда по лемме \ref{lemma1.1.9} выполняется включение $x_k\in X^*_1(\delta_k)$, $k\in K$, а значит, и неравенство (\ref{e1.1.38}). При условии Липшица на функцию $f(x)$ по той же лемме \ref{lemma1.1.9} имеет место и оценка (\ref{e1.1.39}). Теорема доказана.
\end{proof}

Зададим теперь числа $\varepsilon_k$, $k\ge 1$, в методе в виде
\begin{equation}\label{e1.1.40}
\varepsilon_k=\frac{1}{k^p}, 
\end{equation}
где $p>0$. Тогда согласно (\ref{e1.1.38}) для последовательности $\{x_k\}$, $k\in K$, справедливы следующие оценки скорости сходимости:
\begin{equation}\label{e1.1.41}
\|x_k-x^*\|\le \frac{1}{\sqrt{\mu}k^{p/2}},\  k\ge 1.
\end{equation}
Кроме того, если для $f(x)$ выполняется условие Липшица, то наряду, с (\ref{e1.1.41}), имеем также в силу (\ref{e1.1.39}) оценки
\begin{equation}\label{e1.1.42}
|f(x_k)-f^*|\le \frac{L}{\sqrt{\mu}k^{p/2}}, \ k\ge 1.
\end{equation}

\begin{remark}\label{r1.1.9}
При очень большом $p>0$ значение в (\ref{e1.1.40}) $\varepsilon_k$ очень мало. Но при малом $\varepsilon_k$ для построения точки $x_k=y_{i_k}$, для которой с учетом необходимых предположений справедливы оценки типа (\ref{e1.1.41}), (\ref{e1.1.42}),  может понадобится много вспомогательных итераций по $i$, а количество этих вспомогательных итераций с построением точек $y_i$ в оценке никак не учитывается. Отсюда приведенные оценки скорости сходимости (\ref{e1.1.41}), (\ref{e1.1.42}) несколько условны, так как касается только для точек $x_k$ и не отражает полного количества итераций, необходимых для построения точек $x_k$.
\end{remark}

\begin{remark}\label{r1.1.10}
Как уже отмечено в (\ref{e1.1.40}), числа $\varepsilon_k$, $k\ge1$, допустимо выбирать, как и число $\varepsilon_0$, на предварительном шаге метода. Однако, в таком случае последовательность $\{\varepsilon_k\}$, $k\in K$, не будет адаптирована к процессу минимизации. Поэтому в методе предусмотрена возможность выбора чисел $\varepsilon_k$, $k\ge 1$, на шаге 3, то есть  в ходе отыскания приближений. Приведем пример задания последовательности $\{\varepsilon_k\}$, $k\in K$, связанной с процессом построения $\{x_k\}$, $k\in K$.

Не выбирая конкретного значения, считаем $\varepsilon_0$ настолько большим, что $y_0\in D_{\varepsilon_0}$ и $x_0=y_0$. Для всех $k\ge 0$ положим
$$\varepsilon_{k+1}=\sigma_k F(x_k),$$
где $0<\sigma_k<1$. Тогда для $\{\varepsilon_k\}$, $k\in K$, справедливо (\ref{e1.1.19}), а при условии, что $\sigma_k \to 0$, $k\in K$, выполняется (\ref{e1.1.25}).
\end{remark}

\setcounter{equation}{0}

\section[Метод отсечений c отбрасыванием секущих плоскостей  и возможностью построения на его основе смешанных алгоритмов]
{Метод отсечений c отбрасыванием секущих
плоскостей  и возможностью построения
на его основе смешанных алгоритмов}

В данном параграфе предлагается метод отсечений для решения задачи (\ref{e1.1.1}), в которой область ограничений $D$ имеет более общий вид. Предлагаемый метод также характерен тем, что предусматривает возможность обновления аппроксимирующих множеств. Метод отличается от разработанного в $\S$ 1.1 метода, в частности, способами построения аппроксимирующих множеств.

Пусть  для функций $f(x)$, $f_j(x)$, $j\in J=\{1,\ldots, m\}$,  множеств $D_j$, $j\in J$,  выполняются условия, приведенные в $\S \ 1.1$. Пусть
\begin{equation}\label{e1.2.0}
D^\prime = \{x\in {\rm R_n} : f_j(x) \le 0, \ j\in J\},
\end{equation}
$D^{\prime\prime} \subset {\rm R_n}$ -- выпуклое замкнутое множество,
$$D=D^\prime \bigcap D^{\prime\prime}.$$ Решается задача (\ref{e1.1.1}).

Заметим, что внутренность множества $D$ может быть пустой, если, например, ${\rm int\,} D^{\prime} = \emptyset$ или ${\rm int\,} D^{\prime\prime} = \emptyset$.

Определим функцию $F(x)$ как и в $\S$ 1.1, и положим
$$D^\prime_\varepsilon = \{x\in {\rm  R_n} : F(x)\le \varepsilon\},$$
где $\varepsilon\ge 0$. Будем использовать обозначения  $f^*$, $X^*$, $E^*$, $K$ предыдущего параграфа.

Предлагаемый метод решения задачи (\ref{e1.1.1}) вырабатывает последовательности приближений $y_i$, $i\in K$, $x_k$, $z_k$, $k\in K$, и заключается в следующем.

\textbf{Метод 1.2.} Строится выпуклое замкнутое множество $M_0\subset {\rm R_n}$, содержащее хотя бы одну точку множества $X^*$, например, точку $x^*$. Задается число $\varepsilon_0\ge 0$. Полагается $k=0$, $i=0$.

\textbf{1.} Находится точка
\begin{equation}\label{e1.2.1}
y_i\in M_i \bigcap D^{\prime\prime} \bigcap E^*.
\end{equation}
Если выполняется включение
\begin{equation}\label{e1.2.2}
y_i\in D^\prime,
\end{equation}
то $y_i\in X^*$, и процесс завершается.

\textbf{2.} Если $y_i\notin D^\prime_{\varepsilon_k}$, то выбирается выпуклое замкнутое множество $G_i\subset {\rm  R_n}$, содержащее точку $x^*$, полагается
\begin{equation}\label{e1.2.3}
Q_i=M_i\bigcap G_i, 
\end{equation}
\begin{equation}\label{e1.2.4}
u_i=y_i,
\end{equation}
и следует переход к п. 4.

\textbf{3.} Полагается $i_k=i$,
\begin{equation}\label{e1.2.5}
x_k=y_{i_k},
\end{equation}
выбирается выпуклое замкнутое множество $Q_i=Q_{i_k}$ такое, что 
\begin{equation}\label{e1.2.6}
x^*\in Q_i.
\end{equation}
Выбирается точка 
\begin{equation}\label{e1.2.7}
z_k\in M_{i_k} \bigcap D^{\prime\prime},
\end{equation}
удовлетворяющая неравенствам
\begin{equation}\label{e1.2.8}
f(x_k)\le f(z_k)\le f^*.
\end{equation}
Если $z_k\in D^\prime$, то $z_k \in X^*$, и процесс завершается. В противном случае полагается 
\begin{equation}\label{e1.2.9}
u_i=z_k.
\end{equation}
Задается число $\varepsilon_{k+1}\ge 0$, значение $k$ увеличивается на единицу.

\textbf{4.} Выбирается конечное множество
\begin{equation}\label{e1.2.10}
A_i\subset \partial F(u_i)
\end{equation}
или конечные множества
\begin{equation}\label{e1.2.11}
A_i^j\subset \partial f_j(u_i)
\end{equation}
для всех $j\in J_i=\{j\in J : y_i \notin D_j\}$.

\textbf{5.} Полагается 
\begin{equation}\label{e1.2.12}
M_{i+1}=Q_i\bigcap T_i,
\end{equation}
где 
\begin{equation}\label{e1.2.13}
T_i = \{x\in {\rm R_n} : F(u_i) + \langle a, x- u_i\rangle \le 0 \quad \forall a\in A_i\} 
\end{equation}
в случае (\ref{e1.2.10}), или
\begin{equation}\label{e1.2.14}
T_i =\bigcap \limits_{j\in J_i} \{x\in {\rm R_n} : f_j(u_i) + \langle a, x- u_i\rangle \le 0 \quad \forall a\in A_i^j\}
\end{equation}
в случае (\ref{e1.2.11}). Следует переход к п. 1 при $i$, увеличенном на единицу.

Прежде, чем обсуждать свойства строящихся методом последовательностей, приведем некоторые замечания к методу.

\begin{remark}\label{r1.2.1}
Выбор точек $y_i$ из условия (\ref{e1.2.1}) возможен для каждого $i\in K$. Действительно, поскольку множество $M_0$ по построению содержит точку $x^*$, и согласно п. 2 метода $x^*\in G_i$ при всех $i\in K$, то ввиду (\ref{e1.2.3}), (\ref{e1.2.6}), (\ref{e1.2.12}) $x^*\in M_i$ для всех $i\in K$. В то же время $x^*\in D^{\prime\prime}$, и, кроме того, $x^*\in E^*$. Следовательно,
$$M_i\bigcap D^{\prime\prime}\bigcap E^*\ne \emptyset \quad \forall i\in K.$$
\end{remark}
\begin{remark}\label{r1.2.2}
В частности, точку $y_i$  можно отыскивать согласно равенству
\begin{equation}\label{e1.2.15}
f(y_i)=\min \{f(x) : x\in M_i\bigcap D^{\prime\prime}\}.
\end{equation}
Если $f(x)$ -- линейная функция, то точки $y_i$ естественно искать в виде (\ref{e1.2.15}).

Условие (\ref{e1.2.1}) позволяет отыскивать $y_i$, как точки приближенного минимума функции $f(x)$ на множествах $M_i\bigcap D^{\prime\prime}$, $i\in K$.
\end{remark}

\begin{remark}\label{r1.2.3}
Если для точки $y_i$ выполняется включение (\ref{e1.2.2}), то в силу (\ref{e1.2.1}) $y_i\in D$ и $f(y_i)=f^*$, то есть, как отмечено в п. 1 метода, $y_i$ -- решение исходной задачи.
\end{remark}

\begin{remark}\label{r1.2.4}
Множество $M_0$ можно выбирать разными способами. Например, допустимо считать, что $M_0={\rm R_n}$. Если множество $D^{\prime\prime}$ определено системой линейных равенств или неравенств, а функции $f(x)$ и $f_j(x)$, $j\in J^\prime\subset J$, линейны, то множество $M_0$  удобно задать в виде (\ref{e1.1.10}).
\end{remark}

Ввиду (\ref{e1.2.3}), (\ref{e1.2.6}) для задания множеств $Q_i$, участвующих в построении аппроксимирующих множеств $M_{i+1}$, также немало возможностей. Если для всех $i\in K$, независимо от принадлежности точки $y_i$ множеству $D^\prime_{\varepsilon_k}$, положить $Q_i=M_i$, считая $G_i={\rm R_n}$, то согласно (\ref{e1.2.12}) $M_{i+1}=M_i\bigcap T_i$, $i\in K$. В таком случае от итерации к итерации происходит накопление отсекающих плоскостей, и никаких обновлений погружающих множеств $M_i$ в процессе решения задачи не происходит. Покажем теперь аналогично замечанию \ref{r1.1.7} к методу $1.1$, как за счет выбора множеств  $Q_i$ можно проводить в предлагаемом методе $2.1$ вышеупомянутые обновления.

\begin{remark}\label{r1.2.5} Допустим, что  в соотношении (\ref{e1.2.3}) $G_i = {\rm R_n}$ при всех $i\in K$, для некоторого номера $i=i_k$ точка $y_i$ выбрана согласно
\begin{equation}\label{e1.2.17_star1}
y_i\in D_{\varepsilon_k}^\prime,
\end{equation}
 а множество $Q_i$ -- согласно (\ref{e1.1.12}), где $0\le r_i\le i = i_k$. Согласно выбору множества $Q_i=Q_{i_k}$ выполняются включения (\ref{e1.2.6}) при всех $r_i= 0, \ldots, i$. В связи с этим  множество $Q_{i_k}$ можно задать с помощью любых множеств $M_0,\ldots, M_{i_k}$, которые построены к $i_k$-ому шагу. Следовательно, на каждом шаге $i=i_k$ будет инициирован процесс отбрасывания любых накопившихся к шагу $i_k$ секущих плоскостей.
\end{remark}

\begin{remark}\label{r1.2.6}
Для выбора точек $z_k$ из условий (\ref{e1.2.7}), (\ref{e1.2.8}) тоже есть разные возможности. Например, не проводя специальных вычислений, можно положить 
$$z_k=x_k.$$

При дополнительных условиях на $f(x)$  для поиска точки $z_k$ допустимо применить метод внешних штрафных функций. Тогда в качестве $z_k$ можно выбрать ту из точек минимума вспомогательных штрафных функций на $M_{i_k} \bigcap D^{\prime\prime}$, в которой значение целевой функции окажется не меньшим, чем значение $f(x_k)$. В таком случае построится смешанный алгоритм на основе метода 2.1 с привлечением метода штрафов, а поскольку сходимость метода отсечений ниже будет доказана, то тем самым будет обоснована и сходимость смешанного алгоритма.
\end{remark}

Перейдем к исследованию сходимости метода. Будем считать далее в этом параграфе, что последовательность $\{y_i\}$, $i\in K$, и последовательность $\{z_k\}$, $k\in K$,  в случае ее построения, ограничены.

\begin{lemma}\label{lemma1.2.1}
 Для всех $i\in K$ справедливо включение
\begin{equation}\label{e1.2.16}
x^*\in M_i.
\end{equation}
\end{lemma}
\begin{proof}
 При $i=0$ включение (\ref{e1.2.16}) выполняется по условию выбора множества $M_0$ на предварительном шаге метода. Допустим теперь, что (\ref{e1.2.16}) имеет место при $i=l\ge 0$. Покажем справедливость включения (\ref{e1.2.16})  при $i =l+1$, тогда лемма будет доказана.

Отметим, что ввиду (\ref{e1.2.1}), (\ref{e1.2.6}) с учетом индукционного предположения $x^*\in Q_l$. Значит, для обоснования леммы согласно (\ref{e1.2.12}) достаточно доказать, что 
\begin{equation}\label{e1.2.17}
x^*\in T_l.
\end{equation}

Предположим сначала, что множество $T_l$ было выбрано на шаге $i=l$ согласно (\ref{e1.2.13}). Так как функция $F(x)$ выпукла, то для всех $a\in A_l$ выполняется неравенство $F(x^*)-F(u_l)\ge \langle a, x^* - u_l\rangle$. Отсюда с учетом того, что $F(x^*)\le 0$, имеем $F(u_l) + \langle a, x^* - u_l\rangle\le 0 $ для всех $a\in A_l$, и включение (\ref{e1.2.17}) доказано.

Пусть теперь множество $T_l$ имеет вид (\ref{e1.2.14}), где $i=l$. Поскольку согласно (\ref{e1.2.4}), (\ref{e1.2.9}) $u_l\notin D^\prime$, то $J_l \ne \emptyset$. Тогда для всех $a\in A_l^j$, $j\in J_l$  в силу выпуклости функции $f_j(x)$  имеют место неравенства $f_j(x^*)- f_j(u_l)\ge \langle a, x^* - u_l\rangle$. Но $f_j(x^*)\le 0$ для всех $j\in J$, в том числе и для всех $j\in J_l$. Следовательно, $f_j(u_l) + \langle a, x^* - u_l \rangle \le 0$ для всех $a\in A_l^j$, $j\in J_l$. Включение (\ref{e1.2.17}), а с ним и лемма, доказаны.
\end{proof}

\begin{lemma}\label{lemma1.2.2}
 Пусть последовательность $\{u_i\}$, $i\in K$, построена предложенным методом с условием, что для всех $i\in K$, начиная с некоторого номера $\tilde{i} \ge 0$ множества $Q_i$ выбраны согласно (\ref{e1.2.3}). Тогда для любых номеров $i^\prime$, $i^{\prime\prime}\in K$ таких, что 
$$i^{\prime\prime} > i^\prime \ge \tilde{i},$$
выполняется включение 
\begin{equation}\label{e1.2.18}
u_{i^{\prime\prime}}\in T_{i^\prime}.
\end{equation}
\end{lemma}

\begin{proof}
 В силу (\ref{e1.2.1}), (\ref{e1.2.4}) и (\ref{e1.2.7}), (\ref{e1.2.9}) для всех $i\in K$ имеют место включения $u_i\in M_i$. Значит, 
\begin{equation}\label{e1.2.19}
u_{i^{\prime\prime}}\in M_{i^{\prime\prime}}.
\end{equation}

Согласно (\ref{e1.2.12}) и условию леммы $M_{i^{\prime\prime}}\subset M_{i^\prime+1}\subset T_{i^\prime}$. Отсюда и из (\ref{e1.2.19}) следует (\ref{e1.2.18}). Лемма доказана.
\end{proof}

Отметим, что при каждом $i\in K$, независимо от выполнения для точки $y_i$ условия (\ref{e1.2.17_star1}), множество $Q_i$ можно задавать в виде (\ref{e1.2.3}), поскольку в силу леммы \ref{lemma1.2.1} включение (\ref{e1.2.6}) имеет место. Если $Q_i$ для всех $i\ge \tilde{i}\ge 0$ выбирать в виде (\ref{e1.2.3}), то, начиная с шага $\tilde{i}$, будет происходить накопление отсекающих плоскостей и обновлений погружающих множеств не произойдет. С учетом этого замечания сформулируем следующее утверждение.

\begin{lemma}\label{lemma1.2.3}
Пусть последовательность $\{u_i\}$, $i\in K$, построена предложенным методом с условием, что для всех $i\in K$, начиная с некоторого номера $\tilde{i}\ge 0$, множества $Q_i$ выбраны согласно (\ref{e1.2.3}). Тогда любая предельная точка последовательности $\{u_i\}$, $i\in K$, принадлежит множеству $D$.
\end{lemma}

\begin{proof} Допустим, что утверждение леммы неверно. Тогда существует такая сходящаяся подпоследовательность $\{u_i\}$, $i\in K_1\subset K$, последовательности $\{u_i\}$, $i\in K$, что ее предельная точка $\bar{u}$ удовлетворяет неравенству $F(\bar{u})>0$. Положим $F(\bar{u})=\gamma$. Так как функция $F(x)$ непрерывна, то найдется окрестность $\omega$ точки $\bar{u}$ такая, что
$$F(x)\ge \frac{\gamma}{2}$$
для всех $x\in \omega$. Зафиксируем такой номер $\hat{i}\in K_1$, что $\hat{i}\ge \tilde{i}$ и для всех $i\in K_1$, $i\ge \hat{i}$ имеет место включение $u_i\in \omega$. Тогда 
\begin{equation}\label{e1.2.20}
F(u_i)=\max \limits_{j\in J} f_j(u_i) \ge \frac{\gamma}{2} \quad \forall i\in K_1, \quad i\ge \hat{i}.
\end{equation}

В силу (\ref{e1.2.20}) и конечности множества $J$  существует $r\in J$, при котором для бесконечного числа номеров $i\in K_1$, $i\ge \hat{i}$, выполнется неравенство
\begin{equation}\label{e1.2.21}
f_r(u_i)\ge \frac{\gamma}{2}.
\end{equation}

Положим 
$$K_2=\{i\in K_1 : i\ge \hat{i},\ f_r(u_i) \ge \frac{\gamma}{2}\}.$$

Выберем номера $i^\prime, i^{\prime\prime}\in K_2$ так, что $i^{\prime\prime}>i^\prime$. Тогда по лемме \ref{lemma1.2.2} выполняется включение (\ref{e1.2.17}). Если множество $T_i$ при $i=i^\prime$ было выбрано согласно (\ref{e1.2.13}), то ввиду (\ref{e1.2.17})
\begin{equation}\label{e1.2.22}
F(u_{i^\prime})+ \langle a, u_{i^{\prime\prime}} - u_{i^\prime}\rangle \le 0 \quad \forall a\in A_{i^\prime}. 
\end{equation}
Если же $T_{i^\prime}$ имеет вид (\ref{e1.2.14}), то из того же включения (\ref{e1.2.17}) следует неравенство
\begin{equation}\label{e1.2.23}
f_r(u_{i^\prime}) + \langle a, u_{i^{\prime\prime}} - u_{i^\prime}\rangle \le 0
\end{equation}
для всех $a\in A_{i^\prime}^r$, так как $f_r(u_{i^\prime})>0$, и $r\in J_{i^\prime}$. Таким образом, из (\ref{e1.2.22}), (\ref{e1.2.23}) и неравенств (\ref{e1.2.20}), (\ref{e1.2.21}) при $i=i^\prime$ имеем неравенство 
\begin{equation}\label{e1.2.24}
\langle a, u_{i^{\prime\prime}}-u_{i^\prime}\rangle \le -\frac{\gamma}{2}
\end{equation}
для всех $a\in A_{i^\prime}$, если множество $T_{i^\prime}$ выбрано в виде (\ref{e1.2.13}), и имеем то же неравенство (\ref{e1.2.24}) для всех $a\in A_{i^\prime}^r$, если $T_{i^\prime}$ задано равенством (\ref{e1.2.14}).

Далее, ввиду ограниченности последовательности $\{u_i\}$, $i\in K$, найдется (см., напр., \cite[121]{n_polayk_b_t_1}) такое число $\theta>0$, что для всех $a\in \partial F(u_i)$, $a\in \partial f_r(u_i)$, $i\in K_2$, выполняется неравенство
\begin{equation}\label{e1.2.25}
\|a\|\le \theta.
\end{equation}
Тогда из (\ref{e1.2.24}), (\ref{e1.2.25}) следует, что 
\begin{equation}\label{e1.2.26}
\theta \|u_{i^{\prime\prime}}- u_{i^\prime}\|\ge \frac{\gamma}{2}.
\end{equation}
Выберем теперь для каждого $i\in K_2$ такой номер $p_i\in K_2$, что $p_i\ge i+1$. Согласно (\ref{e1.2.26}) $\theta \|u_{p_i} - u_i\|\ge \gamma/2$. Последнее неравенство противоречиво, поскольку $u_i \to \bar{u}$ и $u_{p_i}\to \bar{u}$ при $i\to \infty$, $i\in K_2$. Лемма доказана.
\end{proof}

\begin{theorem}\label{theorem1.2.1}
Если для последовательности $\{u_i\}$, $i\in K$, выполняются условия леммы \ref{lemma1.2.3}, то любая ее предельная точка принадлежит множеству $X^*$.
\end{theorem}

\begin{proof}
Пусть $\{u_i\}$, $i\in K_1\subset K$, -- любая сходящаяся подпоследовательность последовательности $\{u_i\}$, $i\in K$, и $\bar{u}$ -- ее предельная точка. Тогда по лемме \ref{lemma1.2.3} справедливо включение $\bar{u}\in D$, а значит $f(\bar{u}) \ge f^*$. C другой стороны, ввиду (\ref{e1.2.1}), (\ref{e1.2.4}) и (\ref{e1.2.8}), (\ref{e1.2.9}) для всех $i\in K$, $i\ge \tilde{i}$, выполняется неравенство $f(u_i)\le f^*$. Переходя в этом неравенстве к пределу по $i\in K_1$, получим $f(\bar{u})\le f^*$. Таким образом, $f(\bar{u})=f^*$, и утверждение теоремы доказано.
\end{proof}

Заметим,что условия (\ref{e1.2.7})--(\ref{e1.2.9}) выбора точек $z_k$ позволяют задавать точки $u_i=z_k$ в виде (\ref{e1.2.4}), считая $z_k=x_k=y_{i_k}$. Следовательно, допустим случай, когда, начиная с некоторого номера $\tilde{i}\ge 0$, множества $Q_i$ имеют вид (\ref{e1.2.3}), а члены последовательности $\{u_i\}$, $i\in K$, -- вид (\ref{e1.2.4}). С учетом этого замечания приведем следующее утверждение.

\begin{theorem}\label{theorem1.6}
 Пусть выполняются условия леммы \ref{lemma1.2.3}. Кроме того, пусть для всех $i\ge \tilde{i}\ge 0$ выполняются равенства (\ref{e1.2.4}), и точки $y_i$ выбраны согласно условию (\ref{e1.2.15}). Тогда поcледовательность $\{u_i\}$, $i\in K$, сходится к множеству $X^*$.
\end{theorem}

\begin{proof}
В силу (\ref{e1.2.3}), (\ref{e1.2.12}) для всех $i\in K$, $i\ge \tilde{i}$, имеют место включения $M_{i+1}\subset M_i$, а значит, согласно (\ref{e1.2.4}), (\ref{e1.2.15}) $f(u_{i+1})\ge f(u_i)$, $i\in K$, $i\ge \tilde{i}$. Отсюда ввиду ограниченности $\{u_i\}$, $i\in K$, следует, что $\{f(u_i)\}$, $i\in K$, является минимизирующей, и с учетом известной теоремы (напр., \cite[62]{n_vasillev_f_p_1}) утверждение доказано.
\end{proof}

 Перейдем теперь к исследованию последовательностей $\{x_k\}$, $\{z_k\}$, $k\in K$, строящихся описанным методом.

\begin{remark}\label{r1.2.7}
Отметим, что если положить $\varepsilon_0=0$, то при всех $i\in K$ в силу соотношений $D^\prime=D^\prime_{\varepsilon_0}$ и $y_i\in D^\prime$ шаг 3 метода будет пропускаться, и точки $x_k$, $z_k$, $k\in K$, построены не будут. Аналогично, если при некотором $k>0$ положить на шаге 3 метода $\varepsilon_{k+1}=0$, то при всех $i>i_k$ будут выполняться равенства (\ref{e1.2.3}), (\ref{e1.2.4}), значения $k$ перестанут изменяться, и процесс построения точек $x_k$, $z_k$ прекратится. 
\end{remark}

В связи с замечанием \ref{r1.2.7} прежде всего покажем, что при условии положительности всех чисел $\varepsilon_k$, $k\in K$,  вместе с последовательностью $\{y_i\}$, $i\in K$, будут построены и последовательности $\{x_k\}$, $\{z_k\}$, $k\in K$.

\begin{lemma}\label{lemma1.2.4}
 Пусть последовательность $\{y_i\}$, $i\in K$, построена предложенным методом с условием, что числа $\varepsilon_k$, $k\in K$, выбраны согласно (\ref{e1.1.19}). Тогда для каждого $k\in K$ существует такой номер $i=i_k\in K$, что выполняется равенство (\ref{e1.2.5}).
\end{lemma}

\begin{proof}  1) Пусть $k=0$. Если $y_0\in D^\prime_{\varepsilon_0}$, то согласно п. 3 метода $i_0=0$, $x_0=y_0$, и равенство (\ref{e1.2.5}) при $k=0$ выполняется. Поэтому будем считать, что $y_0\notin D^\prime_{\varepsilon_0}$. Покажем тогда существование номера $i=i_0>0$, для которого справедливо включение
\begin{equation}\label{e1.2.27}
y_{i_0}\in D^\prime_{\varepsilon_0}.
\end{equation}

Допустим противное, то есть
\begin{equation}\label{e1.2.28}
y_i\notin D^\prime_{\varepsilon_0} \quad \forall i\in K,\quad  i>0.
\end{equation}
Выделим из последовательности $\{y_i\}$, $i\in K$, $i>0$, сходящуюся подпоследовательность $\{y_i\}$, $i\in K^\prime\subset K$, и пусть $y^\prime$ -- ее предельная точка. Согласно п. 2 метода с учетом допущения (\ref{e1.2.28}) множество $Q_i$ для всех $i\in K$, $i>0$, имеет вид (\ref{e1.2.3}), и для тех же $i$ выполняется равенство (\ref{e1.2.4}). Тогда по лемме \ref{lemma1.2.3} имеем включение $y^\prime\in D$, и 
\begin{equation}\label{e1.2.28.1}
F(y^\prime)\le 0.
\end{equation}

C другой стороны, в силу (\ref{e1.2.28}) $F(y_i)>\varepsilon_0$ для всех $i\in K^\prime$. Переходя в последнем неравенстве к пределу по $i\in K^\prime$, получим $F(y^\prime)\ge\varepsilon_0>0$, что противоречит (\ref{e1.2.28.1}). Таким образом, существование номера $i_0>0$, для которого справедливо (\ref{e1.2.27}), доказано, и равенство (\ref{e1.2.5}) при $k=0$ имеет место.

2) Допустим теперь, что (\ref{e1.2.5}) выполняется при некотором фиксированном $k\ge0$, то есть $x_k=y_{i_k}$ при выбранном $k$. Покажем существование такого номера $i_{k+1}>i_k$, что
\begin{equation}\label{e1.2.29}
y_{i_{k+1}}\in D^\prime_{\varepsilon_{k+1}},
\end{equation}
тогда $x_{k+1}=y_{i_{k+1}}$, и лемма будет доказана.

Предположим противное, то есть
\begin{equation}\label{e1.2.29.1}
y_i\notin D^\prime_{\varepsilon_{k+1}} \quad \forall i\in K, \quad i>i_k.
\end{equation}
Выберем среди точек $y_i$, $i\in K$, $i>i_k$, сходящуюся подпоследовательность $\{y_i\}$, $i\in K^{\prime\prime}$. Пусть $y^{\prime\prime}$ -- ее предельная точка. Согласно предположению (\ref{e1.2.29.1}) для всех $i\in K$, $i>i_k$, выполняются равенства (\ref{e1.2.4}), а множества $Q_i$ задаются в виде (\ref{e1.2.3}). Значит, по лемме \ref{lemma1.2.3} имеет место неравенство
\begin{equation}\label{e1.2.29.2}
F(y^{\prime\prime})\le 0.
\end{equation}
Как и в первой части доказательства из (\ref{e1.2.29.1}) c учетом (\ref{e1.1.19}) легко получается неравенство $F(y^{\prime\prime})>0$, противоречащее (\ref{e1.2.29.2}). Таким образом, доказано существование номера $i_{k+1}$, для которого справедливо (\ref{e1.2.29}). Лемма доказана.

\end{proof}

Из леммы \ref{lemma1.2.4} следует, что в случае (\ref{e1.1.19}) для каждого $k\in K$ зафиксируется в виде (\ref{e1.2.5}) точка $x_k$, а значит, будет выбрана  согласно (\ref{e1.2.8}) и точка $z_k$.

\begin{theorem}\label{theorem1.2.3}
 Пусть последовательности $\{x_k\}$, $\{z_k\}$, $k\in K$, построены предложенным методом с условием, что числа $\varepsilon_k$, $k\in K$, выбраны согласно (\ref{e1.1.19}) и (\ref{e1.1.25}). Тогда любые предельные точки последовательностей $\{x_k\}$, $k\in K$, и $\{z_k\}$, $k\in K$, принадлежат множеству $X^*$.
\end{theorem}

\begin{proof} Напомним, что условие (\ref{e1.1.19}) гарантирует существование последовательностей $\{x_k\}$, $\{z_k\}$, $k\in K$. Отметим, что последовательность $\{x_k\}$, $k\in K$, ограничена в силу (\ref{e1.2.5}) и ограниченности $\{y_i\}$, $i\in K$. Пусть $\{x_k\}$, $k\in K_1\subset K$, и $\{z_k\}$, $k\in K_2\subset K$, -- любые сходящиеся подпоследовательности последовательностей $\{x_k\}$, $k\in K$, и $\{z_k\}$, $k\in K$, а $\bar{x}$ и $\bar{z}$, соответственно, -- их предельные точки. Покажем, что
\begin{equation}\label{e1.2.30}
\bar{x}\in X^*, \quad \bar{z}\in X^*, 
\end{equation}
тогда утверждение будет доказано.

Так как $y_{i_k}\in D^\prime_{\varepsilon_k}$, $k\in K$, то ввиду (\ref{e1.2.5}) $0<F(x_k)\le \varepsilon_k$, $k\in K$, и из (\ref{e1.1.25}) следует, что $\lim \limits_{k\in K} F(x_k)=0$. Тогда
$$\lim \limits_{k\in K} F(x_k)=\lim \limits_{k\in K_1} F(x_k) = F(\bar{x})=0,$$
и $\bar{x}\in D^\prime$. Кроме того, $x_k \in D^{\prime\prime}$ для всех $k\in K$, и в силу замкнутости множества $D^{\prime\prime}$ выполняется включение $\bar{x}\in D^{\prime\prime}$. Таким образом, $\bar{x}\in D$ и потому $f(\bar{x})\ge f^*$. С другой стороны, согласно (\ref{e1.2.1}), (\ref{e1.2.5}) $f(x_k)\le f^*$, $k\in K_1$, а значит, $f(\bar{x})\le f^*$. Следовательно, $f(\bar{x})=f^*$, и первое из включений (\ref{e1.2.30}) доказано.

Далее, последовательности $\{z_k\}$, $k\in K_2$, соответствует подпоследовательность $\{x_k\}$, $k\in K_2$. Выделим из $\{x_k\}$, $k\in K_2$, сходящуюся подпоследовательность $\{x_k\}$, $k\in K^\prime\subset K_2$, и пусть $x^\prime$ -- ее предельная точка. Заметим, что по доказанному $x^\prime \in X^*$, то есть 
\begin{equation}\label{e1.2.31}
f(x^\prime) = f^*.
\end{equation}
Перейдем теперь в неравенствах (\ref{e1.2.8}) к пределу по $k\to \infty$, $k\in K^\prime$, c учетом, что $z_k \to \bar{z}$, $k\in K^\prime$. Тогда $f(x^\prime)\le f(\bar{z})\le f^*$, и ввиду (\ref{e1.2.31}) $f(\bar{z})=f^*$. Cледовательно, второе из включений (\ref{e1.2.31}) также имеет место. Теорема доказана.
\end{proof}

\begin{theorem}\label{theorem1.2.4}
Пусть выполняются условия теоремы \ref{theorem1.2.3}. Если для точек $x_k$, $k\in K$, выполняются неравенства
\begin{equation}\label{e1.2.32}
f(x_{k+1})\ge f(x_k) \quad \forall k\in K, 
\end{equation}
то последовательность $\{x_k\}$, $k\in K$, сходится к множеству $X^*$. Если точки $z_k$, $k\in K$, выбраны с условием, что 
\begin{equation}\label{e1.2.33}
f(z_{k+1})\ge f(z_k) \quad \forall k\in K,
\end{equation}
то последовательность $\{z_k\}$, $k\in K$, сходится к множеству $X^*$.
\end{theorem}

\begin{proof} Пусть имеют место неравенства (\ref{e1.2.32}), (\ref{e1.2.33}). Тогда в силу ограниченности последовательностей $\{x_k\}$, $\{z_k\}$, $k\in K$, последовательности $\{f(x_k)\}$, $\{f(z_k)\}$, $k\in K$, сходятся, и по теореме \ref{theorem1.2.3} выполняются равенства $\lim \limits_{k\in K} f(x_k) = f^*$, $\lim \limits_{k\in K} f(z_k)=f^*$. Поэтому согласно вышеупомянутой теореме (\cite[62]{n_vasillev_f_p_1}) последовательности $\{x_k\}$, $\{z_k\}$, $k\in K$, сходятся к множеству $X^*$. Утверждение доказано.
\end{proof}

При исследовании сходимости метода использовалось условие ограниченности последовательностей $\{y_i\}$, $i\in K$, $\{z_k\}$, $k\in K$. Ясно, что выполнение этого условия можно обеспечить за счет соответствующего выбора множеств $M_0$ и $Q_i$.

Коротко обсудим далее возможности задания последовательности $\{\varepsilon_k\}$, $k\in K$.

\begin{remark}\label{r1.2.8}
Числа $\varepsilon_k$, $k\ge 1$, допустимо выбирать, как и значение $\varepsilon_0$, на предварительном шаге метода. Однако, в таком случае последовательность $\{\varepsilon_k\}$, $k\in K$, не будет адаптирована к процессу  минимизации. Поэтому в методе предусмотрена возможность задания чисел $\varepsilon_k$, $k\ge 1$, на шаге 3, то есть  в ходе построения приближений $x_k$, $z_k$. Приведем пример такого задания последовательности $\{\varepsilon_k\}$, $k\in K$.

Не выбирая конкретного значения, считая $\varepsilon_0$ настолько большим, что $y_0\in D^\prime_{\varepsilon_0}$. Тогда $x_0=y_0$, а точку $z_0$ можно выбрать, например, в виде $z_0=x_0$. Для всех $k\ge 0$ положим $\varepsilon_{k+1} = \gamma_k F(x_k)$ или $\varepsilon_{k+1} = \gamma_k F(z_k)$, где $0<\gamma_k<1$. Для такой последовательности $\{\varepsilon_k\}$, $k\in K$, справедливо (\ref{e1.1.19}), а при условии, что $\gamma_k \to 0$, $k\in K$, выполняется условие (\ref{e1.1.25}).
\end{remark}

Получим теперь для построенного метода при некоторых дополнительных условиях на исходную задачу оценки точности ее решения, а также оценки скорости сходимости последовательностей $\{x_k\}$, $\{z_k\}$, $k\in K$.

Будем предполагать далее в этом параграфе, что $D^{\prime\prime}={\rm R_n}$, то есть $D=D^\prime$. Пусть, кроме того, для функций $f(x)$, $f_j(x)$, $j\in J$, и множества $D$ выполняется предположение \ref{propose1.1.1}. Заметим, что в силу предположения \ref{propose1.1.1} и согласно теореме \ref{theorem1.1.3} решение задачи (\ref{e1.1.1}) единственно, то есть $X^*=\{x^*\}$.

Как и ранее будем считать, что множества $X_1^*(\delta)$, $X_2^*(\gamma)$ заданы в виде (\ref{e1.1.28}), (\ref{e1.1.29}) соответственно.

В силу предположения \ref{propose1.1.1} справедлива лемма \ref{lemma1.1.8}, которая  дает возможность оценивать на каждой итерации близость точек $u_i$, $i\in K$, построенных предложенным методом, к решению исходной задачи. А именно, поскольку для каждого $i\in K$ выполняются включения $u_i\in E^*$, $u_i\in D^\prime_{r_i}$, где $r_i=F(u_i)>0$, и при этом $u_i\ne x^*$, то по лемме \ref{lemma1.1.8} имеем $u_i\in X_1^*(\delta_i)$ или, что то же самое,
$$\|u_i - x^*\|\le \delta_i,$$
где $\delta_i=\sqrt{\frac{F(u_i)}{\mu}}$. Кроме того, если функция $f(x)$ удовлетворяет условию Липшица с константной $L$, то согласно той же лемме \ref{lemma1.1.8} справедливо включение $u_i\in X_1^*(\delta_i)\bigcap X_2^*(L\delta_i)$, а значит,
$$|f(u_i)-f^*|\le L\delta_i.$$

Приведем теперь оценки точности решения для последовательностей $\{x_k\}$, $\{z_k\}$, $k\in K$.

\begin{theorem}\label{theorem1.2.5}
 Для последовательностей $\{x_k\}$, $\{z_k\}$, $i\in K$, построенных предложенным методом с условием   (\ref{e1.1.19}), справедливы следующие три утверждения.

1) Имеет место оценка 
$$\|x_k-x^*\|\le \delta_k \quad \forall k\in K,$$
где $\delta_k=\sqrt{\frac{\varepsilon_k}{\mu}}$.

2) Если функция $f(x)$ удовлетворяет условию Липшица с константой $L$, то 
\begin{equation}\label{e1.2.34}
|f(x_k)-f^*|\le L\delta_k \quad \forall k\in K, 
\end{equation}
\begin{equation}\label{e1.2.35}
|f(z_k)-f^*|\le L\delta_k \quad \forall k\in K.
\end{equation}

3) Если точки $z_k$, $k\in K$, выбраны на шаге 3 метода с дополнительным условием
\begin{equation}\label{e1.2.36}
z_k \in D^\prime_{\varepsilon_k},
\end{equation}
то справедлива оценка 
\begin{equation}\label{e1.2.37}
\|z_k - x^*\|\le \delta_k \quad \forall k\in K.
\end{equation}
\end{theorem}

\begin{proof} Так как согласно пп. 1, 3 метода 
\begin{equation}\label{e1.2.38}
y_{i_k} \in E^*, \quad y_{i_k} \ne x^*, \quad y_{i_k}\in D^\prime_{\varepsilon_k}, \quad k\in K,
\end{equation}
то ввиду леммы \ref{lemma1.1.8} и равенства (\ref{e1.2.5}) для всех $k\in K$ выполняется включение $x_k\in X^*_1(\delta_k)$, где $\delta_k=\sqrt{\frac{\varepsilon_k}{\mu}}$. Следовательно, первое утверждение теоремы доказано.

Далее, при условии Липшица на функцию $f(x)$ по леммы \ref{lemma1.1.8} в силу (\ref{e1.2.38}), (\ref{e1.2.5}) для всех $k\in K$ справедливо включение $x_k\in X^*_2(\gamma_k)$, где $\gamma_k = L\delta_k$, а значит, и неравенство (\ref{e1.2.34}). Из (\ref{e1.2.34}) и (\ref{e1.2.8}) следует (\ref{e1.2.35}).

Пусть, наконец, для точки $z_k$ выполняется условие (\ref{e1.2.36}). Так как, кроме того, $z_k\in E^*$ и $z_k \ne x^*$, $k\in K$, то согласно той же лемме \ref{lemma1.1.8} $z_k \in X^*_1(\delta_k)$, то есть  имеет место оценка (\ref{e1.2.37}). Теорема доказана.
\end{proof}

Приведем при тех же дополнительных условиях на исходную задачу оценки скорости сходимости последовательностей $\{x_k\}$ и $\{z_k\}$, $k\in K$.

Как уже отмечалось, числа $\varepsilon_k$ для всех $k\in K$ допустимо задать на предварительном шаге метода. Пусть числа $\varepsilon_k$, $k\in K$, в методе заданы согласно (\ref{e1.1.40}). Тогда согласно теореме \ref{theorem1.2.5} для последовательности $\{x_k\}$, $k\in K$, справедлива оценка (\ref{e1.1.41}), а для последовательности $\{z_k\}$, $k\in K$, построенной с условием (\ref{e1.2.36}), имеет место оценка
$$\|z_k-x^*\|\le \frac{1}{\sqrt{\mu}k^{p/2}}, \quad  k\ge 1.$$
Кроме того, если для $f(x)$ выполняется условие Липшица, то в силу (\ref{e1.2.34}), (\ref{e1.2.35}) справедливы также оценки
$$|f(x_k) - f^*|\le \frac{L}{\sqrt{\mu}k^{p/2}}, \quad |f(z_k)-f^*|\le \frac{L}{\sqrt{\mu}k^{p/2}}, \quad k\ge 1.$$


\setcounter{equation}{0}

\section
{Метод проектирования точки, использующий аппроксимирующие множества}

В данном параграфе предлагается метод нахождения проекции точки на выпуклое множество. Приближения в алгоритме находятся путем проектирования точки на некоторые множества, которые задаются системами линейных неравенств и аппроксимируют исходное множество в окрестности решения. Метод характерен тем, что не требует вложения каждого из этих множеств в предыдущее. Как и в методах из предыдущих параграфов, такая особенность позволяет при построении аппроксимирующих множеств периодически освобождаться от части задающих их неравенств и тем самым существенно упрощать вспомогательные задачи отыскания итерационных точек.

Пусть множество $D^{\prime\prime}\subset {\rm R_n}$ выпукло и замкнуто, функции $f_j(x)$, $j\in J=\{1,\ldots, m\}$, выпуклы в ${\rm R_n}$, множество $D^\prime$ задано в виде (\ref{e1.2.0}),  и  $D=D^\prime \bigcap D^{\prime\prime}$. Решается задача отыскания проекции точки $y\in {\rm R_n}$ на множество $D$.

 Положим $p^*=Pr(y,D)$ -- проекция точки $y$ на множество $D$, $r^*=\|p^*-y\|$, $D_j = \{x\in {\rm R_n} : f_j(x)\le 0\}$. Будем считать, что множества  $W^1(x, D)$ и $K$ определены как и в $\S \ 1.1$.

Предлагаемый для решения поставленной задачи метод вырабатывает последовательность приближений $\{p_k\}$, $k\in K$, и заключается в следующем.

\textbf{Метод 1.3.} Выбираются точки $v^j\in {\rm int\,} D_j$ для всех $j\in J$. Задается числовая последовательность $\{\varepsilon_k\}$, $k\in K$, такая, что
$$\varepsilon_k>0, \quad k\in K, \quad \varepsilon_k \to 0, \quad k \to \infty.$$ Полагается $k=0$, $i=0$, $M_0={\rm R_n}$, $y_0=Pr(y, D^{\prime \prime})$.

\textbf{1.} Формируется множество $$J_i=\{j\in J\ : \ y_i\notin D_j\}.$$ Если $J_i=\emptyset$, то $y_i$ -- решение задачи, и процесс завершается.

\textbf{2.} Если $y_i\notin D^\prime_{\varepsilon_k}$, то полагается 
\begin{equation}\label{e1.3.1}
Q_i=M_i, 
\end{equation}
и следует переход к п. 3.  В противном случае полагается $i_k=i$, 
\begin{equation}\label{e1.3.2}
p_k=y_{i_k},
\end{equation}
\begin{equation}\label{e1.3.3}
Q_i=Q_{i_k}={\rm R_n},
\end{equation}
значение $k$ увеличивается на единицу, и следует переход к очередному пункту.

\textbf{3.} Для каждого $j\in J_i$ в интервале $(v^j,y_i)$ выбирается точка $z_i^j$ так, чтобы $z^j_i\notin{\rm int\,} D_j$ и при некотором $q_i^j\in [1,q]$, $1\le q<\infty$, для точки $$\bar{y}_i^j=y_i+q_i^j(z_i^j-y_i)$$ выполнялось включение $\bar{y}_i^j \in D_j$. Для всех $j\in J\setminus J_i$ полагается $z_i^j=\bar{y}_i^j=y_i$.

\textbf{4.} Отыскивается номер $j_i\in J_i$  такой, что $$\|y_i-z^{j_i}_i\|=\max \limits_{j\in J_i} \|y_i-z_i^j\|.$$

\textbf{5.} Выбирается вектор $a_i\in W^1(z_i^{j_i},D_{j_i})$, полагается $$M_{i+1}=Q_i\bigcap \{x\in {\rm R_n} \ : \ \langle a_i, x-z_i^{j_i}\rangle \le 0\}.$$

\textbf{6.} Отыскивается точка $$y_{i+1}=Pr(y,M_{i+1}\bigcap D^{\prime \prime}),$$ и следует переход к п. 1 при $i$, увеличенном на единицу.

\begin{remark}\label{r1.3.1}
В силу выбора множеств $M_0$, $Q_i$ и векторов  $a_i$ для всех $i\in K$ выполняется включение $D^\prime \subset M_i$,  а значит,
$$M_i\bigcap D^{\prime \prime}\ne \emptyset \quad \forall i\in K,$$
и построение точки $y_i$ согласно алгоритму возможно.
\end{remark}

\begin{remark}\label{r1.3.2}
Если $D^{\prime \prime}={\rm R_n}$ или множество $D^{\prime \prime}$ задано линейными равенствами или неравенствами, то на каждом шаге алгоритма решается задача проектирования точки $y$ на многогранное множество. Подчеркнем, что эта задача, как задача квадратичного программирования, может быть решена некоторыми известными методами (см., напр., \cite{n_vasillev_f_p_1,n_konnov_i_v_2}) за конечное число шагов.
\end{remark}

\begin{remark}\label{r1.3.3}
Если $\rm {int\,} D^\prime \ne \emptyset$ и известна точка $v\in \rm{int\,} D^\prime$, то в методе удобно положить $v^j=v$ для всех $j\in J$.
\end{remark}

\begin{remark}\label{r1.3.4}
При выполнении включения $y_i\in D^{\prime}_{\varepsilon_k}$ согласно (\ref{e1.3.3}) происходит полное отбрасывание всех накопившихся к шагу $i=i_k$ линейных неравенств, формирующих аппроксимирующие множества, и $M_{i+1}$ задается лишь одним неравенством $\langle a_i,x-z^{j_i}_i\rangle \le 0$. Это существенно упрощает решение  вспомогательных задач построения очередных приближений.
\end{remark}

Теперь перейдем к исследованию сходимости метода. Следующее утверждение показывает, что построенные согласно методу точки $y_i$, $i\in K$, отличаются от $y$ для всех $i>0$.

\begin{lemma}\label{lemma1.3.1}
Для всех $i\ge 1$ выполняется соотношение 
\begin{equation}\label{e1.3.4}
y\notin M_i\bigcap D^{\prime\prime}. 
\end{equation}
\end{lemma}

\begin{proof} Покажем, что 
\begin{equation}\label{e1.3.5}
y\notin M_i \quad \forall i\ge 1,
\end{equation}
тогда (\ref{e1.3.4}) будет доказано.

Пусть $i=1$. Тогда $M_1=Q_0\bigcap\{x\in {\rm R_n}\ : \ \langle a_0, x-z^{j_0}_0 \rangle \le 0\}$. По построению $Q_0={\rm R_n}$ независимо от принадлежности точки $y_0$ множеству $D^{\prime}_{\varepsilon_0}$. Поэтому $M_1$ задается неравенством $\langle a_0,x-z^{j_0}_0\rangle \le 0$, и по выбору вектора $a_0$ имеем $y\notin M_1$, то есть (\ref{e1.3.5}) справедливо при $i=1$.

Пусть теперь утверждение (\ref{e1.3.5}) выполняется при $i=l$, где $l\ge 1$. Покажем, что 
\begin{equation}\label{e1.3.6}
y\notin M_{l+1},
\end{equation}
 тогда утверждение (\ref{e1.3.5}) будет доказано. По построению
 $$M_{l+1}=Q_l \bigcap \{x\in {\rm  R_n} \ : \ \langle a_l, x-z^{j_l}_l\le 0\}.$$ 
Но $Q_l=M_l$ или $Q_l={\rm R_n}$. Если $Q_l={\rm R_n}$, то $M_{l+1}$ задается неравенством $\langle a_l,x-z^{j_l}_l\rangle \le 0$ и согласно выбору $a_l$ имеем (\ref{e1.3.6}). Пусть $Q_l=M_l$. По индукционному предположению $y\notin M_l$. Тогда $y\notin Q_l$, и, следовательно, выполняется (\ref{e1.3.6}). Утверждение доказано.
\end{proof}

Критерий остановки, заложенный в п. 1 алгоритма, обосновывает следующая
\begin{lemma}\label{lemma1.3.2}
Пусть точка $y_i$, построенная методом, такова, что $J_i=\emptyset$. Тогда $y_i=p^*$.
\end{lemma}

\begin{proof}
Поскольку $D^{\prime}\subset M_i$, то $D=D^\prime \bigcap D^{\prime \prime}\subset M_i\bigcap D^{\prime \prime}$. Кроме того, $y_i=Pr(y,M_i\bigcap D^{\prime \prime})$ для всех $i\in K$. Следовательно, $$\|y_i-y\|=\min \limits_{x\in M_i\bigcap D^{\prime \prime}}\|x-y\|\le \min \limits_{x\in D} \|x-y\|=\|p^*-y\|.$$ По условию леммы $y_i\in D^\prime$, а значит, $y_i\in D$, и $\|y_i-y\|\ge \|p^*-y\|$. Отсюда и из предыдущего неравенства следует утверждение леммы.
\end{proof}

\begin{lemma} \label{lemma1.3.3}
Последовательность $\{y_i\}$, $i\in K$, построенная методом, ограничена.
\end{lemma}

\begin{proof} По построению $\|y_i-y\|\le \|x-y\|$ для всех $x\in M_i\bigcap D^{\prime \prime}$, $i\in K$. Поскольку $D\subset M_i\bigcap D^{\prime \prime}$ для всех $i\in K$, то $p^*\in M_i \bigcap D^{\prime \prime}$ для каждого $i\in K$, а значит,
$$\|y_i-y\|\le \|p^*-y\|=r^* \quad \forall i\in K.$$ 
Лемма доказана.
\end{proof}

\begin{lemma}\label{lemma1.3.4}
Пусть последовательность $\{y_i\}$, $i\in K$, построена методом. Тогда для каждого $k\in K$ существует такой номер $i=i_k$, что выполняется (\ref{e1.3.2}).
\end{lemma}

\begin{proof} Пусть $k=0$. Будем считать, что $y_0\notin D^\prime_{\varepsilon_0}$, так как в противном случае равенство (\ref{e1.3.2}) имеет место при $k=0$. Докажем существование номера $i=i_k>0$,  для которого $y_{i_0}\in D^\prime_{\varepsilon_0}$. 

Допустим противное, то есть $y_i\notin D^\prime_{\varepsilon_0}$  для всех $i\in K$, $i>0$. Выделим из последовательности $\{y_i\}$, $i\in K$, $i>0$, сходящуюся подпоследовательность $\{y_i\}$, $i\in K^\prime \subset K$, и пусть $y^\prime$ -- ее предельная точка. Согласно п. 3 метода $Q_i$ для всех $i\in K^\prime$ имеет вид (\ref{e1.3.1}). По методике, используемой при обосновании леммы \ref{lemma1.1.2}, доказывается, что $y^\prime \in D^\prime $ и, следовательно, $F(y^\prime)\le 0$.  C другой стороны, по предположению $F(y_i)>\varepsilon_0$ при всех $k\in K^\prime$, а потому $F(y^\prime)>0$. Полученное противоречие доказывает справедливость равенства (\ref{e1.3.2}) при $k=0$.

Аналогично показывается, что, если $p_k=y_{i_k}$ при некотором $k\ge 0$, то существует такой номер $i_{k+1}>i_k$, что $y_{i_{k+1}}\in D^\prime_{\varepsilon_{k+1}}$ и $p_{k+1}=y_{i_{k+1}}$. Лемма доказана.
\end{proof}
\begin{theorem}\label{theorem1.3.1}
Пусть последовательность $\{p_k\}$, $k\in K$, построена предложенным методом. Тогда любая ее предельная точка совпадает с $p^*$, а если для всех $k\in K$ выполняются неравенства 
\begin{equation}\label{e1.3.7}
\|p_{k+1}-y\|\ge \|p_k-y\|,
\end{equation}
то вся последовательность $\{p_k\}$, $k\in K$, сходится к точке $p^*$.
\end{theorem}
\begin{proof}
Пусть $\{p_k\}$, $k\in K^\prime \subset K$, -- сходящаяся подпоследовательность последовательности $\{p_k\}$, $k\in K$, и $\bar{p}$ -- ее предельная точка. Так как $0<F(p_k)\le \varepsilon_k$, $k\in K$, то с учетом выбора последовательности $\{\varepsilon_k\}$, $k\in K$, имеем $\lim \limits_{k\in K}F(p_k)=\lim \limits_{k\in K^\prime}F(p_k)=F(\bar{p})=0$, и $\bar{p}\in D^\prime$. Кроме того, $p_k\in D^{\prime \prime}$, $k\in K^\prime$, и в силу замкнутости $D^{\prime \prime}$ выполняется включение $\bar{p}\in D^{\prime\prime}$. Таким образом, $\bar{p}\in D$ и $\|\bar{p}-y\|\ge r^*$. С другой стороны, поскольку $p_k=Pr(y,M_{i_k}\bigcap D^{\prime \prime})$, а $D\subset M_{i_k}\bigcap D^{\prime \prime}$, то $\|p_k-y\|\le r^*$, и значит, $\|\bar{p}-y\|\le r^*$. Отсюда следует, что $\|\bar{p}-y\|=r^*=\|p^*-y\|$, и первая часть теоремы доказана.

Пусть теперь для последовательности $\{p_k\}$, $k\in K$, выполняется условие (\ref{e1.3.7}). Тогда ввиду ограниченности $\{p_k\}$, $k\in K$, последовательность $\{\|p_k-y\|\}$, $k\in K$, является сходящейся, и по доказанному выше $\lim \limits_{k\in K} \|p_k-y\|=r^*$. Поэтому в силу известной теоремы \cite[62]{n_vasillev_f_p_1} второе утверждение тоже доказано.
\end{proof}

\begin{remark}\label{r1.3.5}
Условие (\ref{e1.3.7}) для построенной методом последовательности $\{p_k\}$, $k\in K$, имеет место, если, например, для всех $k\in K$ выполняются включения $M_{i_{k+1}}\subset M_{i_k}$. Для выделения из последовательности $\{p_k\}$, $k\in K$, сходящейся к решению подпоследовательности $\{p_{k_l}\}$, $l\in K$, достаточно номера $k_l$ выбрать из условия $\|p_{k_{l+1}}-y\|\ge \|p_{k_l}-y\|$ при всех $l\in K$. 
\end{remark}
          
\setcounter{equation}{0}

\section
{Метод  отсечений, допускающий параллельные вычисления}

В настоящем параграфе предлагается метод отсечений для  решения задачи (\ref{e1.1.1}), позволяющий на каждом шаге распараллеливать процесс вычислений при отыскании очередной итерационной точки. 

Пусть функции $f(x)$, $f_j(x)$, $j\in J$, и множества $D_j$ удовлетворяют условиям $\S \, 1.1$, а множество $D$ имеет вид (\ref{e1.1.0}). Решается задача (\ref{e1.1.1}). Будем пользоваться обозначениями  $\S \, 1.1$.

Предлагаемый метод  вырабатывает последовательность приближений $x_k$, $k\in K$, и заключается в следующем.

\textbf{Метод 1.4.}  Строится выпуклое замкнутое множество $M_0 \subset {\rm R_n}$, содержащее хотя бы одну точку множества $X^*$. Выбираются точки $y^j \in\rm{ int\,} D_j$ для всех $j\in J$, находится точка $x_0 \in M_0$ такая, что $f(x_0)\le f^*$. Пусть уже построена точка $x_k$, $k\in K$.

\textbf{1.} Если
\begin{equation}\label{e1.4.1}
x_k\in D,
\end{equation}
 то процесс заканчивается.

\textbf{2.} Полагается $J_k=\{j \in J \ : \ x_k\notin D_j\}.$

\textbf{3.} Для каждого $j\in J_k$ в интервале $(y^j, x_k)$ выбирается точка $z_k^j \notin \rm{int\,} D_j$ так, чтобы при некотором $1\le q_k^j\le q<\infty$  выполнялось включение 
\begin{equation}\label{e1.4.2}
x_k+q_k^j(z_k^j-x_k)\in D_j,
\end{equation}
и полагается 
\begin{equation}\label{e1.4.3}
\bar{y}_k^j=x_k+q_k^j(z_k^j-x_k).
\end{equation}

\textbf{4.} Для каждого $j\in J_k$ выбирается конечное множество $A_k^j\subset W^1(z_k^j,D_j)$, полагается 
$$M_k^j=M_k\bigcap \{x\in {\rm  R_n} \ : \ \langle a, x-z_k^j\rangle \le 0 \ \forall a\in A_k^j\},$$
 и находится точка $u_k^j$ такая, что
\begin{equation}\label{e1.4.4}
u_k^j\in M_k^j, \quad f(u_k^j) \le f^*.
\end{equation}

\textbf{5.} Отыскивается номер $j_k\in J_k$, удовлетворяющий условию
\begin{equation}\label{e1.4.5}
\|x_k-z_k^{j_k}\|=\max \limits_{j\in J_k} \|x_k-z_k^j\|.
\end{equation}

\textbf{6.} Определяется множество 
$$J^\prime_k=\{j\in J_k \ : \ u_k^j \in M_k^{j_k}\}.$$ Находится номер $l_k\in J^\prime_k$ такой, что $f(u_k^{l_k})=\max \limits_{j\in J^\prime_k}f(u_k^j).$

\textbf{7.} Полагается 
$$M_{k+1}=M_k^{j_k}\bigcap \{x\in {\rm R_n} \ : \ \langle a, x-z_k^{l_k}\rangle \le 0 \ \forall a\in A_k^{l_k}\},$$
\begin{equation}\label{e1.4.6}
x_{k+1}=u_k^{l_k}.
\end{equation}

Сделаем некоторые замечания относительно предложенного метода.

\begin{remark}\label{r1.4.1}
 Если на некотором шаге выполняется включение (\ref{e1.4.1}), то $x_k\in X^*$. Если же в результате работы метода выработана последовательность $\{x_k\}$, $k\in K$, то, как показано ниже, любая предельная точка принадлежит $X^*$. \mbox{На начальном} шаге метода можно положить, в частности, $M_0=\bigcap \limits_{j\in J^\prime}D_j$, где $J^\prime \subset J$. Так как $x_k\in M_0$ для всех $k\in K$, то точки $y^j\in {\rm  int\,} D_j$, $j\in J^\prime$, ни на одной итерации не понадобятся, и в этом случае нет необходимости в их  задании. Такой способ выбора множества $M_0$ удобен, например, когда $\bigcap \limits_{j\in J^\prime} D_j$ -- выпуклый многогранник. При достижении функцией $f(x)$ в ${\rm R_n}$ своего минимального значения можно считать, например, что
 $$M_0={\rm R_n}, \quad x_0={\rm argmin\,} \{f(x) : x\in {\rm R_n}\}.$$
\end{remark}

\begin{remark}\label{r1.4.2}
В случае когда $\rm{int\,} D\neq \emptyset$ и известна точка $y\in\rm{int\,} D$, в методе удобно задать $y^j=y$ для всех $ j\in J$.

Если в (\ref{e1.4.2}) положить $q_k^j=1$, то $z_k^j=\bar{y}_k^j$,  и $z_k^j$ является точкой пересечения отрезка $[x_k, y^j]$ с границей множества $D_j$. Однако условие задания чисел $q_k^j$ позволяет искать указанную точку пересечения приближенно.
\end{remark}

\begin{remark}\label{r1.4.3}
Отметим, что условие (\ref{e1.4.4}) дает возможность отыскивать $u_k^j$ при каждом $j \in J_k$ как точку приближенного минимума $f(x)$ на множестве $M_k^j$. В случае, если для каждого $k\in K$ точки $u_k^j \in M_k^j, \ j \in J_k$, находятся из условия
\begin{equation}\label{e1.4.7}
f(u_k^j)=\min \{f(x) :  x\in M_k^j\},
\end{equation}
то, как доказано ниже, вся последовательность $\{x_k\}$, $k\in K$, сходится к множемтву $X^*$.
\end{remark}

\begin{remark}\label{r1.4.4}
Подчеркнем также, что $J^\prime_k\neq \emptyset$ для всех $k\in K$, так как, по крайней мере, номер $j_k$ <<самого глубокого>> \ отсечения содержится в $J^\prime_k$ согласно включению $u_k^{j_k} \in M_k^{j_k}$. \mbox{Так как} $u_k^{l_k}\in M_k^{j_k}$ и $u_k^{l_k} \in M_k^{l_k}$, то ввиду (\ref{e1.4.6}) $x_{k+1}\in M_{k+1}$.
\end{remark}

Сделаем принципиальное замечание, касающееся проведения параллельных вычислений в методе.
\begin{remark}\label{r1.4.5}
Заметим, наконец, что на каждой итерации задачи отыскания точек $z_k^j$, а затем и точек $u_k^j$ можно решать для всех $j\in J_k$ одновременно и независимо друг от друга. Таким образом, метод действительно позволяет для каждого $k\in K$ проводить на шаге 3 и шаге 4 параллельные вычисления при нахождении точки $x_{k+1}$.
\end{remark}

Перейдем к обоснованию сходимости метода. Будем далее считывать, что последовательность $\{x_k\}$, $k\in K$, построена вышепредложенным  методом является  ограниченной.

\begin{lemma}\label{lemma1.4.1}
Пусть $\{x_k\}$, $k\in K^\prime \subset K$, -- сходящаяся подпоследовательность последовательности $\{x_k\}$, $k\in K$, и номер $r\in J$ таков, что множество 
$K_r=\{k\in K^\prime \ : \ r=j_k\},$ где $j_k$ определен в (\ref{e1.4.5}), состоит из бесконечного числа элементов. Тогда
\begin{equation}\label{e1.4.8}
\lim \limits_{k\in K_r}\|z_k^j-x_k\|=0 \quad \forall j\in J.
\end{equation}
\end{lemma}

\begin{proof}
Согласно п. 3 метода для всех $k\in K$ и $j\in J$ существует такое число $\gamma_k^j\in [0,1)$, что
\begin{equation}\label{e1.4.9}
z_k^j=x_k+\gamma_k^j(y^j-x_k),
\end{equation}
причем $\gamma_k^r>0$ для всех  $k\in K_r$.

Для произвольного $k\in K_r$ зафиксируем номер $p_k\in K_r$ такой, что $p_k>k$. Поскольку  $x_{p_k}\in M_{p_k} , M_{p_k}\subset M_{k+1}, \ M_{k+1} \subset M_k^r,$ а 
$$M_k^r=M_k\bigcap \{x\in {\rm R_n} \ : \ \langle a, x-z_k^r\rangle \le 0 \ \forall a\in A_k^r\},$$
то $A_k^r\subset W^1(z_k^r, M_{p_k}),$ а значит, $\langle a, x_{p_k}-z_k^r\rangle \le0$ для всех $a\in A_k^r$. Отсюда с учетом равенства (\ref{e1.4.9}) при $j=r$ для всех $a\in A_k^r$ имеем 
$$\langle a,x_k-x_{p_k}\rangle \ge \gamma_k^r\langle a,x_k-y^r\rangle.$$
Далее, учитывая включение $y^r \in {\rm int\,} D_r$, а также условия выбора последовательности $\{x_k\}, \ k\in K_r$, и векторов множества $A_k^r$, нетрудно доказать существование такого числа $\delta_r >0$, что $\langle a,x_k-y^r\rangle \ge \delta_r$ для всех $k\in K_r$, $a\in A_k^r$, то
\begin{equation}\label{e1.4.10}
\|x_k - x_{p_k}\|\ge \gamma_k^r\delta_r \quad \forall k, p_k\in K_r, \quad p_k>k.
\end{equation}
 Так как по условию леммы последовательность $\{x_k\}$, $k\in K_r$, является сходящейся , то согласно (\ref{e1.4.10}) $\gamma_k^r \to 0$, $k\to \infty$, $k\in K_r.$ Поэтому из (\ref{e1.4.9}) при $j=r$ с учетом ограниченности величин $\|y^r-x_k\|, \ k \in K_r$, следует равенство
\begin{equation}\label{e1.4.11}
\lim \limits_{k\in K_r} \|z_k^r-x_k\|=0.
\end{equation}
Наконец, в силу условия (\ref{e1.4.5}) выбора номеров $j_k$, $k \in K_r$, при каждом $k\in K_r$ выполняется неравенство $\|z_k^r-x_k\|\ge \|z_k^j-x_k\|$ для всех $j\in J,$ из которого с учетом (\ref{e1.4.11}) следует (\ref{e1.4.8}). Лемма доказана. 
\end{proof}

\begin{theorem}\label{theorem1.4.1}
Любая предельная точка последовательности $\{x_k\}$, $k\in K$, принадлежит множеству $X^*$.
\end{theorem}

\begin{proof}
 Пусть $\{x_k\}$, $k\in K^\prime \subset K$, -- любая сходящаяся подпоследовательность последовательности $\{x_k\}$, $k\in K$,  и $\bar{x}$ -- ее предельная точка. Покажем сначала, что 
\begin{equation}\label{e1.4.12}
\bar{x}\in D.
\end{equation}

Пусть, как и в лемме \ref{lemma1.4.1}, номер $r\in J$ выбран так, что множество $K_r=\{k\in K^\prime  :  j_k=r\}$ бесконечно. Отметим, что последовательность $\{x_k\}$, $k\in K_r$, при каждом $j\in J$ соответствуют последовательности $\{z_k^j\}$, $k\in K_r$, и $\{\bar{y}_k^j\}$, $k\in K_r$. При этом для $\{z_k^j\}$, $k\in K_r$, справедливо (\ref{e1.4.8}), а для точек $\bar{y}_k^j$, $k\in K_r$, согласно п. 3 метода выполняется либо равенство (\ref{e1.4.4}), либо равенство $\bar{y}_k^j=x_k$. Тогда из ограниченности последовательности $\{x_k\}$, $k\in K_r$, следует, что последовательности $\{\bar{y}_k^j\}$, $k\in K_r$, также ограничены для всех $j\in J$.

Выделим теперь для каждого $j\in J$ из $\{\bar{y}_k^j\}$, $k\in K_r$, сходящуюся подпоследовательность $\{\bar{y}_k^j\}$, $k\in K_r^j\subset K_r$, и пусть $v_j$ -- предельная точка. Заметим, что $v_j\in D_j$ для всех $j\in J$ в силу замкнутости множеств $D_j$. Положим для каждого $j\in J$ 
$$P_1^j=\{k\in K^j_r \ : \ j\in J_k\}, \quad P_2^j=K^j_r\setminus P^j_1.$$

Хотя бы одно из множеств $P^j_1$ или $P^j_2$ для каждого $j\in J$ состоит из бесконечного числа номеров. Перейдем теперь при каждом фиксированном $j\in J$ к пределу в равенствах (\ref{e1.4.3}) по $k\to \infty$, $k\in P_1^j$, с учетом (\ref{e1.4.8}), если множество $P_1^j$  бесконечно, или в равенствах $\bar{y}_k^j=x_k$ по $k\to \infty$, $k\in P_2^j$, если $P_2^j$ бесконечно. Тогда получим равенства $\bar{x}=v_j$ для всех $j\in J$, из которых следует включение (\ref{e1.4.12}).

Далее, согласно (\ref{e1.4.12}) выполняется неравенство $f(\bar{x})\ge f^*$. C другой стороны, ввиду второго из условий (\ref{e1.4.4}) и равенства (\ref{e1.4.6}) $f(x_k) \le f^*$ для всех $k\in K$, и потому $\lim \limits_{k\in K^\prime}f(x_k)=f(\bar{x})\le f^*$. Таким образом, $f(\bar{x})=f^*$, и утверждение теоремы доказано. 
\end{proof}

\begin{theorem}\label{theorem1.4.2}
Пусть для последовательности $\{x_k\}$, $k\in K$, приближение $x_0$ является точкой минимума функции $f(x)$  на множестве $M_0$, а точки $u_k^j\in M_k^j$, выбраны для всех $k\in K$ и $j\in J_k$  из условия (\ref{e1.4.7}). Тогда последовательность $\{x_k\}$, $k\in K$, сходится к множеству $X^*$.
\end{theorem}

\begin{proof} Для каждого $k\in K$ согласно пунктам 4, 7 метода выполняется включение $M_{k+1}\subset M_k^{l_k}$, а согласно п. 6 -- включение $u_k^{l_k}\in M_{k+1}$. Кроме того, по условию теоремы $u_l^{l_k}$ -- точка минимума функции $f(x)$ на множестве $M_k^{l_k}$ для каждого $k\in K$. Следовательно, $u_k^{l_k}$ -- точка минимума функции $f(x)$ на множестве $M_{k+1}$. Таким образом, с учетом (\ref{e1.4.6}) и условия выбора точки $x_0$ имеем
\begin{equation}\label{e1.4.13}
x_k={\rm argmin\,} \{f(x) : x\in M_k\} \quad \forall k \in K.
\end{equation}

Далее, для каждого $k\in K$ справедливо включение $M_{k+1}\subset M_k$. Тогда ввиду (\ref{e1.4.13}) $f(x_{k+1})\ge f(x_k)$, $k\in K$. Значит, в силу ограниченности $\{x_k\}$, $k\in K$, последовательность $\{f(x_k)\}$, $k\in K$, сходится, и по предыдущей теореме $\lim \limits_{k\in K}f(x_k)=f^*.$ Отсюда  и из известного утверждения \cite[62]{n_vasillev_f_p_1} следует, что вся последовательность $\{x_k\}$, $k\in K$, сходится к множеству решений задачи (\ref{e1.1.1}). Теорема доказана.
\end{proof}

Утверждение этой теоремы останется в силе, если для каждого $k\in K$ точка $u_k^{l_k}$ выбрана из условия $$f(u_k^{l_k})=\min \{f(x) : x\in M_k^{l_k}\},$$
 а точки $u^j_k$, $j\in J_k\setminus\{l_k\},$ -- из менее жесткого условия (\ref{e1.4.4}).

\setcounter{equation}{0}

\section
{Численные эксперименты}

В этом параграфе приводятся и обсуждаются результаты тестирования методов 1.1 и 1.2. Оба метода программно реализованы на языке программирования С++ в среде Microsoft Visual Studio 2008. Численные эксперименты осуществлялись на компьютере с параметрами OS Windows XP SP3 AMD Athlon(tm) 64 X2 Dual Core Processor 6000+ 3,01 ГГц 2,00 ГБ ОЗУ.

Цели проведенных численных экспериментов состояли в следующем. Требовалось проверить работоспособность каждого из названных методов и выработать рекомендации по применению в этих методах тех или иных процедур обновления аппроксимирующих множеств за счет выбора вспомогательных множеств $Q_i$. Кроме того, предполагалось дать некоторые рекомендации по заданию значений $\varepsilon_k$, которые определяют качество аппроксимации области ограничений погружающими множествами.

Для выработки указанных рекомендаций проведено значительное число экспериментов на тестовых примерах с различным количеством переменных (до 50) и ограничений (до 50). Поскольку каждый тестовый пример решался обоими методами, то тем самым проведено и сравнение методов между собой.

Для выполнения численных экспериментов была выбрана следующая тестовая задача:
$$ f(x) = \langle c,x\rangle \to \min$$
$$f_j(x) = \sum_{i= 1}^n k_{ij} \xi_i^2 -  1\le 0,\quad j=1,\ldots, n,$$
где $n\ge 1$, $c=(-1,-1,\ldots, -1)\in {\rm R_n}$, $x=(\xi_1, \xi_2,\ldots,\xi_n)\in {\rm R_n}$, $k_{ij}=\frac{1}{ij}$ для всех $i=1,\ldots,n$, $j=1,\ldots,n$.

Для начала приведем критерий, на основании которого происходила остановка работы алгоритмов. Процесс отыскания решения задачи прекращался, если на некотором шаге $i^\prime\in K$ для метода 1.1 выполнялось условие $y_{i^\prime}\in D_\varepsilon$, а для метода 1.2 $-$ $y_{i^\prime}\in D^\prime_\varepsilon$, где 
$$\varepsilon = 0.00001.$$
Согласно лемме \ref{lemma1.1.8} выполнение вышеуказанных условий остановки означает справедливость следующих неравенств:

$$\|y_{i^\prime} - x^*\|\le \sqrt{\frac{\varepsilon}{\mu}}, \quad |\langle c, y_{i^\prime} - x^*\rangle|\le \|c\| \sqrt{\frac{\varepsilon}{\mu}}.$$

При решении задач в методе 1.1 точки $v^j$, $j\in J$, задавались в виде $v^j=v=0$, $j\in J$,  для всех $i$ полагалось $H_i=J_i$, $q_i^j=1$, $j\in J_i$. Множество $T_i$ для всех $i$  в методе 1.2 задавалось в виде (\ref{e1.2.13}), а точки $u_i=y_i$, $z_k=x_k$ для всех $i$, $k\in K$. Наконец, множества $G_i$ для всех $i$ выбирались совпадающими с ${\rm R_n}$, множество $M_0$ в обоих методах задавалось в виде куба, содержащего всю допустимую область.

Результаты численных экспериментов представлены в следующей таблице:

\begin{center}
\begin{tabular}{|>{\centering}p{2cm}|>{\centering}p{2cm}|c|c|c|>{\centering}p{4.5cm}|>{\centering}p{3.5cm}|}
\hline
Номер & Метод & $n$  & $\varepsilon_k$ & $Q_i$ &Кол-во  итераций&  Время в мин. \tabularnewline
\hline
1 &1.1 & 10 &$a_1$ &$ b_1$ & 185 & 0,016\tabularnewline
\hline
2 &1.1 & 20 &$a_1$ &$ b_1$ & 618 & 0,1839\tabularnewline
\hline
3 &1.1 & 30 &$a_1$ &$ b_1$ & 1286 & 1,81\tabularnewline
\hline
4 & 1.1 & 40 &$a_1$ &$ b_1$ & 2029 & 9,85\tabularnewline
\hline
5 & 1.1 & 50 &$a_1$ &$ b_1$ & 3282 & 55,56\tabularnewline
\hline
6 & 1.1 & 10 &$a_2$ &$ b_2$ & 305 &0,046\tabularnewline
\hline
7 & 1.1 & 20 &$a_2$ &$ b_2$ & 1479 &0,18\tabularnewline
\hline
8 & 1.1 & 30 &$a_2$ &$ b_2$ & 3176 &0,41\tabularnewline
\hline
9 &1.1 & 40 &$a_2$ &$ b_2$ & 6196 &0,96\tabularnewline
\hline
10 & 1.1 & 50 &$a_2$ &$ b_2$ & 9384 &1,79\tabularnewline
\hline
11 & 1.1 & 10 &$a_3$ &$ b_2$ & 232 &0.02\tabularnewline
\hline
12 & 1.1 & 20 &$a_3$ &$ b_2$ &971  &0,11\tabularnewline
\hline
13 & 1.1 & 30 &$a_3$ &$ b_2$ & 1999 &0,27\tabularnewline
\hline
\end{tabular}
\end{center} 

\begin{center}
\begin{tabular}{|>{\centering}p{2cm}|>{\centering}p{2cm}|c|c|c|>{\centering}p{4.5cm}|>{\centering}p{3.5cm}|}
\hline
14 & 1.1 & 40 &$a_3$ &$ b_2$ & 3444 &0,59\tabularnewline
\hline
15 & 1.1 & 50 &$a_3$ &$ b_2$ & 5281 &1,24\tabularnewline
\hline
16 &1.1 & 10 &$a_4$ &$ b_2$ & 236  &0,02\tabularnewline
\hline
17 & 1.1 & 20 &$a_4$ &$ b_2$ & 827 &0,1\tabularnewline
\hline
18 & 1.1 & 30 &$a_4$ &$ b_2$ &1658  &0,24\tabularnewline
\hline
19 & 1.1 & 40 &$a_4$ &$ b_2$ &2871  &0,54\tabularnewline
\hline
20 & 1.1 & 50 &$a_4$ &$ b_2$ &4231  &1,12\tabularnewline
\hline
21 & 1.1 & 10 &$a_5$ &$ b_2$ &197  &0,02\tabularnewline
\hline
22 & 1.1 & 20 &$a_5$ &$ b_2$ &660  &0,08\tabularnewline
\hline
23 & 1.1 &30 &$a_5$ &$ b_2$ & 1360 & 0,24\tabularnewline
\hline
24 & 1.1 & 40 &$a_5$ &$ b_2$ &2139  & 0,78\tabularnewline
\hline
25 & 1.1 & 50 &$a_5$ &$ b_2$ & 3326  & 2,1\tabularnewline
\hline
26 & 1.1 & 10 &$a_6$ &$ b_2$ &192  &0,02\tabularnewline
\hline
27 & 1.1 & 20 &$a_6$ &$ b_2$ & 694 &0,08\tabularnewline
\hline
28 & 1.1 & 30 &$a_6$ &$ b_2$ &1443  &0,27\tabularnewline
\hline
29 & 1.1 & 40 &$a_6$ &$ b_2$ & 2481  &0,63\tabularnewline
\hline
30 & 1.1 & 50 &$a_6$ &$ b_2$ & 3620 &1,2\tabularnewline
\hline
 31 & 1.1 & 10 &$a_7$ &$ b_2$ & 180 & 0,01\tabularnewline
\hline
32 & 1.1 & 20 &$a_7$ &$ b_2$ & 620 & 0,07\tabularnewline
\hline
33 & 1.1 & 30 &$a_7$ &$ b_2$ & 1336 & 0,26\tabularnewline
\hline
34 & 1.1 & 40 &$a_7$ &$ b_2$ &2207  &0,75\tabularnewline
\hline
35 & 1.1 & 50 &$a_7$ &$ b_2$ & 3297 &2,45\tabularnewline
\hline
36 & 1.1 & 10 &$a_2$ &$ b_3$ &680  &0,69\tabularnewline
\hline
37 & 1.1 & 20 &$a_2$ &$ b_3$ &16689  &2,3\tabularnewline
\hline
38 & 1.1 & 30 &$a_2$ &$ b_3$ &  30390&7,33\tabularnewline
\hline
39 & 1.1 & 40 &$a_2$ &$ b_3$ & 48721 & 26,85\tabularnewline
\hline
40 & 1.1 & 50 &$a_2$ &$ b_3$ & 70194 & 91,42\tabularnewline
\hline
41  & 1.1 & 10 &$a_3$ &$ b_3$ &1502  &0,14\tabularnewline
\hline
42 & 1.1 & 20 &$a_3$ &$ b_3$ &3875  &0,5\tabularnewline
\hline
43 & 1.1 & 30 &$a_3$ &$ b_3$ & 6324 &1,4\tabularnewline
\hline
44 & 1.1 & 40 &$a_3$ &$ b_3$ &  10035& 4,55\tabularnewline
\hline
\end{tabular}
\end{center} 

\begin{center}
\begin{tabular}{|>{\centering}p{2cm}|>{\centering}p{2cm}|c|c|c|>{\centering}p{4.5cm}|>{\centering}p{3.5cm}|}
\hline
45 & 1.1 & 50 &$a_3$ &$ b_3$ & 14685 & 15,07\tabularnewline
\hline
46 & 1.1 & 10 &$a_4$ &$ b_3$ & 745 &0,07\tabularnewline
\hline
47 & 1.1 & 20 &$a_4$ &$ b_3$ & 1946 &0,24\tabularnewline
\hline
48 & 1.1 & 30 &$a_4$ &$ b_3$ & 3946 &0,82\tabularnewline
\hline
49 & 1.1 & 40 &$a_4$ &$ b_3$ &  6249& 2,79\tabularnewline
\hline
50 & 1.1 & 50 &$a_4$ &$ b_3$ & 9864  & 10,4\tabularnewline
\hline
51 & 1.1 & 10 &$a_5$ &$ b_3$ & 375 & 0,04\tabularnewline
\hline
52 & 1.1 & 20 &$a_5$ &$ b_3$ & 920 &0,11\tabularnewline
\hline
53 & 1.1 & 30 &$a_5$ &$ b_3$ &1678  &0,45\tabularnewline
\hline
54 & 1.1 & 40 &$a_5$ &$ b_3$ & 2704 & 1,86\tabularnewline
\hline
55 & 1.1 & 50 &$a_5$ &$ b_3$ & 4303 &8,34\tabularnewline
\hline
56 & 1.1 & 10 &$a_6$ &$ b_3$ & 355 & 0,03\tabularnewline
\hline
57 & 1.1 & 20 &$a_6$ &$ b_3$ & 1162 & 0,15\tabularnewline
\hline
58 & 1.1 & 30 &$a_6$ &$ b_3$ & 2240 & 0,52\tabularnewline
\hline
59 & 1.1 & 40 &$a_6$ &$ b_3$ & 4116 &2,2\tabularnewline
\hline
60 & 1.1 & 50 &$a_6$ &$ b_3$ & 6060 &7,99\tabularnewline
\hline
61 & 1.1 & 10 &$a_7$ &$ b_3$ & 252 & 0,02\tabularnewline
\hline
62 & 1.1 & 20 &$a_7$ &$ b_3$ & 798 &0,11 \tabularnewline
\hline
63 & 1.1 & 30 &$a_7$ &$ b_3$ &1636  &0,49 \tabularnewline
\hline
64 & 1.1 & 40 &$a_7$ &$ b_3$ & 2760 & 2,15\tabularnewline
\hline
65 & 1.1 & 50 &$a_7$ &$ b_3$ & 4160 & 8,23\tabularnewline
\hline
66 & 1.1 & 10 &$a_2$ &$ b_4$ & 17616  & 1,71192\tabularnewline
\hline
67 & 1.1 & 20 &$a_2$ &$ b_4$ &  58301&  14,7\tabularnewline
\hline
68 & 1.1 & 30 &$a_2$ &$ b_4$ & 123597 & 138,78\tabularnewline
\hline
69 & 1.1 & 40 &$a_2$ &$ b_4$ & - &  - \tabularnewline
\hline
70 & 1.1 & 50 &$a_2$ &$ b_4$ & - & -\tabularnewline
\hline
71 & 1.1 & 10 &$a_3$ &$ b_4$ & 15414  & 1,47 \tabularnewline
\hline
72 & 1.1 & 20 &$a_3$ &$ b_4$ & 51442 & 12,975\tabularnewline
\hline
73 & 1.1 & 30 &$a_3$ &$ b_4$ & 109864 &  123,36 \tabularnewline
\hline
74 & 1.1 & 40 &$a_3$ &$ b_4$ & - & -\tabularnewline
\hline
75 & 1.1 & 50 &$a_3$ &$ b_4$ & - & -  \tabularnewline
\hline
\end{tabular}
\end{center}

\begin{center}
\begin{tabular}{|>{\centering}p{2cm}|>{\centering}p{2cm}|c|c|c|>{\centering}p{4.5cm}|>{\centering}p{3.5cm}|}
\hline
76 & 1.1 & 10 &$a_4$ &$ b_4$ &2202  &0,2149 \tabularnewline
\hline
77 & 1.1 & 20 &$a_4$ &$ b_4$ & 6859 &1,7307 \tabularnewline
\hline
78 & 1.1 & 30 &$a_4$ &$ b_4$ & 13733 &15,4242\tabularnewline
\hline
79 & 1.1 & 40 &$a_4$ &$ b_4$ & -  & - \tabularnewline
\hline
80 & 1.1 & 50 &$a_4$ &$ b_4$ & - & - \tabularnewline
\hline
81 & 1.1 & 10 &$a_5$ &$ b_4$ & 1158 & 0,11\tabularnewline
\hline
82 & 1.1 & 20 &$a_5$ &$ b_4$ &  3610 & 0,91 \tabularnewline
\hline
83 & 1.1 & 30 &$a_5$ &$ b_4$ & 7227 & 8,11 \tabularnewline
\hline
84 & 1.1 & 40 &$a_5$ &$ b_4$ & - & - \tabularnewline
\hline
85 & 1.1 & 50 &$a_5$ &$ b_4$ & - & - \tabularnewline
\hline
86 & 1.1 & 10 &$a_6$ &$ b_4$ & 1001 & 0,14\tabularnewline
\hline
87 & 1.1 & 20 &$a_6$ &$ b_4$ & 3117 & 0,77 \tabularnewline
\hline
88 & 1.1 & 30 &$a_6$ &$ b_4$ & 6242 &  7\tabularnewline
\hline
89 & 1.1 & 40 &$a_6$ &$ b_4$ & - & - \tabularnewline
\hline
90 & 1.1 & 50 &$a_6$ &$ b_4$ & - & - \tabularnewline
\hline
91 & 1.1 & 10 &$a_7$ &$ b_4$ & 734  & 0,07 \tabularnewline
\hline
92 & 1.1 & 20 &$a_7$ &$ b_4$ & 2286 & 0,56 \tabularnewline
\hline
93 & 1.1 & 30 &$a_7$ &$ b_4$ &4577  & 5,14 \tabularnewline
\hline
94 & 1.1 & 40 &$a_7$ &$ b_4$ & - & -\tabularnewline
\hline
95 & 1.1 & 50 &$a_7$ &$ b_4$ & - & -\tabularnewline
\hline
96 &1.2 & 10 &$a_1$ &$ b_1$ &328  & 0,0949\tabularnewline
\hline
97 &1.2 & 20 &$a_1$ &$ b_1$ & 992 & 3,3142\tabularnewline
\hline
98 &1.2 & 30 &$a_1$ &$ b_1$ & 2003 & 63,0966 \tabularnewline
\hline
99 & 1.2 & 40 &$a_1$ &$ b_1$ &  -&- \tabularnewline
\hline
100 & 1.2 & 50 &$a_1$ &$ b_1$ & - &-\tabularnewline
\hline
101 & 1.2 & 10 &$a_2$ &$ b_2$ & 474 &0,0617\tabularnewline
\hline
102 & 1.2 & 20 &$a_2$ &$ b_2$ & 2218 &0,3096\tabularnewline
\hline
103 & 1.2 & 30 &$a_2$ &$ b_2$ & 5085 &1,0566\tabularnewline
\hline
104 &1.2 & 40 &$a_2$ &$ b_2$ & 3108 &4,21\tabularnewline
\hline
105 & 1.2 & 50 &$a_2$ &$ b_2$ & 14888 &24,2917\tabularnewline
\hline
106 & 1.2 & 10 &$a_3$ &$ b_2$ &  424&0,0458\tabularnewline
\hline
\end{tabular}
\end{center} 

\begin{center}
\begin{tabular}{|>{\centering}p{2cm}|>{\centering}p{2cm}|c|c|c|>{\centering}p{4.5cm}|>{\centering}p{3.5cm}|}
\hline
107 & 1.2 & 20 &$a_3$ &$ b_2$ & 1748 &0,2039\tabularnewline
\hline
108 & 1.2 & 30 &$a_3$ &$ b_2$ & 0,8995 &3112\tabularnewline
\hline
109 & 1.2 & 40 &$a_3$ &$ b_2$ & 5425 &8,0161\tabularnewline
\hline
110 & 1.2 & 50 &$a_3$ &$ b_2$ & 8091 &43,83\tabularnewline
\hline
111 &1.2 & 10 &$a_4$ &$ b_2$ & 395  &0,0445\tabularnewline
\hline
112 & 1.2 & 20 &$a_4$ &$ b_2$ &  1333&0,2073\tabularnewline
\hline
113 & 1.2 & 30 &$a_4$ &$ b_2$ & 2725 &1,2688\tabularnewline
\hline
114 & 1.2 & 40 &$a_4$ &$ b_2$ & 4548 &8,83\tabularnewline
\hline
115 & 1.2 & 50 &$a_4$ &$ b_2$ & 6945 &65,61\tabularnewline
\hline
116 & 1.2 & 10 &$a_5$ &$ b_2$ & 318 &0,0336\tabularnewline
\hline
117 & 1.2 & 20 &$a_5$ &$ b_2$ & 1065 &0,438\tabularnewline
\hline
118 & 1.2 &30 &$a_5$ &$ b_2$ & 2170 &6,095\tabularnewline
\hline
119 & 1.2 & 40 &$a_5$ &$ b_2$ & 3620 &125,12\tabularnewline
\hline
120 & 1.2 & 50 &$a_5$ &$ b_2$ &  5343 &977,96\tabularnewline
\hline
121 & 1.2 & 10 &$a_6$ &$ b_2$ &356  &0,0318\tabularnewline
\hline
122 & 1.2 & 20 &$a_6$ &$ b_2$ & 1123 &0,19\tabularnewline
\hline
123 & 1.2 & 30 &$a_6$ &$ b_2$ & 2416 &1,87\tabularnewline
\hline
124 & 1.2 & 40 &$a_6$ &$ b_2$ &  4071 &21,25\tabularnewline
\hline
125 & 1.2 & 50 &$a_6$ &$ b_2$ & 6039 &147,69\tabularnewline
\hline
 126 & 1.2 & 10 &$a_7$ &$ b_2$ & 310 &0,0352 \tabularnewline
\hline
127 & 1.2 & 20 &$a_7$ &$ b_2$ & 1052 &0,3575\tabularnewline
\hline
128 & 1.2 & 30 &$a_7$ &$ b_2$ &  2254& 6,1854\tabularnewline
\hline
129 & 1.2 & 40 &$a_7$ &$ b_2$ & 3748 &87,29\tabularnewline
\hline
130 & 1.2 & 50 &$a_7$ &$ b_2$ &5465 &656,0825\tabularnewline
\hline
131 & 1.2 & 10 &$a_2$ &$ b_3$ & 13540 &2,36\tabularnewline
\hline
132 & 1.2 & 20 &$a_2$ &$ b_3$ & 40200 &24,5699\tabularnewline
\hline
133 & 1.2 & 30 &$a_2$ &$ b_3$ & 160789 &378\tabularnewline
\hline
134 & 1.2 & 40 &$a_2$ &$ b_3$ & - &-\tabularnewline
\hline
135 & 1.2 & 50 &$a_2$ &$ b_3$ & - &-\tabularnewline
\hline
136  & 1.2 & 10 &$a_3$ &$ b_3$ &  3168&0,3962\tabularnewline
\hline
137 & 1.2 & 20 &$a_3$ &$ b_3$ & 9266 &5,6791\tabularnewline
\hline
\end{tabular}
\end{center} 

\begin{center}
\begin{tabular}{|>{\centering}p{2cm}|>{\centering}p{2cm}|c|c|c|>{\centering}p{4.5cm}|>{\centering}p{3.5cm}|}
\hline
138 & 1.2 & 30 &$a_3$ &$ b_3$ & 17362&97\tabularnewline
\hline
139 & 1.2 & 40 &$a_3$ &$ b_3$ & 27779 &301\tabularnewline
\hline
140 & 1.2 & 50 &$a_3$ &$ b_3$ & -& -\tabularnewline
\hline
141 & 1.2 & 10 &$a_4$ &$ b_3$ & 1803 &0,3125\tabularnewline
\hline
142 & 1.2 & 20 &$a_4$ &$ b_3$ &5292  &3,62\tabularnewline
\hline
143 & 1.2 & 30 &$a_4$ &$ b_3$ & 10649 &65,38\tabularnewline
\hline
144 & 1.2 & 40 &$a_4$ &$ b_3$ & 24197 &197\tabularnewline
\hline
145 & 1.2 & 50 &$a_4$ &$ b_3$ &   -& -\tabularnewline
\hline
146 & 1.2 & 10 &$a_5$ &$ b_3$ & 665 &0,12 \tabularnewline
\hline
147 & 1.2 & 20 &$a_5$ &$ b_3$ & 1930 &1,72\tabularnewline
\hline
148 & 1.2 & 30 &$a_5$ &$ b_3$ & 5078 &19,15\tabularnewline
\hline
149 & 1.2 & 40 &$a_5$ &$ b_3$ & 13752 & 101\tabularnewline
\hline
150 & 1.2 & 50 &$a_5$ &$ b_3$ &  -&-\tabularnewline
\hline
151 & 1.2 & 10 &$a_6$ &$ b_3$ & 964 & 0,17\tabularnewline
\hline
152 & 1.2 & 20 &$a_6$ &$ b_3$ & 3108 & 2,42\tabularnewline
\hline
153 & 1.2 & 30 &$a_6$ &$ b_3$ & 6576 & 45,49\tabularnewline
\hline
154 & 1.2 & 40 &$a_6$ &$ b_3$ &  19789&178\tabularnewline
\hline
155 & 1.2 & 50 &$a_6$ &$ b_3$ & - &-\tabularnewline
\hline
156 & 1.2 & 10 &$a_7$ &$ b_3$ & 578 & 0,1063 \tabularnewline
\hline
157 & 1.2 & 20 &$a_7$ &$ b_3$ &  1788& 1,75\tabularnewline
\hline
158 & 1.2 & 30 &$a_7$ &$ b_3$ &  3902& 30,99\tabularnewline
\hline
159 & 1.2 & 40 &$a_7$ &$ b_3$ & 8127 &105\tabularnewline
\hline
160 & 1.2 & 50 &$a_7$ &$ b_3$ & - &-\tabularnewline
\hline
161 & 1.2 & 10 &$a_2$ &$ b_4$ &30577  &6,3045 \tabularnewline
\hline
162 & 1.2 & 20 &$a_2$ &$ b_4$ &733848  & 506,27\tabularnewline
\hline
163 & 1.2 & 30 &$a_2$ &$ b_4$ &  -& -\tabularnewline
\hline
164 & 1.2 & 40 &$a_2$ &$ b_4$ & - &-  \tabularnewline
\hline
165 & 1.2 & 50 &$a_2$ &$ b_4$ &  -& -
\tabularnewline
\hline
166 & 1.2 & 10 &$a_3$ &$ b_4$ & 7372 &1,5129 
\tabularnewline
\hline
167 & 1.2 & 20 &$a_3$ &$ b_4$ &  84778& 190,48
\tabularnewline
\hline
168 & 1.2 & 30 &$a_3$ &$ b_4$ &  -& -
\tabularnewline
\hline
\end{tabular}
\end{center} 

\begin{center}
\begin{tabular}{|>{\centering}p{2cm}|>{\centering}p{2cm}|c|c|c|>{\centering}p{4.5cm}|>{\centering}p{3.5cm}|}
\hline
169 & 1.2 & 40 &$a_3$ &$ b_4$ & - &- 
\tabularnewline
\hline
170 & 1.2 & 50 &$a_3$ &$ b_4$ & - &- 
\tabularnewline
\hline
171 & 1.2 & 10 &$a_4$ &$ b_4$ &  4699& 0,968
\tabularnewline
\hline
172 & 1.2 & 20 &$a_4$ &$ b_4$ & 84582 & 112,065
\tabularnewline
\hline
173 & 1.2 & 30 &$a_4$ &$ b_4$ & - & -
\tabularnewline
\hline
174 & 1.2 & 40 &$a_4$ &$ b_4$ & - &- 
\tabularnewline
\hline
175 & 1.2 & 50 &$a_4$ &$ b_4$ &  -&- 
\tabularnewline
\hline
176 & 1.2 & 10 &$a_5$ &$ b_4$ &  1606& 0,3331
\tabularnewline
\hline
177 & 1.2 & 20 &$a_5$ &$ b_4$ &  8883&47,78 
\tabularnewline
\hline
178 & 1.2 & 30 &$a_5$ &$ b_4$ & - & -
\tabularnewline
\hline
179 & 1.2 & 40 &$a_5$ &$ b_4$ & - &- 
\tabularnewline
\hline
180 & 1.2 & 50 &$a_5$ &$ b_4$ &-  &- 
\tabularnewline
\hline
181 & 1.2 & 10 &$a_6$ &$ b_4$ &  2203& 0,4512
\tabularnewline
\hline
182 & 1.2 & 20 &$a_6$ &$ b_4$ & 7710 & 10,518
\tabularnewline
\hline
183 & 1.2 & 30 &$a_6$ &$ b_4$ &69394  & 158,48
\tabularnewline
\hline
184 & 1.2 & 40 &$a_6$ &$ b_4$ & - &- 
\tabularnewline
\hline
185 & 1.2 & 50 &$a_6$ &$ b_4$ &  -&- 
\tabularnewline
\hline
186 & 1.2 & 10 &$a_7$ &$ b_4$ & 893 &0,1837 
\tabularnewline
\hline
187 & 1.2 & 20 &$a_7$ &$ b_4$ & 2232 &2,378 
\tabularnewline
\hline
188 & 1.2 & 30 &$a_7$ &$ b_4$ & 41658 &105,069 
\tabularnewline
\hline
189 & 1.2 & 40 &$a_7$ &$ b_4$ & - &- 
\tabularnewline
\hline
190 & 1.2 & 50 &$a_7$ &$ b_4$ & - &- 
\tabularnewline
\hline
\end{tabular}
\end{center} 

Прокомментируем смысл содержимого колонок. В первой колонке указана нумерация численных экспериментов, во второй колонке $-$ метод, с помощью которого решались тестовые задачи, в третьей $-$ размерность пространства переменных, в четвертной указаны способы задания чисел $\varepsilon_k$, а в пятой колонке указаны способы обновления аппроксимирующих множеств, описание которых будет представлено позже. В  шестой колонке представлены количество итераций, которое понадобилось для достижения заданной точности решения задачи, а в седьмой $-$ время решения. C помощью символа '$-$'  в таблице фиксировались те тестовые примеры, время решения  которых превосходило 24 часов.

Обсудим теперь обозначения  четвертой  и пятой колонок. Отметим, что последовательность $\{\varepsilon_k\}$, $k\in K$,  задавалась как на предварительном шаге, так и во время итерационного процесса. Для обозначения последовательности $\{\varepsilon_k\}$, $k\in K$, справедливы следующие соответствия:

1) $a_1$ $\leftrightarrow$ $\varepsilon_0=0$;

2) $a_2$ $\leftrightarrow$ $\varepsilon_{k+1}=\frac{\varepsilon_k}{1.1}$ для всех $k\in K$;

3) $a_3$ $\leftrightarrow$ $\varepsilon_{k+1}=\frac{\varepsilon_k}{1.5}$ для всех $k\in K$;

4) $a_4$ $\leftrightarrow$ $\varepsilon_{k+1}=\frac{\varepsilon_k}{2}$ для всех $k\in K$;

5) $a_5$ $\leftrightarrow$ $\varepsilon_{k+1}=\frac{\varepsilon_k}{n}$ для всех $k\in K$;

6) $a_6$ $\leftrightarrow$ $\varepsilon_{k+1} = \alpha_k F(x_k)$, $\alpha_k = \frac{1}{1.1^k}$ для всех $k\in K$;

7) $a_7$ $\leftrightarrow$ $\varepsilon_{k+1} = \alpha_k F(x_k)$, $\alpha_k = \frac{1}{2^k}$ для всех $k\in K$.

Для обозначения способов выбора множеств $Q_i$ приняты следующие соответствия:

1) $b_1$ $\leftrightarrow$ $Q_i=M_i$ для всех $i\in K$;

2)  $b_2$ $\leftrightarrow$ $Q_i=M_i$ для всех $i\ne i_k$, множество $Q_{i_k}$ состоит из <<активных>>  в текущей итерационной точке $y_{i_k}$ секущих плоскостей;

3)  a) в методе 1.1: 

 $b_3$ $\leftrightarrow$ $Q_i=M_i$, $0\le i< n+1$,

 $Q_i = \begin{cases}M_i, \ i\ne i_k, \\
 \overset{i_k}{\underset{r=i_k - n - 1}{\cap}}\{x\in{\rm  R_{n}} : \langle a_j, x - z_r^j \rangle\le 0\ \forall j\in J_r, \\
\hfill  a_j\in W^1(z_r^j, D_j) \},  i = i_k,
    \end{cases}$

\quad{} b) в методе 1.2: 

 $b_3$ $\leftrightarrow$ $Q_i=M_i$, $0\le i< n+1$,

 $Q_i = \begin{cases}M_i, \ i\ne i_k, \\
 \overset{i_k}{\underset{r=i_k - n - 1}{\cap}}\{x\in R_{n} : F(u_r) +  \langle a, x - u_r \rangle\le 0,\\
 \hfill \ a\in \partial F(u_r) \},   i = i_k,  \end{cases}$

4)  $b_4$ $\leftrightarrow$ $Q_i = \begin{cases}M_i, & i \ne i_k, \\ M_0,& i = i_k. \end{cases}$.

Перейдем, наконец, к обсуждению полученных результатов. Сравним для начала тестовые примеры, при решении которых аппроксимирующие множества $Q_i$ не обновлялись (см. строки в таблице со значениями $b_1$), с примерами, где множества $Q_i$ обновлялись полностью (см. строки в таблице со значениями $b_4$). На основе приведенных результатов можно сделать вывод о том, что полностью отбрасывать все секущие плоскости, построенные к итерации $i=i_k$, не стоит. Дело в том, что при полном обновлении утрачивается  качество аппроксимации допустимой области, и происходит резский спад значения целевой функции в очередной итерационной точке. В связи с этим для достижения утраченного качества аппроксимации приходится проделать существенной количество итераций. Таким образом, происходит рост трудоемкости решения задачи не только по затраченному времени, но и по количеству проделанных итераций.

Отметим, что если все-таки решено на итерациях $i=i_k$ полностью отбрасывать  все построенные ранее отсекающие гиперплоскости, то, как показали расчеты, последовательность $\{\varepsilon_k\}$, $k\in K$, предпочитительно выбирать быстроубывающей. Например, проанализируем эксперименты с номерами 78 и 83.  В примере 83  последовательность $\{\varepsilon_k\}$, $k\in K$, убывает быстрее по сравнению с последовательностью чисел $\varepsilon_k$, $k\in K$, из примера 78. Обратим внимание на то, что количество проделанных итераций в примере 78 в 3.8 раз больше по сравнению с количеством проделанных итераций из примера 83, а по количеству затраченного времени в 3,42 раза.

Отметим, что решение многих примеров с использованием полного обновления аппроксимирующих множеств было остановлено в связи с истечением лимита времени, который равнялся одному дню.

Сравним теперь примеры, решенные без отбрасывания секущих плоскостей (см. строки в таблице со значениями $b_1$), с примерами, при решении которых в множества $Q_{i_k}$ включаются либо все <<активные>> гиперплоскости (см. записи в таблице со значениями $b_2$),  либо $n+1$ последних секущих плоскостей. Привлекая во внимание результаты численных экспериментов, можно сделать вывод о том, что в момент обновления аппроксимирующего множества любыми  из двух приведенных способов качество аппроксимации допустимой области в окрестности решения сильно не утрачивается. Однако, как показали численные эксперименты, способ обновления, в котором в множество $Q_{i_k}$ включались <<активные>> гиперплоскости, выглядит предпочтительнее, поскольку при таком способе обновления для решения задач затрачивается меньше времени.

Далее обсудим выбор последовательности $\{\varepsilon_k\}$, $k\in K$, которая фактически определяет частоту обновления аппроксимирующих множеств. Как уже отмечалось ранее, в случае полного обновления аппроксимирующих множеств на итерациях $i=i_k$ числа $\varepsilon_k$ следует выбирать быстроубывающими. Однако, при включении в аппроксимирующие множества на итерациях $i=i_k$ всех <<активных>> секущих плоскостей последовательность $\{\varepsilon_k\}$, $k\in K$, рекомендуется выбирать медленноубывающей (см., напр., эксперименты № 6 $-$ 17).

Перейдем, наконец, к сравнению методов 1.1 и 1.2. Учитывая результаты численных экспериментов 1 $-$ 5 и 96 $-$ 100, можно сказать, что метод 1.1  выглядит предпочтительнее метода 1.2., так как для решения указанных задач методом 1.1 затрачено меньше времени и итераций.

\chapter[Методы отсечений с аппроксимацией надграфика целевой функции]
{\quad{} \quad{}\hspace{1mm}Методы отсечений
\newline с аппроксимацией надграфика целевой функции}


В данной главе предлагаются методы отсечений для решения задачи выпуклого программирования, которые используют операцию погружения надграфика целевой функции. Как и методы главы 1, они допускают периодические обновления погружающих множеств за счет отбрасывания отсекающих плоскостей.

Первый параграф  главы посвящен решению задачи условной минимизации дискретной функции максимума. Предлагаемый метод ее решения использует при построении аппроксимирующих надграфик множеств субградиенты целевой функции. В следующем параграфе разрабатывается метод выпуклого программирования, основанный на идеях метода \cite{n_bulatov_v_p_4}. Предлагаемый метод использует для построения отсечений обобщенно-опорные векторы и существенно отличается от упомянутого известного метода возможностью обновлений аппроксимирующих множеств. В $\S\ 2.3$ предлагается модификация известного метода уровней \cite{n_nesterov_yu_e_3,n_lemarechal_c_1}. Главным отличием предлагаемого метода от известного является возможность отбрасывания накапливающихся в процессе счета отсекающих плоскостей.  В $\S$ 2.4, основываясь на способе построения аппроксимирующих множеств из $\S$ 2.2, разрабатывается еще один метод отсечений для задачи выпуклого программирования. Существенным отличием метода из $\S$ 2.4 от методов данной главы является заложенный в нем новый критерий качества аппроксимирующих множеств, а значит, и критерий отбрасывания отсекающих плоскостей. Часть результатов второй главы можно найти в работах авторов \cite{n_yarullin_r_s_1,
n_yarullin_r_s_7,
n_yarullin_r_s_9_1,
n_yarullin_r_s_10,
n_yarullin_r_s_12,
n_yarullin_r_s_19}.

\section
{Метод отсечений  для отыскания дискретного минимакса с отбрасыванием  секущих плоскостей}\label{s2.1}

Пусть $D$ -- выпуклое замкнутое множество из ${\rm R_n}$, $J=\{1,\ldots,m\}$, $m\ge1$, функции $f_j(x)$ для всех $j\in J$ выпуклы в ${\rm R_n}$,
$$f(x)=\max \limits_{j\in J} f_j(x).$$
Решается задача 
\begin{equation}\label{e2.1.0}
\min\{f(x) : x\in D\}.
\end{equation}

Пусть $f^*=\min\{f(x) : x\in D\}$, $X^*=\{x\in D : f(x) = f^*\}$, $\partial f_j(x)$ -- субдифференциал функции $f_j(x)$ в точке $x\in {\rm R_n}$. Положим $K=\{0, 1, ...\}$, ${\rm epi\,} (f,G)=\{(x,\gamma)\in {\rm R_{n+1}} : x\in G, \ \gamma\ge f(x)\}$, где $G\subset {\rm R_n}$.

Метод решения задачи (\ref{e2.1.0}) вырабатывает последовательности $\{y_i\}$, $i\in K$, и $\{x_k\}$, $k\in K$, приближений из множества $D$, и заключается в следующем. 

\textbf{Метод 2.1.} Выбирается выпуклое замкнутое ограниченное  множество $G_0\subseteq D$, содержащее хотя бы одну точку из $X^*$, например, $x^*$. Строится выпуклое замкнутое множество $M_0\subset {\rm R_{n+1}}$ такое, что 

$$(x^*, f^*)\in M_0.$$
Задаются числа $\varepsilon_0\ge0$, $-\infty<\bar{\gamma}_0\le f^*$, и полагается $i=0$, $k=0$. 

\textbf{1.} Отыскивается точка $(y_i, \gamma_i)$, где $y_i\in {\rm R_n}$, $\gamma_i\in {\rm R_1}$, как решение задачи 

$$\gamma \to \min,$$
\begin{equation}\label{e2.1.1}
(x,\gamma)\in M_i, \quad x\in G_i,\quad \gamma\ge \bar{\gamma}_i.
\end{equation}
Если 
\begin{equation}\label{e2.1.2}
f(y_i)=\gamma_i,
\end{equation}
то $y_i\in X^*$, и процесс завершается.

\textbf{2.} Если выполняется неравенство
\begin{equation}\label{e2.1.3}
f(y_i)-\gamma_i>\varepsilon_k,
\end{equation}
то выбирается выпуклое замкнутое множество $S_i\subseteq {\rm R_{n+1}}$, содержащее точку $(x^*, f^*)$, полагается 
\begin{equation}\label{e2.1.4}
Q_i=M_i\bigcap S_i,
\end{equation}
и следует переход к п. 3. В противном случае полагается $i_k=i$, 
\begin{equation}\label{e2.1.5}
x_k=y_{i_k}, \quad \sigma_k=\gamma_{i_k},
\end{equation}
выбирается выпуклое замкнутое множество $Q_i=Q_{i_k}\subseteq {\rm R_{n+1}}$ такое, что 
\begin{equation}\label{e2.1.6}
(x^*, f^*)\in Q_i,
\end{equation}
задается $\varepsilon_{k+1}\ge 0$, значение $k$ увеличивается на единицу.

\textbf{3.}  Полагается
$$J(y_i)=\{j\in J : f(y_i)= f_j(y_i)\}.$$
Выбирается множество $J_i\subset J$ следующим образом. В $J_i$ заносится любой номер $j_i \in J(y_i)$. Кроме $j_i$ в множество $J_i$ включаются любые другие номера из $J$ или не включается более ни одного.

\textbf{4.} Для каждого $j\in J_i$ выбирается конечное подмножество $A_i^j$ множества $\partial f_j(y_i)$ и полагается
\begin{equation}\label{e2.1.7}
M_{i+1}=Q_i\cap\{(x,\gamma)\in {\rm R_{n+1}} : f_j(y_i)+\langle a, x-y_i\rangle \le \gamma \ \forall j\in J_i, \ a\in A_i^j\}.
\end{equation}

\textbf{5.} Выбирается такое выпуклое замкнутое множество $G_{i+1}\subseteq G_0,$ что $x^*\in G_{i+1}$. Задается число $\bar{\gamma}_{i+1}$ из условия $\bar{\gamma}_0 \le \bar{\gamma}_{i+1}\le f^*$. Значение $i$ увеличивается на единицу, и следует переход к п. 1.

Прежде, чем обсуждать свойства последовательностей $\{(y_i$, $\gamma_i)\}$, $i\in K$, и $\{(x_k$, $\sigma_k)\}$, $k\in K$,  сделаем некоторые замечания к методу.

\begin{remark}\label{r2.1.1}
 Если множества $D$ и $G_i$ -- многогранники, то (\ref{e2.1.1}) является задачей линейного программирования. Подчеркнем, что метод пригоден для задачи (\ref{e2.1.0}) и в том случае, если $D={\rm R_n}$. Тогда на подготовительном шаге метода следует построить многогранник $G_0$, содержащий точку безусловного минимума функции $f(x)$.
\end{remark}

\begin{remark} Для выбора множества $M_0$ имеются большие возможности. Его можно задать содержащим ${\rm epi\,} (f, G_0)$, например, положить $M_0={\rm epi\,} (f_{j_0}, {\rm R_n})$, где $j_0\in J$, задать линейным неравенством 
\begin{equation}\label{e2.1.8}
\langle c, x\rangle - \gamma \le \langle c, z\rangle - f_{j_0}(z),
\end{equation}
где $z\in {\rm R_n}$, $c\in \partial f_{j_0}(z)$, или группой подобных неравенств. \mbox{Если положить}
$$M_0={\rm R_n},$$
то пара $(y_0, \gamma_0)$, где $y_0$ -- любая точка из $G_0$, а $\gamma_0 = \bar{\gamma}_0$, может быть принята за решение задачи (\ref{e2.1.1}) при $i=0$.
\end{remark}

\begin{remark}\label{r2.13}
Если нет информации о значении $f^*$, то число $\bar{\gamma}_0$ можно выбрать, например, как решение задачи минимизации переменной $\gamma$ при ограничении (\ref{e2.1.8}) и условии $x\in G_0$. Если $G_0$ -- многогранник, и среди функций $f_j(x)$, $j\in J$, есть линейная функция, то ее минимальное значение на $G_0$ можно также принять за $\bar{\gamma}_0$. Если же множество $M_0$ ограничено, то число $\bar{\gamma}_0$ можно считать  сколь угодно малым, поскольку задача (\ref{e2.1.1}) будет иметь решение при $i=0$ и без ограничения $\gamma \ge \bar{\gamma}_0$. Выбор значений $\bar{\gamma}_i$ для $i>0$ будет обсуждаться позже.
\end{remark}

\begin{remark}\label{r2.1.4}
Для выбора множеств $S_i$ и $G_i$ в пп. 2, 5 также есть немало возможностей. В частности, для всех $i\in K$ допустимо положить $S_i={\rm R_{n+1}}$ и $G_i=D$, если $D$  ограничено. В случае $S_i \ne {\rm R_{n+1}}$ множества $S_i$, как и $M_i$, естественно строить многогранными. Например, если 
$$D=\{x\in {\rm R_n} : F(x)\le 0\},$$
где $F(x)$ -- выпуклая функция, то $S_i$ можно задавать неравенствами вида $F(z)+\langle c(z), x-z\rangle \le \gamma$,
 выбирая разными способами точки $z\in D$ и векторы $c(z)\in \partial F(z)$. Множество $G_{i+1}$, $i\ge 0$, можно, например, задавать в виде
$$G_{i+1}=G_r\bigcap \{x\in {\rm R_n} : \langle b, x-y_i\rangle \le 0 \ \forall b\in B_i\},$$
где $B_i\subset \partial f(y_i)$, $0\le r \le i$. Так как для любых $x^*\in X^*$ и $b\in B_i$ выполняются неравенства $0\ge f(x^*)-f(y_i)\ge \langle b, x^* - y_i\rangle$, и $x^*\in G_r$, то $x^*\in G_{i+1}$, как требуется в п. 5 метода.
\end{remark}

\begin{lemma}\label{lemma2.1.1}
Точка $(x^*, f^*)$ удовлетворяет ограничениям задачи (\ref{e2.1.1}) для всех $i\in K$.
\end{lemma}

\begin{proof}
Согласно выбору множеств $G_i$ и чисел $\bar{\gamma}_i$ для всех $i\in K$ выполняются включение $x^*\in G_i$ и неравенство $f^*\ge \bar{\gamma}_i$. Поэтому для обоснования леммы достаточно показать, что
\begin{equation}\label{e2.1.9}
(x^*, f^*)\in M_i
\end{equation}
для всех $i\in K$.

При $i=0$ включение (\ref{e2.1.9}) выполняется по условию выбора $M_0$. Допустим теперь, что (\ref{e2.1.9}) справедливо при любом фиксированном $i=l\ge 0$. Покажем выполнение (\ref{e2.1.9}) при $i=l+1$, тогда утверждение будет доказано. Действительно, ввиду (\ref{e2.1.4}), (\ref{e2.1.6}), условия выбора $S_l$ и сделанного выше предположения  множество $Q_l$, независимо от (\ref{e2.1.3}), содержит точку $(x^*, f^*)$. Кроме того, $f_j(x^*)-f_j(y_l)\ge \langle a, x^* - y_l\rangle$  для всех $j\in J_l$ и $a\in A_l^j$. Но $f_j(x^*)\le f(x^*)=f^*$. Значит, $f_j(y_l)+\langle a, x^* - y_l\rangle \le f^*$ для всех $ j\in J_l$ и $a\in A_l^j$, и ввиду (\ref{e2.1.7}) $(x^*, f^*)\in M_{l+1}$. Лемма доказана.
\end{proof}

Из этой леммы следует, что задача (\ref{e2.1.1}) разрешима для всех $i\in K$.

Качество аппроксимации надграфика  функции множествами $M_i$ определяется в окрестности точек $y_i$ величиной $f(y_i)-\gamma_i$. На итерациях, где для $y_i$ выполняется неравенство (\ref{e2.1.3}), качество аппроксимации считается недостаточным для фиксирования точки $x_k$. На таких итерациях множества $Q_i$  выбираются в виде (\ref{e2.1.4}), и тогда согласно (\ref{e2.1.7}) при формировании множеств $M_{i+1}$ отсекающие плоскости накапливаются. На некотором шаге $i=i_k$, как будет доказано ниже, для разности $f(y_i)-\gamma_i$ выполнится неравенство, противоположное (\ref{e2.1.3}). В таком случае качество аппроксимации станет удовлетворительным, точка $x_k$ зафиксируется в виде (\ref{e2.1.5}), и за счет практически произвольного выбора множества $Q_i=Q_{i_k}$ произойдет обновление множества $M_{i_k+1}$.

С учетом (\ref{e2.1.6}), (\ref{e2.1.9}) множества $Q_i$ можно задавать в виде (\ref{e2.1.4}), независимо от выполнения условия (\ref{e2.1.3}). Однако, в таком случае от шага к шагу неограниченно растет число отсекающих плоскостей, формирующих аппроксимирующие множества, и вспомогательные задачи (\ref{e2.1.1}) становятся труднорешаемыми. Приведем теперь принципиальное замечание, касающееся возможности периодического обновления аппроксимирующих множеств за счет отбрасывания любого числа отсекающих плоскостей.

\begin{remark}\label{r2.1.5}
Для тех $i$ и $k$, при которых выполняются неравенства 
\begin{equation}\label{e2.1.10}
f(y_i)-\gamma_i\le \varepsilon_k,
\end{equation}
условие (\ref{e2.1.6})  дает большие возможности в выборе множеств $Q_i=Q_{i_k}$, а значит, и в задании множеств $M_{i+1}$. В случае (\ref{e2.1.10}) можно положить, например, 
\begin{equation}\label{e2.1.11}
Q_i=Q_{i_k}={\rm R_{n+1}}.
\end{equation}
Тогда $M_{i+1}$ допустимо задать одним неравенством $f_{j_i}(y_i)+\langle a_i, x-y_i\rangle \le \gamma$, считая, что $J_i=\{j_i\}$, $A_i^{j_i}=\{a_i\}$. При условии (\ref{e2.1.10}) множество $Q_i$ можно задать и в виде 
\begin{equation}\label{e2.1.12}
Q_i=M_{r_i},
\end{equation}
где $0\le r_i\le i=i_k$, поскольку при всех $r_i=0,...,i$ ввиду (\ref{e2.1.9}) включение (\ref{e2.1.6}) выполняется. Согласно (\ref{e2.1.11}), (\ref{e2.1.12}) при каждом $i=i_k$ можно отбрасывать при формировании $M_{i+1}$ любое количество любых накопленных к шагу $i_k$ отсекающих плоскостей. Как будет показано ниже, при определенном способе выбора чисел $\varepsilon_k$ для каждого $k\in K$ найдется решение $(y_i, \gamma_i)=(y_{i_k}, \gamma_{i_k})$ задачи (\ref{e2.1.1}), удовлетворяющее (\ref{e2.1.10}), и, следовательно, представится возможность задания $Q_i$ в виде (\ref{e2.1.11}) или (\ref{e2.1.12}).
\end{remark}
Исследуем  некоторые свойства последовательности $\{(y_i$, $\gamma_i)\}$, $i\in K$,  построенной предложенным методом.

\begin{lemma}\label{lemma2.1.2}
Для построенной методом последовательности $\{\gamma_i\}$, $i\in K$, справедливы неравенства 
\begin{equation}\label{e2.1.13}
\gamma_i \le f^* 
\end{equation}
при всех $i\in K$.
\end{lemma}

\begin{proof}
Согласно (\ref{e2.1.1}) для всех $i\in K$ и всех точек $(x, \gamma)$ таких, что $(x, \gamma)\in M_i$, $x\in G_i$, $\gamma \ge \bar{\gamma}_i$, выполняются неравенства $\gamma_i\le \gamma$. Но по лемме \ref{lemma2.1.1} точка $(x^*, f^*)$ является допустимым решением задачи (\ref{e2.1.1}) при каждом $i\in K$. Отсюда следует (\ref{e2.1.13}). Лемма доказана.
\end{proof}

На основе леммы \ref{lemma2.1.2} обоснуем теперь критерий остановки, заложенный в п. 1 метода.

\begin{theorem}\label{theorem2.1.1}
 Пусть при некотором $i\in K$ выполняется равенство (\ref{e2.1.2}). Тогда $y_i\in X^*$. 
\end{theorem}

\begin{proof}
В силу (\ref{e2.1.2}), (\ref{e2.1.13}) $f(y_i)\le f^*$ для всех $i\in K$. С другой стороны 
\begin{equation}\label{e2.1.14}
f(y_l)\ge f^* \quad \forall l\in K, 
\end{equation}
так как $y_l\in D$  для всех $l\in K$. Таким образом, $f(y_i)=f^*$ при всех $i\in K$, и утверждение доказано.
\end{proof}

\begin{remark}\label{r2.1.6}
Лемма \ref{lemma2.1.2} позволяет выбирать числа $\bar{\gamma}_{i+1}$ в п. 5 метода следующим образом. Можно положить $\bar{\gamma}_{i+1}=\gamma_l$, где $0\le l\le i$, в частности, $\bar{\gamma}_{i+1}=\max \limits_{0\le j\le i} \gamma_j$.
\end{remark}

\begin{remark}\label{r2.1.7}
 В силу (\ref{e2.1.13}), (\ref{e2.1.14}) для каждой из построенных точек $(y_i, \gamma_i)$, $i\in K$, выполняются неравенства
$$\gamma_i\le f^*\le f(y_i).$$
Поэтому на любой итерации есть возможность оценить близость значения $f(y_i)$ к искомому значению $f^*$.
\end{remark}

Как уже было отмечено, согласно п. 2 метода, независимо от условия (\ref{e2.1.3}) множества $Q_i$ для всех $i\in K$ могут иметь вид (\ref{e2.1.4}). С учетом этого сформулируем следующее вспомогательное утверждение.

\begin{lemma}\label{lemma2.1.3}
Пусть последовательность $\{(y_i, \gamma_i)\}$, $i\in K$,  построена с условием, что для всех $i\in K$, начиная с некоторого номера $i^\prime\ge 0$, множества $Q_i$ были выбраны согласно (\ref{e2.1.4}). Тогда выполняется равенство
\begin{equation}\label{e2.1.15}
\lim \limits_{i\in K} (f(y_i)-\gamma_i)=0.
\end{equation}
\end{lemma}

\begin{proof}
Отметим, что ввиду (\ref{e2.1.13}), (\ref{e2.1.14}) $f(y_i)-\gamma_i\ge 0$ для всех $i\in K$. Предположим, что утверждение леммы неверно. Тогда для $\{(y_i, \gamma_i)\}$, $i\in K$, найдется такая подпоследовательность точек $(y_i, \gamma_i)$, $i\in K_1\subset K$, $i\ge i^\prime$, что 
\begin{equation}\label{e2.1.16}
f(y_i)-\gamma_i\ge \delta > 0  \quad \forall i\in K_1, \quad  i\ge i^\prime.
\end{equation}
Пусть $\{(y_i,\gamma_i)\}$, $i\in K_2\subset K_1$ -- сходящаяся подпоследовательность, выделенная из точек $(y_i, \gamma_i)$, $i\in K_1$, $i\ge i^\prime$, а $(\bar{y}, \bar{\gamma})$ -- предельная точка этой подпоследовательности.

По построению множеств $J_i$ существует такой номер $p\in J$, что для бесконечного числа номеров $i\in K_2$ одновременно выполняются включения $p\in J_i$ и $p\in J(y_i)$. Положим

$$K_3=\{i\in K_2 : p\in J_i, p\in J(y_i) \}.$$
Заметим, что согласно п. 4 метода для всех $i\in K_3$ были построены и использованы при нахождении $M_{i+1}$ множества $A_i^p\subset \partial f_p(y_i)$.

Зафиксируем такие номера $r$, $l \in K_3$, что $r\le l-1$. Пусть $c\in A_r^p$. Поскольку $c\in \partial f_p(y_r)$ и при этом $p\in J(y_r)$, а $f_p(x)\le f(x)$ для всех $x\in {\rm R_n}$, то
\begin{equation}\label{e2.1.17}
c\in \partial f(y_r).
\end{equation}
Так как в силу (\ref{e2.1.7}) выполняется включение
$$M_l\subset \{(x,\gamma)\in {\rm R_{n+1}} : f_p(y_r)+\langle a, x-y_r\rangle \le \gamma \ \forall a\in A_r^p\},$$
а $c\in A_r^p$ и $(y_l, \gamma_l)\in M_l$ ввиду (\ref{e2.1.1}), то $\gamma_l \ge f_p(y_r) + \langle c, y_l - y_r\rangle$. Отсюда с учетом равенства $f_p(y_r)=f(y_r)$ имеем
$$\gamma_l \ge f(y_r) + \langle c, y_l-y_r\rangle, $$
а поскольку согласно (\ref{e2.1.16}) $f(y_r)\ge \delta + \gamma_r$, то
\begin{equation}\label{e2.1.18}
\gamma_l\ge \delta +\gamma_r +\langle c, y_l - y_r\rangle . 
\end{equation}
Положим $\theta=\max\{\|a(x)\| : a(x)\in \partial f(x), \ x\in G_0\}$.  Сразу отметим, что $\theta < \infty$ (см., напр., \cite[c. 121]{n_polayk_b_t_1}). Тогда из (\ref{e2.1.18}) и (\ref{e2.1.17}) следует неравенство 
\begin{equation}\label{e2.1.19}
\gamma_{l} - \gamma_r + \theta \|y_{l}-y_r\|\ge \delta.
\end{equation}

Зафиксируем теперь для каждого $i\in K_3$ такой номер $l_i\in K_3$, что $l_i\ge i+1$. Тогда в силу (\ref{e2.1.19})
$$\gamma_{l_i} - \gamma_i + \theta \|y_{l_i}-y_i\|\ge \delta \quad \forall i\in K_3.$$
Последние неравенства противоречивы, так как $y_i \to \bar{y}$, $y_{l_i} \to \bar{y}$ по $i\in K_3$, и $\gamma_i \to \bar{\gamma}$, $\gamma_{l_i} \to\bar{\gamma}, \  i\in K_3$. Лемма доказана.
\end{proof}

\begin{consequence}\label{consequence2.1.1}
Пусть выполняются условия леммы \ref{lemma2.1.3}. Тогда для любой предельной точки $(\bar{y}, \bar{\gamma})$ последовательности $\{(y_i, \gamma_i)\}$, $i\in K$, имеет место равенство $f(\bar{y})=\bar{\gamma}$.
\end{consequence}
\begin{dokaz} утверждения следует непосредственно из леммы \ref{lemma2.1.3}.
\end{dokaz}
\begin{theorem}\label{theorem2.1.2}
Пусть выполняются условия леммы \ref{lemma2.1.3}. Тогда справедливы равенства
\begin{equation}\label{e2.1.20}
\lim \limits_{i\in K} f(y_i)=f^*, \quad \lim \limits_{i\in K}\gamma_i=f^*.
\end{equation}
\end{theorem}

\begin{proof}
В силу (\ref{e2.1.13}), (\ref{e2.1.14}) $0\le f(y_i)-f^*\le f(y_i)-\gamma_i$ для всех $i\in K$. Отсюда и из (\ref{e2.1.15}) следует первое из равенств (\ref{e2.1.20}). Тогда с учетом (\ref{e2.1.15}) имеет место и второе утверждение. Теорема доказана.
\end{proof}

Заметим, что согласно теореме \ref{theorem2.1.2} и известной теореме (\cite[62]{n_vasillev_f_p_1}) при условиях леммы \ref{lemma2.1.3} последовательность $\{y_i\}$, $i\in K$, сходится к множеству $X^*$.

Подчеркнем, что утверждения (\ref{e2.1.20}) доказаны при условии, что в методе для всех $i\in K$ полагается
\begin{equation}\label{e2.1.21}
\begin{split}
M_{i+1}=M_i\bigcap S_i&\bigcap \{(x, \gamma)\in {\rm R_{n+1}} : \\ &  f_j(y_i)+\langle a, x-y_i\rangle \le \gamma 
 \forall j\in J_i, \ a\in A_i^j\}.
\end{split}
\end{equation}
В таком случае при построении $\{(y_i, \gamma_i)\}$, $i\in K$, число отсекающих плоскостей неограниченно растет, обновлений аппроксимирующих множеств не происходит, и потребность в выделении подпоследовательности $\{y_{i_k}\}$, $k\in K$, отпадает. Отметим также, что, если в методе положить $\varepsilon_0=0$, то множества $Q_i$ и $M_{i+1}$ при всех $i$ будут иметь вид (\ref{e2.1.4}) и (\ref{e2.1.21}) соответственно, и ни одна из точек $y_{i_k}=x_k$ построена не будет.

Учитывая эти замечания, перейдем, наконец, к исследованию таких последовательностей $\{(y_i, \gamma_i)\}$, $i\in K$, при построении которых бесконечное число раз выполнялось (\ref{e2.1.5}) и происходило обновление множеств $M_{i+1}$ за счет соответствующего выбора $Q_i$. Другими словами, исследуем свойства последовательностей $\{x_k\}$, $\{\sigma_k\}$, $k\in K$.

Прежде всего сформулируем условие, при котором наряду с $\{(y_i, \gamma_i)\}$, $i\in K$, будет построена согласно (\ref{e2.1.5}) и последовательность $\{(x_k, \sigma_k)\}$, $k\in K$.

\begin{lemma}\label{lemma2.1.4}
Пусть последовательность $\{(y_i, \gamma_i)\}$, $i\in K$, построена предложенным методом с условием, что
\begin{equation}\label{e2.1.22}
\varepsilon_k>0 \quad \forall k\in K.
\end{equation}
Тогда для каждого $k\in K$ существует такой номер $i=i_k$, что выполняются равенства (\ref{e2.1.5}).
\end{lemma}

\begin{proof} 1) Пусть $k=0$. Если $f(y_0)-\gamma_0\le \varepsilon_0$, то согласно п. 2 метода $i_0=0$, $x_0=y_{i_0}=y_0$, $\sigma_0=\gamma_0$, и равенства (\ref{e2.1.5}) при $k=0$ выполняются. Поэтому будем считать, что $f(y_0)-\gamma_0>\varepsilon_0$. Покажем тогда существование номера $i=i_0>0$, для которого справедливо неравенство 
\begin{equation}\label{e2.1.23}
f(y_{i_0})-\gamma_{i_0}\le \varepsilon_0.
\end{equation}

Допустим противное, то есть 
\begin{equation}\label{e2.1.24}
f(y_i)>\varepsilon_0 + \gamma_i \quad \forall i\in K, \quad i>0.
\end{equation}
Выделим из ограниченной последовательности $\{(y_i, \gamma_i)\}$, $i\in K$, $i>0$, сходящуюся подпоследовательность $\{(y_i, \gamma_i)\}$, $i\in K^\prime$, и пусть $(y^\prime, \gamma^\prime)$ -- ее предельная точка. Тогда из неравенств (\ref{e2.1.24}) следует, что 
\begin{equation}\label{e2.1.25}
f(y^\prime)\ge \varepsilon_0+\gamma^\prime.
\end{equation}
С другой стороны, в силу (\ref{e2.1.24}) для всех $i\in K$, $i>0$, выполняется (\ref{e2.1.4}), а значит, по следствию \ref{consequence2.1.1} имеет место равенство $f(y^\prime)=\gamma^\prime$, противоречащее (\ref{e2.1.25}). Таким образом, существование номера $i_0$, для которого справедливо (\ref{e2.1.23}), доказано, и равенства (\ref{e2.1.5}) имеют место для $k=0$.

2) Допустим теперь, что (\ref{e2.1.5}) выполняются при некотором фиксированном $k\ge 0$, то есть  найден номер $i=i_k$, удовлетворяющий условию (\ref{e2.1.10}) при этом $k$. Покажем существование такого $i_{k+1}>i_k$, что
\begin{equation}\label{e2.1.26} 
f(y_{i_{k+1}})-\gamma_{i_{k+1}}\le \varepsilon_{k+1},
\end{equation}
тогда $x_{k+1}=y_{i_{k+1}}$, $\sigma_{k+1}=\gamma_{i_{k+1}}$, и лемма будет доказана.

Допустим противное, то есть 
\begin{equation}\label{e2.1.27}
f(y_i)-\gamma_i>\varepsilon_{k+1} \quad \forall i>i_k.
\end{equation}
Выберем среди точек $(y_i, \gamma_i)$, $i\in K$, $i>i_k$, сходящуюся подпоследовательность $\{(y_i, \gamma_i)\}$, $i\in K^{\prime\prime}$, и пусть $(y^{\prime\prime}, \gamma^{\prime\prime})$ -- ее предельная точка. Тогда из (\ref{e2.1.27}) следует с одной стороны неравенство $f(y^{\prime\prime})\ge \varepsilon_{k+1}+\gamma^{\prime\prime}$, а с другой стороны, как и в первой части доказательства леммы, -- противоречащее ему равенство $f(y^{\prime\prime})=\gamma^{\prime\prime}$. Таким образом, показано наличие номера $i_{k+1}>i_k$, удовлетворяющего (\ref{e2.1.26}). Лемма доказана.
\end{proof}

Итак, в лемме \ref{lemma2.1.4} приведено условие, гарантирующее существование последовательности $\{(x_k, \sigma_k)\}$, $k\in K$. Обоснуем теперь ее сходимость при дополнительном предположении о последовательности $\{\varepsilon_k\}$, $k\in K$.

\begin{theorem}\label{theorem2.1.3}
Пусть последовательность $\{(x_k, \sigma_k)\}$, $k\in K$, построена предложенным методом с условием, что числа $\varepsilon_k$, $k\in K$, выбраны согласно (\ref{e2.1.22}), и 
\begin{equation}\label{e2.1.28}
\varepsilon_k\to 0, \quad k\to \infty.
\end{equation}
Тогда 
\begin{equation}\label{e2.1.29}
\lim \limits_{k\in K} f(x_k)=f^*, \quad \lim \limits_{k \in K}\sigma_k = f^*.
\end{equation}
\end{theorem}

\begin{proof}
Ввиду (\ref{e2.1.5}) $f(y_{i_k})-\gamma_{i_k}\le \varepsilon_k$ или
\begin{equation}\label{e2.1.30}
f(x_k)\le \sigma_k+\varepsilon_k
\end{equation}
для всех $k\in K$. Кроме того, согласно (\ref{e2.1.13}), (\ref{e2.1.14})
\begin{equation}\label{e2.1.31}
\sigma_k\le f^*, \quad f(x_k)\ge f^* \quad \forall k\in K.
\end{equation}
Из (\ref{e2.1.30}), (\ref{e2.1.31}) для всех $k\in K$ имеем неравенства $f^*\le f(x_k)\le f^*+\varepsilon_k$. Тогда с учетом (\ref{e2.1.28}) $f^*\le \lim \limits_{\bar{k\in K}} f(x_k) \le \bar{\lim \limits_{k\in K}} f(x_k)\le f^*$, и первое утверждение доказано.

Далее, в силу тех же неравенств (\ref{e2.1.30}), (\ref{e2.1.31}) $f^*-\varepsilon_k\le \sigma_k\le f^*$ для всех $k\in K$. Отсюда следует второе из равенств (\ref{e2.1.29}). Теорема доказана.
\end{proof}

Теперь обсудим возможности задания последовательности $\{\varepsilon_k\}$, $k\in K$. 
\begin{remark}\label{r2.1.8}
В методе числа $\varepsilon_k$, $k\ge 1$, можно задать, как и число $\varepsilon_0$, на предварительном шаге. Но в этом случае $\{\varepsilon_k\}$, $k\in K$, не будет адаптирована к процессу минимизации. Поэтому в  методе  заложена возможность выбора чисел $\varepsilon_k$, $k\ge 1$, на шаге 2, то есть  в процессе построения приближений $x_k$, $k\in K$. Приведем пример такого способа  задания $\{\varepsilon_k\}$, $k\in K$. 

Пусть $\varepsilon_0$ -- сколь угодно большое положительное число, тогда $x_0=y_0$, $\sigma_0=\gamma_0$. Положим для всех $k\ge 0$ в п. 2 метода 
$$\varepsilon_{k+1} = \alpha_k(f(x_k) - \sigma_k),$$
где $0\le \alpha_k<1$. Отметим, что для такой последовательности $\{\varepsilon_k\}$ справедливо (\ref{e2.1.22}), и выполняется условие (\ref{e2.1.28}), если $\alpha_k \to 0$, $k\in K$.
\end{remark}

В заключение приведем оценки, касающиеся точности решения задачи. Согласно (\ref{e2.1.14}) при всех $i\in K$, в том числе и при $i = i_k$, выполняются неравенства
$$\gamma_i \le f^* \le f(y_i).$$
Следовательно, с учетом (\ref{e2.1.6})
\begin{equation}\label{e2.1.33}
0\le f(x_k) - f^* \le f(x_k) - \sigma_k\le\varepsilon_k, \quad k\in K.
\end{equation}

Пусть теперь $f_j(x)$, $j\in J$, -- сильно выпуклые функции с константами сильной выпуклости $\mu_j$ соответственно . Тогда согласно лемме \ref{lemma1.1.4} функция $f(x)$  сильно выпукла с константой сильной выпуклости $\mu = \min \limits_{j\in J} \mu_j$. C учетом этого условия справедлива следуюшая
\begin{theorem}\label{theorem2.1.4}
Пусть последовательность $\{x_k\}$, $k\in K$, построена методом 2.1. Тогда для каждого $k\in K$ имеет место оценка
\begin{equation}\label{e2.1.34}
\|x_k - x^*\|\le \sqrt{\frac{\varepsilon_k}{\mu}}.
\end{equation}
\end{theorem}
\begin{proof}
Справедливость утверждения следует из известного неравенства $\mu \|x_k - x^*\|^2 \le f(x_k) - f(x^*)$ (напр., \cite[182]{n_vasillev_f_p_2}) и оценки (\ref{e2.1.33}).
\end{proof}

Далее приведем оценки скорости сходимости для последовательности $\{x_k\}$, $k\in K$. Пусть в методе числа $\varepsilon_k$, $k\in K$, заданы следующим образом:
\begin{equation}\label{e2.1.35}
\varepsilon_k = \frac{1}{k^p}, \quad k\ge 1,
\end{equation}
где $p\ge1$. Тогда с учетом (\ref{e2.1.33}), (\ref{e2.1.34}), (\ref{e2.1.35}) справедливы оценки скорости сходимости
\begin{equation}\label{e2.1.36}
\|x_k - x^*\|\le \sqrt{\frac{1}{\mu k ^p}},\quad f(x_k) - f^* \le \frac{1}{k^p}, \quad k\ge 1.
\end{equation}
\setcounter{equation}{0}

\section
{Метод отсечений с обновлением аппроксимирующих множеств, использующий обобщенно-опорные элементы}

Пусть $D$ -- выпуклое замкнутое  множество из $n$-мерного евклидова пространства ${\rm R_n}$, $f(x)$ -- достигающая на множестве $D$ своего минимального значения выпуклая функция. Решается задача (\ref{e2.1.0}).

Пусть $f^*$, $X^*$, ${\rm epi\,} (f, {\rm R_n})$, $\partial f(x)$, $K$ определены как и в $\S \, 2.1$. Положим $X^*_\varepsilon = \{x\in D : f(x)\le f^* + \varepsilon\}$, где $\varepsilon\ge 0$, $W(u, Q)=\{a\in {\rm R_{n+1}} : \langle a, z-u\rangle \le 0 \ \forall z\in Q\}$ -- множество обобщенно-опорных к множеству $Q\subset {\rm R_{n+1}}$ в точке $u\in {\rm R_{n+1}}$ векторов, $W^1(u, Q) =\{a\in W(u, Q) : \|a\|=1\}$,   $J=\{1,\dots,m\}$, где $m\ge 1$. Пусть $x^*\in X^*$.

Предлагаемый метод решения задачи (\ref{e2.1.0}) вырабатывает последовательности $\{y_i\}$, $i\in K$, $\{x_k\}$, $k\in K$, приближений из множества $D$, и заключается в следующем. 

\textbf{Метод 2.2.} Выбираются точки $v^j\in{\rm int\, epi\,}(f, {\rm R_n})$ для всех $j\in J$ и выпуклое ограниченное замкнутое множество $G_0\subset D$, содержащее точку $x^*$. Строится выпуклое замкнутое множество $M_0\subset {\rm R_{n+1}}$ такое, что
$$(x^*, f^*)\in M_0.$$
Задаются числа $\varepsilon_0>0$, $\bar{\gamma}_0\le f^*$, и полагается $i=0$, $k=0$.

\textbf{1.} Отыскивается точка $(y_i,\gamma_i)$, где $y_i\in {\rm R_{n}}$, $\gamma_i\in {\rm R_1}$, как решение задачи
\begin{equation}\label{e2.2.1}
\gamma \to \min
\end{equation}
\begin{equation}\label{e2.2.2}
x\in G_i,\quad  (x,\gamma)\in M_i, \quad \gamma\ge \bar{\gamma}_i.
\end{equation}
Если
\begin{equation}\label{e2.2.3}
f(y_i)=\gamma_i,
\end{equation}
то $y_i\in X^*$, и процесс заканчивается.

\textbf{2.} Если выполняется неравенство
\begin{equation}\label{e2.2.4}
f(y_i)-\gamma_i>\varepsilon_k,
\end{equation}
то выбирается выпуклое замкнутое множество $Q_i\subset {\rm R_{n+1}}$ такое, что 
\begin{equation}\label{e2.2.5}
Q_i\subset M_i,
\end{equation}
\begin{equation}\label{e2.2.6}
(x^*, f^*)\in Q_i,
\end{equation}
полагается
\begin{equation}\label{e2.2.7}
u_i=y_i,
\end{equation}
и следует переход к п. 4. В противном случае выполняется п. 3.

\textbf{3.} Выбирается выпуклое замкнутое множество $Q_i\subset {\rm R_{n+1}}$ так, чтобы выполнялось включение (\ref{e2.2.6}). Выбирается точка $x_k\in {\rm  R_n}$, удовлетворяющая условиям
\begin{equation}\label{e2.2.8}
x_k\in G_i, \quad f(x_k)\le f(y_i).
\end{equation}
Полагается
\begin{equation}\label{e2.2.9}
i_k = i, \quad \sigma_k=\gamma_{i_k},
\end{equation}
\begin{equation}\label{e2.2.10}
u_i=u_{i_k}=x_k,
\end{equation}
задается число $\varepsilon_{k+1}>0$, значение $k$ увеличивается на единицу.

\textbf{4.} Для каждого $j\in J$ в интервале $(v^j, (u_i,\gamma_i))$ выбирается точка $z_i^j\in {\rm R_{n+1}}$ так, чтобы $z_i^j\notin {\rm int\, epi\,}(f, {\rm R_n})$ и при некотором $q_i^j\in [1,q]$, $q<+\infty$, для точки 
\begin{equation}\label{e2.2.11}
\bar{v}_i^j=(u_i,\gamma_i) + q_i^j(z_i^j-(u_i,\gamma_i))
\end{equation}
выполнялось включение $\bar{v}_i^j\in {\rm epi\,} (f, {\rm R_n})$.

\textbf{5.} Для каждого $j\in J$ выбирается конечное множество $A_i^j\subset W^1(z_i^j$, ${\rm epi\,} (f, {\rm R_n}))$, и полагается
\begin{equation}\label{e2.2.12}
M_{i+1} = Q_i\bigcap T_i,
\end{equation}
где $T_i=\bigcap \limits_{j\in J} \{(x,\gamma)\in {\rm R_{n+1}} : \langle a, (x,\gamma) - z_i^j\rangle \le 0 \ \forall a\in A_i^j\}.$

\textbf{6.} Выбирается выпуклое замкнутое множество $G_{i+1}\subset G_0$, содержащее точку $x^*$. Задается число $\bar{\gamma}_{i+1}$ из условия
\begin{equation}\label{e2.2.13}
\bar{\gamma}_0\le \bar{\gamma}_{i+1}\le f^*.
\end{equation}
Значение $i$ увеличивается на единицу, и следует переход к п. 1.

Сделаем некоторые замечания к методу. Прежде всего покажем, что множество ограничений (\ref{e2.2.2}) непусто, то есть задача (\ref{e2.2.1}), (\ref{e2.2.2}) при всех $i\in K$ разрешима.

\begin{lemma}\label{lemma2.2.1}
Точка $(x^*, f^*)$ для всех $i\in K$ принадлежит множеству ограничений задачи (\ref{e2.2.1}), (\ref{e2.2.2}). 
\end{lemma}

\begin{proof}Согласно заданию $\bar{\gamma}_0$ и условию (\ref{e2.2.13}) $f^*\ge \bar{\gamma}_i$, $i\in K$. По выбору множеств $G_i$, $i\in K$, выполняется включение $x^*\in G_i$. Поэтому для обоснования леммы покажем, что 
\begin{equation}\label{e2.2.14}
(x^*, f^*)\in M_i
\end{equation}
для всех $i\in K$.

При $i=0$ включение (\ref{e2.2.14}) выполняется по выбору множества $M_0$. Далее, для множеств $Q_i$ в силу пп. 2, 3 метода имеет место включение (\ref{e2.2.6})  при всех $i\ge 0$. Кроме того, $\langle a, (x^*, f^*) - z_i^j\rangle \le 0$ для всех $a\in A_i^j$, $j\in J$, $i\ge 0$, то есть  $(x^*, f^*)\in T_i$ для всех $i\ge 0$. Таким образом, ввиду (\ref{e2.2.12}) $(x^*, f^*)\in M_{i+1}$ для всех $i\ge 0$. Лемма доказана.
\end{proof}

Обоснуем критерий остановки, заложенный в п. 1 метода.

\begin{lemma}\label{lemma2.2.2}
Пусть последовательность $\{\gamma_i\}$, $i\in K$, построена предложенным методом. Тогда
\begin{equation}\label{e2.2.15}
\gamma_i\le f^* \quad \forall i\in K.
\end{equation}
\end{lemma}

\begin{proof}
Для каждой точки $(x,\gamma)$, удовлетворяющей условиям (\ref{e2.2.2}), и каждого $i\in K$ выполняется неравенство $\gamma_i\le \gamma$. Но по лемме \ref{lemma2.2.1} $(x^*, f^*)$ -- допустимое решение задачи (\ref{e2.2.1}), (\ref{e2.2.2}) при любом $i\in K$. Отсюда следует (\ref{e2.2.15}). Лемма доказана.
\end{proof}

\begin{theorem}\label{theorem2.2.1}
Если при некотором $i\in K$ выполняется равенство (\ref{e2.2.3}), то $y_i\in X^*$.
\end{theorem}

\begin{dokaz}  вытекает из неравенства (\ref{e2.2.15}) и включения $y_i\in D$, $i\in K$, c учетом предположения (\ref{e2.2.3}).
\end{dokaz}

Следующая теорема дает способ построения в точках $z_i^j$ обобщенно-опорных к множеству ${\rm epi\,} (f, {\rm R_n})$ векторов на шаге 5 метода 2.2.

\begin{theorem}\label{ptheorem2.2.2}
Пусть $\varphi(x)$ -- выпуклая в ${\rm R_n}$  функция, точка $\bar{x}\in {\rm R_n}$ и число $\bar{\gamma}$ таковы, что 
\begin{equation}\label{p2.2.22}
\varphi(\bar{x})\ge \bar{\gamma},
\end{equation}
а $c_\varphi$ -- субградиент функции $\varphi(x)$ в точке $\bar{x}$. Тогда вектор $c=(c_\varphi, -1)$ является обобщенно-опорным к множеству ${\rm epi\,} (\varphi, {\rm R_n})$  в точке $\bar{u}=(\bar{x}, \bar{\gamma})$. 
\end{theorem}

\begin{proof}
Для обоснования утверждения зафиксируем точку $u=(x, \gamma)\in {\rm epi\,}(\varphi, {\rm R_n})$ и покажем, что
\begin{equation}\label{p2.2.23}
\langle c, u - \bar{u}\rangle \le 0.
\end{equation}
Так как $\langle c, u - \bar{u}\rangle = \langle (c_\varphi, -1), (x-\bar{x}, \gamma - \bar{\gamma})\rangle = \langle c_\varphi, x - \bar{x}\rangle - \gamma + \bar{\gamma}$, и по условию $\varphi(x) - \varphi(\bar{x}) \ge \langle c_\varphi, x - \bar{x}\rangle$, то
$$\langle c, u - \bar{u}\rangle \le \varphi(x) - \varphi(\bar{x}) - \gamma +\bar{\gamma}.$$ 
Отсюда с учетом неравенства (\ref{p2.2.22}) и неравенства $\varphi(x)\le \gamma$, вытекающего из условия выбора точки $u$, следует (\ref{p2.2.23}). Теорема доказана.
\end{proof}

\begin{remark}\label{r2.2.1}
Выбор в методе чисел $\bar{\gamma_i}$, $i\in K$, не представляет особого труда. Число $\bar{\gamma}_0$ можно выбрать, например, как решение задачи минимизации переменной $\gamma$ при ограничении 
\begin{equation}\label{p17}
\langle c, x\rangle -\gamma\le \langle c, u\rangle - f(u),
\end{equation}
где $u\in {\rm R_n}$, $c\in \partial f(u)$, и условии $x\in G_0$. Если задача (\ref{e2.2.1}), (\ref{e2.2.2}) при $i=0$ имеет решение и без ограничения $\gamma\ge \bar{\gamma}_0$, например, в случае ограниченности $M_0$, то число $\bar{\gamma}_0$ можно считать сколь угодно большим отрицательным. 
\end{remark}

\begin{remark}\label{r2.2.2}
При $i\ge 0$ лемма \ref{lemma2.2.2} позволяет выбирать числа $\bar{\gamma}_{i+1}$ из условия (\ref{e2.2.13}) следующим образом. Как и в замечании \ref{r2.1.6} можно положить $\bar{\gamma}_{i+1}=\gamma_l$, где $0\le l \le i$, в частности, $\bar{\gamma}_{i+1} = \max \limits_{0\le j\le i}\gamma_j$.
\end{remark}

Обсудим теперь способы задания множеств $M_0$, $G_i$, $Q_i$, $i\in K$. Понятно, что эти множества удобно строить так, чтобы задачи (\ref{e2.2.1}), (\ref{e2.2.2}) при всех $i\in K$ были задачами линейного программирования.

\begin{remark}\label{r2.2.3}
Для выбора множества $M_0$ имеется много возможностей. Его можно выбрать содержащим множество ${\rm epi\,}(f, G_0)$ или ${\rm epi\,} (f, {\rm R_n})$, задав, например, линейным неравенством (\ref{p17}) или группой подобных неравенств. Допустимо положить 
$$M_0={\rm R_{n+1}},$$
тогда пара $(y_0,\gamma_0)$, где $y_0$ -- любая точка из $G_0$, а $\gamma_0=\bar{\gamma}_0$, может быть принята за решение задачи (\ref{e2.2.1}), (\ref{e2.2.2}) при $i=0$.
\end{remark}

\begin{remark}\label{r2.2.4}
Множества $G_i$ можно задавать также разными способами. В частности, для всех $i\in K$ допустимо положить 
$$G_i=D,$$
когда $D$ ограничено. Если в задаче (\ref{e2.1.0}) $D={\rm R_n}$, то многогранник $G_0$ необходимо построить на предварительном шаге метода содержащим точку безусловного минимума функции $f(x)$. Множество $G_{i+1}$ при $i\ge 0$, можно строить в следующем виде:
$$G_{i+1}=G_r\bigcap \{x\in {\rm R_n} :\langle b, x - y_i\rangle \le 0 \ \forall b\in B_i\},$$
где $0\le r\le i$, $B_i\subset \partial f(y_i)$. Так как для любых $b\in B_i$ выполняются неравенства $0\ge f(x^*) - f(y_i)\ge \langle b, x^* - y_i\rangle$, и $x^*\in G_r$, то $x^*\in G_{i+1}$, как требуется в п. 6 метода.
\end{remark}

Подробнее остановимся на способах задания множеств $Q_i$, поскольку именно за счет их выбора удается периодически обновлять аппроксимирующие надграфик множества $M_i$.

Обратим внимание на то, что согласно пп. 2, 3 метода, независимо от выполнения условия (\ref{e2.2.4}), множества $Q_i$  можно выбирать для всех $i\in K$ из условий (\ref{e2.2.5}), (\ref{e2.2.6}), например, $Q_i=M_i$ или $Q_i=M_i\bigcap S_i$, где $S_i\subset {\rm R_{n+1}}$, $(x^*, f^*)\in S_i$. Однако, в таком случае от шага к шагу ввиду (\ref{e2.2.12}) неограниченно растет число отсекающих плоскостей, формирующих аппроксимирующие множества, и, как следствие, растет трудоемкость решения задач (\ref{e2.2.1}), (\ref{e2.2.2}) построения приближений.

Покажем теперь, как на итерациях с номерами $i=i_k$ за счет выбора множеств $Q_i$ проводить упомянутые обновления.

\begin{remark}\label{r2.2.5}
Пусть для точки $(y_i,\gamma_i)$ выполняется неравенство
\begin{equation}\label{e2.2.17}
f(y_i)-\gamma_i\le \varepsilon_k.
\end{equation}
В таком случае согласно п. 3 метода множество $Q_i=Q_{i_k}$ выбирается удовлетворяющим лишь условию (\ref{e2.2.6}). Положим, например, $Q_i={\rm R_{n+1}}$ или $Q_i=M_0$. Тогда, соответственно, $M_{i+1}=T_i$ или $M_{i+1}=M_0\bigcap T_i$, и в формировании $M_{i+1}$ не участвуют ни одна из построенных ранее отсекающих плоскостей. При условии (\ref{e2.2.17})  множество $Q_i=Q_{i_k}$ можно задать и в виде
$$Q_i=M_{r_i},$$
где $0<r_i\le i = i_k$, так как при всех $r_i=1,\dots, i$ в силу (\ref{e2.2.14}) включение (\ref{e2.2.6}) выполняется. В таком случае при построении множества $M_{i+1}$ отбрасывается только часть накопленных к шагу $i = i_k$ отсечений. Как будет доказано ниже, для каждого $k\in K$ найдется решение $(y_i, \gamma_i)=(y_{i_k}, \gamma_{i_k})$ задачи (\ref{e2.2.1}), (\ref{e2.2.2}), удовлетворяющее условию (\ref{e2.2.17}), а значит, представится возможность обновления аппроксимирующих множеств.
\end{remark}

Отметим, что условие (\ref{e2.2.17}) фактически определяет качество аппроксимации ${ \rm epi\,} (f, {\rm R_n})$ множеством $M_i$ в окрестности $y_i$. При выполнении (\ref{e2.2.17}) качество аппроксимации считается достаточным для фиксирования точки $y_{i_k}$, а значит, и основной итерационной точки $x_k$.

Попутно подчеркнем также, что на итерациях с номерами $i=i_k$ отсечения множеств $M_{i_k}$ проводятся и c использованием точек $u_{i_k}$ вида (\ref{e2.2.10}), что удобно с практической точки зрения.

Перейдем к исследованию сходимости метода.

Докажем, что вместе с последовательностью $\{(y_i, \gamma_i)\}$, $i\in K$, построенной методом, будет построена и последовательность $\{(x_k, \sigma_k)\}$, $k\in K$.

\begin{lemma}\label{lemma2.2.3}
Пусть $\{(y_i, \gamma_i)\}$, $i\in K$, построена предложенным методом. Тогда для каждого $k\in K$ существует номер $i = i_k\in K$, для которого выполняется неравенство
\begin{equation}\label{e2.2.18}
f(y_{i_k}) - \gamma_{i_k}\le \varepsilon_k.
\end{equation}
\end{lemma}

\begin{proof}
Зафиксируем $k\in K$, и докажем существование номера $i_k\in K$, удовлетворяющего (\ref{e2.2.18}). Предположим противное, то есть, что для всех $i\in K$ имеет место неравенство (\ref{e2.2.4}). Сразу отметим, что в таком случае согласно п. 2 метода множества $Q_i$ для всех $i\in K$ выбраны из условия (\ref{e2.2.5}), и в силу (\ref{e2.2.12}) справедливы включения
\begin{equation}\label{e2.2.19}
M_{i+1}\subset M_i \quad \forall i\in K.
\end{equation}

Выделим из последовательности $\{(y_i, \gamma_i)\}$, $i\in K$, сходящуюся подпоследовательность $\{(y_i, \gamma_i)\}$, $i\in K^\prime\subset K$, и пусть $(\bar{y},\bar{\gamma})$ -- ее предельная точка. Тогда ввиду (\ref{e2.2.4})
\begin{equation}\label{e2.2.20}
f(\bar{y}) - \bar{\gamma}\ge \varepsilon_k.
\end{equation}

Покажем, что при сделанном предположении выполняется равенство 
\begin{equation}\label{e2.2.21}
\lim \limits_{i\in K^\prime} \|z_i^j - (y_i, \gamma_i)\| = 0 \quad \forall j\in J.
\end{equation}
Зафиксируем $j\in J$. Заметим, что для каждого $i\in K$ точки $u_i$ имеют вид (\ref{e2.2.7}), а значит, согласно выбору  точек $z_i^j$  существуют такие числа $\tau_i^j\in (0, 1)$, что 
\begin{equation}\label{e2.2.22}
z_i^j=(y_i, \gamma_i) + \tau_i^j (v^j - (y_i, \gamma_i)) \quad \forall i\in K.
\end{equation}
Выберем номера $i$, $p_i\in K^\prime$ с условием, что $p_i>i$. Тогда ввиду (\ref{e2.2.19}) $M_{p_i}\subset M_i$. Кроме того, выполняется включение $(y_{p_i}, \gamma_{p_i})\in M_{p_i}$, а любой элемент $a\in A_i^j$ является обобщенно-опорным к множеству $M_{p_i}$  в точке $z_i^j$. Следовательно,
$$\langle a, (y_{p_i}, \gamma_{p_i}) - z_i^j \rangle \le 0 \quad \forall a\in A_i^j,$$
а с учетом (\ref{e2.2.22})
$$\langle a, (y_i, \gamma_i) - (y_{p_i}, \gamma_{p_i})\rangle \ge \tau_i^j \langle a, (y_i, \gamma_i) - v^j\rangle \quad \forall a\in A_i^j.$$
Согласно лемме \ref{lemma1.1.1} найдется такое число $\delta^j>0$, что $\langle a, (y_i, \gamma_i) - v^j\rangle \ge \delta^j$ для всех $a\in A_i^j$. Поэтому 
$\langle a, (y_i, \gamma_i)-(y_{p_i},\gamma_{p_i}\rangle \ge \tau_i^j \delta^j$ для всех $a\in A_i^j$, а так как $\|a\|=1$, то 
\begin{equation}\label{e2.2.23}
\|(y_i, \gamma_i) - (y_{p_i}, \gamma_{p_i})\|\ge \tau_i^j \delta^j \quad \forall i, p_i\in K^\prime, \quad p_i>i.
\end{equation}
Поскольку $\{(y_i, \gamma_i)\}$, $i\in K^\prime$, сходится, то в силу (\ref{e2.2.23}) $\tau_i^j\to 0$, $i\to \infty$, $i\in K^\prime$. Значит, из (\ref{e2.2.22}) с учетом ограниченности последовательности $\{\|v^j - (y_i,\gamma_i)\|\}$, $i\in K^\prime$, вытекает доказываемое утверждение (\ref{e2.2.21}).

Далее, так как последовательность $\{(y_i, \gamma_i)\}$, $i\in K^\prime$, ограничена, то из (\ref{e2.2.11}), (\ref{e2.2.7}), (\ref{e2.2.21}) следует ограниченность последовательностей $\{\bar{v}_i^j\}$, $i\in K^\prime$, для всех $j\in J$. Выделим теперь для некоторого $r\in J$ из последовательности $\{\bar{v}_i^r\}$, $i\in K^\prime$, сходящуюся подпоследовательность $\{\bar{v}_i^r\}$, $i\in K_r\subset K^\prime$. Пусть $\tilde{v}^r$ -- ее предельная точка. Отметим, что 
\begin{equation}\label{e2.2.24}
\tilde{v}^r\in {\rm epi\,} (f, {\rm R_n})
\end{equation}
в силу замкнутости множества ${\rm epi\,} (f, {\rm R_n})$. Перейдем в равенствах (\ref{e2.2.11}) при $j=r$ к пределу по $i\to \infty$, $i\in K_r$, с учетом (\ref{e2.2.7}), (\ref{e2.2.21}). Тогда $(\bar{y}, \bar{\gamma}) = \tilde{v}^r$, и ввиду (\ref{e2.2.24}) выполняется включение $(\bar{y}, \bar{\gamma})\in {\rm epi\,} (f, {\rm R_n})$. Cледовательно, $\bar{\gamma}\ge f(\bar{y})$. Но с другой стороны, $\bar{\gamma}\le f(\bar{y})$, поскольку для всех $i\in K^\prime$ выполняются неравенства (\ref{e2.2.15}) и $f^*\le f(y_i)$. Таким образом, получено равенство $f(\bar{y})=\bar{\gamma}$, противоречащее (\ref{e2.2.20}). Лемма доказана.
\end{proof}

Из леммы \ref{lemma2.2.3} следует, что для каждого $k\in K$ согласно п. 3 метода будут зафиксированы в виде (\ref{e2.2.9}) номер $i_k$ и число $\sigma_k$, а также построена точка $x_k$, удовлетворяющая условиям (\ref{e2.2.8}), где $i=i_k$.

\begin{theorem}\label{theorem2.2.2}
Пусть последовательность $\{(x_k, \sigma_k)\}$, $k\in K$, построена методом с условием, что
числа $\varepsilon_k$, $k\in K$, выбраны согласно (\ref{e2.1.28})
Тогда 
\begin{equation}\label{e2.2.26}
\lim \limits_{k\in K} f(x_k)=f^*, \quad \lim \limits_{k\in K} \sigma_k = f^*.
\end{equation}
\end{theorem}

\begin{proof}
Для каждого $k\in K$ выполняется неравенство (\ref{e2.2.18}). Кроме того, ввиду (\ref{e2.2.8}) $f(x_k)\le f(y_{i_k})$, $k\in K$. Следовательно, с учетом (\ref{e2.2.9})
\begin{equation}\label{e2.2.27}
f(x_k)\le \sigma_k + \varepsilon_k \quad \forall k\in K.
\end{equation}
Согласно лемме \ref{lemma2.2.2} и включению $x_k\in D$ справедливо
\begin{equation}\label{e2.2.28}
\sigma_k\le f^*,\quad \quad f(x_k)\ge f^* \quad \forall k\in K.
\end{equation}
Из (\ref{e2.2.27}), (\ref{e2.2.28}) для всех $k\in K$ имеем $f^*\le f(x_k)\le f^* + \varepsilon_k$. Тогда с учетом (\ref{e2.1.28}) $f^*\le \lim \limits_{\overline{k\in K}} f(x_k)\le \overline{\lim \limits_{k\in K}} f(x_k) \le f^*$, и первое из равенств (\ref{e2.2.26}) доказано. Далее, в силу тех же неравенств (\ref{e2.2.27}), (\ref{e2.2.28}) $f^*\le \sigma_k + \varepsilon_k \le f^* + \varepsilon_k$, $k\in K$. Отсюда следует второе из утверждений (\ref{e2.2.26}). Теорема доказана.
\end{proof}

\begin{remark}\label{r2.2.6}
Числа $\varepsilon_k$ при $k>0$ могут быть выбраны, как и число $\varepsilon_0$, на предварительном шаге метода. Однако в методе заложена возможность их задания в процессе счета, что, по-видимому, более естественно. При проведении численных экспериментов последовательность $\{\varepsilon_k\}$, $k\in K$, задавалось, в частности, следующим образом. Число $\varepsilon_0$ считалось сколь угодно большим, то есть полагалось $x_0=y_0$, $\sigma_0=\gamma_0$, а затем при всех $k\ge 0$ значение $\varepsilon_{k+1}$ вычислялось в виде $\varepsilon_{k+1} = \alpha_k (f(x_k)-\sigma_k)$, где $\alpha_k\in (0,1)$.
\end{remark}

Обсудим теперь возможность смешивания предложенного метода с другими методами математического программирования.

\begin{remark}\label{r2.2.6.1}
Условие (\ref{e2.2.8}) дает большую свободу в выборе точек $x_k$, отличных от $y_{i_k}$. Оно позволяет привлекать для получения приближений $x_k$ известные алгоритмы условной минимизации выпуклых функций. Например, считая $y_{i_k}$ точкой начального приближения, можно для минимизации $f(x)$ на множестве $G_{i_k}$ проделать, начиная с $y_{i_k},$ определенное число шагов известным градиентным или субградиентным алгоритмом (напр., \cite{n_vasillev_f_p_1}, \cite{n_konnov_i_v_2},  \cite{n_polayk_b_t_1}) в зависимости от свойств целевой функции. Таким образом, на основе предложенного метода отсечений будет построен новый алгоритм, причем сходимость такого смешанного алгоритма уже является обоснованной в силу теоремы \ref{theorem2.2.2}.
\end{remark}

\begin{remark}\label{r2.2.6.2}
Общность условия (\ref{e2.2.8}) задания точек $x_k$ позволяет также применить на этапе их нахождения параллельные вычисления. А именно, можно на шаге $i=i_k$ любыми алгоритмами параллельно построить некоторое число вспомогательных точек $w^1_k,\dots, w_k^l$, $l\ge 1$, таких, что $w^j_k\in G_{i_k}$, $f(w^j_k)<f(y_{i_k})$, $j=1,\dots, l$, а затем искомую точку $x_k$ выбрать из условия $f(x_k)=\min \{f(w^1_k),\dots, f(w^l_k)\}$.
\end{remark}

В заключение получим для предложенного метода 2.2 оценку точности решения.

Заметим, что ввиду леммы \ref{lemma2.2.2}, включения $x_k\in D$ и равенств (\ref{e2.2.9}) для всех $k\in K$ справедливы неравенства
\begin{equation}\label{e2.2.29}
\sigma_k\le f^*\le f(x_k).
\end{equation}
Следовательно, при каждом $k\in K$ можно оценить близость значения $f(x_k)$ к оптимальному.

Пусть теперь $f(x)$ -- сильно выпуклая функция с константой сильной выпуклости $\mu$. Тогда с учетом теоремы \ref{theorem2.1.4} для последовательности $\{x_k\}$, $k\in K$, справедливы оценки точности решения 
\begin{equation}\label{e2.2.30}
\|x_k - x^*\|\le \sqrt{\frac{\varepsilon_k}{\mu}}, \quad k\in K.
\end{equation}

Теперь пусть числа $\varepsilon_k$ заданы в виде (\ref{e2.1.35}). Значит, ввиду (\ref{e2.2.29}) и (\ref{e2.2.30})  для последовательности $\{x_k\}$, $k\in K$, выполняются   оценки скорости сходимости (\ref{e2.1.36}).
          
\setcounter{equation}{0}


\section
{Модифицированный метод уровней}

Одним из известных методов отсечений, использующих аппроксимацию надграфика целевой функции, является так называемый метод уровней, предложенный в работе \cite{n_lemarechal_c_1} (см., также  \cite{n_nesterov_yu_e_3}).  Как и в большинстве методов, исследуемого класса, в нем отсутствует процедура обновления аппроксимирующих множеств, что затрудняет его практичекое использование.  \mbox{В настоящем} параграфе строится метод, обобщающий метод уровней, в который такую процедуру обновления удалось включить.

Пусть $D$ -- выпуклое ограниченное замкнутое множество из  ${\rm R_n}$, $f(x)$ -- выпуклая в ${\rm R_n}$ функция. Решается задача (\ref{e2.1.0}).

Пусть $f^*$, $X^*$, ${\rm epi\,} (f, D)$, $\partial f(x)$,  $K$ определены как и в $\S\ 2.1$, $x^*\in X^*$.

 Предлагаемый метод решения задачи (\ref{e2.1.0}) вырабатывает последовательности приближений $\{x_i\}$, $i\in K$, $\{z_k\}$, $k\in K$, и заключается в следующем.

\textbf{Метод 2.3.} Выбирается выпуклое замкнутое множество $M_0\subset {\rm  R_{n+1}}$ такое, что ${ \rm epi\,}(f,D)\subset M_0$. Задаются числа $\varepsilon_0$, $\alpha_0$, $\beta_{-1}$, удовлетворяющие условиям
$$\varepsilon_0\ge 0,\quad \alpha_0 \le f^*\le \beta_{-1}.$$
Полагается $\delta_0=+\infty$, $i=0$, $k=0$.

\textbf{1.} Отыскивается решение $(y_i, \gamma_i)$, где $y_i\in {\rm R_n}$, $\gamma_i\in {\rm R_1}$, следующей задачи:
\begin{equation}\label{e2.3.1}
\min\{\gamma : (x,\gamma)\in M_i, \ x\in D, \ \gamma\ge\alpha_i\}.
\end{equation}
Если 
\begin{equation}\label{e2.3.2}
f(y_i)=\gamma_i,
\end{equation}
то $y_i$ -- решение задачи (\ref{e2.1.0}), и процесс завершается.

\textbf{2.} Выбирается число $\lambda_i$ согласно условию
\begin{equation}\label{e2.3.3}
0<\lambda_i\le \bar{\lambda}<1,
\end{equation}
полагается
\begin{equation}\label{e2.3.4}
\beta_i = \min \{\beta_{i-1}, \delta_i\},
\end{equation}
\begin{equation}\label{e2.3.5}
l_i = (1-\lambda_i)\gamma_i + \lambda_i\beta_i.
\end{equation}
Строится множество $U_i$ следующим образом. В $U_i$ включается каждая из точек $x\in D$, для которой при некотором $\gamma_x\le l_i$ выполняется включение $(x,\gamma_x)\in M_i$.

\textbf{3.} Выбирается точка $x_i\in U_i$.

\textbf{4.} Если выполняется неравенство
\begin{equation}\label{e2.3.6}
f(x_i)-\gamma_i > \varepsilon_k,
\end{equation}
то полагается 
\begin{equation}\label{e2.3.7}
Q_i=M_i,
\end{equation}
и следует переход к п. 5. В противном случае полагается $i_k=i$, 
\begin{equation}\label{e2.3.8}
z_k=x_{i_k}, \quad \sigma_k=\gamma_{i_k},
\end{equation}
выбирается выпуклое замкнутое множество $Q_i\subset {\rm  R_{n+1}}$ такое, что 
\begin{equation}\label{e2.3.9}
(x^*, f^*)\in Q_i,
\end{equation}
задается $\varepsilon_{k+1}\ge 0$, значение $k$ увеличивается на единицу.

\textbf{5.} Полагается
\begin{equation}\label{e2.3.10}
M_{i+1}=Q_i\bigcap \{(x,\gamma)\in {\rm  R_{n+1}} : f(x_i) + \langle a_i, x - x_i\rangle \le \gamma\},
\end{equation}
где $a_i\in \partial f(x_i)$.

\textbf{6.} Полагается
\begin{equation}\label{e2.3.11}
\alpha_{i+1}=\gamma_i,
\end{equation}
\begin{equation}\label{e2.3.12}
\delta_{i+1} = f(x_i).
\end{equation}
Значение $i$ увеличивается на единицу, и следует переход к п. 1.

Прежде, чем обсуждать свойства строящихся методом последовательностей, сделаем некоторые замечания.

Сразу отметим, что непустота допустимого множества задачи (\ref{e2.3.1}) и множества $U_i$, а также критерий остановки, заложенный в п. 1 метода, будут обоснованы позже.

Если $D$ -- многогранное множество, и $M_0$,  $Q_i$, $i\in K$, выбирать так, чтобы множества $M_{i+1}$, $i\in K$, вида (\ref{e2.3.10}) также были многогранными, то (\ref{e2.3.1}) для всех $i\in K$ является задачей линейного программирования. Обсудим способы задания множеств $M_0$,  $Q_i$.

\begin{remark}\label{r2.3.1}
Для выбора множества $M_0$ имеются большие возможности. Прежде всего отметим, что можно положить
$$M_0={\rm  R_{n+1}}.$$
В таком случае пара $(y_0, \gamma_0)$, где $y_0$ -- любая точка из $D$, а $\gamma_0=\alpha_0$, может быть принята за решение задачи (\ref{e2.3.1}) при $i=0$. Можно задать множество $M_0$ одним линейным неравенством
\begin{equation}\label{e2.3.13}
\langle c, x\rangle - \gamma \le \langle c, u\rangle - f(u),
\end{equation}
где $u\in {\rm R_n}$, $c\in \partial f(u)$, или группой подобных неравенств. Если 
\begin{equation}\label{e2.3.14}
f(x)=\max \limits_{j\in J} f_j(x),
\end{equation}
где $J$ -- конечное множество номеров, $f_j(x)$, $j\in J$, -- выпуклые в ${\rm R_n}$ функции, то допустимо положить 
\begin{equation}\label{e2.3.15}
M_0={\rm epi\,}(f_{j_0},{\rm R_n}),
\end{equation}
где $j_0\in J$, или аналогично (\ref{e2.3.13}) $M_0$ можно задать неравенством $\langle a_{j_0}, x\rangle - \gamma \le \langle a_{j_0}, u\rangle - f_{j_0}(u)$, где $a_{j_0} \in \partial f_{j_0}(u)$. Если функция $f_{j_0}(x)$ линейна, то выбор (\ref{e2.3.15}) удобен с практической точки зрения.
\end{remark}

\begin{remark}\label{r2.3.2}
Число $\alpha_0$ можно выбрать, например, как решение задачи минимизации переменной $\gamma$ при ограничении (\ref{e2.3.13}), где $u\in D$, и условии $x\in D$. За $\alpha_0$ можно также принять минимальное значение на множестве $D$ одной из функций $f_{j_0}(x)$, $j_0\in J$, если $f(x)$ имеет вид (\ref{e2.3.14}). Такое значение нетрудно отыскать, если $D$ -- многогранник, а функция $f_{j_0}(x)$ линейна. 
\end{remark}

Как и ранее, множества $M_i$ будем называть  аппроксимирующими надграфик целевой функции, а плоскости $f(x_i) + \langle a_i, x- x_i\rangle = \gamma$ или $\langle a_i, x\rangle - \gamma = \langle a_i, x_i\rangle - f(x_i)$ -- отсекающими плоскостями.

Прежде, чем перейти к обсуждению способов задания множеств $Q_i$, которые позволяют обновлять аппроксимирующие множества, докажем одно вспомогательное утверждение.

\begin{lemma} \label{lemma2.3.1}Для всех $i\in K$ выполняется включение 
\begin{equation}\label{e2.3.16}
(x^*, f^*)\in M_i.
\end{equation}
\end{lemma}

\begin{proof}
При $i=0$ включение (\ref{e2.3.16}) выполняется по условию выбора $M_0$. Допустим теперь, что (\ref{e2.3.16}) справедливо при любом фиксированном $i=l\ge 0$. Покажем выполнение (\ref{e2.3.16}) при $i=l+1$, тогда утверждение леммы будет доказано. Действительно, ввиду (\ref{e2.3.7}), (\ref{e2.3.9}),  а также сделанного индукционного предположения имеет место включение
\begin{equation}\label{e2.3.17}
(x^*, f^*)\in Q_l.
\end{equation}
Кроме того, $f(x^*)-f(x_l)\ge \langle a_l, x^* - x_l\rangle$. Значит, $f(x_l) + \langle a_l, x^* - x_l\rangle \le f^*$, и согласно (\ref{e2.3.10}), (\ref{e2.3.17}) $(x^*, f^*)\in M_{l+1}$. Лемма доказана.
\end{proof}

Далее, ввиду условия (\ref{e2.3.9}) и включения (\ref{e2.3.16}) множества $Q_i$ при всех $i\in K$ можно задавать в виде (\ref{e2.3.7}), независимо от выполнения неравенств (\ref{e2.3.6}). Однако, в таком случае от шага к шагу неограниченно растет число отсекающих плоскостей, формирующих аппроксимирующие множества, и вспомогательные задачи (\ref{e2.3.1}) становятся труднорешаемыми. Приведем теперь принципиальное замечание относительно задания множеств $Q_i$, на основе которого выявится возможность периодического обновления аппроксимирующих множеств за счет отбрасывания любого числа отсекающих плоскостей.

\begin{remark} \label{r2.3.3}
Для тех пар номеров $i$, $k$, при которых выполняются неравенства
\begin{equation}\label{e2.3.18}
f(x_i)-\gamma_i \le \varepsilon_k,
\end{equation}
условие (\ref{e2.3.9}) дает большие возможности в выборе множеств  $Q_i=Q_{i_k}$, а значит, и в задании множеств $M_{i+1}$. В случае (\ref{e2.3.18}) можно положить, например,
\begin{equation}\label{e2.3.19}
Q_i=Q_{i_k}={\rm  R_{n+1}}.
\end{equation}
Тогда $M_{i+1}$ задается лишь неравенством $f(x_i)+\langle a_i, x - x_i\rangle \le \gamma$. То есть в случае (\ref{e2.3.19}) все полученные к шагу $i=i_k$ отсекающие плоскости в построении $M_{i+1}$ не участвуют. Далее, при выполнении условия (\ref{e2.3.18}) множества $Q_i=Q_{i_k}$ можно задавать также в виде
\begin{equation}\label{e2.3.20}
Q_i=M_{r_i},
\end{equation}
 где $0\le r_i\le i=i_k$, поскольку при всех $r_i=0,\ldots , i$ ввиду (\ref{e2.3.16})   включение (\ref{e2.3.9}) выполняется. Согласно (\ref{e2.3.19}), (\ref{e2.3.20}) при каждом $i=i_k$ можно отбрасывать при формировании $M_{i_k+1}$ любое количество любых накопленных к шагу $i_k$ отсекающих плоскостей.
\end{remark}

Как будет показано ниже, для каждого $k\in K$ при выборе чисел $\varepsilon_k$ положительными найдется точка $(x_i,\gamma_i)$, удовлетворяющая (\ref{e2.3.18}), и, следовательно, представится возможность задания $Q_i$ в виде (\ref{e2.3.19}) или (\ref{e2.3.20}).

Для исследования сходимости метода изучим некоторые свойства строящихся методом последовательностей.

\begin{lemma}\label{lemma2.3.2}
Для последовательности $\{\alpha_i\}$, $i\in K$, построенной предложенным методом, справедливы неравенства
\begin{equation}\label{e2.3.21}
\alpha_i\le f^*
\end{equation}
при всех $i\in K$.
\end{lemma}
\begin{proof}
Неравенство (\ref{e2.3.21}) следует из вида ограничений задачи (\ref{e2.3.1}) и равенства (\ref{e2.3.11}) с учетом включения (\ref{e2.3.16}).
\end{proof}

Из лемм {\ref{lemma2.3.1}}, {\ref{lemma2.3.2}}  непосредственно вытекает
\begin{lemma}\label{lemma2.3.3}
Точка $(x^*, f^*)$ удовлетворяет ограничениям задачи (\ref{e2.3.1}) при всех $i\in K$.
\end{lemma}

На основе леммы \ref{lemma2.3.3} по схеме обоснования леммы \ref{e2.1.2}  доказывается

\begin{lemma}\label{lemma2.3.4}
Для построенной методом последовательности $\{\gamma_i\}$, $i\in K$, справедливы неравенства
\begin{equation}\label{e2.3.22}
\gamma_i\le f^* \quad \forall i\in K.
\end{equation}
\end{lemma}

Отметим, что ниже будет доказано существование последовательности $\{\sigma_k\}$, $k\in K$, при условии положительности чисел $\varepsilon_k$, $k\in K$, причем, по лемме \ref{lemma2.3.4} для нее выполняется неравенство 
\begin{equation}\label{e2.3.23}
\sigma_k \le f^* \quad \forall k\in K.
\end{equation}

На основе леммы {\ref{lemma2.3.4}} обоснуем теперь критерий остановки, заложенный в п. 1 метода. 

\begin{theorem}\label{theorem2.3.1}
Пусть при некотором $i\in K$ выполняется равенство (\ref{e2.3.2}). Тогда $y_i\in X^*$.
\end{theorem}
\begin{proof}
В силу (\ref{e2.3.2}), (\ref{e2.3.22}) $f(y_i)\le f^*$. С другой стороны, $f(y_i)\ge f^*$, так как $y_i\in D$. Таким образом, $f(y_i)=f^*$, и утверждение доказано.
\end{proof}

\begin{lemma}\label{lemma2.3.5}
Для каждого $i\in K$ множество $U_i$ непусто.
\end{lemma}
\begin{proof}
Зафиксируем произвольно номер $i=r\in K$. Для обоснования утверждения достаточно доказать существование точки $\bar{x}$ и числа $\bar{\gamma}$ таких, что
\begin{equation}\label{e2.3.24}
\bar{x}\in D, \quad \bar{\gamma}\le l_r, \quad (\bar{x}, \bar{\gamma})\in M_r.
\end{equation}

Cогласно (\ref{e2.3.4}), (\ref{e2.3.12}) $\beta_r\ge f^*$. Cледовательно, с учетом (\ref{e2.3.22}), (\ref{e2.3.5}), (\ref{e2.3.3}) справедливо неравенство
\begin{equation}\label{e2.3.25}
l_r\ge \gamma_r.
\end{equation}
Далее, в силу (\ref{e2.3.1}) выполняются включения $y_r\in D$ и $(y_r, \gamma_r)\in M_r$. Таким образом, положив $\bar{x}=y_r$, $\bar{\gamma}=\gamma_r$, получаем отсюда и из (\ref{e2.3.25})  соотношения (\ref{e2.3.24}). Лемма доказана.
\end{proof}

Как уже было отмечено, согласно п. 4 метода множества $Q_i$ для всех $i\in K$ могут иметь вид (\ref{e2.3.7}) независимо от выполнения условия (\ref{e2.3.6}). C учетом этого замечания сформулируем три следующих вспoмогательных  утверждения.

\begin{lemma}\label{lemma2.3.6}
Пусть последовательность $\{x_i\}$, $i\in K$, построена предложенным методом с условием, что для всех $i\in K$, начиная с некоторого номера $\tilde{i}>0$, множества $Q_i$ выбраны согласно (\ref{e2.3.7}). Тогда для любых номеров $i^\prime$, $i^{\prime\prime}\in K$ таких, что
\begin{equation}\label{e2.3.26}
i^{\prime\prime}>i^\prime\ge \tilde{i},
\end{equation}
справедливо неравенство 
\begin{equation}\label{e2.3.27}
f(x_{i^\prime})+\langle a_{i^\prime}, x_{i^{\prime\prime}}-x_{i^\prime}\rangle \le l_{i^{\prime\prime}}.
\end{equation}
\end{lemma}
\begin{proof}
Пусть $i^\prime$, $i^{\prime\prime}\in K$ выбраны с условием (\ref{e2.3.26}), и 
$$T_{i^\prime}=\{(x,\gamma)\in {\rm  R_{n+1}} : f(x_{i^\prime}) + \langle a_{i^\prime}, x-x_{i^\prime}\rangle \le \gamma\}.$$
Так как согласно п. 3 метода $x_{i^{\prime\prime}}\in U_{i^{\prime\prime}}$, то существует такое число $\gamma^{\prime\prime}$, что выполняются неравенство
\begin{equation}\label{e2.3.28}
\gamma^{\prime\prime}\le l_{i^{\prime\prime}}
\end{equation}
и включение $(x_{i^{\prime\prime}}, \gamma^{\prime\prime})\in M_{i^{\prime\prime}}$. Но в силу (\ref{e2.3.7}), (\ref{e2.3.10}), (\ref{e2.3.26}) $M_{i^{\prime\prime}}\subset T_{i^\prime}$, а значит, $(x_{i^{\prime\prime}}, \gamma^{\prime\prime}) \in T_{i^\prime}$. Отсюда и из (\ref{e2.3.28}) следует (\ref{e2.3.27}). Лемма доказана.
\end{proof}

\begin{lemma}\label{lemma2.3.7}
Пусть последовательности $\{x_i\}$, $\{\gamma_i\}$, $\{\beta_i\}$, $i\in K$, построены методом с условием, что, начиная с некоторого номера $\tilde{i}>0$, множества $Q_i$ выбраны согласно (\ref{e2.3.7}). Тогда выполняется равенство
$$\lim \limits_{i\in K} (\beta_i - \gamma_i)=0.$$
\end{lemma}
\begin{proof}
Отметим, что ввиду условий (\ref{e2.3.4}), (\ref{e2.3.12}) задания чисел $\beta_i$ и неравенств (\ref{e2.3.22}) $\beta_i - \gamma_i \ge 0$ для всех $i\in K$. Предположим, что утверждение леммы неверно. Тогда найдутся бесконечное подмножество номеров $K_1\subset K$ и число $\varepsilon > 0$ такие, что для всех $i\in K_1$, $i\ge \tilde{i}$, выполняется неравенство
\begin{equation}\label{e2.3.29}
\beta_i - \gamma_i \ge \varepsilon.
\end{equation}
Выделим из последовательности $\{x_i\}$, $i\in K_1\subset K$, сходящуюся подпоследовательность $\{x_i\}$, $i\in K_2\subset K_1$, и пусть $\bar{x}$ -- ее предельная точка.

Зафиксируем номера $i^\prime$, $i^{\prime\prime}\in K_2$ согласно (\ref{e2.3.26}). Поскольку в силу (\ref{e2.3.4}), (\ref{e2.3.12}) справедливо $\beta_{i^\prime + 1} \le \delta_{i^\prime+1}= f(x_{i^\prime})$, а последовательность $\{\beta_i\}$, $i\in K$, монотонно убывает, то
\begin{equation}\label{e2.3.30}
f(x_{i^\prime})\ge \beta_{i^{\prime\prime}}.
\end{equation}

Следуя \cite{n_nesterov_yu_e_3}, получим с учетом (\ref{e2.3.30}), (\ref{e2.3.5}) соотношения
$$f(x_{i^\prime}) - (1 - \lambda_{i^{\prime\prime}})(\beta_{i^{\prime\prime}}-\gamma_{i^{\prime\prime}})\ge \beta_{i^{\prime\prime}}  - (1 - \lambda_{i^{\prime\prime}})(\beta_{i^{\prime\prime}}-\gamma_{i^{\prime\prime}}) = l_{i^{\prime\prime}}.$$
Отсюда и из леммы {\ref{lemma2.3.6}} имеем неравенство
\begin{equation}\label{e2.3.31}
f(x_{i^\prime})-(1-\lambda_{i^{\prime\prime}})(\beta_{i^{\prime\prime}}-\gamma_{i^{\prime\prime}})\ge f(x_{i^\prime}) + \langle a_{i^\prime}, x_{i^{\prime\prime}}-x_{i^\prime}\rangle.
\end{equation}
Так как $x_{i}\in D$, $i\in K$, а множество $D$ ограничено, то найдется (см., напр., \cite[121]{n_polayk_b_t_1}) такое $\theta < +\infty$, что 
\begin{equation}\label{e2.3.32}
\|a\|\le \theta \quad \forall a\in \partial f(x_i), \quad i\in K.
\end{equation}
Тогда из (\ref{e2.3.31}), (\ref{e2.3.32}) следует, что $(1-\lambda_{i^{\prime\prime}})(\beta_{i^{\prime\prime}}-\gamma_{i^{\prime\prime}})\le \theta\|x_{i^{\prime\prime}} - x_{i^\prime}\|$, а поскольку $\beta_{i^{\prime\prime}}-\gamma_{i^{\prime\prime}}\ge \varepsilon$ согласно (\ref{e2.3.29}), то с учетом (\ref{e2.3.3})
\begin{equation}\label{e2.3.33}
\varepsilon \le \theta \frac{\|x_{i^{\prime\prime}}-x_{i^\prime}\|}{1-\bar{\lambda}}.
\end{equation}

Выберем теперь для каждого $i\in K_2$ такой номер $p_i\in K_2$, что $p_i\ge i+1$. Тогда в силу (\ref{e2.3.33}) $\varepsilon \le \theta\frac{ \|x_{p_i} - x_i\|}{1-\bar{\lambda}}$, $i\in K_2$. Последнее неравенство противоречит условию выбора $\varepsilon$, так как $x_i\to \bar{x}$, $x_{p_i} \to \bar{x}$ по $i\in K_2$. Лемма доказана.
\end{proof}

\begin{lemma}\label{lemma2.3.8}
Пусть выполняются условия леммы {\ref{lemma2.3.7}}. Тогда имеет место равенство $\lim \limits_{i\in K} (f(x_i)-\gamma_{i})=0$.
\end{lemma}
\begin{proof}
Отметим, что согласно (\ref{e2.3.22}) и включению $x_i \in D$ для всех $i\in K$ справедливы неравенства $f(x_i)-\gamma_i\ge 0$. Допустим, что утверждение леммы неверно. Тогда найдутся подпоследовательности $\{x_i\}$, $\{\gamma_i\}$, $i\in K_1\subset K$, последовательностей $\{x_i\}$, $\{\gamma_i\}$, $i\in K$, соответственно и число $\varepsilon>0$ такие, что 
$$f(x_i)-\gamma_i\ge \varepsilon$$
для всех $i\in K_1$, $i\ge\tilde{i}$.

Пусть $\{x_i\}$, $i\in K_2\subset K_1$, -- сходящаяся подпоследовательность последовательности $\{x_i\}$, $i\in K_1$, а $\bar{x}$ -- ее предельная точка. Зафиксируем номера $i^\prime$, $i^{\prime\prime}\in K_2$ из условия (\ref{e2.3.26}). Поскольку $f(x_{i^\prime})\ge \gamma_{i^\prime} + \varepsilon$, то из (\ref{e2.3.27}) следует, что $l_{i^{\prime\prime}}\ge \gamma_{i^\prime} + \varepsilon + \langle a_{i^\prime}, x_{i^{\prime\prime}}-x_{i^\prime}\rangle$. Тогда, учитывая (\ref{e2.3.32}), получаем неравенство 
$$l_{i^{\prime\prime}}\ge \gamma_{i^\prime} + \varepsilon - \theta\|x_{i^{\prime\prime}}-x_{i^\prime}\|.$$
C другой стороны, $l_{i^{\prime\prime}}=\gamma_{i^{\prime\prime}}+\lambda_{i^{\prime\prime}}(\beta_{i^{\prime\prime}}-\gamma_{i^{\prime\prime}})\le \beta_{i^{\prime\prime}}+\lambda_{i^{\prime\prime}}(\beta_{i^{\prime\prime}}-\gamma_{i^{\prime\prime}})\le \beta_{i^\prime} + \lambda_{i^{\prime\prime}}(\beta_{i^\prime}-\gamma_{i^{\prime\prime}})$. Значит,

\begin{equation}\label{e2.3.34}
\varepsilon - \theta \|x_{i^{\prime\prime}}-x_{i^\prime}\|\le \lambda_{i^{\prime\prime}}(\beta_{i^\prime}-\gamma_{i^{\prime\prime}}) + \beta_{i^\prime}-\gamma_{i^\prime}.
\end{equation}

Далее, в силу (\ref{e2.3.1}) $\alpha_{i+1}\le \gamma_{i+1}$, а согласно (\ref{e2.3.11}) $\alpha_{i+1}=\gamma_i$ для всех $i\in K$, то есть последовательность $\{\gamma_i\}$, $i\in K$,  монотонно возрастает. Следовательно, $\gamma_{i^\prime}\le \gamma_{i^{\prime\prime}}$, и из (\ref{e2.3.34}) с учетом (\ref{e2.3.3}) имеем неравенство
\begin{equation}\label{e2.3.35}
\varepsilon \le \theta \|x_{i^{\prime\prime}}-x_{i^\prime}\| + (1+\bar{\lambda})(\beta_{i^\prime}-\gamma_{i^\prime}).
\end{equation}

Зафиксируем теперь для каждого $i\in K_2$ такой номер $p_i\in K_2$, что $p_i\ge i + 1$. Тогда согласно (\ref{e2.3.35})
\begin{equation}\label{e2.3.36}
\varepsilon\le \theta \|x_{p_i}-x_i\|+(1+\bar{\lambda})(\beta_i-\gamma_i) \quad \forall i\in K_2.
\end{equation}
Так как $x_i\to \bar{x}$ и $x_{p_i}\to \bar{x}$ при $i\to \infty$, $i\in K_2$, а  по лемме {\ref{lemma2.3.7}} выполняется равенство $\lim \limits_{i\in K_2} (\beta_i-\gamma_i)=0$, то из (\ref{e2.3.36}) следует неравенство $\varepsilon\le 0$, противоречащее выбору числа $\varepsilon$. Лемма доказана.
\end{proof}

На основе полученных выше вспомогательных результатов перейдем к исследованию сходимости метода. Прежде всего приведем результат, касающийся случая, когда при построении $\{x_i\}$, $i\in K$, начиная с некоторого номера, множества $Q_i$ выбираются в виде (\ref{e2.3.7}).

\begin{theorem}\label{theorem2.3.2}
Пусть выполняются условия леммы \ref{lemma2.3.7}. Тогда справедливы равенства
$$\lim \limits_{i\in K} f(x_i)=f^*, \quad \lim \limits_{i\in K} \gamma_i=f^*.$$
\end{theorem}

\begin{dokaz}
 теоремы \ref{theorem2.3.2} аналогично обоснованию теоремы \ref{theorem2.1.2}, и  следует из неравенства (\ref{e2.3.22}) и леммы \ref{lemma2.3.8}.
\end{dokaz}

Заметим, что согласно теореме {\ref{theorem2.3.2}} и известной теореме (\cite[62]{n_vasillev_f_p_1}) при условиях леммы {\ref{lemma2.3.7}} последовательность $\{x_i\}$, $i\in K$, сходится к множеству $X^*$.

Напомним, что утверждения теоремы {\ref{theorem2.3.2}} доказаны при условии, что в методе полагается для всех $i\in K$
\begin{equation}\label{e2.3.37}
M_{i+1}=M_i\bigcap \{(x,\gamma)\in {\rm  R_{n+1}} : f(x_i) + \langle a_{i}, x-x_i\rangle \le \gamma\}.
\end{equation}
В таком случае при построении последовательностей $\{x_i\}$, $\{\gamma_i\}$, $i\in K$,  число отсекающих плоскостей неограниченно растет и обновлений аппроксимирующих множеств не происходит. Отметим также, что, если в методе положить $\varepsilon_0=0$, то множества $Q_i$ и $M_{i+1}$ при всех $i\in K$ будут иметь вид (\ref{e2.3.7}) и (\ref{e2.3.37}) соответственно, и ни одна из точек $z_k=x_{i_k}$ зафиксирована не будет.

Учитывая эти замечания, перейдем к исследованию последовательностей $\{x_i\}$, $\{\gamma_i\}$, $i\in K$, при построении которых бесконечное число раз выполнялись присваивания (\ref{e2.3.8}) и происходило обновление множеств $M_{i+1}$ за счет соответствующего выбора множеств $Q_i$. Другими словами, исследуем свойства последовательностей $\{z_k\}$, $\{\sigma_k\}$, $k\in K$.

Прежде всего сформулируем условие выбора чисел $\varepsilon_k$, при котором наряду с $\{x_i\}$, $\{\gamma_i\}$, $i\in K$, будут построены в виде (\ref{e2.3.8}) и последовательности $\{z_k\}$, $\{\sigma_k\}$, $k\in K$.

\begin{lemma}\label{lemma2.3.9}
Пусть последовательности $\{x_i\}$, $\{\gamma_i\}$, $i\in K$, построены методом с условием, что
числа $\varepsilon_k$, $k\in K$, выбраны согласно условию (\ref{e2.1.22}).
Тогда для каждого $k\in K$ существует номер $i=i_k$, при котором выполняются равенства (\ref{e2.3.8}).
\end{lemma}
\begin{proof}
1) Пусть $k=0$. Если $f(x_0)-\gamma_0\le \varepsilon_0$, то согласно п. 4 метода $i_0=0$, $z_0=x_{i_0}=x_0$, $\sigma_0=\gamma_0$, и равенства (\ref{e2.3.8}) при $k=0$ выполняются. Поэтому будем считать, что $f(x_0)-\gamma_0 > \varepsilon_0$. Покажем тогда существование номера $i=i_0>0$, для которого справедливо неравенство 
\begin{equation}\label{e2.3.39}
f(x_{i_0})-\gamma_{i_0}\le \varepsilon_0.
\end{equation}

Допустим противное, то есть $f(x_i)-\gamma_i > \varepsilon_0$ для всех $i\in K$, $i>0$. Тогда множества $Q_i$, $i\in K$, $i>0$, имеют вид (\ref{e2.3.7}), и из леммы {\ref{lemma2.3.8}} следует неравенство $\varepsilon_0\le 0$, противоречащее условию (\ref{e2.1.22}). Таким образом, существование номера $i_0>0$, для которого справедливо (\ref{e2.3.39}), доказано, и равенства (\ref{e2.3.8}) имеют место при $k=0$.

2) Допустим теперь, что (\ref{e2.3.8}) выполняются при некотором фиксированном $k\ge 0$, то есть имеется номер $i=i_k$, удовлетворяющий условию (\ref{e2.3.18}) при выбранном $k$. Покажем существование такого номера $i_{k+1}>i_k$, что 
\begin{equation}\label{e2.3.40}
f(x_{i_{k+1}})-\gamma_{i_{k+1}}\le \varepsilon_{k+1},
\end{equation}
тогда $z_{k+1}=x_{i_{k+1}}$, $\sigma_{k+1}=\gamma_{i_{k+1}}$, и лемма будет доказана. Если предположим противное, то есть, что при всех $i\in K$, $i>i_k$ выполняются неравенства $f(x_i)-\gamma_i>\varepsilon_{k+1}$, то, как и в первой части доказательства, привлекая лемму {\ref{lemma2.3.8}}, получим противоречие с выбором числа $\varepsilon_{k+1}$. Таким образом, номер $i_{k+1}>i_k$, при котором справедливо (\ref{e2.3.40}), действительно существует. Лемма доказана.
\end{proof}

Итак, условие (\ref{e2.1.22}) гарантирует существование последовательностей $\{z_k\}$, $\{\sigma_k\}$, $k\in K$. Обоснуем теперь сходимость $\{f(z_k)\}$, $\{\sigma_k\}$, $k\in K$, к решению при еще одном дополнительном предположении о последовательности $\{\varepsilon_k\}$, $k\in K$.

\begin{theorem}\label{theorem2.3.3}
Пусть в методе числа $\varepsilon_k$, $k\in K$, выбраны согласно (\ref{e2.1.22}), (\ref{e2.1.28}).
Тогда для последовательностей $\{z_k\}$, $\{\sigma_k\}$, $k\in K$, справедливы равенства
\begin{equation}\label{e2.3.42}
\lim \limits_{k\in K} f(z_k)=f^*, \quad \lim \limits_{k\in K} \sigma_k = f^*.
\end{equation}
\end{theorem}
\begin{proof}
Ввиду (\ref{e2.3.8}) $f(x_{i_k})-\gamma_{i_k}\le \varepsilon_k$ или 
\begin{equation}\label{e2.3.43}
f(z_k)\le \sigma_k + \varepsilon_k \quad \forall k\in K.
\end{equation}
Кроме того, согласно (\ref{e2.3.23}) и включению $z_k\in D$, $k\in K$, имеют место неравенства
\begin{equation}\label{e2.3.44}
\sigma_k\le f^*\le f(z_k) \quad \forall k\in K.
\end{equation}
Поэтому $f^*\le f(z_k)\le f^* + \varepsilon_k$, $k\in K$. Тогда с учетом (\ref{e2.1.28}) 
$$f^*\le \lim \limits_{\overline{k\in K}} f(z_k)\le \overline{\lim \limits_{k\in K}} f(z_k)\le f^*,$$
и первое из равенств (\ref{e2.3.42}) доказано.

Далее, в силу тех же неравенств (\ref{e2.3.43}), (\ref{e2.3.44}) $f^*- \varepsilon_k\le \sigma_k \le f^*$, $k\in K$. Отсюда и из (\ref{e2.1.28}) следует второе из равенств (\ref{e2.3.42}). Теорема доказана.
\end{proof}

Относительно последовательности $\{\varepsilon_k\}$, $k\in K$, можно сделать замечание, близкое к замечанию \ref{r2.1.8}. В данном методе удобно положить, например,
$$\varepsilon_{k+1}=\tau_k (f(z_k)-\sigma_k), \quad k\ge 0,$$
считая $\varepsilon_0$ сколь угодно большим.

В заключение параграфа приведем оценки, касающиеся точности решения задачи. Согласно (\ref{e2.3.22}) при всех $i\in K$, в том числе и при $i=i_k$, выполняются неравенства $\gamma_i\le f^*\le f(x_i)$. Следовательно, с учетом (\ref{e2.3.8})
\begin{equation}\label{e2.3.45}
0\le f(z_k)-f^*\le f(z_k)-\sigma_k\le \varepsilon_k, \quad k\in K.
\end{equation}

При дополнительном условии сильной выпуклости функции $f(x)$ с константой сильной выпуклости $\mu$ для всех $k\in K$ справедлива оценка
\begin{equation}\label{e2.3.46}
\|z_k - x^*\|\le \sqrt{\frac{\varepsilon_k}{\mu}},
\end{equation}
которая следует из известного неравенства $\mu\|z_k-x^*\|^2\le f(z_k)-f(x^*)$ (напр., \cite[207]{n_vasillev_f_p_1}) и оценки (\ref{e2.3.45}).

Как и в $\S$ 2.1, выбирая числа $\varepsilon_k$ в виде (\ref{e2.1.35}), из (\ref{e2.3.45}), (\ref{e2.3.46}) получаются оценки скорости сходимости (\ref{e2.1.36}), где \mbox{$x_k=z_k$}.
         
\setcounter{equation}{0}

\section
{Метод отсечений и построение на его основе смешанных алгоритмов минимизаций}\label{s2.4}

В данном параграфе предлагается метод отсечений с другим критерием качества аппроксимации, то есть с другим условием отбрасывания дополнительных ограничений. За конечное число шагов методом строится $\varepsilon$-решение исходной задачи. Метод характерен также тем, что его можно смешивать с любыми другими релаксационными методами выпуклого программирования. При этом сходимость смешанных алгоритмов гарантирована сходимостью предложенного метода отсечений. Кроме того, предлагаемый метод допускает возможность использования параллельных вычислений на этапе выбора основной итерационной точки.

Пусть $D$ -- выпуклое замкнутое множество, $f(x)$ -- выпуклая достигающая на множестве $D$ своего минимального значения функция. Пусть $f^*$, $X^*$ определены как и в $\S \ 2.1$,  $X^*(\varepsilon)=\{x\in D : f(x)\le f^* + \varepsilon\}$, где $\varepsilon\ge 0$.

Предлагаемый ниже метод решения задачи (\ref{e2.1.0}) позволяет в случае $\varepsilon>0$ за конечное число шагов получить итерационную точку, принадлежащую множеству $X^*(\varepsilon)$, то есть найти $\varepsilon$-решение задачи (\ref{e2.1.0}). Если же после нахождения $\varepsilon$-решения процесс построения приближений не заканчивается, то, как показано ниже, предельные точки построенной методом последовательности приближений принадлежат множеству $X^*$.

Введем вспомогательную функцию
\begin{equation}\label{p2.2.2}
g(x) = f(x) + \varepsilon,
\end{equation}
где $\varepsilon \ge 0$. Так как по условию $X^*\ne \emptyset$, то легко видеть, что и функция $g(x)$ достигает на множестве $D$ минимального значения. Положим $g^*=\min \{g(x) : x\in D\}$, $X^*_g= \{x\in D : g(x) = g^*\}$. Пусть, кроме того, ${\rm epi\,} (g, G)= \{(x, \gamma)\in {\rm R_{n+1}} : x\in G ,\ \gamma \ge g(x)\}$, где $G\subset {\rm R_n}$, ${\rm int\,}Q$ -- внутренность множества $Q\subset {\rm R_{n+1}}$, $W(u, Q)= \{a\in {\rm R_{n+1}}:\langle a, v - x\rangle \le 0\ \forall v\in Q\}$ -- множество обобщенно-опорных к $Q$ в точке $u\in {\rm R_{n+1}}$ векторов, $W^1(u, Q)= \{a\in W(u, Q) : \|a\|=1\}$, $K=\{0,1,\ldots\}$, $J=\{1, \ldots , m\}$, где $m\ge 1$, $x^* \in X^*$.

Заметим, что с учетом (\ref{p2.2.2}) легко доказать равенства
\begin{equation}\label{p2.2.3}
X^*=X^*_g,
\end{equation}
\begin{equation}\label{p2.2.4}
g^*=f^*+ \varepsilon.
\end{equation}

Сформулируем критерий принадлежности точки множеству $X^*(\varepsilon)$.

\begin{lemma}\label{plemma2.2.1}
Пусть $\varepsilon>0$, $(y, \gamma)\in {\rm epi\,} (f, D)$, причем
\begin{equation}\label{p2.2.5}
\gamma\le g^*.
\end{equation}
Тогда имеет место включение
$$y\in X^*(\varepsilon).$$
\end{lemma}
\begin{dokaz}
Так как по условию $f(y)\le \gamma$, то из (\ref{p2.2.4}), (\ref{p2.2.5}) следует, что $f(y)\le f^*+\varepsilon$. Отсюда с учетом включения $y\in D$ вытекает утверждение леммы.
\end{dokaz}

Предлагаемый метод решения задачи (\ref{e2.1.0}) вырабатывает вспомогательную последовательность приближений $y_i$, $i\in K$, и основную последовательность $\{x_k\}$, $k\in K$, и заключается в следующем.

\textbf{Метод 2.4.}
Выбираются точки $v^j\in {\rm int\, epi\,}(g, {\rm R_n})$ для всех $j\in J$ и выпуклое ограниченное замкнутое множество $G_0\subset D$, содержащее точку $x^*$. Строится выпуклое замкнутое множество $M_0\subseteq {\rm R_{n+1}}$ такое, что ${\rm epi\,} (g, {\rm R_n})\subset M_0$. Задаются числа $\varepsilon>0$, $\delta_0>0$, $\bar{\gamma}_0\le g^*$. Полагается $i=0$, $k=0$.

\textbf{1.} Отыскивается точка $u_i=(y_i, \gamma_i)$, где $y_i\in {\rm R_n}$, $\gamma_i\in {\rm R_1}$, как решение задачи
\begin{equation}\label{p2.2.6}
\gamma \to \min 
\end{equation}
\begin{equation}\label{p2.2.7}
(x, \gamma)\in M_i, \quad x\in G_i, \quad \gamma \ge \bar{\gamma}_i.
\end{equation}
Если
\begin{equation}\label{p2.2.8}
\gamma_i=g(y_i),
\end{equation}
то $y_i\in X^*$,  и процесс завершается.

\textbf{2.} Если
\begin{equation}\label{p2.2.9}
\gamma_i \ge f(y_i),
\end{equation}
то $y_i\in X^*(\varepsilon)$.

\textbf{3.} Для каждого $j\in J$ в интервале $(v^j, u_i)$ выбирается точка $\bar{u}_i^j\in {\rm R_{n+1}}$ так, чтобы $\bar{u}_i^j\notin {\rm int\, epi\,} (g, {\rm R_n})$ и при некотором $q_i^j\in [1, q]$, $q<+\infty$, для точки 
\begin{equation}\label{p2.2.10}
\bar{v}_i^j = u_i + q_i^j (\bar{u}_i^j - u_i)
\end{equation}
выполнялось включение $\bar{v}_i^j \in {\rm epi}(g, {\rm R_n})$.

\textbf{4.} Если
\begin{equation}\label{p2.2.11}
\|\bar{u}_i^j - u_i\|>\delta_k \quad \forall j\in J,
\end{equation}
то выбирается выпуклое замкнутое множество $Q_i\subset {\rm R_{n+1}}$ такое, что 
\begin{equation}\label{p2.2.12}
Q_i\subseteq M_i,
\end{equation}
\begin{equation}\label{p2.2.13}
{\rm epi} (g, {\rm R_n}) \subset Q_i,
\end{equation}
и следует переход к п. 6. В противном случае выполняется п. 5.

\textbf{5.} Выбирается выпуклое замкнутое множество $Q_i\subseteq {\rm R_{n+1}}$ так, чтобы выполнялось включение (\ref{p2.2.13}). Полагается $i_k=i$,
\begin{equation}\label{p2.2.14}
\sigma_k = \gamma_{i_k},
\end{equation}
выбирается точка $x_k\in G_{i_k}$, удовлетворяющая неравенству
\begin{equation}\label{p2.2.15}
g(x_k)\le g(y_{i_k}).
\end{equation}
Задается $\delta_{k+1}>0$, значение $k$ увеличивается на единицу.
\linebreak

\textbf{6.} Для каждого $j\in J$ выбирается конечное множество $A_i^j\subset W^1(\bar{u}_i^j$, ${\rm epi\,} (g, {\rm R_n}))$, и полагается
\begin{equation}\label{p2.2.16}
M_{i+1}=Q_i\bigcap T_i,
\end{equation}
где $T_i=\bigcap \limits_{j\in J} \{u\in {\rm R_{n+1}} : \langle a, u - \bar{u}_i^j \rangle \le 0 \ \forall a\in A_i^j\}$.

\textbf{7.} Выбирается выпуклое замкнутое множество $G_{i+1}\subseteq G_0$, содержащее точку $x^*$. Задается число $\bar{\gamma}_{i+1}$ из условия
\begin{equation}\label{p2.2.17}
\bar{\gamma}_0\le \bar{\gamma}_{i+1}\le g^*.
\end{equation}
Значение $i$ увеличивается на единицу, и следует переход к п. 1.

Сделаем относительно этого метода некоторые замечания. Покажем сначала, что множество ограничений вспомогательной задачи (\ref{p2.2.6}), (\ref{p2.2.7}) непусто.
\begin{lemma}\label{plemma2.2.2}
Точка $(x^*, g^*)$ для всех $i\in K$ удовлетворяет условиям (\ref{p2.2.7}).
\end{lemma}
\begin{dokaz}
Согласно выбору множеств $G_i$ и условию (\ref{p2.2.17}) для всех $i\in K$ выполняются включения $x^*\in G_i$ и неравенство $g^*\ge \bar{\gamma}_i$. Поэтому для обоснования леммы достаточно показать, что 
\begin{equation}\label{p2.2.18}
(x^*, g^*)\in M_i
\end{equation}
для всех $i\in K$.

Заметим, что при $i=0$ включение (\ref{p2.2.18}) выполняется по условию выбора множества $M_0$. Зафиксируем произвольно номер $i=l\ge 0$, и покажем, что 
\begin{equation}\label{p2.2.19}
(x^*, g^*)\in M_{l+1}.
\end{equation}
Действительно, во-первых, $(x^*, g^*)\in Q_l$ ввиду условия (\ref{p2.2.13}) выбора множеств $Q_i$, $i\in K$. Во-вторых, так как $(x^*, g^*)\in {\rm epi\,} (g, {\rm R_n})$, то $\langle a, (x^*, g^*) - \bar{u}_l^j\rangle \le 0$ для всех $a\in A_l^j$, $j\in J$,  а значит, $(x^*, g^*)\in T_l$. Отсюда и из (\ref{p2.2.16}) имеем (\ref{p2.2.19}). Лемма доказана.
\end{dokaz}

Обоснуем далее критерий остановки, заложенный в п. 1 метода.
\begin{lemma}\label{plemma2.2.3}
При каждом $i\in K$ справедливо неравенство
\begin{equation}\label{p2.2.20}
\gamma_i\le g^*.
\end{equation}
\end{lemma}
\begin{dokaz}
Согласно (\ref{p2.2.6}), (\ref{p2.2.7}) для каждой точки $(x, \gamma)$, удовлетворяющей условиям (\ref{p2.2.7}), выполняется неравенство $\gamma_i\le \gamma$. Но по лемме \ref{plemma2.2.2}  точка $(x^*, \gamma^*)$ является допустимым решением задачи (\ref{p2.2.6}), (\ref{p2.2.7}) при любом $i\in K$. Отсюда следует (\ref{p2.2.20}). Лемма доказана.
\end{dokaz}

\begin{remark}\label{r2.4.1}
Лемма \ref{plemma2.2.3} позволяет выбирать числа $\bar{\gamma}_{i+1}$ в п. 7 метода следующим образом. Можно положить $\bar{\gamma}_{i+1}=\gamma_l$, где $0\le l\le i$, в частности, $\bar{\gamma}_{i+1}=\max \limits_{0\le j\le i} \gamma_j$.
\end{remark}

\begin{theorem}\label{ptheorem2.2.1}
Пусть при некотором $i\in K$ выполняется равенство (\ref{p2.2.8}). Тогда $y_i\in X^*$.
\end{theorem}
\begin{dokaz}
В силу (\ref{p2.2.8}), (\ref{p2.2.20}) $g(y_i)\le g^*$. С другой стороны,
\begin{equation}\label{p2.2.21}
g(y_i)\ge g^* \quad \forall i\in K,
\end{equation}
поскольку $y_i\in D$ для всех $i\in K$. Таким образом, $g(y_i)=g^*$, а значит, $y_i\in X^*_g$. Отсюда с учетом (\ref{p2.2.3}) следует утверждение теоремы.
\end{dokaz}

Далее, предположим, что для точки $(y_i,\gamma_i)$ выполнилось неравенство (\ref{p2.2.9}). Тогда по лемме \ref{plemma2.2.1} с учетом (\ref{p2.2.20}), как и отмечено в п. 2 метода, точка $y_i$ -- является $\varepsilon$-решением исходной задачи.

\begin{remark}\label{r2.4.2}
Задание точки $\bar{u}_i^j$, $j\in J$, не представляет особого труда. Их можно выбирать в виде точки пересечения лучей $v^j + \alpha (u_i- v^j)$, $\alpha\ge 0$, $j\in J$,  с границей
\linebreak
\linebreak
множества ${\rm epi\,} (g, {\rm R_n})$. Условие (\ref{p2.2.10}) позволяет находить эти точки пересечения приближенно.
\end{remark}

\begin{remark}\label{r2.4.3}
Задание начального аппроксимирующего множества $M_0$ также не вызывает затруднений. Его легко задать линейными неравенствами так, чтобы выполнялось включение ${\rm epi\,} (g, {\rm R_n})\subset M_0$. Если положить $M_0={\rm R_{n+1}}$, то точка $(y_0, \gamma_0)$, где $y_0\in G_0$, $\gamma_0=\bar{\gamma}_0$, может быть принята за решение задачи (\ref{p2.2.6}), (\ref{p2.2.7}) при $i=0$.
\end{remark}

Как было подчеркнуто выше, предлагаемый метод обладает важной практической особенностью -- возможностью периодического обновления аппроксимирующих надграфик множеств $M_i$ за счет отбрасывания отсекающих плоскостей. Поясним, как такие обновления можно проводить за счет выбора множеств $Q_i$ на итерациях с номерами $i=i_k$.

\begin{remark}\label{r2.4.3.1}
Пусть $i$, $k$ таковы, что для точки $u_i=(y_i, \gamma_i)$ хотя бы при одном $j\in J$ выполняется неравенство
\begin{equation}\label{p2.2.24}
\|\bar{u}_i^j - u_i\|\le \delta_k.
\end{equation}
Тогда согласно п. 5 метода выпуклое замкнутое множество $Q_i=Q_{i_k}$ выбирается в ${\rm R_{n+1}}$ лишь с одним условием -- (\ref{p2.2.13}). Если положить, например, $Q_{i_k}={\rm R_{n+1}}$ или $Q_{i_k}=M_0$, то согласно (\ref{p2.2.16}) $M_{i_k+1}=T_{i_k}$ или $M_{i_k+1}=M_0\bigcap T_{i_k}$ соответственно. В таком случае ни одна из построенных к шагу $i=i_k$ отсекающих плоскостей принимать участие в построении $M_{i_k+1}$ не будет, то есть произойдет полное обновление множества $M_{i_k+1}$. Если при условии (\ref{p2.2.24}) положить 
$$Q_{i_k}=M_{r_i},$$
где $0<r_i<i_k$, то при построении $M_{i_k+1}$ отбросится лишь часть накопленных к шагу $i_k$ отсечений, и будет иметь место частичное обновление аппроксимирующего множества.
\end{remark}

Как показано ниже, для каждого $k\in K$ найдется точка $u_i=u_{i_k}$ строящейся последовательности $\{u_i\}$, для которой выполнится условие (\ref{p2.2.24}), а значит, представится возможность полного или частичного обновления аппроксимирующего множества $M_{i_k+1}$. Условия (\ref{p2.2.11}), (\ref{p2.2.24}) фактически определяют качество аппроксимации надграфика ${\rm epi\,} (g, {\rm R_n})$ множеством $M_i$  в окрестности точки $(y_i, \gamma_i)$. При выполнении для текущего номера $k\in K$ условия (\ref{p2.2.24}) качество аппроксимации считается достаточным для фиксирования числа $\sigma_k$ в виде (\ref{p2.2.14}) и выбора точки $x_k$ из условия (\ref{p2.2.15}), в частности, $x_k=y_{i_k}$.

Прокомментируем далее возможность построения смешанных алгоритмов на основе предложенного метода отсечений и других известных или новых алгоритмов выпуклого программирования.

\begin{remark}\label{r2.4.3.2}
Возможность смешивания с другими методами обеспечивается условием (\ref{p2.2.15}) выбора точек $x_k$, если
\begin{equation}\label{p2.2.25}
x_k\ne y_{i_k}.
\end{equation}
А именно, в зависимости от свойств целевой функции задачи (\ref{e2.1.0}) на итерации $i=i_k$ выбирается некоторый релаксационный метод $\Gamma_i$, пригодный для минимизации функции $g(x)$ на множестве $G_i$ (см., напр., \cite{n_vasillev_f_p_1}). Методом $\Gamma_i$ делается $p_i$ шагов, считая $y_i$ точкой начального приближения, и найденная на шаге $p_i$ итерационная точка принимается за $x_k$. Таким образом, будет построен новый алгоритм решения задачи (\ref{e2.1.0}), причем сходимость этого смешанного алгоритма уже гарантирована обоснованной ниже сходимостью предложенного метода отсечений.
\end{remark}

\begin{remark}\label{r2.4.3.3}
Точку $x_k\in G_i$ с условиями (\ref{p2.2.15}), (\ref{p2.2.25}) допустимо находить, используя  параллельные вычисления. Для этого можно на шаге $i=i_k$ одновременно  алгоритмами $\Gamma_i^1, \ldots, \Gamma_i^l$, $l\ge 1$, построить соответственно вспомогательные точки $w_i^1, \ldots, w_i^l$ такие, что $w_i^r\in G_i$, $g(w_i^r)< g(y_i)$, $r=1,\ldots, l$, и в качестве $x_k$ выбрать одну из них, например, в которой значение функции $g(x)$ минимально.
\end{remark}

Перейдем теперь к исследованию сходимости метода.

Так как выполняются включения $y_i\in G_i\subset G_0$, $i\in K$, множество $G_0$ по условию ограничено, а для чисел $\gamma_i$, $i\in K$, имеют место неравенства $\bar{\gamma}_0\le \gamma_i\le g^*$, то последовательность $\{u_i\}$, $i\in K$, ограничена. Кроме того, поскольку $\|\bar{u}_i^j - u_i\|\le \|v^j - u_i\|$ для всех $i\in K$, $j\in J$, то с учетом (\ref{p2.2.10}) и последовательность $\{\bar{v}_i^j\}$, $i\in K$, является ограниченной при каждом $j\in J$.

Заметим также, что ввиду пп. 4, 5 метода, независимо от условий (\ref{p2.2.11}) или  (\ref{p2.2.24}), множества $Q_i$ при всех $i\in K$ могут выбираться согласно (\ref{p2.2.12}), например, $Q_i=M_i$. В таком случае в силу (\ref{p2.2.16})
\begin{equation}\label{p2.2.26}
M_{i+1}\subset M_i,
\end{equation}
$i\in K$. Кроме того, отметим, что наряду с последовательностью $\{y_i\}$, $i\in K$, будет построена и последовательность $\{x_k\}$, $k\in K$, так как $y_i\notin X_g^*$, $i\in K$. С учетом этих замечаний сформулируем следующее вспомогательное утверждение.

\begin{lemma}\label{plemma2.2.4}
Пусть последовательность $\{u_i\}$, $i\in K$, построена с условием, что для всех $i\in K$, начиная с некоторого номера $i^\prime\ge 0$, множества $Q_i$ были выбраны согласно (\ref{p2.2.12}). Тогда для каждого $j\in J$ выполняется равенство
\begin{equation}\label{p2.2.27}
\lim \limits_{i\in K} \|\bar{u}_i^j - u_i\|=0.
\end{equation}
\end{lemma}
\begin{proof}
Зафиксируем номер $j\in J$ и докажем для него равенство (\ref{p2.2.27}). Предположим, что (\ref{p2.2.27}) не выполняется. Тогда найдутся бесконечное подмножество $K^\prime$ множества $K$ и число $\Delta>0$ такие, что
\begin{equation}\label{p2.2.28}
\|\bar{u}_i^j - u_i\|\ge \Delta \quad \forall i\in K^\prime.
\end{equation}
Выделим из последовательности $\{u_i\}$, $i\in K^\prime$, сходящуюся подпоследовательность $\{u_i\}$, $i\in K^{\prime\prime}\subset K^\prime$.

Согласно выбору точек $\bar{u}_i^j$ для каждого $i\in K$ найдется такое число $\tau_i^j\in (0,1)$, что
\begin{equation}\label{p2.2.29}
\bar{u}_i^j=u_i + \tau_i^j(v^j - u_i).
\end{equation}
Заметим также, что по условию леммы для всех $i\ge i^\prime$, $i\in K$, справедливо включение (\ref{p2.2.26}). Зафиксируем теперь номера $i$, $p_i\in K^{\prime\prime}$, удовлетворяющие неравенствам $p_i> i\ge i^\prime$. Тогда $M_{p_i}\subset M_i$. Так как по построению $u_{p_i}\in M_{p_i}$, а любой элемент $a\in A_i^j$ является обобщенно-опорным к множеству $M_{p_i}$  в точке $\bar{u}_i^j$, то
$$\langle a, u_{p_i}-\bar{u}_i^j\rangle \le 0 \quad \forall a\in A_i^j,$$
а с учетом (\ref{p2.2.29})
$$\langle a, u_i - u_{p_i} \rangle \ge \tau_i^j \langle a, u_i - v^j\rangle \quad \forall a\in A_i^j.$$
Согласно лемме \ref{lemma1.1.1} найдется такое число $\beta^j > 0$, что $\langle a, u_i - v^j\rangle \ge \beta ^j$ для всех $a\in A_i^j$. В связи с этим $\langle a, u_i - u_{p_i}\rangle \ge \tau_i^j \beta^j$ для всех $a\in A_i^j$, а так как $\|a\|=1$, то для выбранных номеров $i$, $p_i\in K^{\prime\prime}$ имеем
\begin{equation}\label{p2.2.30}
\|u_i - u_{p_i}\|\ge \tau_i^j \beta^j.
\end{equation}
Поскольку последовательность $\{u_i\}$, $i\in K^{\prime\prime}$, сходится, то в силу (\ref{p2.2.30}) $\tau_i^j\to 0$, $i\to \infty$, $i\in K^{\prime\prime}$. Значит, из (\ref{p2.2.29})  с учетом ограниченности последовательности $\{\|v^j - u_i\|\}$, $i\in K^{\prime\prime}$, следует, что $\|\bar{u}_i^j - u_i\|\to 0$, $i\to \infty$, $i\in K^{\prime\prime}$. Это предельное соотношение противоречит неравенствам (\ref{p2.2.28}). Лемма доказана.
\end{proof}

\begin{lemma}\label{plemma2.2.5}
Пусть выполняются условия леммы \ref{plemma2.2.4}, и, кроме того, $\{u_i\}$, $i\in K^\prime\subset K$, -- сходящаяся  подпоследовательность последовательности $\{u_i\}$, $i\in K$, а $\bar{u}$=$(\bar{y}, \bar{\gamma})$ -- предельная точка этой подпоследовательности. Тогда
\begin{equation}\label{p2.2.32}
g(\bar{y})=\bar{\gamma}.
\end{equation}
\end{lemma}
\begin{proof}
Зафиксируем $r\in J$, и выделим из последовательности $\{\bar{v}_i^r\}$, $i\in K^\prime$, сходящуюся подпоследовательность  $\{\bar{v}_i^r\}$, $i\in K_r\subset K^\prime$. Пусть $\tilde{v}^r$ -- ее предельная точка. Отметим, что 
\begin{equation}\label{p2.2.33}
\tilde{v}^r \in {\rm epi\,} (g, {\rm R_n})
\end{equation}
в силу замкнутости множества ${\rm epi} (g, {\rm R_n})$. Перейдем в равенстве (\ref{p2.2.10}), где $j=r$, к пределу по $i\to \infty$, $i\in K_r$, с учетом утверждения (\ref{p2.2.27}). Тогда $\bar{u}=\tilde{v}^r$, и ввиду (\ref{p2.2.33}) выполняется включение $\bar{u}\in {\rm epi\,}(g, {\rm R_n})$. Следовательно,
\begin{equation}\label{p2.2.34}
\bar{\gamma}\ge g(\bar{y}).
\end{equation}
С другой стороны, для всех $i\in K$ согласно (\ref{p2.2.20}) выполняются неравенства $\gamma_i\le g(y_i)$, а значит, $\bar{\gamma}\le g(\bar{y})$. Отсюда и из (\ref{p2.2.34}) вытекает равенство (\ref{p2.2.32}). Лемма доказана.
\end{proof}

\begin{lemma}\label{plemma2.2.6}
Пусть выполняются условия леммы \ref{plemma2.2.4}. Тогда 
\begin{equation}\label{p2.2.35}
\lim \limits_{i\in K}(g(y_i) - \gamma_i) = 0.
\end{equation}
\end{lemma}
\begin{proof}
Отметим, что $g(y_i)-\gamma_i > 0$ для всех $i\in K$. Предположим, что утверждение леммы неверно. Тогда найдется  подпоследовательность $\{u_i\}$, $i\in K^\prime \subset K$, $i\ge i^\prime$, последовательности $\{u_i\}$, $i\in K$, и число $\Delta>0$ такие, что 
\begin{equation}\label{p2.2.36}
g(y_i)-\gamma_i\ge \Delta \quad \forall i\in K^\prime, \quad i\ge i^\prime.
\end{equation}
Пусть $\{u_i\}$, $i\in K^{\prime\prime}\subset K^\prime$, -- сходящаяся подпоследовательность, выделенная из точек $u_i$, $i\in K^\prime$, $i\ge i^\prime$, а $\bar{u}=(\bar{y}, \bar{\gamma})$ -- предельная точка этой подпоследовательности. Тогда с одной стороны в силу предположения (\ref{p2.2.36}) выполняется неравенство $g(\bar{y}) - \bar{\gamma} \ge \Delta$. Полученное противоречие доказывает утверждение леммы.
\end{proof}

\begin{theorem}\label{ptheorem2.2.3}
Если для последовательности  $\{u_i\}$, $i\in K$, выполняются условия леммы \ref{plemma2.2.4}, то справедливы равенства
\begin{equation}\label{p2.2.37}
\lim \limits_{i\in K} g(y_i)= g^*, \quad \lim \limits_{i\in K} \gamma_i = g^*.
\end{equation}
\end{theorem}
\begin{proof}
В силу (\ref{p2.2.20}) $0\le g(y_i)-g^*\le g(y_i) - \gamma_i$, $i\in K$. Отсюда и из (\ref{p2.2.35}) следует первое из равенств (\ref{p2.2.37}). Тогда в силу того же предельного соотношения (\ref{p2.2.35}) справедливо и второе из утверждение (\ref{p2.2.37}).
\end{proof}

Заметим, что согласно теореме \ref{ptheorem2.2.3} и известной теореме (напр., \cite[62]{n_vasillev_f_p_1}) при выполнении условий  леммы \ref{plemma2.2.4} последовательность $\{y_i\}$, $i\in K$, сходится с учетом (\ref{p2.2.3}) к множеству $X^*$.

Перейдем к исследованию последовательности $\{(x_k, \sigma_k)\}$, $k\in K$. Прежде всего докажем, что наряду с $\{(y_i, \gamma_i)\}$, $i\in K$, будет построена и последовательность $\{(x_k, \sigma_k)\}$, $k\in K$.

\begin{lemma}\label{plemma2.2.7}
Если последовательность $\{(y_i, \gamma_i)\}$, $i\in K$, построена предложенным методом, то для каждого $k\in K$ существует точка $(x_k, \sigma_k)$, полученная согласно (\ref{p2.2.14}), (\ref{p2.2.15}).
\end{lemma}
\begin{proof}
1) Пусть $k=0$. Если допустить, что $\|\bar{u}_0^j - u_0\|\le \delta_0$ для некоторого $j\in J$, то согласно п. 5 метода $i_0=0$, $\sigma_0=\gamma_0$, и найдена точка $x_0\in G_0$, удовлетворяющая неравенству $g(x_0)\le g(y_0)$, то есть пара $(x_0, \sigma_0)$ построена. Поэтому будем предполагать, что $\|\bar{u}_0^j - u_0\|>\delta_0$ для всех $j\in J$. Докажем тогда существование номера $i=i_0>0$, для которого выполняется неравенство
\begin{equation}\label{p2.2.38}
\|\bar{u}_{i_0}^j - u_{i_0}\|\le\delta_0
\end{equation}
при некотором $j\in J$.

Предположим противное, то есть
\begin{equation}\label{p2.2.39}
\|\bar{u}_i^j - u_i\|>\delta_0 \quad \forall i>0, \quad j\in J.
\end{equation}
Пусть $\{u_i\}$, $i\in K^\prime \subset K$, сходящаяся подпоследовательность, выделенная из точек $u_i$, $i\in K$, $i>0$, а $\tilde{u}$ -- ее предельная точка. Выберем теперь из множества $K^\prime$ такое бесконечное подмножество $K^{\prime\prime}$, что для каждого $j\in J$ последовательность $\{\bar{u}_i^j\}$, $i\in K^{\prime\prime}$,  окажется сходящейся. Пусть при каждом $j\in J$ точка $\tilde{u}^j$ является предельной для последовательности $\{\bar{u}_i^j\}$, $i\in K^{\prime\prime}$. Тогда из неравенств (\ref{p2.2.39}) следует, что
\begin{equation}\label{p2.2.40}
\|\tilde{u}^j - \tilde{u}\|\ge \delta_0 \quad \forall j\in J.
\end{equation}
C другой стороны, в силу предположения (\ref{p2.2.39}) для всех $i>0$ множества $Q_i$ при построении последовательности $\{u_i\}$, $i\in K$, $i>0$, были выбраны согласно условию (\ref{p2.2.12}). Значит, по лемме \ref{plemma2.2.4} для всех $j\in J$ имеют место равенства $\|\tilde{u}^j - \tilde{u}\|=0$, противоречащее (\ref{p2.2.40}). Таким образом, доказано существование номера $i_0>0$, для которого справедливо (\ref{p2.2.38}), а следовательно, доказано наличие точки $(x_0, \sigma_0)$, полученной в ходе построения последовательности $\{(y_i, \gamma_i)\}$, $i\in K$.

2) Допустим теперь, что при некотором  фиксированном $k\ge 0$ точка $(x_k, \sigma_k)$ построена, то есть существует номер $i=i_k$, для которого выполняется неравенство $\|\bar{u}_{i_k}^j - u_{i_{k}}\|\le \delta_{k+1}$ хотя бы для одного $j\in J$. Покажем при этом предположении существование такого номера $i_{k+1}> i_k$, что
\begin{equation}\label{p2.2.41}
\|\bar{u}^j_{i_{k+1}} - u_{i_{k+1}}\|\le \delta_{k+1}
\end{equation}
при некотором $j\in J$. Тогда $\sigma_{k+1}=\gamma_{i_{k+1}}$, точка $x_{k+1}\in G_{i_{k+1}}$ с условием, что $g(x_{k+1})\le g(y_{i_{k+1}})$, сушествует, и лемма будет доказана.

Допустим противное, то есть для всех $i>i_k$ и всех $j\in J$
\begin{equation}\label{p2.2.42}
\|\bar{u}_i^j - u_i\|> \delta_{k+1}.
\end{equation}
Как и в первой части доказательства леммы, выделим из точек $u_i$, $\bar{u}_i^j$, $j\in J$, с номерами $i\in K$, $i>i_k$ сходящиеся подпоследовательности $\{u_i\}$, $\{\bar{u}_i^j\}$, $j\in J$, $i\in \bar{K}$, и пусть $\bar{u}$, $\bar{u}^j$,  соответственно, их предельные точки. Тогда из неравенств (\ref{p2.2.42}) следует, что $\|\bar{u}^j - \bar{u}\|\ge \delta_{k+1}$ для всех $j\in J$, а с другой стороны с учетом леммы \ref{plemma2.2.4} для каждого $j\in J$ имеет место равенства $\|\bar{u}^j - \bar{u}\|=0$. Полученное противоречие доказывает наличие номера $i_{k+1}>i_k$, удовлетворяющего при некотором $j\in J$ неравенству (\ref{p2.2.41}). Лемма доказана.
\end{proof}

Согласно лемме \ref{plemma2.2.7} для каждого $k\in K$ найдутся номера $i=i_k$ и $j_k\in J$ такие, что выполняется неравенство $\|\bar{u}_{i_k}^{j_k} - u_{i_k}\|\le \delta_k$. Значит, наряду с последовательностью $\{u_i\}$, $i\in K$, методом будет построена и последовательность $\{u_{i_k}\}$, $k\in K$. С учетом этого замечания сформулируем следующее утверждение.

\begin{lemma}\label{plemma2.2.8}
Пусть последовательность $\{u_i\}$, $i\in K$, построена методом с условием, что
\begin{equation}\label{p2.2.43}
\delta_k\to 0, \quad k \to \infty.
\end{equation}
Тогда для последовательности $\{u_{i_k}\}$, $k\in K$, выполняется равенство
\begin{equation}\label{p2.2.44_0}
\lim \limits_{k\in K} (g(y_{i_k}) - \gamma_{i_k}) = 0.
\end{equation}
\end{lemma}
\begin{proof}
Заметим, что $g(y_{i_k})- \gamma_{i_k} >0$ для всех $k\in K$. Предположим, что утверждение леммы неверно. Тогда существует подпоследовательность $\{(y_{i_k}, \gamma_{i_k})\}$, $k\in K_1\subset K$, последовательности $\{(y_{i_k}, \gamma_{i_k})\}$, $k\in K$, и число $\Delta > 0$ такие, что
\begin{equation}\label{p2.2.44}
g(y_{i_k}) - \gamma_{i_k} \ge \Delta \quad  \forall k\in K_1.
\end{equation}
Пусть $\{(y_{i_k}, \gamma_{i_k})\}$, $k\in K_2\subset K_1$, -- сходящаяся подпоследовательность последовательности $\{(y_{i_k}, \gamma_{i_k})\}$, $k\in K_1$, и $\bar{u} = (\bar{y}, \bar{\gamma})$ -- ее предельная точка, причем согласно (\ref{p2.2.44})
\begin{equation}\label{p2.2.45}
g(\bar{y}) - \bar{\gamma} \ge \Delta.
\end{equation}

Далее, существует такой номер $r\in J$, что для бесконечного числа номеров $k\in K_2$ имеет место неравенство 
$$\|\bar{u}_{i_k}^r - u_{i_k}\|\le \delta_k.$$
Занесем в множество $K_3$ все номера $k\in K_2$, для которых выполняется последнее неравенство. Тогда из него с учетом (\ref{p2.2.43}) следует, что
\begin{equation}\label{p2.2.46}
\lim \limits_{k\in K_3} \|\bar{u}_{i_k}^r - u_{i_k}\|=0.
\end{equation}

Выделим теперь из последовательности $\{\bar{v}_{i_k}^r\}$, $k\in K_3$, сходящуюся подпоследовательность $\{\bar{v}_{i_k}^r\}$, $k\in K_4\subset K_3$, и пусть $\tilde{v}^r$ -- ее предельная точка. Сразу отметим, что для нее справедливо включение (\ref{p2.2.33}) ввиду замкнутости множества ${\rm epi\,} (g, {\rm R_n})$. Поскольку согласно (\ref{p2.2.10})
$$\bar{v}_{i_k}^r = u_{i_k} + q_{i_k}^r (\bar{u}_{i_k}^r - u_{i_k}), \quad k\in K_4,$$
то с учетом равенства (\ref{p2.2.46}) и ограниченности $\{q_{i_k}^r\}$, $k\in K_4$, отсюда имеем $\tilde{v}^r = \bar{u}$. Тогда в силу (\ref{p2.2.33}) $\bar{u}\in {\rm epi\,} (g, {\rm R_n})$, а значит, справедливо неравенство (\ref{p2.2.34}). Но с другой стороны $\bar{\gamma}\le g(\bar{y})$, так как $\gamma_{i_k}\le g(y_{i_k})$, $k\in K$. Таким образом, получено равенство (\ref{p2.2.32}), противоречащее (\ref{p2.2.45}). Лемма доказана.
\end{proof}

\begin{theorem}\label{ptheorem2.2.4}
Пусть последовательность $\{(x_k, \sigma_k)\}$, $k\in K$, построена предложенным методом с условием (\ref{p2.2.43}). Тогда
\begin{equation}\label{p2.2.47}
\lim \limits_{k\in K} g(x_k) = f^* + \varepsilon, \quad \lim \limits_{k\in K} \sigma_k = f^* + \varepsilon.
\end{equation}
\end{theorem}
\begin{proof}
Так как $x_k\in D$, $k\in K$, и согласно (\ref{p2.2.14}), (\ref{p2.2.20}) $\sigma_k =\gamma_{i_k} \le g^*$, $k\in K$, то $0\le g(x_k) - g^* \le g(y_{i_k}) - g^* \le g(y_{i_k}) -\sigma_k = g(y_{i_k}) - \gamma_{i_k} $ для всех $k\in K$. Отсюда с учетом (\ref{p2.2.44_0}), (\ref{p2.2.4}) следуют равенства (\ref{p2.2.47}). Теорема доказана.
\end{proof}

Согласно теореме \ref{ptheorem2.2.4} и упомянутой известной теореме \mbox{\cite[62]{n_vasillev_f_p_1}} последовательность $\{x_k\}$, $k\in K$, сходится к множеству $X^*$.

Отметим также, что в силу (\ref{p2.2.44_0}) найдется такой номер $l\in K$, что $g(y_l) - \gamma_l\le \varepsilon$. Тогда ввиду (\ref{p2.2.2}) $f(y_l)\le \gamma_l$, то есть выполняется неравенство (\ref{p2.2.9}) при $i=l$, и точка $y_l$ -- $\varepsilon$-решение задачи (\ref{e2.1.0}).
         
\setcounter{equation}{0}


\section
{Численные эксперименты}

В настоящем параграфе приводятся и обсуждаются результаты численных экспериментов, которые проводились в связи с исследованиями предложенных выше методов 2.1 и 2.2. Цель этих экспериментов заключалась в следующем.  Необходимо было проверить работоспособность вышеупомянутых методов и предложить рекомендации по применению в них тех или иных процедур обновления аппроксимирующих множеств $Q_i$. Кроме того,  предполагалось предложить  некоторые способы  задания чисел $\varepsilon_k$, которые определяют уровень аппроксимации надграфика целевой функции аппроксимирующими множествами. 

Данные методы программно реализованы на языке программирования С++. Тестирование методов осуществлялись на персональном компьютере, конфигурация которого была описана в $\S \ 1.5$.

При формировании вышеуказанных рекомендаций сделано значительное число экспериментов на тестовых примерах с различным числом переменных (до 50). Поскольку для каждого примера  численные эксперименты осуществлялись обоими методами, то тем самым проведено и сравнение методов между собой.

Численные эксперименты проводились на следующей тестовой задаче:
$$f(x)=\sum_{i=1}^n i^2 \xi^2_i \to \min$$
$$a_i \le \xi_i \le b_i, \quad i=1,\ldots, n,$$
где $n\ge 1$, $x=(\xi_i,\ldots,\xi_n)\in {\rm R_n}$,  $a_i=-50$, $b_i = 50$ для всех $i=1,\ldots, n$.

В первую очередь приведем критерий, на основании которого происходила остановка работы алгоритмов. Процесс отыскания решения задачи прекращался, если на некотором шаге $i^\prime\in K$ для методов 2.1 и 2.3 выполнялось условие $f(y_{i^\prime}) - \gamma_{i^\prime}\le \varepsilon$, где 
$$\varepsilon = 0.00001.$$
Согласно (\ref{e2.1.13}), (\ref{e2.2.15}) и  включению $y_{i^\prime}\in D$  для методов 2.1 и 2.3 выполнение вышеуказанного условия остановки означает справедливость следующего  неравенства:
$$f(y_{i^\prime}) - f^*\le \varepsilon.$$

При решении тестовых задач в методе 2.1 множества $S_i$ для каждого $i\in K$ выбирались совпадающими с ${\rm R_{n+1}}$, для всех $i$ полагалось $J_i=J(y_i)$. В методе 2.2 точки $v^j$, $j\in J$, строились в виде $v^j=v=(0,\ldots, 0, 100)$, для каждого $i$ задавалось $q_i^j=1$, $j\in J$. Наконец, в обоих методах при любом $i$ множества $G_i$ выбирались совпадающими с $D$, полагалось $M_0={\rm R_{n+1}}$, $\bar{\gamma}_0=-1000000$,  $\bar{\gamma}_{i+1}=\gamma_i$, $i\ge0$.

Результаты численных экспериментов представлены в следующей таблице:

\begin{center}
\begin{tabular}{|>{\centering}p{2cm}|>{\centering}p{2cm}|c|c|c|>{\centering}p{4.5cm}|>{\centering}p{3.5cm}|}
\hline
Номер & Метод & $n$  & $\varepsilon_k$ & $Q_i$ &Кол-во  итераций&  Время в мин. \tabularnewline
\hline
1& 2.1 &30 & $a_1$&$b_1$&528&1,41\tabularnewline
\hline
2& 2.1 &30 & $a_1$&$b_2$&1556&2,35\tabularnewline
\hline
3& 2.1 &30 & $a_1$&$b_3$&87151&69,25\tabularnewline
\hline
4& 2.1 &30 & $a_2$&$b_1$&907&2,12\tabularnewline
\hline
5& 2.1 &30 & $a_2$&$b_2$&2217&4,71\tabularnewline
\hline
6& 2.1 &30 & $a_2$&$b_3$&60067&51,60\tabularnewline
\hline
7& 2.1 &30 & $a_3$&$b_1$&1098&3,76\tabularnewline
\hline
8& 2.1 &30 & $a_3$&$b_2$&2580&9,41\tabularnewline
\hline
9& 2.1 &30 & $a_3$&$b_3$&34785&21,69\tabularnewline
\hline
10& 2.1 &30 & $a_4$&$b_1$&3125&11,18\tabularnewline
\hline
11& 2.1 &30 & $a_4$&$b_2$&3415&11,76\tabularnewline
\hline
12& 2.1 &30 & $a_4$&$b_3$&7917&11,90\tabularnewline
\hline
13& 2.1 &50 & $a_1$&$b_1$&875&2,00\tabularnewline
\hline
\end{tabular}
\end{center}

\begin{center}
\begin{tabular}{|>{\centering}p{2cm}|>{\centering}p{2cm}|c|c|c|>{\centering}p{4.5cm}|>{\centering}p{3.5cm}|}
\hline
14& 2.1 &50 & $a_1$&$b_2$&2520&4,98\tabularnewline
\hline
15& 2.1 &50 & $a_1$&$b_3$&-&-\tabularnewline
\hline
16& 2.1 &50 & $a_2$&$b_1$&1549&4,40\tabularnewline
\hline
17& 2.1 &50 & $a_2$&$b_2$&3926&8,6\tabularnewline
\hline
18& 2.1 &50 & $a_2$&$b_3$&-&-\tabularnewline
\hline
19& 2.1 &50 & $a_3$&$b_1$&1827&6,80\tabularnewline
\hline
20& 2.1 &50 & $a_3$&$b_2$&4350&17,8\tabularnewline
\hline
21& 2.1 &50 & $a_3$&$b_3$&-&-\tabularnewline
\hline
22& 2.1 &50 & $a_4$&$b_1$&5289&20,26\tabularnewline
\hline
23& 2.1 &50 & $a_4$&$b_2$&5913&22,94\tabularnewline
\hline
24& 2.1 &50 & $a_4$&$b_3$&-&-\tabularnewline
\hline
25& 2.2 &30 & $a_1$&$b_1$&672&0,94\tabularnewline
\hline
26& 2.2 &30 & $a_1$&$b_2$&1271&1,76\tabularnewline
\hline
27& 2.2 &30 & $a_1$&$b_3$&83976&64,11\tabularnewline
\hline
28& 2.2 &30 & $a_2$&$b_1$&903&2,59\tabularnewline
\hline
29& 2.2 &30 & $a_2$&$b_2$&2126&5,88\tabularnewline
\hline
30& 2.2 &30 & $a_2$&$b_3$&60254&45,52\tabularnewline
\hline
31& 2.2 &30 & $a_3$&$b_1$&1104&3,29\tabularnewline
\hline
32& 2.2 &30 & $a_3$&$b_2$&2788&8,24\tabularnewline
\hline
33& 2.2 &30 & $a_3$&$b_3$&41303&28,73\tabularnewline
\hline
34& 2.2 &30 & $a_4$&$b_1$&3310&11,76\tabularnewline
\hline
35& 2.2 &30 & $a_4$&$b_2$&3186&12,94\tabularnewline
\hline
36& 2.2 &30 & $a_4$&$b_3$&7827&11,7\tabularnewline
\hline
37& 2.2 &50 & $a_1$&$b_1$&804&1,60\tabularnewline
\hline
38& 2.2 &50 & $a_1$&$b_2$&2905&3,9\tabularnewline
\hline
39& 2.2 &50 & $a_1$&$b_3$&-&-\tabularnewline
\hline
40& 2.2 &50 & $a_2$&$b_1$&1592&3,60\tabularnewline
\hline
41& 2.2 &50 & $a_2$&$b_2$&3657&9,39\tabularnewline
\hline
42& 2.2 &50 & $a_2$&$b_3$&-&-\tabularnewline
\hline
43& 2.2 &50 & $a_3$&$b_1$&1844&6,00\tabularnewline
\hline
44& 2.2 &50 & $a_3$&$b_2$&4834&14,22\tabularnewline
\hline
\end{tabular}
\end{center} 

\begin{center}
\begin{tabular}{|>{\centering}p{2cm}|>{\centering}p{2cm}|c|c|c|>{\centering}p{4.5cm}|>{\centering}p{3.5cm}|}
\hline
45& 2.2 &50 & $a_3$&$b_3$&-&-\tabularnewline
\hline
46& 2.2 &50 & $a_4$&$b_1$&5856&24,58\tabularnewline
\hline
47& 2.2 &50 & $a_4$&$b_2$&5879&20,46\tabularnewline
\hline
48& 2.2 &50 & $a_4$&$b_3$&-&-\tabularnewline
\hline
\end{tabular}
\end{center} 

Сделаем комментарии к содержимому колонок. Нумерация численных экcпериментов указана в первой колонке, во второй колонке $-$ метод, с помощью которого решались тестовые задачи. Размерность пространтсва переменных приведена в третьей колонке, в четвертой $-$ способы задания чисел $\varepsilon_k$, а в пятой колонке приведены приемы обновления аппроксимирующих множеств. Количество итераций, которое необходимо для достижения заданной точности решения задачи, представлено в шестой колонке, а в седьмой $-$ время решения. В таблице символом <<->> обозначены примеры, время работы которых превысило 1 день.

Перейдем к обсуждению обозначений четвертой и пятой колонок. Для задания в таблице последовательности $\{\varepsilon_k\}$, $k\in K$, выбраны следующие соответствия:

1) $a_1$ $\leftrightarrow$ $\varepsilon_{k+1}=\frac{\varepsilon_k}{1.1}$ для всех $k\in K$;

2) $a_2$ $\leftrightarrow$ $\varepsilon_{k+1}=\frac{\varepsilon_k}{2}$ для всех $k\in K$;

3) $a_3$ $\leftrightarrow$ $\varepsilon_{k+1}=\frac{\varepsilon_k}{n}$ для всех $k\in K$;

4) $a_4$ $\leftrightarrow$ $\varepsilon_{k+1} = \alpha_k (f(x_k)-\sigma_k)$, $\alpha_k = \frac{1}{2^k}$ для всех $k\in K$.

Способы обновления множеств $Q_i$ обозначались следующим образом:

1)  $b_1$ $\leftrightarrow$$Q_i=M_i$ для всех $i\ne i_k$, а множество $Q_{i_k}$ состоит из <<активных>>  в текущей итерационной точке $(y_{i_k}, \gamma_{i_k})$ секущих плоскостей;

2)
 a) в методе 2.1:

 $b_2$ $\leftrightarrow$ $Q_i=M_i$, $0\le i< n+1$,

  $Q_i = \begin{cases}M_i, \ i\ne i_k, \\
 \overset{i_k}{\underset{r=i_k - n - 1}{\cap}}\{(x, \gamma)\in {\rm R_{n+1}} : f_j(y_r) + \langle a_j, x - y_r \rangle\le \gamma,  \\
 \hfill \forall j\in J, a_j\in \partial f_j(y_r) \},  i = i_k,  \end{cases}$

$i\ge n+1$;

b) в методе 2.2: 

 $b_2$ $\leftrightarrow$ $Q_i=M_i$, $0\le i< n+1$,

 $Q_i = \begin{cases}M_i, \ i\ne i_k, \\
 \overset{i_k}{\underset{r=i_k - n - 1}{\cap}}\{s\in {\rm R_{n+1}} :\langle a_j, s- z_r^j\rangle \le 0,  \forall j\in J, \\
\hfill a_j\in W^1(z_r^j, {\rm epi\,}(f, {\rm R_n})) \}, i = i_k,  \end{cases}$

$i\ge n+1$;

3)  $b_3$ $\leftrightarrow$ $Q_i = \begin{cases}M_i, & i \ne i_k, \\ M_0,& i = i_k. \end{cases}$.

Обратим внимание на то, что все результаты, приведенные в вышеуказанной таблице, получены методами 2.1 и 2.2 с привлечением тех или иных способов обновления аппроксимирующих множеств. Указанный пример  решались (при $n=30, 50$) методами 2.1, 2.2 и без использования каких-либо обновлений, то есть  при всех $i$ полагалось $Q_i=M_i$. В этих случаях указанная выше точность решения была достигнута методом 2.1 при $n=30$ за 342 итераций и 45 минут, в случае $n=50$ $-$ за 761 итераций и 69 минут, а методом 2.2 при $n=30$ $-$  за 402 итераций и 39 минут, в случае $n=50$ $-$ за 807 итераций и 76 минут.

Обсудим, наконец, полученные результаты. Для начала сравним тестовые примеры, при решении которых  аппроксимирующие множества $Q_i$ обновлялись полностью, с примерами, где множества $Q_i$ не обновлялись. Учитывая результаты тестирования примеров можно сделать вывод о том, что  отбрасывание всех секущих плоскостей, построенных к итерации $i=i_k$, не рекомендуется. Это объясняется тем, что при полном обновлении утрачивается приемлимый уровень аппроксимации надграфика целевой функции, и тем самым происходит резкий рост значения целевой функции в очередной итерационной точке. Следовательно, приходится делать  значительное количество итераций для достижения утраченного уровня аппроксимации. Таким образом, одновременно ростет трудоемкость решения задачи как по затраченному времени, так и по количеству проделанных итераций. 

Заметим, что если принято решение, при котором на итерациях $i=i_k$ требуется полностью отбросить все построенные ранее отсекающие плоскости, то последовательность $\{\varepsilon_k\}$, $k\in K$, рекомендуется выбирать быстроубывающей.

Перейдем к сравнению примеров, которые решались без отбрасывания секущих плоскостей, с примерами, при решении которых в множества $Q_{i_k}$ включались $n+1$ последние отсекающие гиперплоскости, либо все <<активные>> плоскости. Учитывая результаты численных экспериментов, можно сделать вывод о том, что указанные приемы обновления аппроксимирующих множеств существенно ускоряют по времени  процесс отыскания решения задачи. Выяснилось, что в момент обновления очередного аппроксимирующего множества $Q_{i_k}$ лучше выбирать прием с включением всех <<активных>> гиперплоскостей. Предпочтительность выбора этого приема обусловлена меньшими затратами по времени для решения тестовых задач.

Теперь приведем рекомендации по выбору последовательности $\{\varepsilon_k\}$, $k\in K$, с помощью которой задается частота обновления аппроксимирующих множеств. В начале параграфа было отмечено, что при полном обновлении аппроксимирующих множеств на итерациях $i=i_k$ последовательность $\{\varepsilon_k\}$, $k\in K$, следует задавать быстроубывающей. Но в случае включения в аппроксимирующие множества на итерациях $i=i_k$ всех <<активных>> отсекающих гиперплоскостей числа $\varepsilon_k$ следует выбирать медленноубывающей.

Заметим, что в результате  проделанных численных экспериментов  не выявлено преимуществ одного метода над другим.

\chapter[Методы отсечений с одновременной аппроксимацией допустимой области и надграфика целевой функции]{
\quad{} \quad{}\hspace{1mm}Методы отсечений
\newlineс одновременной аппроксимацией 
\newline  \hspace*{1.7cm}допустимой области 
\newline  и надграфика целевой функции}
\setcounter{equation}{0}

В главах 1, 2 были построены методы отсечений, в которых при построении приближений, использовалась аппроксимация либо допустимой области, либо надграфика целевой функции.  \mbox{В настоящей} главе разрабатываются методы отсечений,  в которых при построении приближений аппроксимируется и допустимая область, и надграфик целевой функции. Такая одновременная аппроксимация позволяет независимо от вида функций ограничений и целевой функции исходной задачи решать на каждом шаге при построении итерационной точки вспомогательную задачу линейного программирования. Материал этой главы опубликован авторами в  \cite{n_yarullin_r_s_5,n_yarullin_r_s_2018_LJM,n_yarullin_r_s_2016_IOP}.

\section
[Метод с аппроксимацией области ограничений и надграфика на основе обобщенно-опорных отсекающих плоскостей]
{Метод с аппроксимацией области ограничений 
\linebreak и надграфика на основе обобщенно-опорных отсекающих плоскостей}\label{s3.1}

Предлагается общий метод решения задачи выпуклого программирования с вложением погружающих множеств, в котором используется  одновременная аппроксимация надграфика и области ограничений.

Пусть $f(x)$, $F(x)$ $-$ выпуклые в ${\rm R_n}$ функции,
$$D=\{x\in {\rm R_n} : F(x)\le 0\},$$
функция $f(x)$ достигает на множестве $D$ минимального значения, внутренность множества $D$ непуста. Решается задача
\begin{equation}\label{e3.1.1}
\min\{f(x) : x\in D\}.
\end{equation}

Как обычно будем использовать обозначения $f^*$, $X^*$, $K$, определенные  в $\S \, 1.1$. Пусть  ${\rm epi\,}(f, G)=\{u=(x,\gamma)\in {\rm R_{n+1}} : x\in G,\ \gamma\ge f(x)\}$, где $G\subset {\rm R_n}$,  $W^\prime(z, Q)=\{b\in {\rm R_{n+1}} : \|b\|=1,\ \langle b, u- z\rangle \le 0 \ \forall u\in Q\}$ и  $W^{\prime\prime}(t, L) = \{a\in {\rm R_n} : \|a\| = 1, \ \langle a, y - t\rangle \le 0 \ \forall y\in L\}$ $-$ множества нормированных обобщенно-опорных векторов для  множества $Q\subset {\rm R_{n+1}}$ в точке $z\in {\rm R_{n+1}}$ и множества $L\subset {\rm R_n}$ в точке $t\in {\rm R_n}$ соответственно, $x^*\in X^*$.

Предлагаемый метод  вырабатывает последовательность приближений $\{(y_i$, $\gamma_i)\}$, $i\in K$, и заключается в следующем.

\textbf{Метод 3.1.} Строятся выпуклое ограниченное замкнутое множество $G_0\subset {\rm R_n}$ и выпуклое замкнутое множество $M_0\subset {\rm R_{n+1}}$ такие, что
\begin{equation}\label{e3.1.2}
x^*\in G_0, \quad {\rm epi\,} (f, {\rm R_n})\subset M_0.
\end{equation}
Выбираются точки $v^\prime \in{\rm int\, epi\,}(f, {\rm R_n})$, $v^{\prime\prime} \in{\rm int\,} D$. Задается такое число $\alpha$, что $\alpha\le \bar{f}=\min\{f(x):x\in G_0\}$. Полагается $i=0$.

\textbf{1.} Находится решение $u_i=(y_i,\gamma_i)$, где $y_i\in {\rm R_n}$, $\gamma_i\in {\rm R_1}$, следующей задачи:
\begin{equation}\label{e3.1.3}
\gamma \to \min
\end{equation}
\begin{equation}\label{e3.1.4}
(x,\gamma)\in M_i, \quad x\in G_i, \quad \gamma\ge \alpha.
\end{equation}
Если 
\begin{equation}\label{e3.1.5}
y_i\in D, \quad f(y_i)=\gamma_i,
\end{equation}
то $y_i$ $-$ решение задачи (\ref{e3.1.1}), и  процесс построения итерационных точек завершается.

\textbf{2.} В интервале $(v^\prime, u_i)$ выбирается точка $z_i^\prime$ таким образом, чтобы $z_i^\prime\notin {\rm int\,epi\,}(f, {\rm R_n})$ и при некотором $q_i^\prime\in [1,q^\prime]$, $q^\prime<+\infty$, для точки 
\begin{equation}\label{e3.1.6}
\bar{v}_i^\prime=u_i+q_i^\prime(z_i^\prime-u_i)
\end{equation}
выполнялось включение $\bar{v}_i^\prime\in {\rm epi\,} (f, {\rm R_n})$.

\textbf{3.} Выбирается конечное множество $B_i\subset W^\prime(z_i^\prime, {\rm epi\,} (f, R_{n}))$, и полагается
\begin{equation}\label{e3.1.7}
M_{i+1}=M_i\bigcap T_i,
\end{equation}
где $T_i=\{u\in {\rm R_{n+1}} : \langle b, u-z_i^\prime\rangle\le 0 \ \forall \ b\in B_i\}$.

\textbf{4.}
 Если $y_i\notin D$, то в интервале $(v^{\prime\prime}, y_i)$ выбирается точка $z_i^{\prime\prime}$ так, чтобы $z_i^{\prime\prime}\notin{\rm int\,} D$ и при некотором $q_i^{\prime\prime} \in [1, q^{\prime\prime}]$, $q^{\prime\prime}<+\infty$, для точки
\begin{equation}\label{e3.1.8.0}
\bar{v}_i^{\prime\prime} = y_i + q_i^{\prime\prime}(z_i^{\prime\prime} - y_i)
\end{equation}
выполнилось включение $\bar{v}_i^{\prime\prime} \in D$. Следует переход к очередному пункту. В противном случае  полагается
\begin{equation}\label{e3.1.8}
G_{i+1}=G_i,
\end{equation}
$z^{\prime\prime}_i =\bar{v}^{\prime\prime}_i =  y_i$, и следует переход к п. 6.

\textbf{5.}  Выбирается конечное множество $A_i\subset W^{\prime\prime}(z_i^{\prime\prime}, D)$, и полагается
\begin{equation}\label{e3.1.9}
G_{i+1}=G_i\bigcap P_i,
\end{equation}
где $P_i = \{x\in {\rm R_n} : \langle a, x - z_i^{\prime\prime}\rangle \le 0 \ \forall a\in A_i\}$.

\textbf{6.} Значение $i$ увеличивается на единицу, и следует переход к п. 1.

Сделаем относительно предложенного алгоритма некоторые замечания. Покажем прежде всего, что задача (\ref{e3.1.3}), (\ref{e3.1.4}) разрешима при любом $i\in K$.

\begin{lemma}\label{lemma3.1.1}
 Для каждого $i\in K$ выполняются следующие включения:
\begin{equation}\label{e3.1.10}
x^*\in G_i,
\end{equation} 
\begin{equation}\label{e3.1.11}
{\rm epi\,} (f, {\rm R_n})\subset M_i.
\end{equation}
\end{lemma}

\begin{proof}
Поскольку множества $M_0$, $G_0$ выбраны согласно (\ref{e3.1.2}), то  справедливы включения (\ref{e3.1.10}), (\ref{e3.1.11}) при $i=0$. Допустим теперь, что (\ref{e3.1.10}), (\ref{e3.1.11}) справедливы при любом фиксированном $i=l\ge 0$. Покажем, что
\begin{equation}\label{e3.1.12}
x^*\in G_{l+1},
\end{equation}
\begin{equation}\label{e3.1.13}
{\rm epi\,} (f, {\rm R_n})\subset M_{l+1},
\end{equation}
тогда утверждение леммы будет доказано.

Cогласно алгоритму множество $G_{l+1}$ задается  в виде (\ref{e3.1.8}) или (\ref{e3.1.9}). Если $G_{l+1}=G_l$, то включение (\ref{e3.1.12}) выполняется ввиду сделанного индукционного предположения. Поэтому будем считать, что $G_{l+1}=G_l\bigcap P_l$. Cогласно выбору точки $z_l^{\prime\prime}$ и множества $A_l$ выполняется включение $D\subset P_l$.  Отсюда и из (\ref{e3.1.9}) c учетом того, что $x^*\in G_l$, следует  (\ref{e3.1.12}). 

Далее, по индукционному предположению ${ \rm epi\,} (f, {\rm R_n})\subset M_l$. Кроме того, по выбору точки $z_l^{\prime}$ и множества $B_l$ выполняется включение ${\rm epi\,} (f, {\rm R_n})\subset T_l$. Тогда ввиду (\ref{e3.1.7}) справедливо соотношение (\ref{e3.1.13}). Лемма доказана.
\end{proof}

В силу выбора значений $\alpha$, $\bar{f}$ и множества $G_0$
\begin{equation}\label{e3.1.14}
f^*\ge \alpha.
\end{equation}
Тогда из леммы \ref{lemma3.1.1} и неравенства (\ref{e3.1.14}) вытекает

\begin{lemma}\label{lemma3.1.2}
 Для всех $i\in K$ точка $(x^*, f^*)$  удовлетворяет ограничениям задачи (\ref{e3.1.3}), (\ref{e3.1.4}).
\end{lemma}

Леммой \ref{lemma3.1.2} обоснована разрешимость задачи (\ref{e3.1.3}), (\ref{e3.1.4}).  Докажем теперь критерий остановки, заложенный в п. 1 алгоритма.

Поскольку для всех точек $(x,\gamma)$, удовлетворяющих условиям (\ref{e3.1.4}), при всех $i\in K$ выполняется неравенство $\gamma_i\le \gamma$, то из леммы \ref{lemma3.1.2} следует
\newline

\begin{lemma}\label{lemma3.1.3}
Пусть последовательность $\{\gamma_i\}$, $i\in K$, построена предложенным алгоритмом. Тогда справедливы неравенства
\begin{equation}\label{e3.1.15}
\gamma_i\le f^* \quad \forall i\in K.
\end{equation}
\end{lemma}

\begin{theorem}\label{theorem3.1.1}
Пусть при некотором $i\in K$ выполняются соотношения (\ref{e3.1.5}). Тогда $y_i\in X^*$.
\end{theorem}

\begin{proof}
Поскольку выполняются неравенства (\ref{e3.1.15}) и второе из условий (\ref{e3.1.5}), то  $f(y_i)\le f^*$. С другой стороны, $f(y_i)\ge f^*$, так как  по условию выполняется включение $y_i\in D$. Таким образом, $f(y_i)=f^*$, и теорема доказана.
\end{proof}

\begin{remark}\label{r3.1.1}
Аппроксимирующие множества $M_0$, $G_0$  удобно задавать многогранными  с практической точки зрения. В таком случае в силу (\ref{e3.1.7}) $-$ (\ref{e3.1.9}) задачи (\ref{e3.1.3}), (\ref{e3.1.4}) при всех $i\in K$ будут задачами линейного программирования.
\end{remark}

\begin{remark}\label{r3.1.2}
Множество $M_0$ можно задать совпадающим со всем пространством ${\rm R_{n+1}}$. Этим объясняется наличие неравенства $\gamma\ge \alpha$ в ограничениях (\ref{e3.1.4}).
\end{remark}

\begin{remark}\label{remark3.1.3_2019}
В случае, если допустимое множество $D$ задано системой линейных неравенств, то естественно положить $G_0=D$, поскольку нет необходимости в построении аппроксимирующих $D$ множеств. Тогда согласно (\ref{e3.1.8}) будут выполняться равенства $G_i=D$ для всех $i\in K$.
\end{remark}

\begin{remark}\label{remark3.1.4_2019}
Если 
$$f(x)=\max \limits_{1\le j\le m} f_j(x),$$
и функции $f_j(x)$ линейны, то удобно положить $M_0={\rm epi\,} (f, {\rm R_n})$. Тогда точками $z_i^\prime$ можно считать точки пересечения интервала $(v^\prime, u_i)$  с границей множества ${\rm epi\,} (f, {\rm R_n})$, чтобы  полагать $M_{i+1}=M_i$ для всех $i\in K$.
\end{remark}

\begin{lemma}\label{lemma3.1.4}
Пусть последовательность $\{u_i\}$,  $i\in K$, построена предложенным методом 3.1, и $\{u_i\}$, $i\in K'\subset K$, $-$ ее сходящаяся подпоследовательность. Тогда 
\begin{equation}\label{e3.1.16}
\lim \limits_{i\in K^\prime} \|z_i^\prime - u_i\|=0.
\end{equation}
\end{lemma}
\begin{proof}
Заметим, что для каждого $i\in K$ 
\begin{equation}\label{e3.1.17}
z_i^\prime = u_i + \tau_i (v^\prime  - u_i),
\end{equation}
где $0<\tau_i < 1$. Для произвольного $i\in K^\prime$ зафиксируем номер $p_i\in K^\prime$ такой, что $p_i > i$. \mbox{В силу} (\ref{e3.1.7}) $M_{p_i}\subset M_i$. Кроме того, ввиду (\ref{e3.1.4}) $u_{p_i}\in M_{p_i}$, а любой элемент множества $B_i$ является обобщенно-опорным для множества $M_{p_i}$ в точке $z^\prime_i$. Следовательно,
$$\langle b, u_{p_i} - z^\prime_i\rangle \le 0 \quad b\in B_i,$$
а с учетом (\ref{e3.1.17})
$$\langle b, u_i - u_{p_i}\rangle \ge \tau_i \langle b, u_i - v^\prime \rangle \quad \forall b\in B_i.$$

Согласно лемме \ref{lemma1.1.1} найдется такое число $\lambda >0$, что 
$$\langle b, u_i - v\rangle \ge \lambda \quad b\in B_i,$$
а поскольку $\|b\|=1$ для всех $b\in B_i$, то
\begin{equation}\label{e3.1.18}
\|u_i - u_{p_i}\|\ge \tau_i \lambda \quad \forall i, p_i \in K^\prime, \quad p_i > i.
\end{equation}
Так как подпоследовательность $\{u_i\}$, $i\in K^\prime$, по условию сходится, то в силу (\ref{e3.1.18}) $\tau_i \to 0$, $i\to \infty$, $i\in K^\prime$. Поэтому из (\ref{e3.1.17}) с учетом ограниченности последовательности $\{\|v^\prime - u_i\|\}$, $i\in K^\prime$, следует (\ref{e3.1.16}). Лемма доказана.
\end{proof}

\begin{lemma}\label{lemma3.1.5}
Пусть последовательность $\{y_i\}$, $i\in K$, построена предложенным методом 3.1. Тогда любая предельная точка последовательности $\{y_i\}$, $i\in K$, принадлежит множеству $D$.
\end{lemma}
\begin{proof}
Выделим из ограниченной последовательности $\{y_i\}$, $i\in K$, сходящуюся подпоследовательность $\{y_i\}$, $i\in K^\prime \subset K$. Пусть $\bar{y}$ $-$ предельная точка последовательности $\{y_i\}$, $i\in K^\prime$. Покажем справедливость включения
\begin{equation}\label{e3.1.19}
\bar{y} \in D.
\end{equation}

Положим
$$ K^\prime_1 = \{i\in K^\prime : y_i\in D\}, \quad K^\prime_2 = K^\prime \setminus K^\prime_1.$$
Хотя бы одно из множеств $K^\prime_1$, $K^\prime_2$ состоит из бесконечного числа номеров. Рассмотрим 2 случая. 

1) Если  множество $K^\prime_1$ бесконечно, то включение  (\ref{e3.1.19}) имеет место в силу построения $K^\prime_1$.

2) Пусть теперь множество $K_2^\prime$ состоит из бесконечного числа номеров.
Докажем сначала равенство
\begin{equation}\label{e3.1.20.0}
\lim \limits_{i\in K^\prime_2} \|z^{\prime\prime}_i - y_i\|=0,
\end{equation}
учитывая, что вместе с последовательностью $\{y_i\}$, $i\in K^\prime_2$,  построены и последовательности $\{z_i^{\prime\prime}\}$, $\{\bar{v}_i^{\prime\prime}\}$, $i\in K^\prime_2$.

Заметим, что для всех $i\in K$
\begin{equation}\label{e3.1.20}
z^{\prime\prime}_i = y_i + \gamma_i ( v^{\prime\prime} - y_i),
\end{equation}
где $\gamma_i \in [0, 1)$, причем $\gamma_i > 0$ для всех $i\in K_2^\prime$. Для произвольного $i\in K^\prime_2$ зафиксируем номер $p_i\in K_2^\prime$ такой, что $p_i>i$. \mbox{В силу} (\ref{e3.1.8}), (\ref{e3.1.9}) выполняется включение $G_{p_i}\subset G_i$. Кроме того, ввиду (\ref{e3.1.4}) $y_{p_i}\in G_{p_i}$, а любой элемент множества $A_i$ является обобщенно-опорным и для множества $G_{p_i}$ в точке $z_i^{\prime\prime}$. Следовательно,
$$\langle a, y_{p_i} - z_i^{\prime\prime}\rangle \le 0 \quad \forall a\in A_i.$$
Отсюда с учетом (\ref{e3.1.20}) для всех $a\in A_i$ имеем
$$\langle a, y_i - y_{i_p}\rangle \ge \gamma_i \langle a, y_i - v^{\prime\prime}\rangle.$$
По лемме \ref{lemma1.1.1} найдется такое число $\delta>0$, что $\langle a, y_i - v^{\prime\prime}\rangle \ge \delta$ для всех $i\in K^\prime_2$, $a\in A_i$. Значит, $\langle a, y_i - y_{i_p}\rangle \ge \gamma_i \delta$ для всех $a\in A_i$, а поскольку $\|a\|=1$ для всех $a\in A_i$, то
\begin{equation}\label{e3.1.21}
\|y_i - y_{p_i}\|\ge \gamma_i \delta \quad \forall i, p_i\in K^\prime_2, \quad p_i >i.
\end{equation}
Так как последовательность $\{y_i\}$, $i\in K^\prime_2$, является сходящейся, то согласно (\ref{e3.1.21}) $\gamma_i \to 0$, $i\to \infty$, $i\in K^\prime_2$. Поэтому из (\ref{e3.1.20}) с учетом ограниченности последовательности $\{\|v^{\prime\prime} - y_i\|\}$, $i\in K^\prime_2$, следует равенство (\ref{e3.1.20.0}).

Cогласно п. 4 метода для каждого $i\in K^\prime_2$ точка $\bar{v}^{\prime\prime}_i$  имеет вид (\ref{e3.1.8.0}).  Так как последовательность $\{y_i\}$, $i\in K^\prime_2$, ограничена, то отсюда с учетом (\ref{e3.1.20}) следует ограниченность последовательности $\{\bar{v}^{\prime\prime}_i\}$, $i\in K^\prime_2$. Выделим теперь из последовательности $\{\bar{v}^{\prime\prime}_i\}$, $i\in K^\prime_2$, сходящуюся подпоследовательность $\{\bar{v}^{\prime\prime}_i\}$, $i\in K^\prime_3\subset K^\prime_2$, и пусть $\bar{v}$ $-$ ее предельная точка. Заметим, что $\bar{v}\in D$ в силу замкнутости множества $D$. Перейдем теперь к пределу в равенстве (\ref{e3.1.8.0}) по $i\to \infty$, $i\in K^\prime_3$, с учетом (\ref{e3.1.20.0}). Тогда получим равенство $\bar{y}=\bar{v}$, из которого следует включение (\ref{e3.1.19}).  Лемма доказана.
\end{proof}

\begin{lemma}\label{lemma3.1.6}
Пусть последовательность $\{u_i\}$, $i\in K$, построена методом 3.1. Тогда любая ее предельная точка принадлежит множеству ${\rm epi\,}$ $(f, {\rm R_n})$.
\end{lemma}
\begin{proof}
Пусть подпоследовательность $\{u_i\}$, $i\in K^\prime\subset K$, сходится, и $\bar{u}$ $-$ ее предельная точка. Покажем, что 
\begin{equation}\label{e3.1.24}
\bar{u}\in {\rm epi\,}(f, {\rm R_n}).
\end{equation}
Так как $\{u_i\}$, $i\in K$, ограничена, то  из (\ref{e3.1.6}), (\ref{e3.1.16}) следует ограниченность последовательности $\{\bar{v}^\prime_i\}$, $i\in K^\prime$. Выделим из последовательности $\{\bar{v}_i^\prime\}$, $i\in K^\prime$, сходящуюся подпоследовательность $\{\bar{v}^\prime_i\}$, $i\in K^{\prime\prime}\subset K^\prime$, и пусть $\bar{v}$ $-$  ее предельная точка. Заметим, что 
\begin{equation}\label{e3.1.25}
\bar{v}\in {\rm epi\,}(f, {\rm R_n})
\end{equation}
в силу замкнутости множества ${\rm epi\,}(f, {\rm R_n})$. Но согласно (\ref{e3.1.6}) и равенству $\lim \limits_{i\in K^{\prime\prime}} \|z^\prime_i - u_i\|=0$, вытекающему из леммы \ref{lemma3.1.4}, имеем $\bar{v}=\bar{u}$. Отсюда и из (\ref{e3.1.25}) следует (\ref{e3.1.24}). Лемма доказана.
\end{proof}

\begin{theorem}
Пусть последовательность $\{(y_i, \gamma_i)\}$, $i\in K$, построена предложенным методом. Тогда для любой ее предельной точки $\bar{u}=(\bar{y},\bar{\gamma})$ справедливо включение $\bar{y}\in X^*$. 
\end{theorem}

\begin{proof}
 Пусть $\bar{u}=(\bar{y},\bar{\gamma})$ $-$ предельная точка сходящейся подпоследовательности $\{(y_i, \gamma_i)\}$, $i\in K^\prime\subset K$. Сразу отметим, что $\bar{y}\in D$ согласно лемме \ref{lemma3.1.5}, и, значит,
\begin{equation}\label{e3.1.26}
f(\bar{y})\ge f^*.
\end{equation}
Кроме того, в силу (\ref{e3.1.24})
\begin{equation}\label{e3.1.27}
\bar{\gamma}\ge f(\bar{y}).
\end{equation}

Далее, поскольку $(y_i, f(y_i))\in {\rm epi\,}(f, {\rm R_n})$, $i\in K$, то ввиду (\ref{e3.1.11})
\begin{equation}\label{e3.1.28}
(y_i, f(y_i))\in M_i \quad \forall i\in K.
\end{equation}
Так как $y_i\in G_i\subset G_0$, $i\in K$, то $f(y_i)\ge \bar{f}$, $i\in K$, и, следовательно,
\begin{equation}\label{e3.1.29}
f(y_i)\ge \alpha \quad \forall i\in K.
\end{equation}
Согласно (\ref{e3.1.3}), (\ref{e3.1.4}) для любой точки $(x, \gamma)\in {\rm R_{n+1}}$, удовлетворяющей условиям (\ref{e3.1.4}), выполняется неравенство $\gamma\ge \gamma_i$ при каждом $i\in K$. Но для точки $(y_i, f(y_i))$ эти условия имеют место согласно включениям $y_i\in G_i$ и (\ref{e3.1.28}), а также неравенству (\ref{e3.1.29}). Таким образом, 
\begin{equation}\label{e3.1.30}
f(y_i)\ge \gamma_i
\end{equation}
для всех $i\in K$. Переходя в (\ref{e3.1.30}) к пределу по $i\to \infty$, $i\in K^\prime$, получим $f(\bar{y})\ge \bar{\gamma}$. Отсюда  и из (\ref{e3.1.27}) следует равенство $f(\bar{y})=\bar{\gamma}$, а значит, из (\ref{e3.1.15}) имеем $f(\bar{y})\le f^*$. Тогда с учетом (\ref{e3.1.26}) $f(\bar{y})=f^*$, и теорема доказана.

\end{proof}
\setcounter{equation}{0}
\section
{Двухэтапный метод отсечений}\label{s3.2}

В данном параграфе предлагается метод условной минимизации выпуклой негладкой функции. Предлагаемый метод отличается от метода \ref{s3.1} возможностью периодического обновления аппроксимирующих надграфик целевой функции множеств. 

Пусть функция $f(x)$ выпукла в $n$-мерном евклидовом пространстве ${\rm R_n}$,  множество $D\subset {\rm R_n}$ выпукло и замкнуто. Решается задача (\ref{e3.1.1}). 

Пусть обозначения $f^*$, $X^*$,  $K$ определены как в $\S \ \ref{s2.1}$, а ${\rm int}\, Q$, $W^1(u, Q)$ -- как в $\S \ \ref{s2.4}$. Положим $x^*\in X^*$,  ${\rm epi}\,(f, {\rm R_n}) = \{(x, \gamma)\in {\rm R_{n+1}} : x\in {\rm R_n}, \gamma \ge f(x)\}$.

\textbf{Метод 3.2.} Предлагаемый метод решения задачи (\ref{e3.1.1}) вырабатывает вспомогательную последовательность $\{u_i\}$, $i\in K$, и основную последовательность приближений $\{x_k\}$, $k\in K$, и заключается в следующем. Выбираются точки $v^\prime \in {\rm int}\, D$ и $v^{\prime\prime} \in {\rm int\, epi}\, (f, {\rm R_n})$. Выбираются выпуклое ограниченное замкнутое множество $M_0\subset {\rm R_n}$ и выпуклое замкнутое множество $G_0\subset {\rm R_{n+1}}$ такие, что $x^*\in M_0$, ${\rm epi}\,(f, {\rm R_n})\subset G_0$. Задаются числа $\bar{\gamma}$, $\varepsilon_k$, $k\in K$,  с условиями, что $\bar{\gamma}\le \min \{f(x) : x \in M_0\}$, 
\begin{equation}\label{o2}
\varepsilon_k>0, \quad k\in K,\quad \varepsilon_k\to 0,\quad k\to \infty.
\end{equation}
Полагается $i=0$, $k=0$.

\textbf{1.} Находится решение $u_i=(y_i, \gamma_i)$, где $y_i\in {\rm R_n}$, $\gamma_i\in {\rm R_1}$, следующей задачи:
\begin{equation}\label{o3}
\min\{\gamma :(x, \gamma)\in G_i, \ x\in M_k, \ \gamma \ge \bar{\gamma}\}.
\end{equation}

\textbf{2.} Отыскивается $\bar{u}_i=(\bar{y}_i, \bar{\gamma}_i)\in {\rm R_{n+1}}$, как точка пересечения отрезка  $[v^{\prime\prime}, u_i]$ с границей множества ${\rm epi}\, (f, {\rm R_n})$. Если
\begin{equation}\label{o4}
\bar{u}_i=u_i, \quad y_i\in D,
\end{equation}
то $y_i\in X^*$, и процесс решения задачи (\ref{e3.1.1}) окончен.

\textbf{3.} Если 
\begin{equation}\label{o5}
\|\bar{u}_i - u_i\|>\varepsilon_k,
\end{equation}
то полагается
\begin{equation}\label{o6}
G_{i+1}=S_i \bigcap \{u\in {\rm R_{n+1}} : \langle b_i, u - \bar{u}_i\rangle \le 0\},
\end{equation} 
где
\begin{equation}\label{o7}
b_i \in W^1(\bar{u}_i, {\rm epi}\, (f, {\rm R_n})),
\end{equation}
\begin{equation}\label{o8}
S_i=G_i,
\end{equation}
и следует переход к п. 1 при $i$, увеличенном на единицу. Иначе выполняется п. 4.

\textbf{4.} Полагается $i_k=i$,
\begin{equation}\label{o9}
x_k=y_{i_k}, \quad \sigma_k=\gamma_{i_k},
\end{equation}
и строится множество $G_{i+1}$ согласно (\ref{o6}), (\ref{o7}), где $S_i$ -- выпуклое замкнутое множество, выбранное с условием
\begin{equation}\label{o10}
{\rm epi}\, (f, {\rm R_n})\subset S_i.
\end{equation}

\textbf{5.} Если $x_k\in D$, то полагается $M_{k+1}=M_k$. В противном случае отыскивается         $\bar{x}_k$, как точка пересечения отрезка  $[v^\prime, x_k]$ с границей множества $D$, и полагается 
$$M_{k+1}=M_k\bigcap \{x\in {\rm R_n} : \langle a_k, x - \bar{x}_k\rangle \le 0\},$$
где $a_k\in W^1(\bar{x}_k, D)$.

\textbf{6.} Значения $i$ и $k$ увеличиваются на единицу, и следует переход к п. 1.

Сделаем некоторые замечания к предложенному методу. Прежде всего обоснуем критерий остановки, заложенный в п. 2 метода.

\begin{theorem} Пусть при некотором $i\in K$ выполняются соотношения (\ref{o4}). Тогда $y_i$ -- решение задачи (\ref{e3.1.1}).
\end{theorem}

\begin{proof}  C учетом условий выбора множеств $M_0$, $G_0$, условия (\ref{o10}) и способа построения секущих плоскостей легко доказать по индукции, что для всех $i\in K$, $k\in K$ точка  $(x^*, f^*)$ удовлетворяет ограничениям задачи (\ref{o3}). Значит, для решения $(y_i, \gamma_i)$ этой задачи при любых $i\in K$, $k\in K$ справедливо неравенство
\begin{equation}\label{o11}
\gamma_i\le f^*.
\end{equation}

Пусть номер $i\in K$ такой, что выполняется (\ref{o4}). Тогда $\bar{y}_i=y_i$, $\bar{\gamma}_i=\gamma_i=f(\bar{y}_i)$, и следовательно, $f(y_i)=\gamma_i$. Но по условию $y_i\in D$, и поэтому $f(y_i)\ge f^*$. С другой стороны, ввиду (\ref{o11})  $f(y_i)\le f^*$. Таким образом, $f(y_i) = f^*$, и утверждение теоремы доказано.
\end{proof}

Далее, обратим внимание на то, что на итерациях с номерами $i=i_k$ можно производить обновления аппроксимирующих множеств $G_{i+1}$. Эти обновления проводятся за счет выбора множеств $S_{i_k}$. А именно, ниже будет показано, что для каждого $k\in K$ при построении  $\{u_i\}$, $i\in K$, зафиксируется такой номер $i=i_k$, что выполнится неравенство
\begin{equation}\label{o12}
\|\bar{u}_i - u_i\|\le \varepsilon_k.
\end{equation}
Тогда с учетом (\ref{o10}) можно положить, например, $S_{i_k}=G_{r_i}$, где $0\le r_i\le i_k-1$. В таком случае будет отброшено некоторое число построенных к шагу $i_k$ последних отсечений.  Можно, напротив, $S_{i_k}$ формировать на основе последних отсекающих плоскостей, в частности, привлекая лишь активные в точке $u_{i_k}$ плоскости.

Отметим, что согласно (\ref{o10}) допустимо положить $S_{i_k}=G_{i_k}$ при каждом $k\in K$. Тогда независимо от условий (\ref{o5}), (\ref{o12})  в (\ref{o6}) полагается $S_i=G_i$, и обновлений аппроксимирующих множеств не происходит.

\begin{lemma}\label{lem_3.2.1}
 Пусть последовательность $\{u_i\}$, $i\in K$, построена с условием, что $S_i$ имеют вид (\ref{o8}) для всех $i\in K$. Тогда $\lim \limits_{i\in K} \| \bar{u}_i - u_i\|=0.$
\end{lemma}

\begin{proof} 
Допустим, что утверждение неверно. Тогда найдется подпоследовательность  $\{u_i\}$, $i\in K^\prime\subset K$,  такая, что
\begin{equation}\label{o13}
\|\bar{u}_i - u_i\|\ge \Delta>0\quad \forall i\in K^\prime.
\end{equation}
Выделим из последовательности $\{u_i\}$, $i\in K^\prime$, сходящуюся подпоследовательность $\{u_i\}$, $i\in K^{\prime\prime}\subset K^\prime$. Пусть номера $i$, $p_i\in K^{\prime\prime}$ такие, что $p_i>i$. В силу (\ref{o6}), (\ref{o8}) $G_{p_i}\subset G_i$. Но $u_{p_i}\in G_{p_i}$, и, кроме того, $b_i\in W^1(\bar{u}_i, G_{p_i})$. Следовательно, $\langle b_i, u_{p_i} - \bar{u}_i\rangle \le 0$. Так как для каждого $l\in K$ 
\begin{equation}\label{o14}
\bar{u}_l=u_l +\alpha_l(v^{\prime\prime}-u_l),
\end{equation}
где $\alpha_l\in (0, 1)$, то отсюда и из последнего неравенства имеем
$\langle b_i, u_i - u_{p_i}\rangle \ge \alpha_i \langle b_i, u_i - v^{\prime\prime}\rangle.$
Поскольку $v^{\prime\prime}\in {\rm int \, epi}\,(f, {\rm R_n})$, то нетрудно доказать существование такого $\sigma>0$, что $\langle b_l, v^{\prime\prime} - u_i\rangle \le -\sigma$ для всех $l\in K$. В связи с этим $\langle b_i, u_i - u_{p_i}\rangle \ge \alpha_i\sigma$, или $\|u_i - u_{p_i}\|\ge \alpha_i\sigma$. Отсюда и из сходимости последовательности  $\{u_i\}$, $i\in K^{\prime\prime}$, имеем $\alpha_i\to 0$, $i\in K^{\prime\prime}$. Значит, из равенства (\ref{o14}), где $l=i$, следует, что 
$\|\bar{u}_i - u_i\|\to 0$, $i\in K^{\prime\prime}.$ Это предельное соотношение противоречит (\ref{o13}). Лемма доказана.
\end{proof}

На  основе леммы \ref{lem_3.2.1} покажем теперь, что при построении методом последовательности $\{u_i\}$, $i\in K$, будет построена согласно (\ref{o9}) и последовательность $\{(x_k, \sigma_k)\}$, $k\in K$.

\begin{lemma} Пусть $\{u_i\}$, $i\in K$, построена предложенным методом. Тогда для каждого $k\in K$ существует такой номер $i=i_k$, что выполняется (\ref{o12}).
\end{lemma}

\begin{proof}
Пусть $k\in K$. Докажем справедливость неравенства (\ref{o12}) при некотором $i=i_k$. Допустим противное, т.е. для всех $i\in K$ выполняется (\ref{o5}). Тогда для всех $i\in K$ множества $S_i$ имеют вид (\ref{o8}), и утверждение леммы \ref{lem_3.2.1} противоречит неравенству (\ref{o5}).
\end{proof}

\begin{lemma}\label{lem_3.2.3}
Пусть $\{x_k\}$, $k\in K^\prime \subset K$, -- сходящаяся подпоследовательность последовательности  $\{x_k\}$, $k\in K$, и $\bar{x}$ -- ее предельная точка. Тогда выполняется включение
\begin{equation}\label{o15}
\bar{x}\in D.
\end{equation}
\end{lemma}

\begin{proof}
Если предположить, что для бесконечного числа номеров $k\in K^\prime$ выполняется включение $x_k\in D$, то ввиду замкнутости множества $D$ утверждение (\ref{o15}) очевидно. Поэтому будем предполагать, что $x_k\notin D$ для всех номеров $k\in K^\prime$, начиная с номера $N\in K^\prime$.

Отметим, что
\begin{equation}\label{o16}
\bar{x}_k=x_k +\theta_k(v^\prime - x_k), \quad k\in K^\prime, \quad k\ge N,
\end{equation}
где $\theta_k\in (0, 1)$. Пусть номера $l$, $p_l\in K^\prime$ таковы, что $p_l>l\ge N$. Поскольку $M_{p_l}\subset M_l$, то $a_l\in W^1(\bar{x}_l, M_{p_l})$.  Но согласно  (\ref{o9}), (\ref{o13}) $x_{p_l}\in M_{p_l}$, следовательно, $\langle a_l, x_{p_l}-\bar{x}_l\rangle \le 0,$ и с учетом (\ref{o16}) 
\begin{equation}\label{o17}
\langle a_l, x_l-x_{p_l}\rangle \ge \theta_l\langle a_l, x_l - v^\prime\rangle.
\end{equation}
Так как $x_k\notin D$, $k\ge N$, $k\in K^\prime$, а $v^\prime\in {\rm int}\, D$, то найдется такое  $\delta >0$, что $\langle a_k, v^\prime - x_k\rangle\le -\delta$ для всех  $k\ge N$, $k\in K^\prime$. Отсюда и из  (\ref{o17}) с учетом, что  $\|a_k\|=1$, имеем для любых $l$, $p_l\in K^\prime$ таких, что $p_l>l\ge N$, неравенство  $\|x_l-x_{p_l}\|\ge \theta_l\delta$. Из этого неравенства ввиду сходимости последовательности $\{x_k\}$, $k\in K^\prime$, следует, что  $\theta_k \to 0$, $k\in K^\prime$. Тогда согласно  (\ref{o16}) и ограниченности $\{\|v^\prime - x_k\|\}$, $k\in K^\prime$, справедливо
\begin{equation}\label{o18}
\|\bar{x}_k -x_k\|\to 0, \quad k\in K^\prime.
\end{equation}
Выделим теперь из последовательности $\{\bar{x}_k\}$, $k\in K^\prime$, сходящуюся подпоследовательность $\{\bar{x}_k\}$, $k\in K^{\prime\prime}\subset K^\prime$, и пусть  $\tilde{x}$  -- ее предельная точка. Тогда в силу (\ref{o18}) $\tilde{x}=\bar{x}$, и из замкнутости множества  $D$ следует включение   (\ref{o15}). Лемма доказана.
\end{proof}

Сформулируем, наконец, теорему сходимости предложенного метода.

\begin{theorem} Для любой предельной точки $(\bar{x}, \bar{\sigma})$  последовательности $\{(x_k, \sigma_k)\}$, $k\in K$, построенной предложенным методом, имеют место следующие соотношения:  $\bar{x}\in X^*$, $\bar{\sigma}=f^*$.
\end{theorem}

\begin{proof} Пусть $(\bar{x}, \bar{\sigma})$ -- предельная точка сходящейся подпоследовательности  $\{(x_k, \sigma_k)\}$, $k\in K^\prime\subset K$, последовательности  $\{(x_k, \sigma_k)\}$, $k\in K$. По лемме \ref{lem_3.2.3}  $\bar{x}\in D$, т.е.
\begin{equation}\label{o19}
f(\bar{x})\ge f^*.
\end{equation}
Выделим из последовательности  $\{\bar{u}_{i_k}\}$, $k\in K^{\prime\prime}\subset K$ сходящуюся подпоследовательность $\{\bar{u}_{i_k}\}$, $k\in K^\prime$. Пусть  $\bar{u}$ -- ее предельная точка. Так как множество   ${\rm epi}\, (f, {\rm R_n})$ замкнуто, то
\begin{equation}\label{o20}
\bar{u}\in {\rm epi}\, (f, {\rm R_n}).
\end{equation}
С учетом условия  (\ref{o2}) выбора чисел   $\varepsilon_k$, а также равенств  $u_{i_k}=(x_k, \sigma_k)$ из неравенств   $\|\bar{u}_{i_k}-u_{i_k}\|\le \varepsilon_k$, $k\in K^{\prime\prime}$, следует, что $\bar{u}=(\bar{x}, \bar{\sigma}).$ Тогда ввиду  (\ref{o20}) выполняется включение  $(\bar{x}, \bar{\sigma})\in {\rm epi}\, (f, {\rm R_n})$, т.е. 
\begin{equation}\label{o21}
f(\bar{x})\le \bar{\sigma}.
\end{equation}
Но согласно (\ref{o11}) $\sigma_k\le f^*$, $k\in K^{\prime\prime}$, и $\bar{\sigma}\le f^*$. Отсюда и из неравенств (\ref{o21}), (\ref{o19}) следуют равенства $f(\bar{x})=\bar{\sigma}=f^*$. Теорема доказана.
\end{proof}

\setcounter{equation}{0}


\section
[Метод с аппроксимацией допустимой области 
 и надграфика, использующий обновление погружающих множеств]
{Метод с аппроксимацией допустимой области
\linebreak и надграфика, использующий обновление погружающих множеств}

В данном параграфе предлагается метод отсечений для задачи выпуклого программирования, характерный тем, что допускает возможность периодического  обновления множеств, аппроксимирующих как надграфик целевой функции, так и область ограничений.

Пусть $f(x)$, $F(x)$ -- выпуклые в  ${\rm R_n}$ функции,
$$D^\prime = \{x\in {\rm R_n} : F(x) \le 0\},$$
$D^{\prime \prime} \subset {\rm R_n}$ -- выпуклое замкнутое множество,
$$D=D^\prime \bigcap D^{\prime\prime},$$
функция $f(x)$ достигает на множестве $D$ минимального значения. Решается задача (\ref{e3.1.1}).

Отметим, что как и в $\S\ 1.2$ внутренность множества $D$ может быть пустой, например, за счет ${\rm int\,} D^\prime = \emptyset$ или ${\rm int\,}  D^{\prime\prime}=\emptyset$.

Пусть $f^*$, $X^*$, $K$ определены как и в $\S \ 1.2$, a ${\rm epi\,} (f, G)$, $\partial f(x)$, $\partial F(x)$ -- как в $\S \ 2.1$. Положим  $D^\prime_\varepsilon = \{x\in {\rm R_n} : F(x) \le \varepsilon \}$, $X^*(\varepsilon) = \{x\in D: f(x) \le f^* + \varepsilon \}$, где $\varepsilon \ge 0$. Пусть, кроме того, $x^*\in X^*$.

\textbf{Метод 3.3.} Метод заключается в следующем. Выбираются выпуклые замкнутые множества $M_0\subset {\rm R_n}$, $G_0\subset {\rm R_{n+1}}$ такие, что
\begin{equation}\label{e3.2.2}
x^*\in M_0, \quad {\rm epi\,} (f, {\rm R_n})\subset G_0,
\end{equation}
и для любой точки $(x, \gamma)\in G_0$, где $x\in M_0$, выполняется неравенство
\begin{equation}\label{e3.2.3}
\gamma \ge \bar{\gamma} > -\infty.
\end{equation}
Задаются числа $\varepsilon_0> 0$, $\delta_0>0$, и полагается $i=0$, $k=0$.

\textbf{1.} Отыскивается точка $(y_i, \gamma_i)$, где $y_i\in {\rm R_n}$, $\gamma_i\in {\rm R_1}$, как решение задачи
\begin{equation}\label{e3.2.4}
\gamma \to \min
\end{equation}
\begin{equation}\label{e3.2.5}
x\in M_i\bigcap D^{\prime\prime}, \quad (x, \gamma) \in G_i.
\end{equation}
Если
\begin{equation}\label{e3.2.6}
y_i\in D^\prime,
\end{equation}
и
\begin{equation}\label{e3.2.7}
f(y_i)=\gamma_i,
\end{equation}
то $y_i\in X^*$, и задача (\ref{e3.1.1}) решена.

\textbf{2.} Если выполняется включение (\ref{e3.2.6}) и 
\begin{equation}\label{e3.2.8}
f(y_i) - \gamma_i \le \delta_k,
\end{equation}
то $y_i\in X^*(\delta_k)$. С нахождением $\delta_k$-решения задачи (\ref{e3.1.1}) процесс завершается, либо следует переход к п. 3.

\textbf{3.} Строятся множества $M_{i+1}$, $G_{i+1}$ способами, зависящими от выполнения для точки $(y_i, \gamma_i)$ следующих условий.

а) Пусть $(y_i, \gamma_i)$ такова, что
\begin{equation}\label{e3.2.9}
y_i\notin D^\prime_{\varepsilon_k}.
\end{equation}
Тогда выбираются конечные множества
\begin{equation}\label{e3.2.10}
A_i\subset \partial F(y_i),
\end{equation}
\newline
\begin{equation}\label{e3.2.11}
B_i\subset \partial f(y_i),
\end{equation}
и полагается 
\begin{equation}\label{e3.2.12}
Q_i=M_i,
\end{equation}
\newline
\begin{equation}\label{e3.2.13}
S_i=G_i,
\end{equation}
\newline
\begin{equation}\label{e3.2.14}
M_{i+1} = Q_i \bigcap \{x\in {\rm R_n} : F(y_i) + \langle a, x - y_i\rangle \le 0 \ \forall a\in A_i\},
\end{equation}
\newline
\begin{equation}\label{e3.2.15}
G_{i+1} = S_i\bigcap \{(x, \gamma)\in {\rm R_{n+1}} : f(y_i) + \langle b, x - y_i\rangle \le \gamma \ \forall b\in B_i\}.
\end{equation}

б) Пусть точка $(y_i, \gamma_i)$ такова, что
\begin{equation}\label{e3.2.16}
y_i \in D^\prime_{\varepsilon_k},
\end{equation}
но при этом выполняется неравенство
\begin{equation}\label{e3.2.17}
f(y_i) - \gamma_i > \delta_k.
\end{equation}
Тогда, во-первых, полагается
\begin{equation}\label{e3.2.18}
M_{i+1} = Q_i = M_i,
\end{equation}
если $y_i\in D^\prime$, либо множество $M_{i+1}$ строится в виде (\ref{e3.2.14}), (\ref{e3.2.12}), (\ref{e3.2.10}), если $y_i\in D^\prime_{\varepsilon_k}\setminus D^\prime$. Во-вторых, строится множество $G_{i+1}$ согласно (\ref{e3.2.15}), (\ref{e3.2.13}), (\ref{e3.2.11}).

в) Пусть, наконец, для точки $(y_i, \gamma_i)$ одновременно выполняются включение (\ref{e3.2.16}) и неравенство (\ref{e3.2.8}). Во-первых, если $y_i\in D^\prime$, то множества $Q_i$ и $M_{i+1}$ строятся в виде (\ref{e3.2.18}), а в противном случае  выбирается выпуклое замкнутое множество $Q_i\subset {\rm R_n}$ с условием, что 
\begin{equation}\label{e3.2.19}
Q_i\subseteq M_0, \quad x^* \in Q_i,
\end{equation}
и $M_{i+1}$ строится в виде (\ref{e3.2.14}), где $A_i$ задано согласно (\ref{e3.2.10}), а $Q_i$ -- согласно (\ref{e3.2.19}). Во-вторых, выбирается выпуклое замкнутое множество $S_i\subset {\rm R_{n+1}}$ такое, что 
\begin{equation}\label{e3.2.20}
S_i\subseteq G_0, \quad {\rm epi\,} (f, {\rm R_n}) \subset S_i,
\end{equation}
 и строится множество $G_{i+1}$  в виде равенства (\ref{e3.2.15}), в котором $B_i$ задано согласно (\ref{e3.2.11}), а $S_i$ -- согласно (\ref{e3.2.20}). В-третьих, полагается $i_k = i$,
\begin{equation}\label{e3.2.21}
x_k = y_{i_k}, \quad \sigma_k=\gamma_{i_k},
\end{equation}
задаются $\varepsilon_{k+1}>0$, $\delta_{k+1}>0$, и значение $k$ увеличивается на единицу.

\textbf{4.} Следует переход к п. 1 при $i$, увеличенном на единицу.

Сделаем относительно предложенного метода некоторые замечания. Прежде всего покажем, что при любом $i\in K$ задача (\ref{e3.2.4}), (\ref{e3.2.5}) разрешима. Для этого обоснуем два следующих вспомогательных утверждения.

\begin{lemma}\label{lemma3.2.1}
Для каждого $i\in K$ имеют место следующие включения:
\begin{equation}\label{e3.2.22}
x^*\in M_i,
\end{equation}

\begin{equation}\label{e3.2.23}
{\rm epi\,} (f, {\rm R_n})\subset G_i.
\end{equation}
\end{lemma}

\begin{proof} При $i=0$ включение (\ref{e3.2.22}), (\ref{e3.2.23}) выполняются ввиду условиям (\ref{e3.2.2}). Допустим теперь, что (\ref{e3.2.22}), (\ref{e3.2.23}) справедливы при любом фиксированном $i=l\ge 0$. Покажем, что
\begin{equation}\label{e3.2.24}
x^*\in M_{l+1},
\end{equation}

\begin{equation}\label{e3.2.25}
{\rm epi\,} (f, {\rm R_n}) \subset G_{l+1}.
\end{equation}
Тогда утверждение леммы будет доказано.

Докажем сначала (\ref{e3.2.24}). Положим
\begin{equation}\label{e3.2.26}
P_i= \{x\in {\rm R_n} : F(y_i) + \langle a, x - y_i \rangle \le 0 \ \forall a\in A_i\}, \quad i\in K.
\end{equation}
В силу п. 3 метода множество $M_{l+1}$ построено согласно либо равенству (\ref{e3.2.18}), либо равенству (\ref{e3.2.14}). В случае $M_{l+1} = M_l$ включение (\ref{e3.2.24}) выполняется ввиду сделанного индукционного предположения. Поэтому пусть $M_{l+1} = Q_l\bigcap P_l$ согласно (\ref{e3.2.14}). Заметим, что в связи с условиями (\ref{e3.2.12}), (\ref{e3.2.19}) выбора множества $Q_l$ справедливо включение
\begin{equation}\label{e3.2.27}
x^*\in Q_l.
\end{equation}
Далее, для каждого $a\in A_l$ имеет место неравенство $F(x^*) - F(y_l)\ge \langle a, x^* - y_l\rangle$, а значит, и неравенство $F(y_l) + \langle a, x^* - y_l\rangle \le 0$. Следовательно, $x^* \in P_l$. Отсюда и из (\ref{e3.2.27}) вытекает (\ref{e3.2.24}).

Докажем теперь (\ref{e3.2.25}). Положим 
\begin{equation}\label{e3.2.28}
T_i = \{(x, \gamma) \in {\rm R_{n+1}} : f(y_i) + \langle b , x - y_i \rangle \le \gamma \ \forall b\in B_i\}, \quad i\in K.
\end{equation}
Согласно п. 3 метода $G_{l+1} = S_l \bigcap T_l$, независимо от условий (\ref{e3.2.16}), (\ref{e3.2.8}) для точки $(y_l, \gamma_l)$. Заметим, что в силу (\ref{e3.2.13}), (\ref{e3.2.20}) c учетом индукционного предположения выполняется включение
\begin{equation}\label{e3.2.29}
{\rm epi\,} (f, {\rm R_n}) \subset S_l.
\end{equation}
Покажем, что 
\begin{equation}\label{e3.2.30}
{\rm epi\,} (f, {\rm R_n})\subset T_l.
\end{equation}
Зафиксируем произвольно точку $(\bar{x}, \bar{\gamma}) \in {\rm epi\,} (f, {\rm R_n})$. Тогда с учетом неравенства $\bar{\gamma} \ge f(\bar{x})$ для всех $b\in B_l$ имеем $\bar{\gamma} - f(y_l) \ge \langle b, \bar{x} - y_l\rangle$. Следовательно, $(\bar{x}, \bar{\gamma})\in T_l$, и утверждение (\ref{e3.2.30}) доказано. Таким образом, из (\ref{e3.2.29}), (\ref{e3.2.30}) следует (\ref{e3.2.25}). Лемма доказана.
\end{proof}

\begin{lemma}\label{lemma3.2.2}
Точка $(x^*, f^*)$ для всех $i\in K$ удовлетворяет ограничениям задачи (\ref{e3.2.4}), (\ref{e3.2.5}).
\end{lemma}
\begin{proof}
Так как $x^*\in D^{\prime\prime}$, и справедливо (\ref{e3.2.22}), то
\begin{equation}\label{e3.2.31}
x^*\in M_i\bigcap D^{\prime\prime} \quad \forall i\in K.
\end{equation}
Кроме того, поскольку $(x^*, f^*)\in {\rm epi\,}(f, {\rm R_n})$, то из (\ref{e3.2.23}) имеем
\begin{equation}\label{e3.2.32}
(x^*, f^*)\in G_i \quad \forall i\in K.
\end{equation}
Тогда из (\ref{e3.2.31}), (\ref{e3.2.32}) следует утверждение леммы.
\end{proof}

Лемма \ref{lemma3.2.2}, а также условия (\ref{e3.2.2}), (\ref{e3.2.3}) выбора множеств $M_0$, $G_0$ говорят о том, что задача (\ref{e3.2.4}), (\ref{e3.2.5}) разрешима при всех $i\in K$. Отметим также, что в силу леммы \ref{lemma3.2.1} для всех $i\in K$ множества $M_{i+1}$, $G_{i+1}$, построенные в виде (\ref{e3.2.14}), (\ref{e3.2.15}), (\ref{e3.2.18}), непусты.

Обоснуем теперь критерий остановки, заложенный в п. 1 метода.

\begin{lemma}\label{lemma3.2.3}
При каждом $i\in K$ справедливо неравенство 
\begin{equation}\label{e3.2.33}
\gamma_i \le f^*.
\end{equation}
\end{lemma}

\begin{proof}
Согласно (\ref{e3.2.4}), (\ref{e3.2.5}) для каждой точки $(x, \gamma)$, удовлетворяющей условиям (\ref{e3.2.5}), выполняется неравенство $\gamma_i\le \gamma$. Но по лемме \ref{lemma3.2.2} точка $(x^*, f^*)$ является допустимым решением задачи (\ref{e3.2.4}), (\ref{e3.2.5}) при любом $i\in K$. Отсюда следует (\ref{e3.2.33}). Лемма доказана.
\end{proof}

\begin{theorem}\label{theorem3.2.1}
Пусть при некотором $i\in K$  для точки $(y_i, \gamma_i)$ выполняется (\ref{e3.2.6}), (\ref{e3.2.7}). Тогда $y_i$ -- решение задачи (\ref{e3.1.1}).
\end{theorem}
\begin{proof}
В силу (\ref{e3.2.33}) и предположения (\ref{e3.2.7}) $f(y_i) \le f^*$. С другой стороны, поскольку справедливо (\ref{e3.2.6}) и ввиду (\ref{e3.2.5}) имеет место включение $y_i\in D^{\prime\prime}$, то $y_i\in D$, а значит, $f(y_i) \ge f^*$. Таким образом, $f(y_i) = f^*$, и теорема доказана.
\end{proof}

Далее, предположим, что для некоторой точки $(y_i, \gamma_i)$ выполнились включение (\ref{e3.2.6}) и неравенство (\ref{e3.2.8}). Тогда $y_i\in D$, и cогласно (\ref{e3.2.8}), (\ref{e3.2.33}) $f(y_i)\le f^*+\delta_k$. Таким образом, как и отмечено в п. 2 метода, $y_i\in X^*(\delta_k)$, то есть $y_i$ является $\delta_k$-решением исходной задачи.

Множества $M_i$ будем называть аппроксимирующими область ограничений исходной задачи в окрестности решения, а множества $G_i$ -- аппроксимирующими надграфик целевой функции. Плоскости $F(y_i) + \langle a, x - y_i\rangle = 0$, где $a\in \partial F(y_i)$, и $f(y_i) + \langle b, x - y_i\rangle = \gamma$, где $b\in \partial f(y_i)$, которые формируют согласно (\ref{e3.2.14}), (\ref{e3.2.15}) аппроксимирующие множества, будем называть отсекающими плоскостями.

Аппроксимирующие множества естественно строить для всех $i\in K$ многогранными. Тогда задачи (\ref{e3.2.4}), (\ref{e3.2.5}) будут задачами линейного программирования. Приведем некоторые замечания по поводу построения аппроксимирующих множеств.

\begin{remark}
При определенных дополнительных условиях на задачу (\ref{e3.1.1}) нет необходимости в специальном построении множеств $M_i$, $G_i$. Предположим сначала, что множество $D^\prime$ является многогранным, и покажем, что метод позволяет выбрать множества $M_i$ при всех $i\in K$ в виде
\begin{equation}\label{e3.2.34}
M_i = D^\prime.
\end{equation}
Действительно, при $i=0$ равенство (\ref{e3.2.34}) возможно в связи с условиями задания множества $M_0$ на предварительном шаге метода. Пусть при некотором $i=l\ge 0$ уже выполняется равенство $M_l=D^\prime$. Тогда в силу первого из включений (\ref{e3.2.5}) $y_l\in M_l$, а значит, $y_l\in D^\prime$. В таком случае согласно пп. 3б), 3в) метода $M_{l+1}=Q_l=M_l=D^\prime$, и возможность выбора множеств $M_i$  в виде (\ref{e3.2.34}) доказана.
\end{remark}

\begin{remark}
Предположим теперь, что целевая функция задачи (\ref{e3.1.1}) имеет вид
\begin{equation}\label{e3.2.35}
f(x) = \max \limits_{j\in J} f_j(x),
\end{equation}
где $J=\{1,\ldots,m\}$, $m\ge 1$. Пусть все функции $f_j(x)$ линейны. Тогда нет необходимости в построении специальных многогранных множеств, аппроксимирующих надграфик. Покажем, что метод допускает выбор множеств $G_i$ в виде
\begin{equation}\label{e3.2.36}
G_i = {\rm epi\,} (f, {\rm R_n})
\end{equation}
для всех $i\in K$. Действительно, при $i=0$  это утверждение справедливо согласно условиям задания множества $G_0$. Пусть на некотором шаге $l\ge 0$ уже построено множество $G_l$ в виде
\begin{equation}\label{e3.2.37}
G_l = {\rm epi\,} (f, {\rm R_n}).
\end{equation}
Тогда $G_{l+1}$ должно быть построено согласно (\ref{e3.2.15}), (\ref{e3.2.13}), если $y_l \notin D^\prime_{\varepsilon_k}$, либо согласно (\ref{e3.2.15}), (\ref{e3.2.20}), если $y_l\in D^\prime_{\varepsilon_k}$. \mbox{В первом} случае
\begin{equation}\label{e3.2.38}
S_l = {\rm epi\,} (f, {\rm R_n})
\end{equation}
ввиду (\ref{e3.2.13}), (\ref{e3.2.37}). Во втором случае условие (\ref{e3.2.20}) также позволяет задать $S_l$ в виде (\ref{e3.2.38}). Таким образом, в обоих случаях
$$G_{l+1} = {\rm epi\,}(f, {\rm R_n})\bigcap T_l,$$
а поскольку справедливо включение (\ref{e3.2.30}), то $G_{l+1} = {\rm epi\,} (f, {\rm R_n})$. Следовательно, возможность задания множеств $G$ в виде (\ref{e3.2.36}) доказана.
\end{remark}

Отметим, что для выбора множеств $M_0$, $G_0$ есть много возможностей. В частности, допустимо положить $M_0={\rm R_n}$. Однако, в этом случае при задании множества $G_0$ необходимо учесть условие (\ref{e3.2.3}). Если функция $f(x)$ имеет вид (\ref{e3.2.35}), и среди функций	 $f_j(x)$ есть линейные, например, с номерами $j\in J^\prime \subset J$, то удобно выбрать  $G_0={\rm epi\,} (g, {\rm R_n})$, где $g(x) = \max \limits_{j\in J^\prime} f_j(x)$.

Далее, обратим внимание на то, что условия (\ref{e3.2.19}), (\ref{e3.2.20}) дают большие возможности в выборе множеств $Q_i$, $S_i$, а значит, и соответствующих им множеств $M_{i+1}$, $G_{i+1}$. В частности, для всех $i\in K$, при которых имеют место соотношения (\ref{e3.2.16}), (\ref{e3.2.8}), множества $Q_i$, $S_i$ можно задать в виде (\ref{e3.2.12}), (\ref{e3.2.13}), поскольку условия (\ref{e3.2.19}), (\ref{e3.2.20}) для них будут выполнены в силу леммы \ref{lemma3.2.1}. Таким образом, с учетом (\ref{e3.2.18}) множества $Q_i$, $S_i$ при всех $i\in K$ либо имеют вид (\ref{e3.2.12}), (\ref{e3.2.13}), либо могут быть заданы в виде (\ref{e3.2.12}), (\ref{e3.2.13}).

\begin{remark}
Предположим теперь, что для всех $i\in K$, начиная с некоторого номера $i^\prime$, выполняется (\ref{e3.2.12}), (\ref{e3.2.13}). Тогда для всех  $i\ge i^\prime$ множества $M_{i+1}$ имеют вид (\ref{e3.2.18}), либо 
\begin{equation}\label{e3.2.39}
M_{i+1} = M_i \bigcap P_i,
\end{equation}
где $P_i$ определены в (\ref{e3.2.26}), а множества $G_{i+1}$ имеют вид
\begin{equation}\label{e3.2.40}
G_{i+1} = G_i\bigcap T_i,
\end{equation}
где $T_i$ определены в (\ref{e3.2.28}). Заметим, что в таком случае происходит накопление отсекающих плоскостей при построении множеств $G_i$, а если $y_i\notin D^\prime$ для бесконечного числа номеров $i\ge i^\prime$, то идет процесс накопления и плоскостей, формирующих множества $M_i$, то есть никаких обновлений погружающих множеств не происходит.
\end{remark}

Покажем теперь, как за счет выбора множеств $Q_i$, $S_i$ из условий (\ref{e3.2.19}), (\ref{e3.2.20}) можно проводить упомянутые обновления путем отбрасывания накапливающихся отсекающих плоскостей.

\begin{remark}
Пусть при некотором $i=i_k$ для точки $(y_i, \gamma_i)$ выполнились одновременно условия (\ref{e3.2.16}), (\ref{e3.2.8}). Положим, например,
\begin{equation}\label{e3.2.41}
Q_i = Q_{i_k} = M_r, \quad S_i = S_{i_k} = G_l,
\end{equation}
где $0\le r\le i_k$, $0\le l \le i_k$. Легко доказать, что $M_r \subset M_0$, $G_l \subset G_0$ для всех $r=0,\ldots, i_k$, $l=0,\ldots, i_k$. Кроме того, ввиду (\ref{e3.2.22}), (\ref{e3.2.23}) $x^*\in M_r$, ${\rm epi\,} (f, {\rm R_n})\subset G_l$ для всех $r=0, \ldots, i_k$, $l=0,\ldots,i_k$. Таким образом, при всех указанных $r$, $l$ условия п. 3в) метода для множеств (\ref{e3.2.41}) выполняются. Следовательно, в качестве $Q_{i_k}$ можно выбрать любое из построенных к $i_k$-ому шагу множеств $M_0, \ldots,M_{i_k}$, а в качестве $S_{i_k}$ -- любое из множеств $G_0,\ldots, G_{i_k}$. В частности, если положить 
$$Q_{i_k}=M_0, \quad S_{i_k}=G_0,$$
то произойдет отбрасывание всех накопившихся к шагу $i=i_k$ отсекающих плоскостей, и множества $M_{i+1}$, $G_{i+1}$ будут иметь, соответственно, вид $M_{i+1}=M_0\bigcap P_i$, $G_{i+1}=G_0\bigcap T_i$.
\end{remark}

Понятно, что в случае (\ref{e3.2.16}), (\ref{e3.2.8}) имеются и другие, отличные от (\ref{e3.2.41}), способы задания множеств $Q_i=Q_{i_k}$ и $S_i=S_{i_k}$.

Перейдем теперь к изучению свойств последовательностей, построенных предложенным методом и к исследованию сходимости метода. Будем предполагать далее, что построенная методом последовательность $\{(y_i, \gamma_i)\}$, $i\in K$, ограничена.

Как уже отмечено, независимо от условий (\ref{e3.2.16}), (\ref{e3.2.8}), множества $Q_i$ для всех $i\in K$ могут быть заданы в виде (\ref{e3.2.12}). С учетом этого сформулируем следующее вспомогательное утверждение.

\begin{lemma}\label{lemma3.2.4}
Пусть последовательность $\{(y_i, \gamma_i)\}$, $i\in K$, построена методом с условием, что для всех $i\in K$, начиная с некоторого номера $i^\prime \ge 0$ множества $Q_i$ имеют вид (\ref{e3.2.12}). Тогда любая предельная точка последовательности $\{y_i\}$, $i\in K$, принадлежит множеству $D$.
\end{lemma}
\begin{proof} Пусть $\{y_i\}$, $i\in K_1\subset K$, -- сходящаяся подпоследовательность последовательности $\{y_i\}$, $i\in K$, а $\bar{y}$ -- ее предельная точка. Если $\bar{y}\in D^\prime$, то в силу включений $y_i\in M_i \bigcap D^{\prime\prime}$, $i\in K_1$, и замкнутости $D^{\prime\prime}$ имеем $\bar{y}\in D^{\prime\prime}$, и утверждение доказано. Поэтому предположим далее, что $\bar{y}\notin D^\prime$, то есть 
$$F(\bar{y})>0.$$

Пусть $F(\bar{y}) = \beta$. Так как функция $F(x)$ непрерывна, то найдется окрестность $\omega$ точки $\bar{y}$ такая, что 
\begin{equation}\label{e3.2.42}
F(x) \ge \frac{\beta}{2}
\end{equation}
для всех $x\in \omega$. Зафиксируем такой номер $\tilde{i}\in K_1$, что $\tilde{i}\ge i^\prime$, и для всех $i\in K_1$, $i\ge \tilde{i}$ имеет место включение $y_i\in \omega$. Тогда с учетом (\ref{e3.2.42})
\begin{equation}\label{e3.2.43}
F(y_i)\ge \frac{\beta}{2} \quad \forall i\in K_1, \quad i\ge \tilde{i}.
\end{equation}
Отметим, что поскольку $y_i\notin D^\prime$ для всех $i\in K_1$, $i\ge \tilde{i}$, то согласно п. 3 метода и условию леммы множества $M_{i+1}$ для всех $i\in K_1$, $i\ge \tilde{i}$ построены в виде (\ref{e3.2.39}).

Зафиксируем номера $l$, $r\in K_1$ с условием, что $r>l\ge \tilde{i}$. Так как $M_{l+1}=M_l\bigcap P_l$, то $M_r \subset M_{l+1} \subset P_l$,  а поскольку $y_r\in M_r$ ввиду (\ref{e3.2.5}), то $y_r\in P_l$. Значит,  $F(y_l)+\langle a, y_r - y_l\rangle \le 0$ для всех $a\in A_l$. Отсюда и из (\ref{e3.2.43}) имеем
\begin{equation}\label{e3.2.44}
\langle a, y_r - y_l\rangle \le - \frac{\beta}{2}
\end{equation}
для всех $a\in A_l$.

Далее, ввиду ограниченности последовательности $\{y_i\}$, $i\in K$, найдется (cм., напр., \cite[121]{n_polayk_b_t_1}) такое число $\theta > 0$, что для всех $a\in \partial F(y_i)$, $b\in \partial f(y_i)$, $i\in K$, выполняются неравенства
\begin{equation}\label{e3.2.45}
\|a\|\le \theta, \quad \|b\|\le \theta.
\end{equation}
Тогда из (\ref{e3.2.44}), (\ref{e3.2.45}) следует, что
\begin{equation}\label{e3.2.46}
\theta \|y_r - y_l\|\ge \frac{\beta}{2}.
\end{equation}

Выберем теперь для каждого $i\in K_1$, $i\ge \tilde{i}$, такой номер $p_i\in K_1$, что $p_i\ge i+1$. Согласно (\ref{e3.2.46}) для выбранных номеров $i$, $p_i$ имеем $\theta\|y_{p_i}-y_i\|\ge \beta/ 2$. Последнее неравенство противоречиво, поскольку $y_i\to \bar{y}$ и $y_{p_i}\to \bar{y}$ при $i\to \infty$, $i\in K_1$. Лемма доказана.
\end{proof}

\begin{lemma}\label{lemma3.2.5}
Для точки $(y_i, \gamma_i)$ при каждом $i\in K$ справедливо неравенство
\begin{equation}\label{e3.2.47}
f(y_i)\ge \gamma_i.
\end{equation}
\end{lemma}
\begin{proof}
Так как $(y_i, f(y_i))\in {\rm epi\,} (f, {\rm R_n})$, $i\in K$, то в силу (\ref{e3.2.23})
\begin{equation}\label{e3.2.48}
(y_i, f(y_i))\in G_i \quad \forall i\in K.
\end{equation}
Согласно (\ref{e3.2.4}), (\ref{e3.2.5}) для любой точки $(y, \gamma)\in G_i$ такой, что $y\in M_i \bigcap D^{\prime\prime}$, выполняется неравенство $\gamma \ge \gamma_i$. Тогда отсюда и из (\ref{e3.2.48}) с учетом включения $y_i\in M_i \bigcap D^{\prime\prime}$ следует (\ref{e3.2.47}). Лемма доказана.
\end{proof}

Как уже сказано выше, множества $S_i$ для всех $i\in K$ могут иметь вид (\ref{e3.2.13}). C учетом этого сформулируем следующее утверждение.
\begin{lemma}\label{lemma3.2.6}
Пусть последовательность $\{(y_i, \gamma_i)\}$, $i\in K$, построена предложенным методом с условием, что для всех $i\in K$, начиная с некоторого номера $i^\prime\ge 0$ множества $S_i$ имеют вид (\ref{e3.2.13}). Тогда выполняется равенство
\begin{equation}\label{e3.2.49}
\lim \limits_{i\in K} (f(y_i) - \gamma_i) = 0.
\end{equation}
\end{lemma}
\begin{proof} Отметим, что ввиду (\ref{e3.2.47}) $f(y_i) - \gamma_i\ge 0$ для всех $i\in K$. Предположим, что утверждение леммы неверно. Тогда для $\{(y_i, \gamma_i)\}$, $i\in K$, найдется такая подпоследовательность $\{(y_i, \gamma_i)\}$, $i\in K_1\subset K$, $i\ge i^\prime$, что 
\begin{equation}\label{e3.2.50}
f(y_i) - \gamma_i \ge \alpha > 0 \quad \forall i\in K_1, \quad i\ge i^\prime.
\end{equation}
Пусть $\{(y_i, \gamma_i)\}$, $i\in K_2\subset K_1$, -- сходящаяся подпоследовательность, выделенная из точек $(y_i, \gamma_i)$, $i\in K_1$, $i\ge i^\prime$, а $(\bar{y}, \bar{\gamma})$ -- предельная точка этой подпоследовательности.

Зафиксируем такие номера $l$, $r\in K_2$, что $r>l$. Согласно п. 3 метода и условию леммы о том, что $S_i=G_i$, $i\in K_2$, множества $G_{i+1}$ для всех $i\in K_2$ имеют вид (\ref{e3.2.40}). Значит,
$$G_r \subset G_{l+1} = G_l \bigcap T_l \subset T_l,$$
а поскольку $(y_r, \gamma_r)\in G_r$, то $(y_r, \gamma_r)\in T_l$. Отсюда следует, что $f(y_l) + \langle b, y_r - y_l \rangle \le \gamma_r$ для всех $b\in B_l$. Но согласно (\ref{e3.2.50}) $f(y_l)\ge \alpha + \gamma_l$. Поэтому $\gamma_r - \gamma_l + \langle b, y_l - y_r\rangle \ge \alpha$ для всех $b\in B_l$. Из последнего неравенства и из (\ref{e3.2.45}) имеем
\begin{equation}\label{e3.2.51}
\gamma_r - \gamma_l + \theta \|y_l - y_r\|\ge \alpha.
\end{equation}

Зафиксируем теперь для каждого $i\in K_2$ такой номер $p_i\in K_2$, что $p_i\ge i+1$. Тогда в силу (\ref{e3.2.51}) для всех $i\in K_2$ выполняются неравенства
$$\gamma_{p_i} - \gamma_i + \theta \|y_i - y_{p_i}\|\ge \alpha.$$
Эти неравенства противоречивы, так как $y_i \to \bar{y}$, $y_{p_i}\to \bar{y}$ по $i\in K_2$, и $\gamma_i \to \bar{\gamma}$, $\gamma_{p_i} \to \bar{\gamma}$ по $i\in K_2$. Лемма доказана.
\end{proof}

\begin{consequence}\label{consequence3.2.1}
Пусть выполняются условия леммы \ref{lemma3.2.6}. Тогда для любой предельной точки $(\bar{y}, \bar{\gamma})$  последовательности $\{(y_i, \gamma_i)\}$, $i\in K$, имеет место равенство
$$f(\bar{y}) = \bar{\gamma}.$$
\end{consequence}
\begin{dokaz} утверждения непосредственно вытекает из леммы \ref{lemma3.2.6}.
\end{dokaz}
\begin{theorem}\label{theorem3.2.2}
Пусть последовательность $\{(y_i, \gamma_i)\}$, $i\in K$, построена методом с условием, что для всех $i\in K$, начиная с некоторого номера $i^\prime\ge 0$ множества $Q_i$, $S_i$ имеют, соответственно,  вид (\ref{e3.2.12}), (\ref{e3.2.13}). Тогда для любой предельной точки $(\bar{y}, \bar{\gamma})$ этой последовательности выполняется включение $\bar{y}\in X^*$.
\end{theorem}
\begin{proof}
Пусть $(\bar{y}, \bar{\gamma})$ -- предельная точка сходящейся подпоследовательности $\{(y_i, \gamma_i)\}$, $i\in K_1 \subset K$. По лемме \ref{lemma3.2.4} имеем $\bar{y}\in D$, то есть
\begin{equation}\label{e3.2.52}
f(\bar{y})\ge f^*.
\end{equation}
С другой стороны, согласно (\ref{e3.2.33}) $\gamma_i \le f^*$, $i\in K_1$, а значит, $\bar{\gamma} \le f^*$. Тогда по следствию \ref{consequence3.2.1} выполняется неравенство $f^*\ge f(\bar{y})$. Отсюда и из (\ref{e3.2.52}) следует, что $f(\bar{y}) = f^*$, и теорема доказана.
\end{proof}

Подчеркнем, что утверждение теоремы \ref{theorem3.2.2} доказано при условиях, что в п. 3 метода для всех $i\in K$ множества $M_{i+1}$ строятся в виде (\ref{e3.2.18}) или (\ref{e3.2.39}), а множества $G_{i+1}$ -- в виде (\ref{e3.2.40}). В таком случае при построении $\{(y_i, \gamma_i)\}$, $i\in K$, число отсекающих плоскостей неограниченно растет, обновлений аппроксимирующих множеств не происходит, и потребность в выделенных подпоследовательностях (\ref{e3.2.21}) отпадает.

Учитывая это замечание, перейдем к исследованию последовательностей $\{(y_i, \gamma_i)\}$, $i\in K$, при построении которых бесконечное число раз происходили обновления аппроксимирующих множеств за счет соответствующим условиям (\ref{e3.2.19}), (\ref{e3.2.20}) выбора множеств $Q_i$, $S_i$.

Покажем сначала, что при каждом $k\in K$ найдется точка $(y_i, \gamma_i)=(y_{i_k}, \gamma_{i_k})$ строящейся последовательности, для которой выполнятся условия п. 3в) метода, то есть условия (\ref{e3.2.16}), (\ref{e3.2.8}).

\begin{lemma}\label{lemma3.2.7}
Если последовательность $\{(y_i, \gamma_i)\}$, $i\in K$, построена согласно методу, то для каждого $k\in K$ существует такой номер $i=i_k$, что выполняются равенства (\ref{e3.2.21}).
\end{lemma}
\begin{proof} 1) Пусть $k=0$. Если $y_0\in D^\prime_{\varepsilon_0}$, и $f(y_0) - \gamma_0\le \delta_0$, то для выбранного $k$ равенства (\ref{e3.2.21}) имеют место при $i=i_0=0$. Поэтому будем считать, что $y_0\notin D^\prime_{\varepsilon_0}$ или $f(y_0) - \gamma_0 > \delta_0$. Покажем существование номера $i=i_0>0$, для которого одновременно выполняются включение
\begin{equation}\label{e3.2.53}
y_{i_0} \in D^\prime_{\varepsilon_0}
\end{equation}
и неравенство
\begin{equation}\label{e3.2.54}
f(y_{i_0}) - \gamma_{i_0} \le \delta_0.
\end{equation}

Предположим противное, то есть для всех $i>0$
\begin{equation}\label{e3.2.55}
y_i \notin D^\prime_{\varepsilon_0}
\end{equation}
или 
\begin{equation}\label{e3.2.56}
f(y_i) - \gamma_i > \delta_0.
\end{equation}
Cразу отметим, что в таком случае согласно п. 3 метода для всех $i\ge 1$ множества $Q_i$, $S_i$ имеют, соответственно, вид (\ref{e3.2.12}), (\ref{e3.2.13}). Выделим из последовательности $\{(y_i, \gamma_i)\}$, $i\in K$, $i>0$, сходящуюся подпоследовательность $\{(y_i, \gamma_i)\}$, $i\in K^\prime \subset K$,  и пусть $(y^\prime, \gamma^\prime)$ -- ее предельная точка. Тогда согласно лемме \ref{lemma3.2.4} и лемме \ref{lemma3.2.6} имеем
\begin{equation}\label{e3.2.57}
y^\prime\in D, \quad f(y^\prime) = \gamma^\prime.
\end{equation}

Допустим, что в множестве $K^\prime$ существует бесконечное число номеров $i$, при которых выполняется (\ref{e3.2.55}). Положим $\bar{K}^\prime = \{i\in K^\prime : y_i \notin D^\prime_{\varepsilon_0}\}$. Тогда $F(y_i) > \varepsilon_0$ для всех $i\in \bar{K}^\prime$, а значит, $F(y^\prime)\ge \varepsilon_0>0$, что невозможно ввиду первого из соотношений (\ref{e3.2.57}). Следовательно, найдется такой номер $i^\prime \in K^\prime$, что для всех $i\ge i^\prime$, $i\in K^\prime$, выполняется включение $y_i\in D^\prime_{\varepsilon_0}$. Но в таком случае согласно предположению (\ref{e3.2.55}), (\ref{e3.2.56}) для всех $i\ge i^\prime$, $i\in K^\prime$, справедливо неравенство \ref{e3.2.56}, и $f(y^\prime) - \gamma^\prime\ge \delta_0$. Последнее неравенство противоречит второму из соотношений (\ref{e3.2.57}). Таким образом, доказано существование номера $i_0>0$, для которого одновременно имеют место соотношения (\ref{e3.2.53}), (\ref{e3.2.54}),  а значит, и равенства $x_0 = y_{i_0}$, $\sigma_0 = \gamma_{i_0}$.

2) Предположим теперь, что при некотором фиксированном $k\ge 0$ уже выполняются равенства $x_k = y_{i_k}$, $\sigma_k = \gamma_{i_k}$, то есть  для точки $(y_i, \gamma_i) = (y_{i_k}, \gamma_{i_k})$ одновременно справедливы (\ref{e3.2.16}), (\ref{e3.2.8}). Покажем существование такого номера $i=i_{k+1}> i_k$, что
\begin{equation}\label{e3.2.58}
y_{i_{k+1}}\in D^\prime_{\varepsilon_{k+1}}, \quad f(y_{i_{k+1}}) - \gamma_{i_{k+1}} \le \delta_{k+1},
\end{equation}
тогда $x_{k+1}=y_{i_{k+1}}$, $\sigma_{k+1} = \gamma_{i_{k+1}}$, и лемма будет доказана.

Допустим противное, то есть для всех $i\in K$, $i>i_k$ выполняется $y_i \notin D^\prime_{\varepsilon_{k+1}}$ или $f(y_i) - \gamma_i > \delta_{k+1}$. Выделим из последовательности $\{(y_i, \gamma_i)\}$, $i\in K$, $i>i_k$,  сходящуюся подпоследовательность $\{(y_i, \gamma_i)\}$, $i\in K^{\prime\prime}\subset K$, и пусть $(y^{\prime\prime}, \gamma^{\prime\prime})$ -- ее предельная точка. Поскольку согласно сделанному
\linebreak
\newline
предположению множества $Q_i$, $S_i$ для всех $i\in K$, $i>i_k$ имеют вид (\ref{e3.2.12}), (\ref{e3.2.13}), то в силу леммы \ref{lemma3.2.4}, \ref{lemma3.2.6} имеем $y^{\prime\prime}\in D$,
\begin{equation}\label{e3.2.59}
f(y^{\prime\prime}) = \gamma^{\prime\prime}.
\end{equation}
Далее, как и в п. 1) доказательства, обосновывается существование такого номера $i^{\prime\prime}\in K^{\prime\prime}$, что $y_i\in D^\prime_{\varepsilon_{k+1}}$ для всех $i\in K^{\prime\prime}$, $i\ge i^{\prime\prime}$,  а затем получается неравенство $f(y^{\prime\prime}) -\gamma^{\prime\prime}\ge \delta_{k+1}$, противоречащее (\ref{e3.2.59}). Таким образом, номер $i=i_{k+1}$, при котором справедливы соотношения (\ref{e3.2.58}), существует. Лемма доказана.
\end{proof}

Леммой \ref{lemma3.2.7} обосновано существование наряду с $\{(y_i, \gamma_i)\}$, $i\in K$, последовательности $\{(x_k, \sigma_k)\}$, $k\in K$.

\begin{theorem}\label{theorem3.2.3}
Пусть последовательность $\{(x_k, \sigma_k)\}$, $i\in K$, построены методом 3.2 с условием, что
\begin{equation}\label{e3.2.60}
\varepsilon_k \to 0, \quad k\to \infty, \quad \delta_k \to 0, \quad k\to \infty.
\end{equation}
Тогда для любой ее предельной точки $(\bar{x}, \bar{\sigma})$ выполняется
\begin{equation}\label{e3.2.61}
\bar{x}\in X^*, \quad \bar{\sigma}=f^*.
\end{equation}
\end{theorem}
\begin{proof}
Отметим, что последовательность $\{(x_k, \sigma_k)\}$, $k\in K$, ограничена ввиду ограниченности $\{(y_i, \gamma_i)\}$, $i\in K$. Пусть $(\bar{x}, \bar{\gamma})$ -- предельная точка сходящейся подпоследовательности $\{(x_k, \sigma_k)\}$, $k\in K^\prime \subset K$, выделенной из последовательности  $\{(x_k, \sigma_k)\}$, $k\in K$. Так как по построению $x_k = y_{i_k}\in D^\prime_{\varepsilon_k}$, $k\in K$, то $F(x_k) \le \varepsilon_k$, $k\in K$, и с учетом (\ref{e3.2.60}) $F(\bar{x})\le 0$ или $\bar{x}\in D^\prime$. Кроме того, согласно (\ref{e3.2.5}) $x_k\in D^{\prime\prime}$, $k\in K$, и в силу замкнутости множества $D^{\prime\prime}$ имеет место включение $\bar{x}\in D^{\prime\prime}$. Таким образом, $\bar{x}\in D$, и потому
\begin{equation}\label{e3.2.62}
f(\bar{x})\ge f^*.
\end{equation}

Далее, ввиду выполнения для номеров $i=i_k$, $k\in K$, условия (\ref{e3.2.8}) справедливы неравенства $f(y_{i_k}) - \gamma_{i_k}\le \delta_k$ или $f(x_k) \le \delta_k + \sigma_k$, $k\in K$.
Отсюда и из второго предельного сотношения (\ref{e3.2.60}) имеем
$$f(\bar{x})\le \bar{\sigma}.$$
С другой стороны, $f(x_k) \ge \sigma_k$, $k\in K$, согласно (\ref{e3.2.47}), а значит, $f(\bar{x})\le \bar{\sigma}$. Следовательно,
\begin{equation}\label{e3.2.63}
f(\bar{x})=\bar{\sigma}.
\end{equation}
По лемме \ref{lemma3.2.3} для всех $k\in K$ выполняются неравенства $\sigma_k \le f^*$. Поэтому $\bar{\sigma} \le f^*$. Отсюда и из (\ref{e3.2.62}), (\ref{e3.2.63}) следуют равенства $f(\bar{x}) = f^* = \bar{\sigma}$, то есть утверждения (\ref{e3.2.61}) справедливы. Теорема доказана.\end{proof}

Пусть при условиях теоремы \ref{theorem3.2.3} последовательность $\{f(x_k)\}$, $k\in K$, является монотонной. Тогда в силу ограниченности $\{(x_k, \sigma_k)\}$, $k\in K$, последовательность $\{f(x_k)\}$, $k\in K$, сходится, и по теореме \ref{theorem3.2.3} выполняется равенство $\lim \limits_{k\in K} f(x_k) = f^*$. Отсюда и из известной теоремы (напр., \mbox{\cite[62]{n_vasillev_f_p_1}}) следует сходимость последовательности $\{x_k\}$, $k\in K$, к множеству $X^*$.

При исследовании сходимости метода использовалось условие ограниченности последовательности $\{(y_i, \gamma_i)\}$, $i\in K$. Ясно, что выполнение этого условия можно обеспечить, например, соответствующим выбором множества $M_0$.
         
\setcounter{equation}{0}


\section
{Численные эксперименты}

Данный параграф посвящен обсуждению результатов численных экспериментов, которые получены при исследовании метода 3.3. При проведении численных экспериментов была поставлена следующая цель. Нужно было проверить работоспособность метода 3.3 и дать рекомендации по использованию тех или иных приемов обновления аппроксимируюших множеств $M_i$, $G_i$. Кроме того, требовалось предложить способы задания чисел $\delta_k$, $\varepsilon_k$, которые определяют качества аппроксимации надграфика целевой функции и  допустимой области соответственно.

Метод 3.3 программно реализован на языке программирования С++. Численные эксперименты проводились на персональном компьютере, конфигурация которого была описана в $\S \ 1.5$.

Сделано значительное число экспериментов на тестовых примерах с различным чиcлом переменных (до 50) для выработки вышеуказанных рекомендаций.

Для проведения численных экспериментов была выбрана следующая тестовая задача:
$$f(x) = \sum_{i=1}^ni (\xi_i - a)^2 \to \min,$$
$$f_1(x) = \sum_{i=1}^ni(\xi_i - b)^2-1600\le 0,$$
где $n \ge 1$, $a=10$, $b=5$, $x=(\xi_1,\ldots,\xi_n)\in {\rm R_n}$.

Приведем критерий, при выполнении которого прекращалась работа алгоритма. Процесс отыскания решения задачи останавливался, если на некотором шаге $i^\prime\in K$ для метода 3.3 выполнялось условие $y_i\in X^*(\varepsilon)$, где множество $X^*(\varepsilon)$ определено в $\S \ 3.3$, а
$$\varepsilon=0.00001.$$
Таким образом, по определению множества $X^*(\varepsilon)$ имеет место следующее соотношение:
$$f(y_{i^\prime}) - f^*\le \varepsilon.$$

В методе 3.3 при решении тестовых примеров  множество $M_0$  задавалось в виде куба, содержащего всю допустимую область,   $G_0$ выбиралось совпадающим со  всем пространством ${\rm R_{n+1}}$, и полагалось $\bar{\gamma}= -1000000$.

Результаты численных экпериментов приведены в следующей таблице:

\begin{center}
\begin{tabular}{|c|c|c|c|c|c|c|c|}
\hline
Номер  &  $n$  & $\delta_k$ & $\varepsilon_k$ &  $Q_i$ & $S_i$ &Кол-во  итераций&  Время в мин. \\
\hline
1& 30 & $a_1$ & $c_1$&$b_1$&$d_1$&426&3,53\\
\hline
2& 30 & $a_1$ & $c_1$&$b_2$&$d_2$&1016&10,59\\
\hline
3& 30 & $a_1$ & $c_1$&$b_3$&$d_3$&984791&664,14\\
\hline
4& 30 & $a_2$ & $c_2$&$b_1$&$d_1$&758&8,00\\
\hline
5& 30 & $a_2$ & $c_2$&$b_2$&$d_2$&2216&17,06\\
\hline
6& 30 & $a_2$ & $c_2$&$b_3$&$d_3$&602881&436,89\\
\hline
7& 30 & $a_3$ & $c_3$&$b_1$&$d_1$&1054&16,00\\
\hline
8& 30 & $a_3$ & $c_3$&$b_2$&$d_2$&2798&29,41\\
\hline
9& 30 & $a_3$ & $c_3$&$b_3$&$d_3$&385678&237,65\\
\hline
10& 30 & $a_4$ & $c_4$&$b_1$&$d_1$&4632&66,47\\
\hline
11& 30 & $a_4$ & $c_4$&$b_2$&$d_2$&4465&62,94\\
\hline
12& 30 & $a_4$ & $c_4$&$b_3$&$d_3$&59288&112,37\\
\hline
13& 50 & $a_1$ & $c_1$&$b_1$&$d_1$&719&7,60\\
\hline
14& 50 & $a_1$ & $c_1$&$b_2$&$d_2$&1774&19,42\\
\hline
15& 50 & $a_1$ & $c_1$&$b_3$&$d_3$&-&-\\
\hline
16& 50 & $a_2$ & $c_2$&$b_1$&$d_1$&1471&12,40\\
\hline
17& 50 & $a_2$ & $c_2$&$b_2$&$d_2$&3706&34,29\\
\hline
18& 50 & $a_2$ & $c_2$&$b_3$&$d_3$&-&-\\
\hline
19& 50 & $a_3$ & $c_3$&$b_1$&$d_1$&1631&24,00\\
\hline
\end{tabular}
\end{center}

\begin{center}
\begin{tabular}{|c|c|c|c|c|c|c|c|}
\hline
Номер  &  $n$  & $\delta_k$ & $\varepsilon_k$ &  $Q_i$ & $S_i$ &Кол-во  итераций&  Время в мин. \\
\hline
20& 50 & $a_3$ & $c_3$&$b_2$&$d_2$&4368&69,45\\
\hline
21& 50 & $a_3$ & $c_3$&$b_3$&$d_3$&-&-\\
\hline
22& 50 & $a_4$ & $c_4$&$b_1$&$d_1$&7588&118,59\\
\hline
23& 50 & $a_4$ & $c_4$&$b_2$&$d_2$&7585&111,10\\
\hline
24& 50 & $a_4$ & $c_4$&$b_3$&$d_3$&-&-\\
\hline
\end{tabular}
\end{center}

Перейдем к обсуждению обозначения четвертой и пятой колонок. Для удобства записи в таблице  последовательностей $\{\varepsilon_k\}$, $k\in K$, выбраны следующие соответствия:

1) $a_1$ $\leftrightarrow$ $\varepsilon_{k+1}=\frac{\varepsilon_k}{1.1}$ для всех $k\in K$;

2) $a_2$ $\leftrightarrow$ $\varepsilon_{k+1}=\frac{\varepsilon_k}{2}$ для всех $k\in K$;

3) $a_3$ $\leftrightarrow$ $\varepsilon_{k+1}=\frac{\varepsilon_k}{n}$ для всех $k\in K$;

4) $a_4$ $\leftrightarrow$ $\varepsilon_{k+1} = \alpha_k F(x_k)$, $\alpha_k = \frac{1}{2^k}$ для всех $k\in K$,

а для последовательности $\{\delta_k\}$, $k\in K$:

1) $c_1$ $\leftrightarrow$ $\delta_{k+1}=\frac{\delta_k}{1.1}$ для всех $k\in K$;

2) $c_2$ $\leftrightarrow$ $\delta_{k+1}=\frac{\delta_k}{2}$ для всех $k\in K$;

3) $c_3$ $\leftrightarrow$ $\delta_{k+1}=\frac{\delta_k}{n}$ для всех $k\in K$;

4) $c_4$ $\leftrightarrow$ $\delta_{k+1} = \alpha_k (f(x_k)-\sigma_k)$, $\alpha_k = \frac{1}{2^k}$ для всех $k\in K$.

Способы обновления множеств $S_i$ в таблице обозначены следующим образом:

1)  $b_1$ $\leftrightarrow$  $S_i=G_i$ для всех $i\ne i_k$, а $S_{i_k}$ состоит из   <<активных>>  в текущей итерационной точке $(y_{i_k}, \gamma_{i_k})$ секущих плоскостей;

2) $b_2$ $\leftrightarrow$ $S_i=G_i$, $0\le i< n+1$,

  $S_i = \begin{cases}G_i, \ i\ne i_k, \\
 \overset{i_k}{\underset{r=i_k - n - 1}{\cap}}\{(x, \gamma)\in R_{n+1} : f(y_r) + \langle a, x - y_r \rangle\le \gamma,  \\
 \hfill a\in \partial f(y_r) \},  i = i_k,  \end{cases}$

$i\ge n+1$;

3)  $b_3$ $\leftrightarrow$ $S_i = \begin{cases}G_i, & i \ne i_k,  \\ G_0,& i = i_k. \end{cases}$

Для  множеств $Q_i$ приняты следующие соответствия:

1)  $d_1$ $\leftrightarrow$  $Q_i=M_i$ для всех $i\ne i_k$, а  множество $Q_{i_k}$ выбирается как пересечение  $M_0$ с некоторым множеством $T_{i_k}$, которое задается  <<активными>>  в текущей итерационной точке $y_{i_k}$ секущими плоскостями;

2)  $d_2$ $\leftrightarrow$ $Q_i=M_i$, $0\le i< n+1$,

  $Q_i = \begin{cases}M_i, \ i\ne i_k, \\
 \overset{i_k}{\underset{r=i_k - n - 1}{\cap}}\{x\in R_{n} : f_1(y_r) + \langle a, x - y_r \rangle\le 0,  \\ 
 \hfill a\in \partial f_1(y_r) \}, i = i_k,  \end{cases}$

$i\ge n+1$;

3)    $d_3$ $\leftrightarrow$ $S_i = \begin{cases}M_i, & i \ne i_k,  \\ M_0,& i = i_k. \end{cases}$

Обратим внимание на то, что все результаты, приведенные в таблице этого параграфа, получены методом 3.3 с привлечением тех или иных способов обновления аппроксимирующих множеств. Указанный пример (при $n=30, 50$) решался методом 3.3 и без использования каких-либо обновлений, то есть при всех $i$ полагалось  $Q_i=M_i$, $S_i=G_i$. В этих случаях указанная выше точность решения была получена за 5128 итераций, 174,15 минут при $n=30$, а в случае $n=50$ за $12923$ итераций и $431,72$ минут.

Для начала объясним смысл содержимого колонок. В первой колонке указана нумерация численных экспериментов, во второй $-$ размерность пространства переменных, в третьей колонке приведены способы  выбора чисел $\delta_k$, а в четвертой $-$ задание чисел $\varepsilon_k$. В пятой и шестой колонках определены приемы обновления погружающих множеств, которые аппроксимируют надграфик целевой функции и допустимую область соответственно. В седьмой колонке представлены количество итераций, которое потребовалось для отыскания решения задачи с  заданной точностью, а в восьмой $-$ время решения.

Приступим к обсуждению полученных результатов. \mbox{При сравнении} результатов тестовых примеров можно сделать вывод о том, что на итерациях $i=i_k$ не рекомендуется полностью обновлять все погружающие множества, аппроксимирующие как надграфик целевой функции, так и допустимую область. Вызвано это тем, что в момент полного обновления аппроксимирующих множеств  полностью утрачивается качество аппроксимации надграфика и допустимой области. Отсюда значение целевой функции в очередной итерационной точке может резко либо вырасти, либо уменьшится. В связи с этим для достижения утраченного качества аппроксимации надграфика и допустимой области требуется значительное число дополнительных итераций. Если решено делать  в процессе решения полное обновление аппроксимирующих множеств, то последовательности чисел $\delta_k$, $\varepsilon_k$ следует выбирать быстроубывающими.

Перейдем теперь к обсуждению других способов обновления аппроксимирующих множеств. Пусть на  всех итерациях $i=i_k$ в множества $Q_{i_k}$ и $S_{i_k}$ включаются все  <<активные>> гиперплоскости, либо во множества $Q_{i_k}$ и $S_{i_k}$ входят $n+1$ последних отсекающих плоскостей. Применение этих  процедур  обновления позволяет существенно уменьшить время решения тестовых задач. Следовательно, при решении задач нецелесообразно копить все отсекающие гиперплоскости. При проведении численных экспериментов выяснилось, что более предпочтительным  является прием, на основе которого  в аппроксимирующие множества $Q_{i_k}$ и $S_{i_k}$ включаются только <<активные>> гиперплоскости.

Завершим данный параграф обсуждением выбора чисел $\delta_k$, $\varepsilon_k$, $k\in K$.  Ранее отмечалось, что при полном обновлении аппроксимирующих множеств последоваетльности чисел $\delta_k$, $\varepsilon_k$, $k\in K$, следует выбирать быстроубывающей. Если же при обновлении аппроксимирующих множеств выбираются только <<активные>> гиперплоскости, то числа $\delta_k$, $\varepsilon_k$, $k\in K$, рекомендуется выбирать медленноубывающими.

\printbibliography[title=Литература]
}
\begin{center}
\textit{Научное издание}
\end{center}

\begin{center}
\textbf{Заботин} Игорь Ярославич\linebreak
\textbf{Яруллин} Рашид Саматович
\end{center}

\begin{center}
\textbf{МЕТОДЫ ОТСЕЧЕНИЙ БЕЗ ВЛОЖЕНИЯ 
АППРОКСИМИРУЮЩИХ МНОЖЕСТВ ДЛЯ ЗАДАЧ МАТЕМАТИЧЕСКОГО ПРОГРАММИРОВАНИЯ}
\end{center}

\vspace{7cm}

\begin{center}
Подписано в печать 16.10.2019.\linebreak
Бумага офсетная. Печать цифровая.\linebreak
Формат 60х84 1/16. Гарнитура «Times New Roman». Усл. печ. л. 10,5.\linebreak
Уч.-изд. л. 43,2. Тираж 200 экз. Заказ 201/10\linebreak
\end{center}

\begin{center}
Отпечатано с готового оригинал-макета\linebreak
в типографии Издательства Казанского университета
\end{center}

\begin{center}
420008, г. Казань, ул. Профессора Нужина, 1/37\linebreak
тел. (843) 233-73-59, 233-73-28 
\end{center}

\end{document}